\documentclass[11pt]{amsart}
\usepackage[hmargin=3cm,vmargin=3cm]{geometry}

\usepackage[latin1]{inputenc}
\usepackage{xspace,amssymb,amsfonts,euscript}
\usepackage{amsthm,amsmath,amscd,stmaryrd,latexsym}
\usepackage{palatino, euscript}
\usepackage{mathtools}
\usepackage{graphicx}
\usepackage[all]{xy}
\usepackage{tikz}
\usetikzlibrary{calc, matrix, arrows, decorations.pathmorphing}
\usepackage[dvips]{epsfig}
\usepackage{pinlabel}
\usepackage{psfrag}

\usepackage{float}
\usepackage{blkarray}

\IfFileExists{comments.sty}{\usepackage{comments}}

\RequirePackage{color}
\definecolor{myred}{rgb}{0.75,0,0}
\definecolor{mygreen}{rgb}{0,0.5,0}
\definecolor{myblue}{rgb}{0,0,0.65}

\RequirePackage{ifpdf}
\ifpdf
  \IfFileExists{pdfsync.sty}{\RequirePackage{pdfsync}}{}
  \RequirePackage[pdftex,
   colorlinks = true,
   urlcolor = myblue, % \href{...}{...} external (URL)
   citecolor = mygreen,  % \cite{...}
   linkcolor = myred, % \ref{...} and \pageref{...}
   pagebackref,
   bookmarksopen=true]{hyperref}
\else
  \RequirePackage[hypertex]{hyperref}
\fi

\newtheorem{thm}{Theorem}[section]
\newtheorem{lemma}[thm]{Lemma}

\newtheorem{prop}[thm]{Proposition}
\newtheorem{cor}[thm]{Corollary}

\newtheorem*{prop*}{Proposition}

\theoremstyle{definition}
\newtheorem{defn}[thm]{Definition}

\newtheorem{notation}[thm]{Notation}

\newtheorem{example}[thm]{Example}

\theoremstyle{remark}
\newtheorem{remark}[thm]{Remark}
\newtheorem{rmk}[thm]{Remark}

\numberwithin{equation}{section}

     \def\ak{{\mathbbm{k}}}

\let\phi=\varphi

\usepackage{bbm}

\def\1{\mathbbm{1}}

\newcommand{\ot}{\otimes}

\newcommand{\symp}{\mathfrak{sp}}

\newcommand{\refequal}[1]{\xy {\ar@{=}^{#1}
(-1,0)*{};(1,0)*{}};
\endxy}

\DeclareMathOperator{\Hom}{{\rm Hom}}

\DeclareMathOperator{\End}{{\rm End}}

\DeclareMathOperator{\Ext}{{\rm Ext}}

\DeclareMathOperator{\id}{{\rm id}}

\DeclareMathOperator{\Fund}{\bf{Fund}}
\DeclareMathOperator{\Rep}{\bf{Rep}}

\DeclareMathOperator{\Tilt}{\bf{Tilt}}
\DeclareMathOperator{\wt}{\rm{wt}}
\DeclareMathOperator{\eval}{\Xi}
\DeclareMathOperator{\cupa}{\textbf{cup}_1}
\DeclareMathOperator{\cupb}{\textbf{cup}_2}
\DeclareMathOperator{\capa}{\textbf{cap}_1}
\DeclareMathOperator{\capb}{\textbf{cap}_2}
\DeclareMathOperator{\pp}{\textbf{p}}
\DeclareMathOperator{\ii}{\textbf{i}}
\newcommand{\VV}{\text{$V$}}
\newcommand{\DD}{\text{$\mathcal{D}_{\mathfrak{sp}_4}$}}
\newcommand{\DDk}{\text{$\mathcal{D}_{\mathfrak{sp}_4}^{\mathbbm{k}}$}}

\newcommand{\blues}{\textcolor{blue}{1}}
\newcommand{\greent}{\textcolor{green}{2}}

\begin{document}

\begin{abstract} Using Kuperberg's web calculus \cite{Kupe}, and following Elias and Libedinsky, we describe a ``light leaves" algorithm to construct a basis of morphisms between arbitrary tensor products of fundamental representations for $\mathfrak{sp}_4$ (and the associated quantum group). Our argument has very little dependence on the base field. As a result, we prove that when $[2]_q\ne 0$, the Karoubi envelope of the $C_2$ web category is equivalent to the category of tilting modules for the divided powers quantum group $U_q^{\mathcal{A}}(\mathfrak{sp}_4)$.

%Finally, we define a new diagrammatic category which we argue is equivalent to both the monoidal category $\Rep(\mathbb{C}^{\times}\rtimes 
%\mathbb{Z}/2)$ and the semisimplification of the category of tilting modules for $U_q^{\mathbb{Z}}(\mathfrak{sp}_4)$ at a primitive third root %of unity. 

%Then, following the philosophy of object adapted cellular categories, we derive "triple clasp" formulas for idempotents projecting to the top summand in each tensor product of fundamental representations and derive explicit formulas for the coefficients occurring in the clasp formulas (these are computed as "local intersection forms"). Our formulas provide further evidence for Elias's clasp conjecture, which was given for type $A$ webs, and suggests how to generalize the entire story to non-simply laced types.

\end{abstract}

\title{Web Calculus and Tilting Modules in Type $C_2$}

\author{Elijah Bodish} \address{University of Oregon, Eugene} \email{ebodish@uoregon.edu}

\maketitle

\tableofcontents

%%%%%%%%%%%%%%%%%%%%%%%%%
\section{Introduction}
\label{sec-intro}
%%%%%%%%%%%%%%%%%%%%%%%%%

Let $\mathfrak{g}$ be a complex semisimple Lie algebra and let $\Rep(\mathfrak{g})$ denote the category of finite dimensional modules for $\mathfrak{g}$.  By Weyl's theorem on complete reducibility $\Rep(\mathfrak{g})$ is a semisimple category, so as an abelian category $\Rep(\mathfrak{g})$ is determined by the number of its simple objects. Since isomorphism classes of finite dimensional irreducible $\mathfrak{g}$-modules are in bijection with the countably infinite set of dominant integral weights $X_+$, $\Rep(\mathfrak{g})\cong \Rep(\mathfrak{g}')$ as abelian categories, for any two semisimple Lie algebras.  

A Lie algebra acts on the tensor product of two representations, so $\Rep(\mathfrak{g})$ is a monoidal category. Viewing $\Rep(\mathfrak{g})$ as a monoidal semisimple category, we capture much more information about $\mathfrak{g}$ (the amount of information can be made precise through Tannaka--Krein duality). One then may ask for a presentation by generators and relations of the monoidal category $(\Rep(\mathfrak{g}), \ot)$. A modern point of view on this problem is to find a combinatorial replacement for $\Rep(\mathfrak{g})$ and then use planar diagrammatics to describe the combinatorial replacement by generators and relations. 

By combinatorial replacement, we mean a full subcategory of $\Rep(\mathfrak{g})$ monoidally generated by finitely many objects, such that all objects in $\Rep(\mathfrak{g})$ are direct sums of summands of objects in the subcategory. We will focus on the combinatorial replacement $\Fund(\mathfrak{g})$, which is the full subcategory of $\Rep(\mathfrak{g})$ monoidally generated by the irreducible modules $\VV(\varpi)$ of highest weight $\varpi$ for all fundamental weights $\varpi$. Note that $\Fund(\mathfrak{g})$ is not an additive category. 

We use the terminology $\mathfrak{g}$-webs to refer to a diagrammatic category equivalent to $\Fund(\mathfrak{g})$. The history of $\mathfrak{g}$-webs begins with the Temperley--Lieb algebra \cite{RTW, TemLie} for $\mathfrak{sl}_2$ and Kuperberg's ``rank two spiders" \cite{Kupe} for $\mathfrak{sl}_3$, $\mathfrak{sp}_4\cong \mathfrak{so}_5$, and $\mathfrak{g}_2$. D. Kim gave a conjectural presentation for $\mathfrak{sl}_4$-webs \cite{Kim03}, and then Morrison gave a conjectural description of $\mathfrak{sl}_n$-webs \cite{morr07}. Proving that the diagrammatic category was equivalent to $\Fund(\mathfrak{sl}_n)$ proved difficult, but was eventually carried out by Cautis, Kamnitzer, and Morrison using skew Howe duality \cite{CKM}. Recently a conjectural description of $\mathfrak{sp}_6$-webs has appeared in a preprint by Rose and Tatham for $\mathfrak{sp}_6$ \cite{rose2020webs}. 

The category of $\mathfrak{g}$-webs has a $q$ deformation, and an integral form, which we denote by $\mathcal{D}_{\mathfrak{g}}$, over $\mathbb{Z}[q, q^{-1}]$ (or some localization). On the representation theory side we have Lusztig's divided powers form of the quantum group, denoted $U_q^{\mathbb{Z}}(\mathfrak{g})$. This algebra has modules $V^{\mathbb{Z}}(\varpi)$, which are lattices inside $\VV(\varpi)$, for each fundamental weight. One should keep in mind that these lattices may not be irreducible after scalar extension to a field. The full subcategory monoidally generated by the modules $V^{\mathbb{Z}}(\varpi)$ will be denoted $\Fund(U_q^{\mathbb{Z}}(\mathfrak{g}))$. Taking all sums of summands of objects in $\Fund(\ak \ot U_q^{\mathbb{Z}}(\mathfrak{g}))$, one obtains the category of tilting modules $\Tilt(\ak \ot U_q^{\mathbb{Z}}(\mathfrak{g}))$. 

Let $\ak$ be a field and let $q\in \ak^{\times}$. We can specialize the integral versions of both the diagrammatic category and the combinatorial replacement category to $\ak$. It is natural to ask if these two categories are equivalent \cite[5A.4]{tiltcellular}. 

For $\mathfrak{g}= \mathfrak{sl}_n$ an answer to this question appears in a paper of Elias \cite{elias2015light}. Using ideas from Libedinsky's work \cite{LibRA} on constructing bases for maps between Soergel bimodules, Elias constructs a set of diagrams, denoted $\mathbb{LL}$ and referred to as double ladders, in the $\mathbb{Z}[q, q^{-1}$-linear category $\mathcal{D}_{\mathfrak{sl}_n}$. There are two main arguments in \cite{elias2015light}. First, a diagrammatic argument shows that $\mathbb{LL}$ spans the category over $\mathbb{Z}[q, q^{-1}]$. Second, Elias describes a functor $\Gamma: \mathcal{D}_{\mathfrak{sl}_n}\rightarrow \Fund(U_q^{\mathbb{Z}}(\mathfrak{sl}_n))$ and proves that $\Gamma(\mathbb{LL})$ is linearly independent. After observing that the ranks of homomorphism spaces in $\Fund(\ak \ot U_q^{\mathbb{Z}}(\mathfrak{sl}_n))$ are equal to $\#\mathbb{LL}$ \cite{QWeylDPS}, it follows that the diagrams $\ak \ot\mathbb{LL}$ are a basis for $\ak\ot \mathcal{D}_{\mathfrak{sl}_n}$ and the functor $\ak\ot\Gamma$ is an equivalence. 

Kuperberg proved \cite{Kupe} there is a monoidal equivalence $\ak \ot \mathcal{D}_{\mathfrak{sp}_4}\rightarrow\Fund(\ak \ot U_q(\mathfrak{sp}_4))$, when $\ak =\mathbb{C}(q)$and when $\ak= \mathbb{C}$ and $q=1$. Our goal is to prove this equivalence with as few restrictions on $\ak$ and $q$ as possible.  

The present work is completely indebted to Elias's approach, and the basis we construct for Kuperberg's $\mathcal{D}_{\mathfrak{sp}_4}$ webs is the analogue of Elias's light ladder basis for $\mathfrak{sl}_n$-webs in \cite{elias2015light}. However, our arguments take less effort, since we can use Kuperberg's result \cite{Kupe} that non-elliptic webs span $\mathcal{D}_{\mathfrak{sp}_4}$ over $\mathbb{Z}[q, q^{-1}]$, and are a basis for $\mathcal{D}_{\mathfrak{sp}_4}$ over $\mathbb{C}$, when $q=1$. Most of our work is to carefully construct an explicit functor $\eval: \mathcal{D}_{\mathfrak{sp}_4}\rightarrow U_q^{\mathbb{Z}}(\mathfrak{sp}_4)-\text{mod}$. 

The following theorem is the main result of the paper.
\begin{thm}\label{mainthm}
If $\ak$ is a field and $q\in \ak^{\times}$ is such that $q+ q^{-1} \ne 0$, then the functor 
\[
\eval: \ak \ot \mathcal{D}_{\mathfrak{sp}_4} \longrightarrow \Fund(\ak\ot U_q^{\mathbb{Z}}(\mathfrak{sp}_4)).
\]
is a monoidal equivalence, and therefore induces a monoidal equivalence between the Karoubi envelope of $\ak \ot \mathcal{D}_{\mathfrak{sp}_4}$ and the category $\textbf{Tilt}(\ak \ot U_q^{\mathbb{Z}}(\mathfrak{sp}_4))$. 
\end{thm}

\begin{remark}
The reader who is already well acquainted with \cite{Kupe} may wonder why we are talking about type $C_2$ and $\mathfrak{sp}_4$, instead of type $B_2$ and $\mathfrak{so}_5$. This certainly makes no difference classically, since $\mathfrak{sp}_4(\mathbb{C}) \cong \mathfrak{so}_5(\mathbb{C})$. For the purposes of this paper there is no difference integrally either. Under our hypothesis that $q+q^{-1} \ne 0$, there is an isomorphism $\ak \ot U^{\mathbb{Z}}_q(\mathfrak{sp}_4)\cong \ak \ot U^{\mathbb{Z}}_q(\mathfrak{so}_5)$, as well as an equivalence between $\ak \ot \mathcal{D}_{\mathfrak{sp}_4}$ and the base change from $\mathbb{Z}[q, q^{-1}]$ to $\ak$ of Kuperberg's $B_2$ spider category. 

We chose $C_2$ over $B_2$ hoping it would prevent confusion, since the defining relations in $\mathcal{D}_{\mathfrak{sp}_4}$ are slightly different than the relations in Kuperberg's $B_2$ spider.
\end{remark}

\begin{remark}
If $q+ q^{-1} = 0$, then the fundamental representation $\ak \ot V^{\mathbb{Z}}(\varpi_2)$ is not tilting. So if one is interested in tilting objects the category $\Fund(\mathfrak{g})$ is not the correct category to study. Also, the category $\mathcal{D}_{\mathfrak{sp}_4}$ is not defined when $q+q^{-1}= 0$, because some relations have coefficients with $q+ q^{-1}$ in the denominator. One could clear denominators in the relations and obtain a category which is defined when $q+ q^{-1} = 0$. However, we do not know what this diagrammatic category would describe. 
\end{remark}

The following result is a consequence of our main theorem, and is new even if $\ak = \mathbb{C}$ and $q=1$ or if $\ak =\mathbb{C}(q)$.
\begin{thm}
Let $\ak$ be a field and let $q\in \ak^{\times}$ so that $q+ q^{-1}\ne 0$. The double ladder diagrams defined in section \eqref{doubleladslabel} form a basis for the morphism spaces in $\ak \ot\mathcal{D}_{\mathfrak{sp}_4}$. 
\end{thm}

\begin{remark}
As we have already mentioned, Kuperberg's $B_2$ web category is spanned by the same non-elliptic diagrams over $\mathbb{Z}[q, q^{-1}]$. The work of Sikora--Westbury \cite{confluencegraphs} proves that these diagrams are linearly independent whenever $q+q^{-1}\ne 0$. Although their techniques are quite different than ours and certainly are worth studying, their result is a consequence of ours.

Suppose that one could show that either double ladder diagrams span or are linearly independent. Since the number of double ladders is equal to the number of non-elliptic webs, the result from \cite{confluencegraphs} would imply that the double ladder diagrams are a basis. 

However, it is not possible to obtain our main theorem with just their result. Even though their paper and some basic representation theory implies the dimensions of homomorphism spaces in $\ak \ot \mathcal{D}_{\mathfrak{sp}_4}$ and $\Fund(\ak \ot U_q^{\mathbb{Z}}(\mathfrak{sp}_4))$ are equal, it is not enough to deduce that $\ak \ot \eval$ is an equivalence. The difficulty is best illustrated via analogy: the lattice $\mathbb{Z}$ becomes a one-dimensional vector space after base change to any field, but the map $\mathbb{Z}\xrightarrow{x\mapsto 2x} \mathbb{Z}$ is not an isomorphism after tensoring with a field of characteristic two. We really need to know that the map $\ak \ot \eval$ is an isomorphism and to do this we must explicitly construct and analyze the functor $\eval$.
\end{remark}

\begin{remark}
It remains an open problem to adapt the arguments in \cite{elias2015light} to prove that double ladder diagrams span $\mathcal{D}_{\mathfrak{sp}_4}$ without using Kuperberg's results about non-elliptic webs. 
\end{remark}

\begin{remark}
It is work in progress of Victor Ostrik and Noah Snyder to find the precise relationship between Kuperberg's $G_2$ webs and tilting modules. 
\end{remark}

%%%%%%%%%%%%%%%%%%%%%%%%%
\subsection{Potential Applications}
\label{sec-appappapp}
%%%%%%%%%%%%%%%%%%%%%%%%%

Let $\ak= \mathbb{C}$ and let $q= e^{2\pi i/2\ell}$. Soergel conjectured \cite{SoergelCombinatoric} and then proved \cite{SoergelTilt} a formula for the character of a tilting module for $\ak\ot U_q^{\mathbb{Z}}(\mathfrak{g})$ when $\ell> h$, where $h$ is the Coxeter number of $\mathfrak{g}$.

The results of this paper \eqref{Dissoacc} imply that the category $\mathcal{D}_{\mathfrak{sp}_4}$ is a strictly object adapted cellular category \cite{ELauda}. Thus, the discussion in \cite[11.5]{soergelbook} allows one to adapt the algorithm in \cite{jensen2015pcanonical} from the context of Soergel bimodules to $\mathfrak{sp}_4$-webs. So one can compute tilting characters for the quantum group at a root of unity as long as $\ell\ge 3$ (the $\ell= 2$ case is ruled out by the assumption in our theorem that $q+ q^{-1} \ne 0$). The Coxeter number of $\mathfrak{sp}_4$ is $h= 4$. This means that when $\ell= 3$, Soergel's conjecture for tilting characters does not apply but the diagrammatic category $\ak \ot\mathcal{D}_{\mathfrak{sp}_4}$ does still describe tilting modules. 

There may be a conjecture for the characters of tilting modules of quantum groups that includes $\ell\le h$, along the lines of \cite{eliaslosev} and \cite[Theorem 1.6]{riche2020smithtreumann}. Ideally, the conjecture would relate tilting characters for the quantum group at a root of unity to singular, antispherical Kazhdan-Lusztig polynomials. One could use $\mathfrak{sp}_4$ webs to check such a conjecture for small weights.

There are other open questions related to tilting modules when $\ell$ is large enough for the diagrammatic category to be equivalent to the category of tilting modules, but $\ell$ is still less than the Coxeter number. For example, what is the semisimplification of the category of tilting modules for such $\ell$? The solution to this problem when $\ell> h$ is very well known, and provides a wealth of examples of finite tensor categories. When $\ell>h$ satisfies certain congruence conditions based on the root system of $\mathfrak{g}$ (for $\mathfrak{sp}_4$ the condition is $\ell$ is even) the semisimplification of the category of tilting modules is a modular category \cite{rowellqmod} which gives rise to a Reshetikhin--Turaev $3$-manifold invariant \cite{Tur}. Theorem \eqref{mainthm} implies that $\mathcal{D}_{\mathfrak{sp}_4}$ can be used to aid in the calculation of these three manifold invariants. 

By interpreting $\mathfrak{sl}_n$-webs in terms of the Schur algebra, Brundan, Entova-Aizenbud, Etingof, and Ostrik \cite{brundan2020semisimplification} were able to use results of Donkin to reprove that $\ak \ot \mathcal{D}_{\mathfrak{sl}_n}$ is equivalent to $\ak \ot \Fund(U_q^{\mathbb{Z}}(\mathfrak{sl}_n))$ for any field when $q= 1$. The main result of \cite{brundan2020semisimplification} is that when $q= 1$ and char $\ak< h$, the semisimplification of $\Tilt(\ak \ot U_q(\mathfrak{gl}_n))$ is a semisimple monoidal category, which may have infinitely many objects, and is related to Kazhdan--Lusztig cells the affine Hecke algebra. 

When Soergel's results on tilting characters of the quantum group are known to hold, Ostrik proved \cite{ostrik_1997} that there is a bijection between cells in the antispherical module for the Langlands dual affine Hecke algebra and thick monoidal ideals in the category of tilting modules $\Tilt(\ak \ot U_q^{\mathbb{Z}}(\mathfrak{g}))$. On the other hand, a deep theorem of Lusztig \cite{Lusztigcells_1989} is that there is also a bijection between cells in the antispherical module for affine Weyl group and orbits in the nilpotent cone of $G$. Note that the bijection between nilpotent orbits and thick monoidal ideals no longer appears to involve Langlands duality.

The maximal thick monoidal ideal in the category of tilting modules corresponds to the ``highest" cell in the antispherical module which in turn corresponds to the regular nilpotent orbit. This maximal ideal coincides with the ideal of negligible morphisms,  denoted by $\mathcal{N}$, and therefore the quotient is what is referred to as the semisimplification of the category of tilting modules. Since Soergel's methods of proof don't apply when $\ell< h$, it follows that Ostrik's results also do not apply. There may still be a non-trivial negligible ideal, but it might be that the objects in it now correspond to a different cell in the antispherical module and correspondingly a different nilpotent orbit.

When $\mathfrak{g}= \mathfrak{sp}_4$ and $\ell=3$, we still have a nontrivial semisimplification (this is not the case when $\ell=2$) and now the ``highest" cell is replaced by the unique reduced expression cell. The unique reduced expression cell corresponds via Lusztig's bijection to the sub regular nilpotent orbit $\mathcal{O}_{\text{subreg}}$. The group $\text{Sp}_4(\mathbb{C})$ acts on this orbit by conjugation. Now, fix a point $u\in \mathcal{O}_{\text{subreg}}$ in the orbit. The stabilizer of $u$ is an algebraic group with maximal reductive quotient, denoted $G_u$, a two component disconnected group with a one dimensional torus for the identity component. As an abstract group $G_u$ is an extension of $\mathbb{Z}/2$ by $\mathbb{C}^{\times}$. We conjecture that $G_u$ is a split but nontrivial extension. 

Motivated by these observations, we expect the following. Let $\ak= \mathbb{C}$ and let $q\in \mathbb{C}$ be a primitive $2\ell$-th root of unity for $\ell =3$ or $4$. There is an equivalence of monoidal categories $\Tilt(\ak \ot U_q^{\mathbb{Z}}(\mathfrak{sp}_4))/\mathcal{N}\longrightarrow \Rep(\mathbb{C}^{\times}\rtimes\mathbb{Z}/2)$. In order to prove this we will certainly need to use the results of this paper, as well as develop something like webs for the group $\mathbb{C}^{\times}\rtimes\mathbb{Z}/2$. Other work in progress of the author which stems from the results in this paper is adapting Elias's clasp conjectures \cite{elias2015light} to $\mathfrak{sp}_4$ webs. Work in progress of Ben Elias and Geordie Williamson uses $\mathcal{D}_{\mathfrak{sp}_4}$ to extend the quantum algebraic Satake equivalence \cite{Elias_2017} to type $B_2/C_2$. 

%As well as finding a diagrammatic description of the category of tilting modules when $p=2$, and using the diagrammatic category to study the monoidal ideals in $\Tilt(\text{Sp}_4(\overline{\mathbb{F}_2}))$. 

%%%%%%%%%%%%%%%%%%%%%%%%%
\subsection{Structure of the Paper}
\label{sec-struckture}
%%%%%%%%%%%%%%%%%%%%%%%%%

Section $2$: We discuss how to decompose tensor products of representations for $\mathfrak{sp}_4$. Then use the plethysm patterns to describe an algorithm for light ladder diagrams. Finally we define the double ladder diagrams. Section $3$: We define an evaluation functor from the diagrammatic category to the representation theoretic category. After reviewing some of the theory of tilting modules for quantum groups/reductive algebraic groups, we interpret the image of the evaluation functor as an integral form of the category of tilting modules. Then we argue that the main theorem follows from linear independence of the image of the double ladder diagrams. Section $4$: We argue that the double ladder diagrams are linearly independent. 

%Chapter $4$: We recall some results from the theory of monoidal categories about semisimplification. Then we introduce a new diagrammatic category $%\mathcal{D}_{G_u}$. We then argue this new diagrammatic category is equivalent to both $\Rep(\mathbb{C}^{\times}\rtimes\mathbb{Z}/2)$ and the %semisimplification of $\Tilt$ at a third root of unity. 

%%%%%%%%%%%%%%%%%%%%%%%%%
\subsection{Acknowledgements}
\label{sec-ackackack}
%%%%%%%%%%%%%%%%%%%%%%%%%

I want to thank Ben Elias for teaching me the philosophy of light leaves, which this work is guided by, and helping me prepare this document for mass consumption. I also want to thank Victor Ostrik and Noah Snyder for some very helpful discussions about webs and tilting modules. Finally, I am very thankful to both referees for giving me substantial comments to help improve the exposition.

%%%%%%%%%%%%%%%%%%%%%%%%%
\section{Light Ladders in Type $C_2$}
\label{sec-webs}
%%%%%%%%%%%%%%%%%%%%%%%%%

%===========
\subsection{$C_2$-Webs}
\label{subsec-webs}
%===========

We use the convention that the quantum integers in $\mathbb{Z}[q, q^{-1}]$ are defined as 
\begin{equation}
[n]_q= \dfrac{q^n- q^{-n}}{q-q^{-1}}, \ \ \ \ \ \text{for} \ n\in \mathbb{Z}.
\end{equation}
Let $\mathcal{A} = \mathbb{Z}[q, q^{-1}, [2]_q^{-1}]$, the ring $\mathbb{Z}[q, q^{-1}]$ localized at $[2]_q$.

\begin{defn}
Let $\mathcal{D}$ be the $\mathcal{A}$-linear monoidal category defined by generators and relations. The generating objects  are $\blues$ and $\greent$, the generating morphisms are the following diagrams. 
\begin{figure}[H]
\centering
\includegraphics[width=0.1\textwidth]{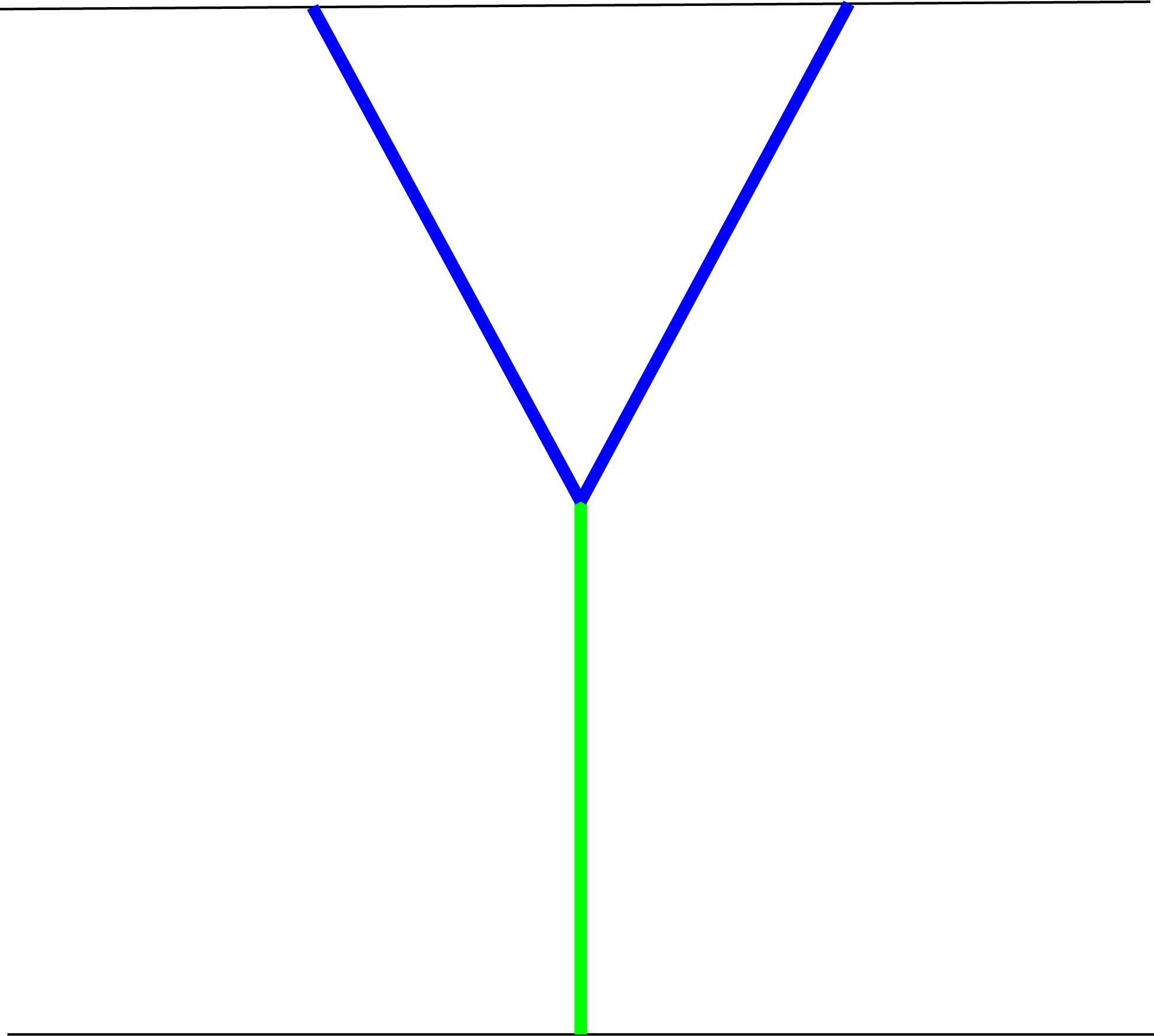} \ \ \ \ \ 
\includegraphics[width=0.1\textwidth]{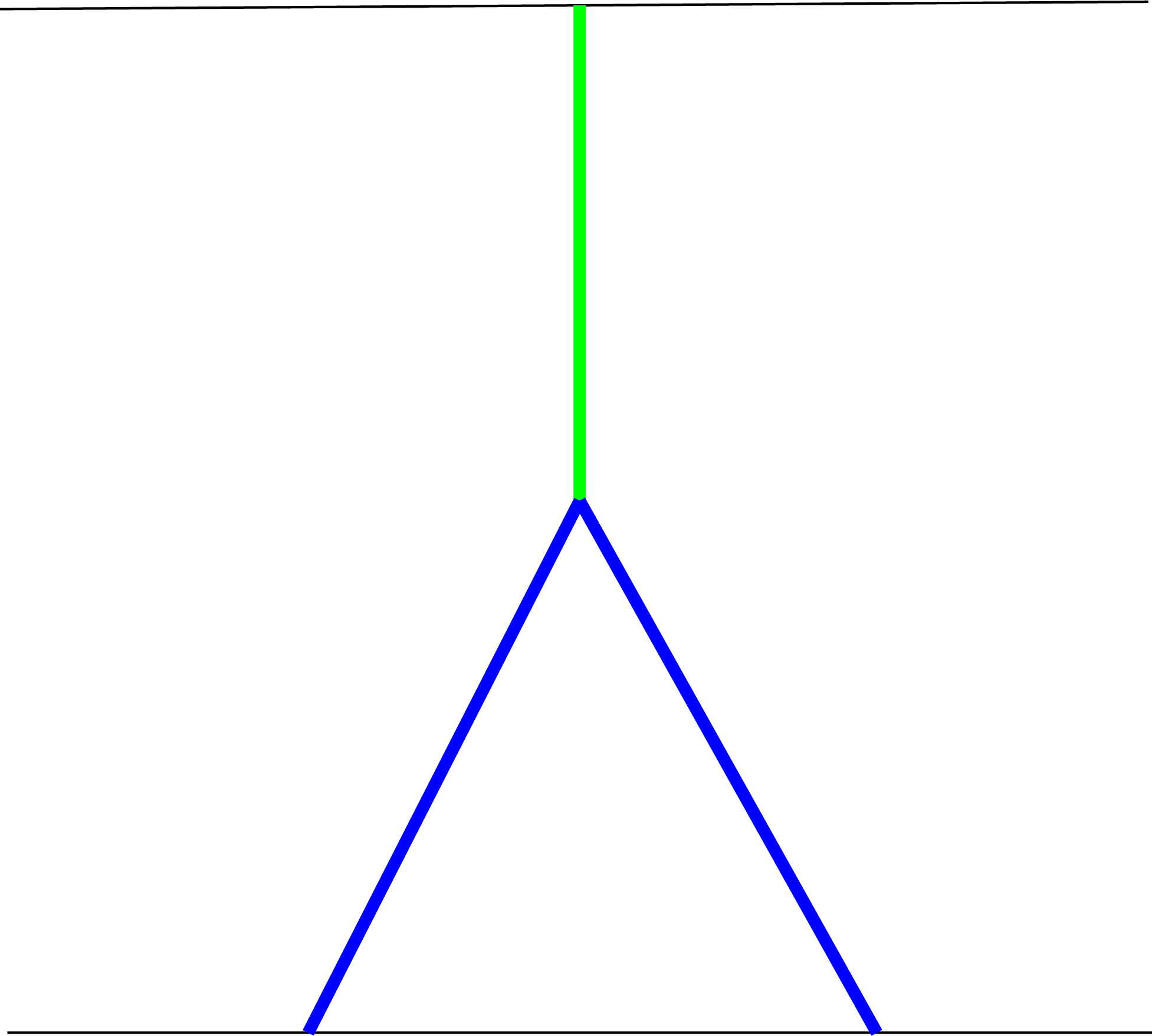} \ \ \ \ \ 
\includegraphics[width=0.1\textwidth]{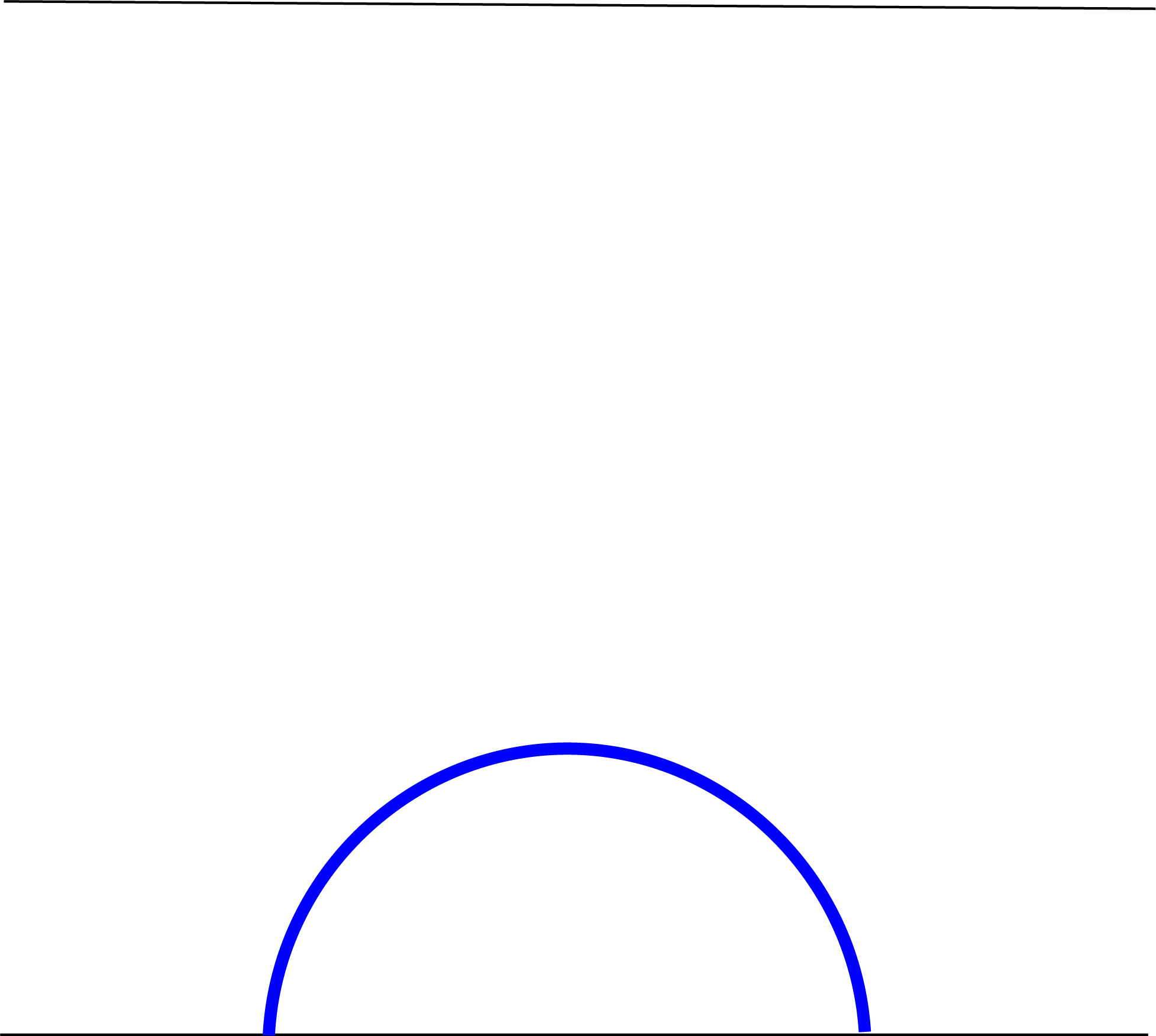} \ \ \ \ \ 
\includegraphics[width=0.1\textwidth]{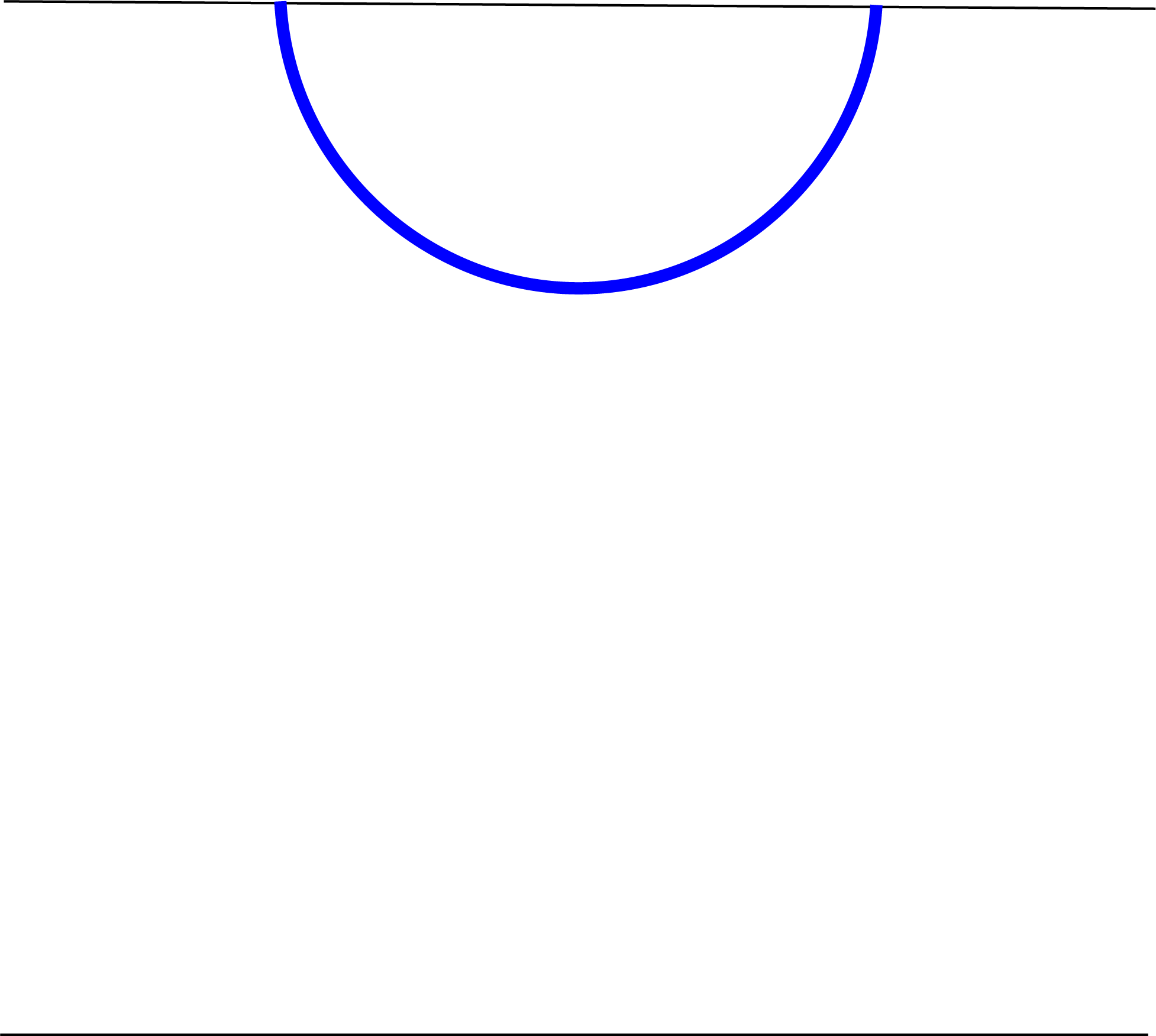} \ \ \ \ \ 
\includegraphics[width=0.1\textwidth]{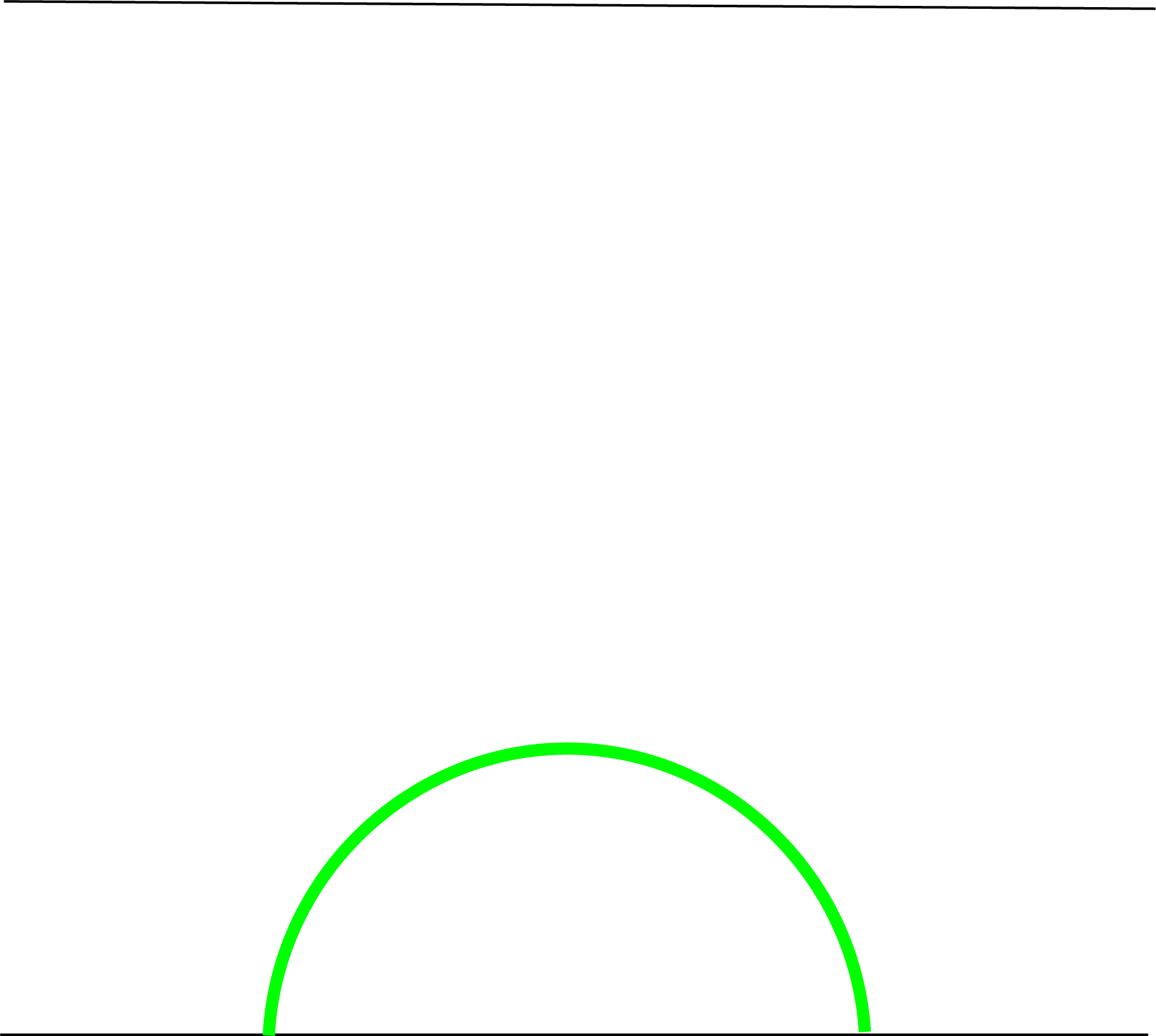} \ \ \ \ \ 
\includegraphics[width=0.1\textwidth]{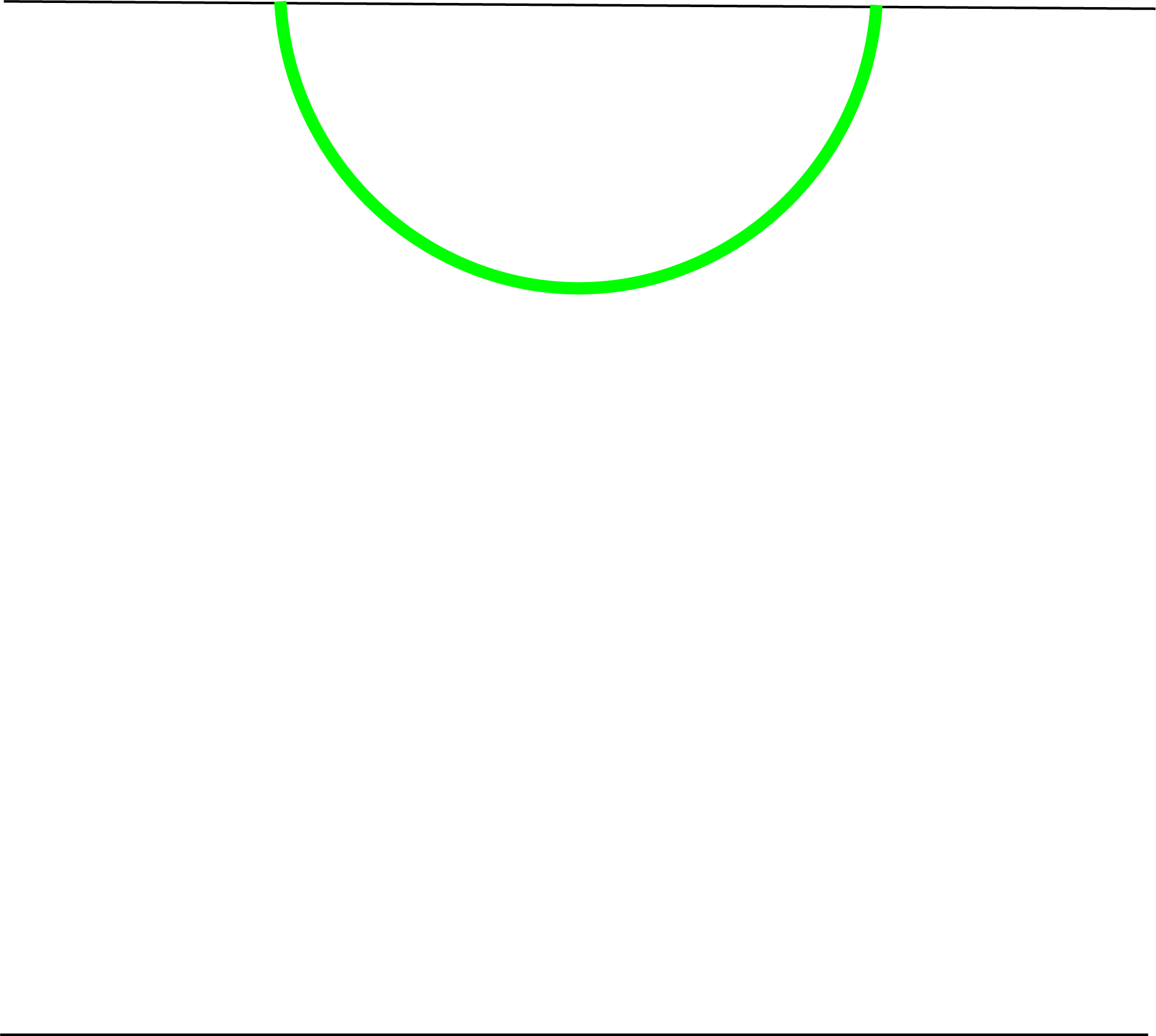}
\end{figure}
\noindent The relations are the following local relations on diagrams. 
\begin{figure}[H]
\centering
\includegraphics[width=.10\textwidth]{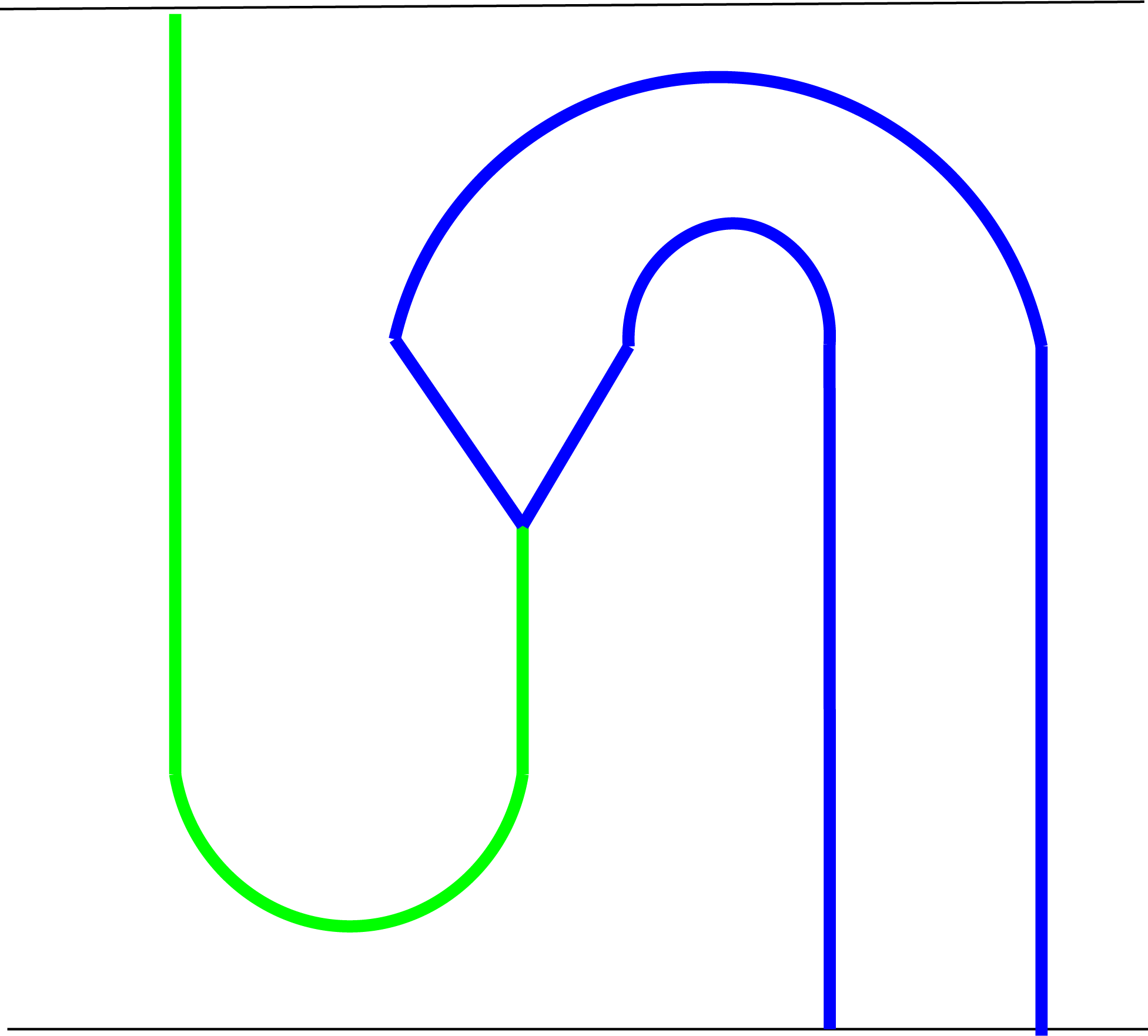}\put(8, 17){$=$} \ \ \ \ \ \ \ \ \ \includegraphics[width=.10\textwidth]{figs1/psst}\put(8, 17){$= $} \ \ \ \ \ \ \ \ \ 
\includegraphics[width=.10\textwidth]{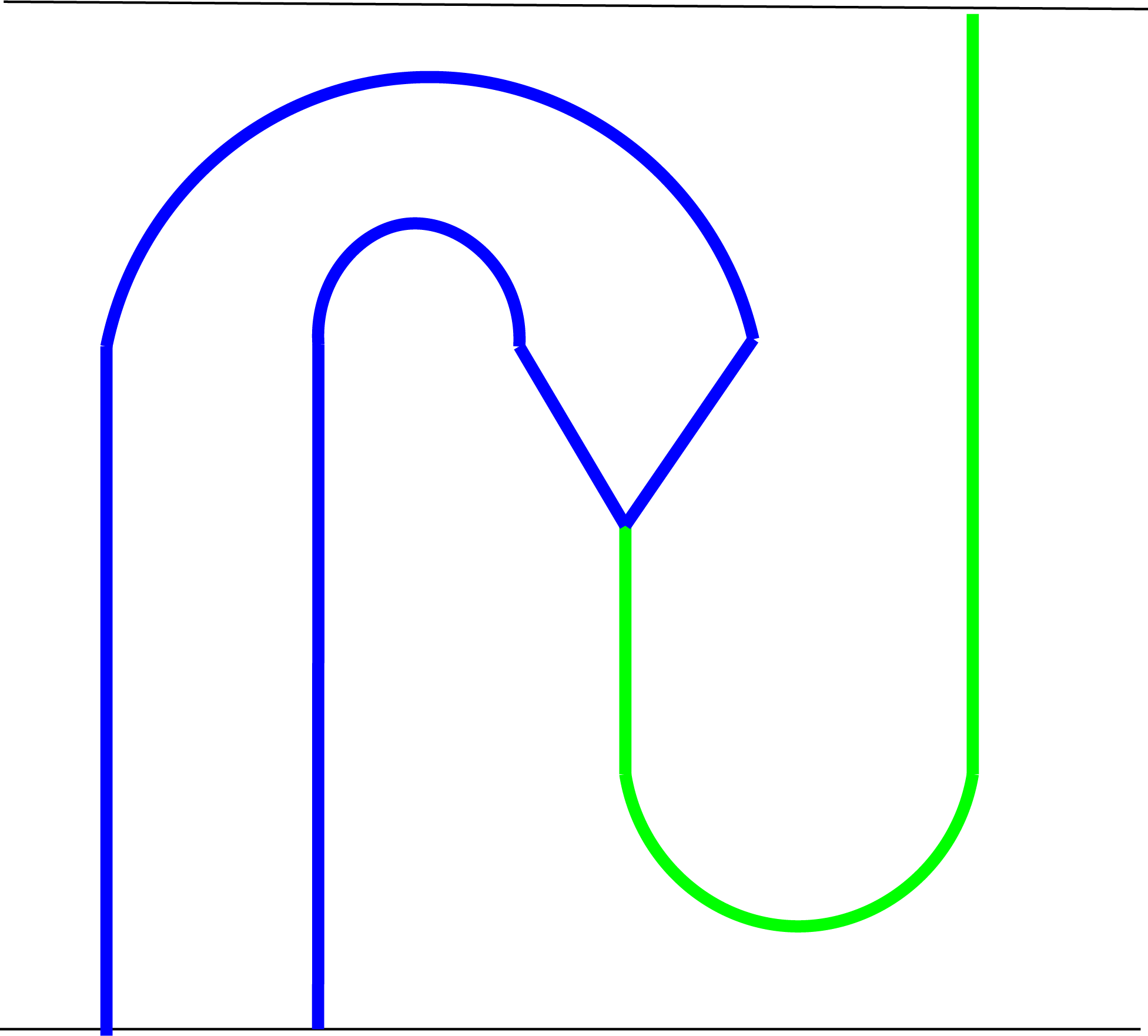}
\end{figure}
\begin{figure}[H]
\centering
\includegraphics[width=.10\textwidth]{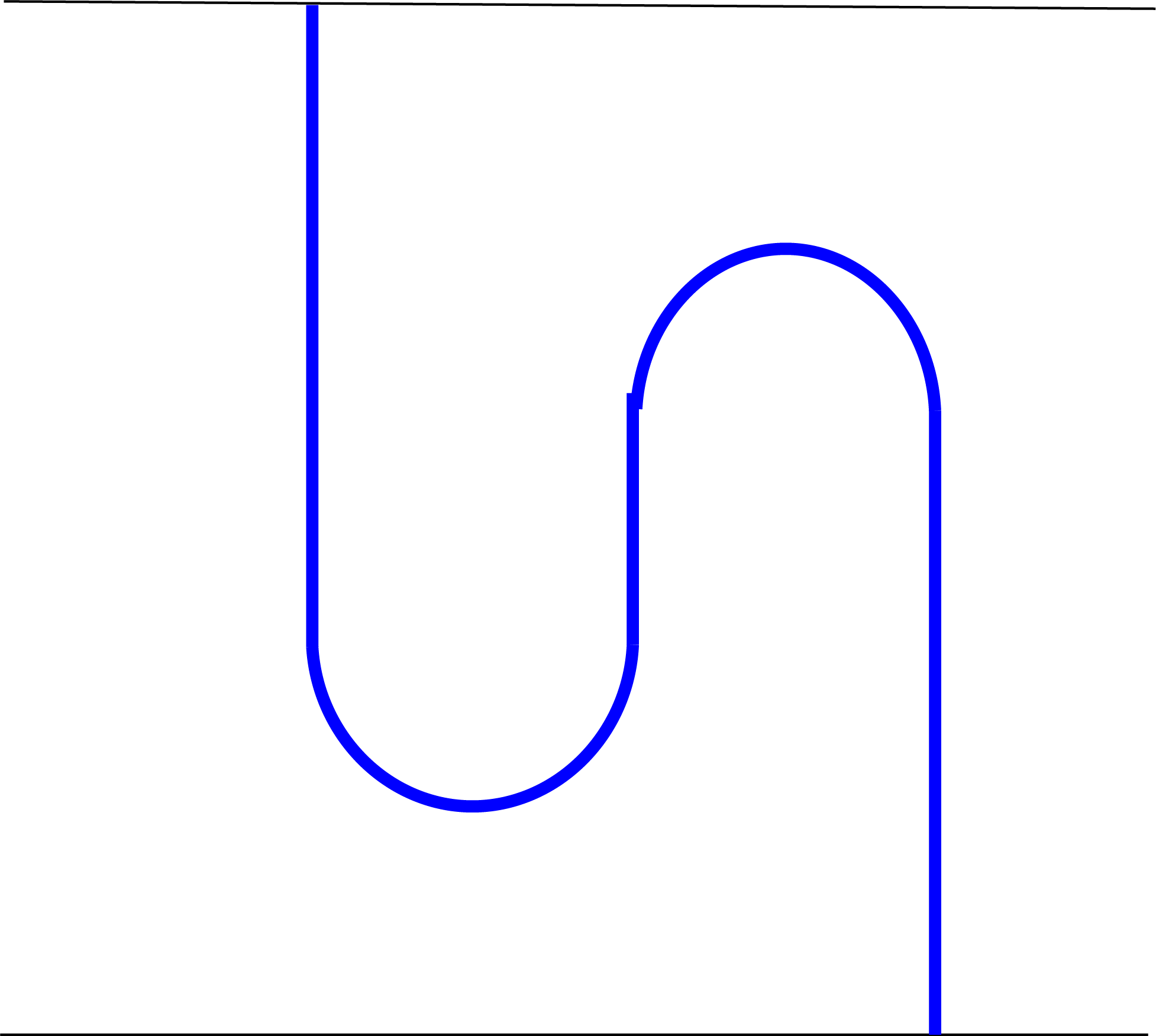}\put(8, 17){$=$} \ \ \ \ \ \ \ \ \ \includegraphics[width=.10\textwidth]{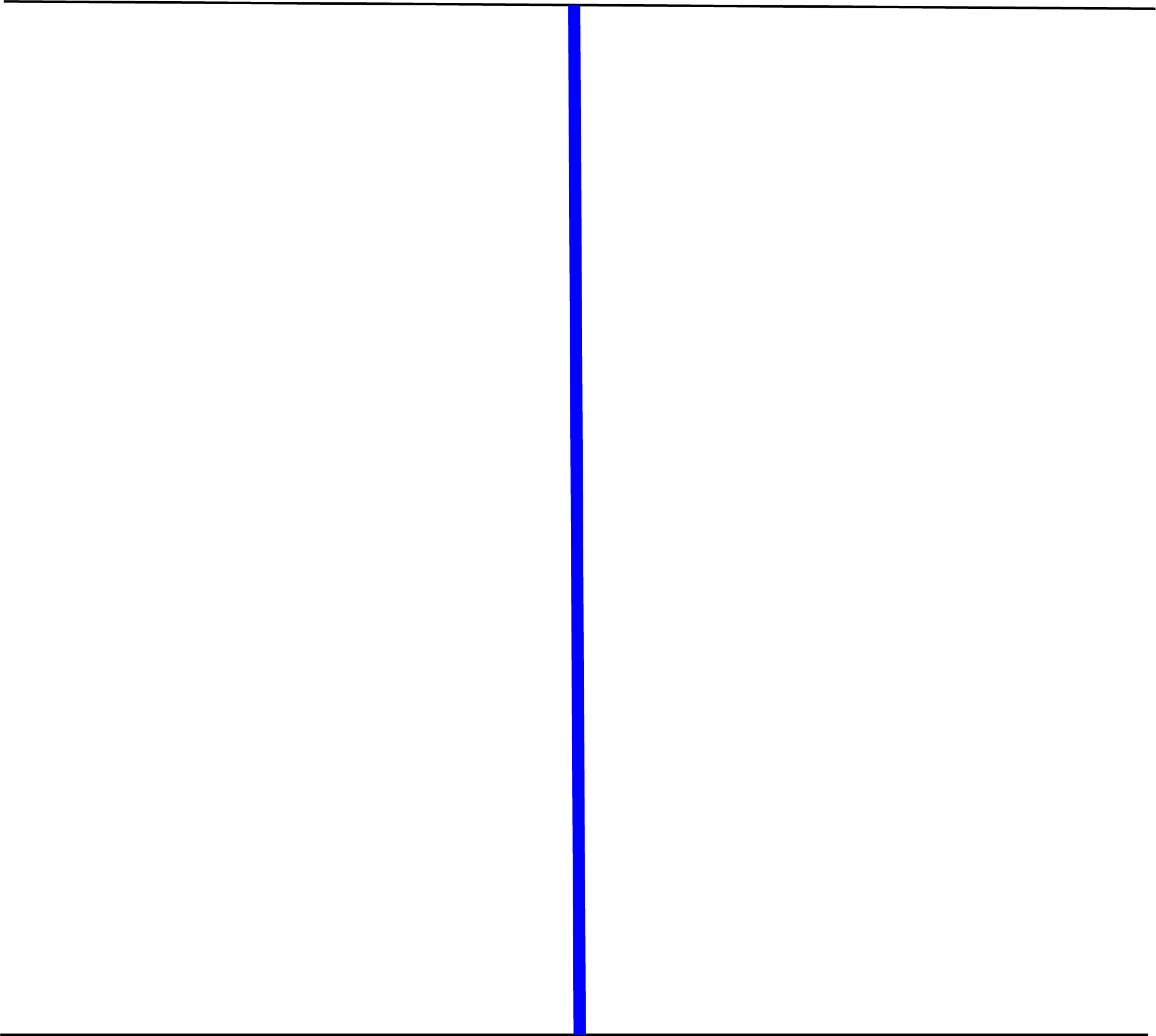}\put(8, 17){$= $} \ \ \ \ \ \ \ \ \ 
\includegraphics[width=.10\textwidth]{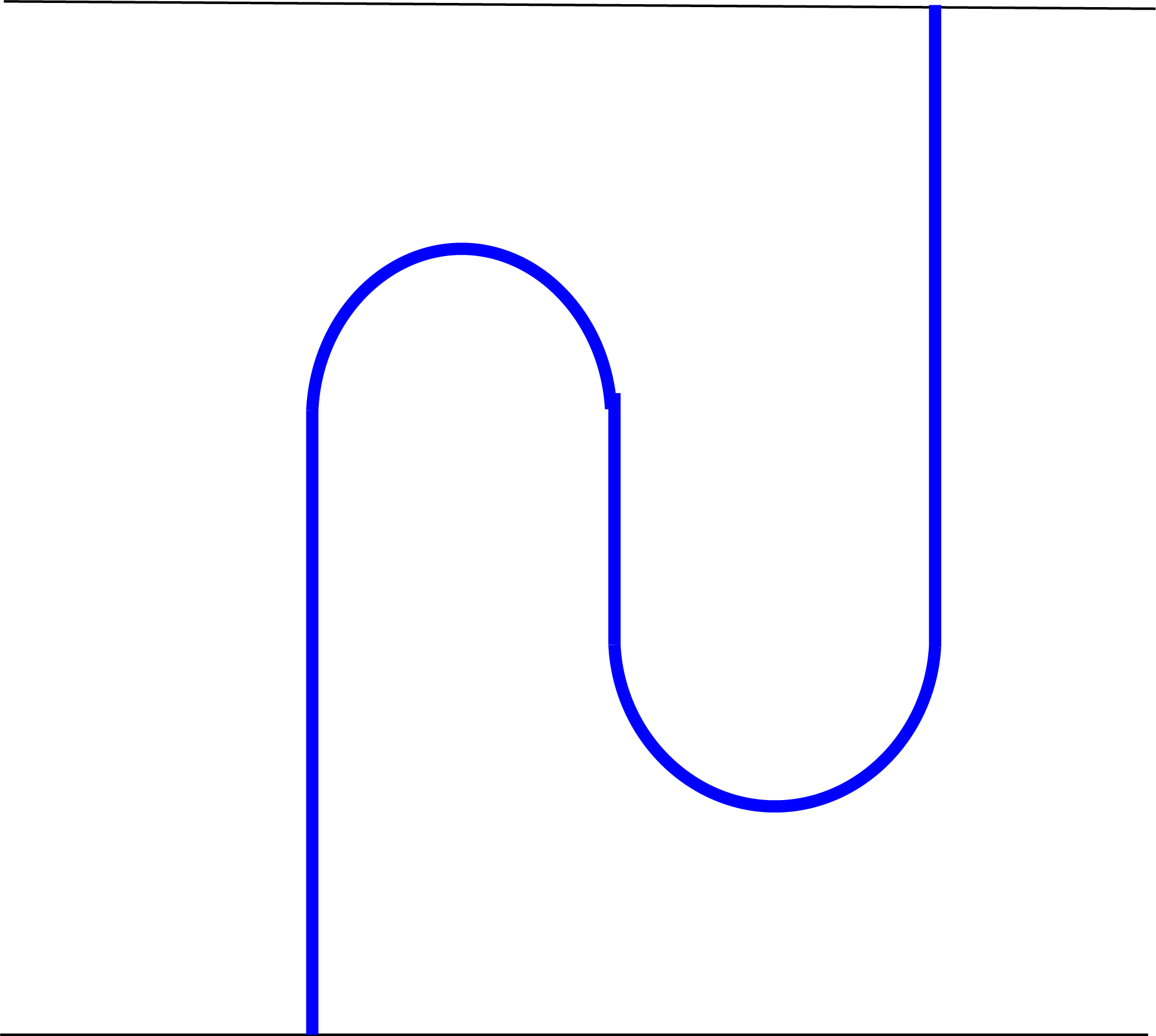}
\end{figure}
\begin{figure}[H]
\centering
\includegraphics[width=.10\textwidth]{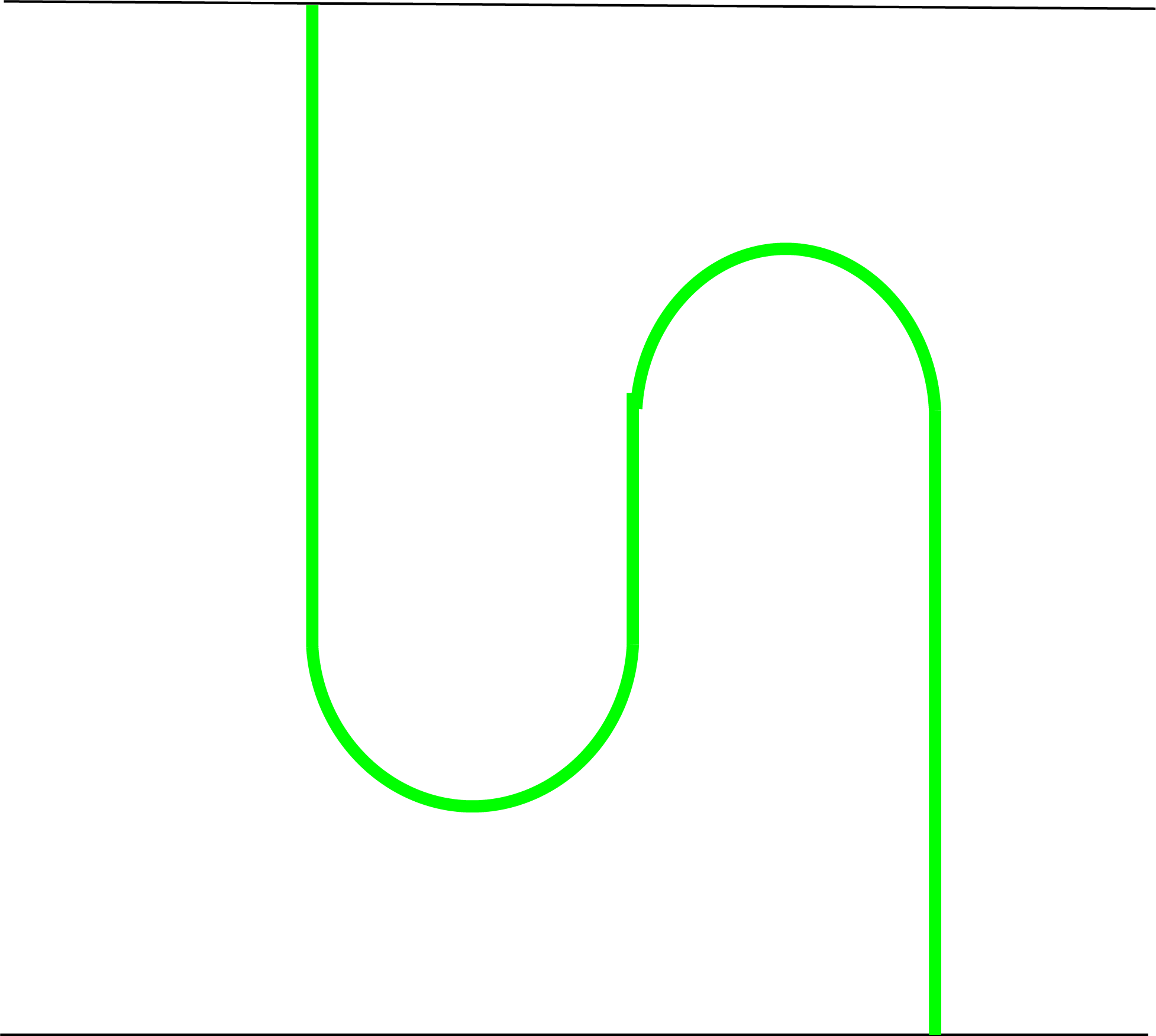}\put(8, 17){$=$} \ \ \ \ \ \ \ \ \ \includegraphics[width=.10\textwidth]{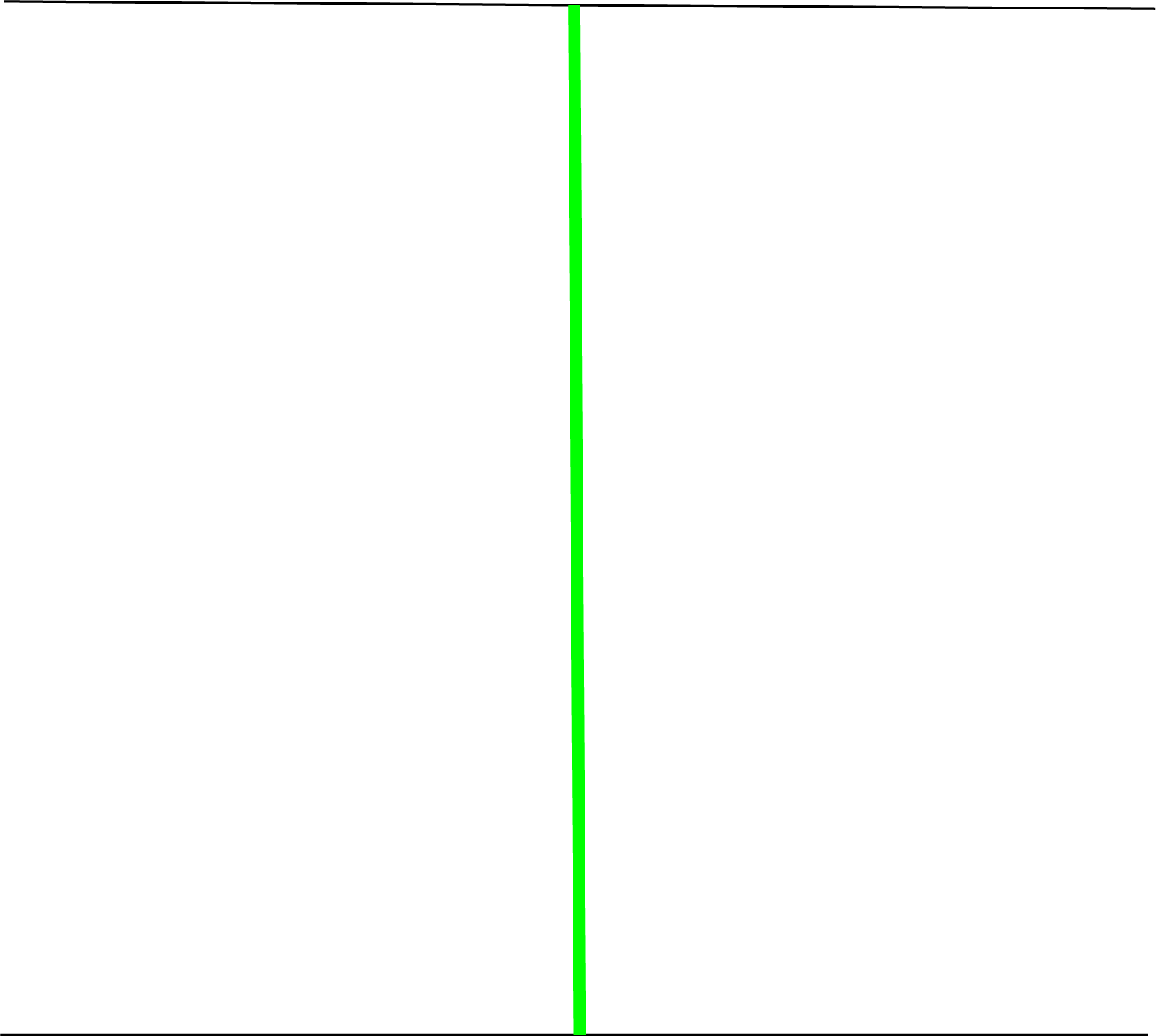}\put(8, 17){$= $} \ \ \ \ \ \ \ \ \ 
\includegraphics[width=.10\textwidth]{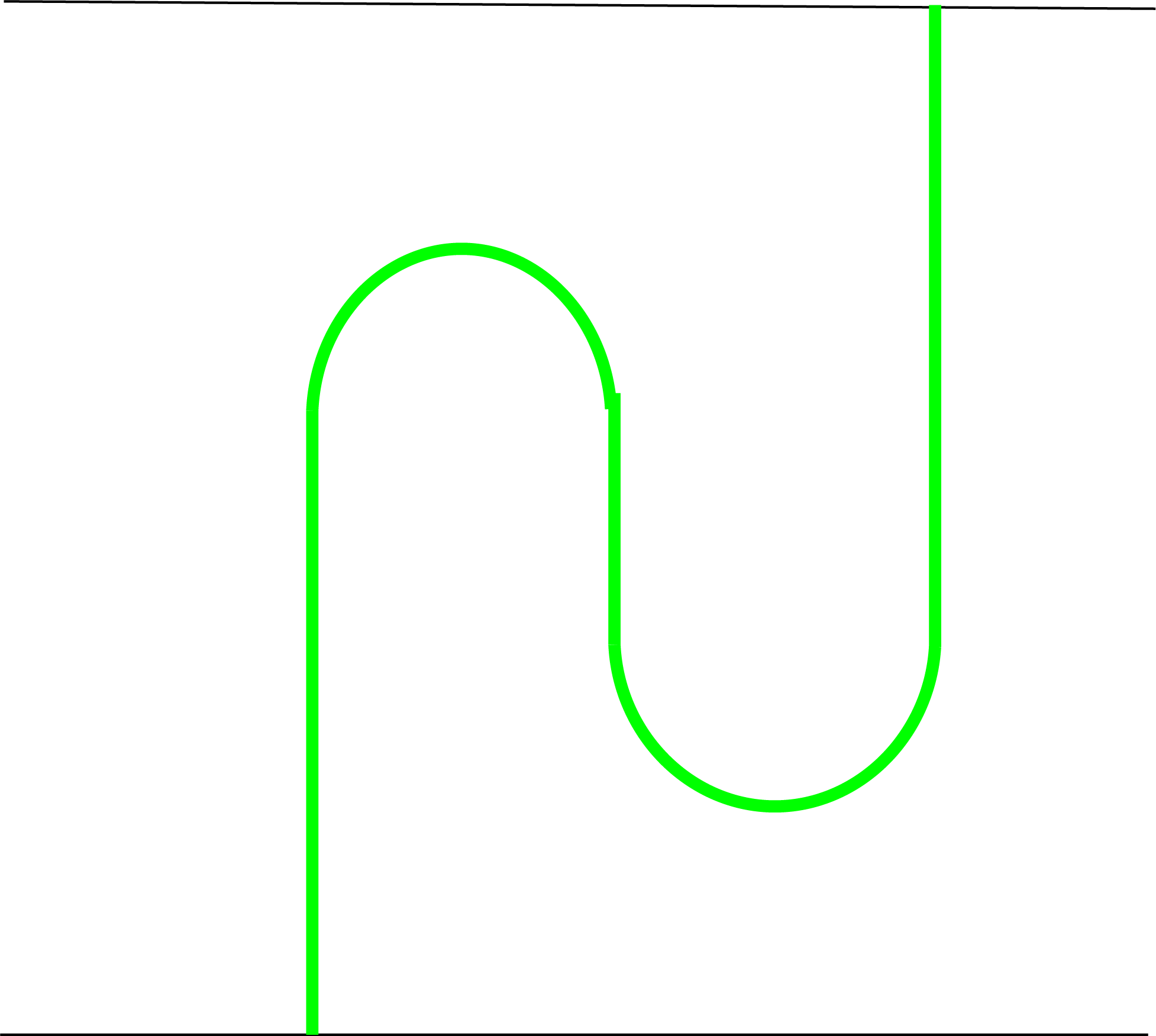}
\end{figure}
\end{defn}

\begin{remark}
Our convention is that diagrams are read as morphisms from the bottom boundary to the top boundary. Composition of morphisms is vertical stacking. The monoidal structure on objects is concatenation of words and the monoidal unit is the empty word. The monoidal product on morphisms is horizontal concatenation of diagrams, and the identity morphism of the empty word is the empty diagram.  
\end{remark}

\begin{figure}[H]
\caption{The identity morphism of $\blues\greent\blues\blues$ and a morphism from $\blues\greent\blues\blues\blues$ to $\blues\blues\greent\greent$.}
\centering
\includegraphics[width=0.1\textwidth]{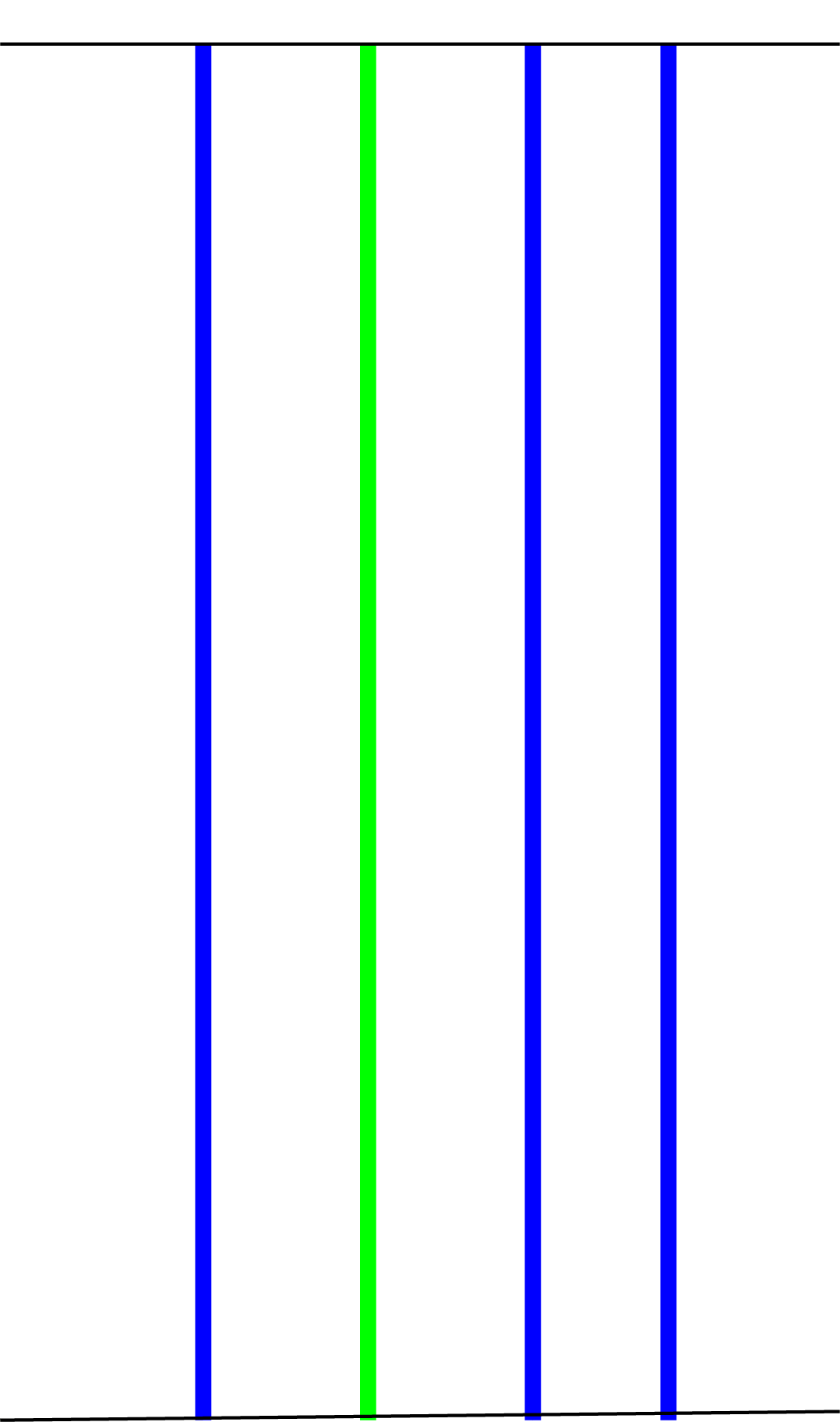} \ \ \ \ \ \ \ \ \ \ \ \ \ \ \ \ \ \ \ \ 
\includegraphics[width=0.11\textwidth]{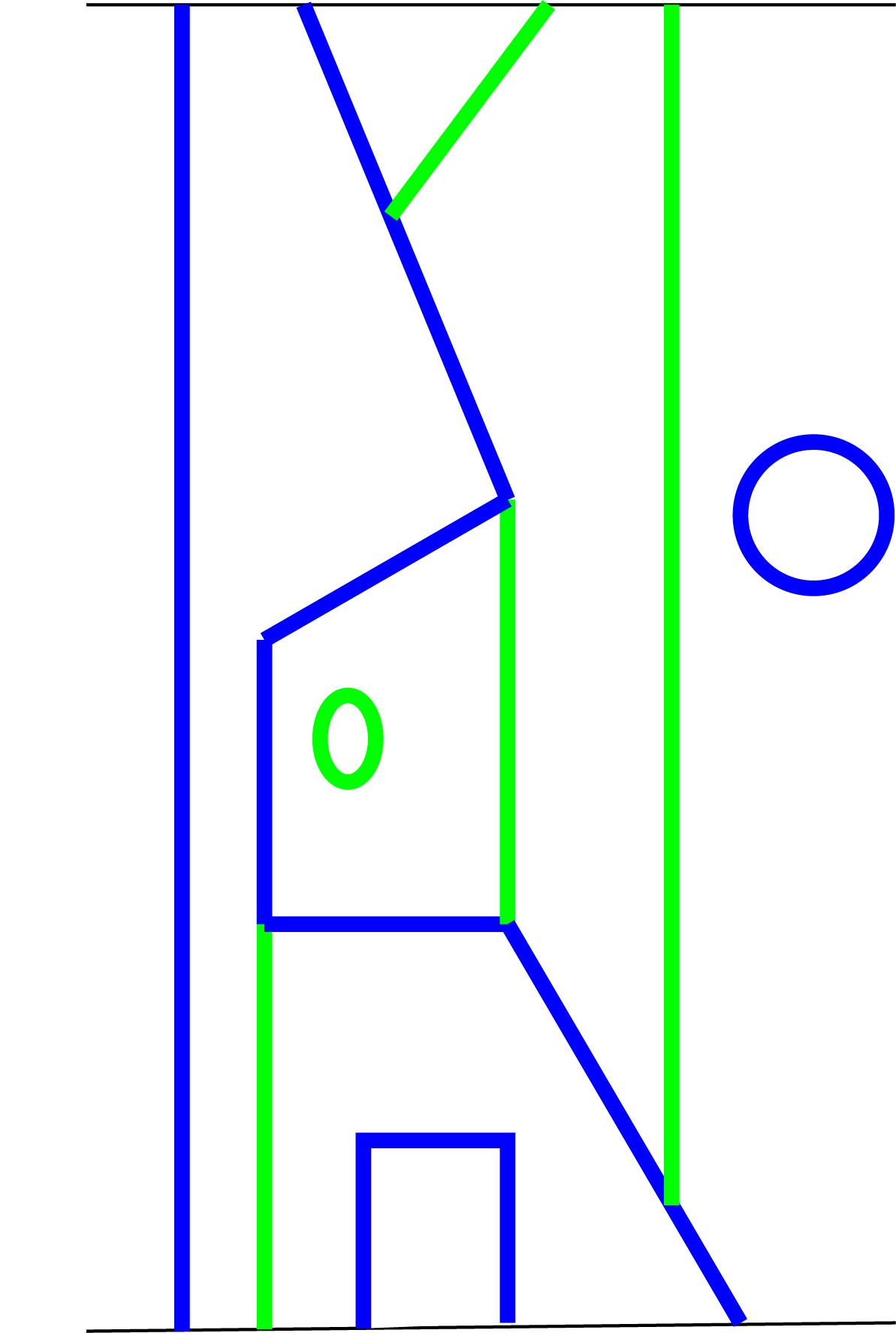}
\end{figure}

\begin{notation}\label{isotopyrmk}
The defining relations in $\mathcal{D}$ imply the following equalities of morphisms in $\Hom_{\mathcal{D}}(\blues\greent, \blues)$. 
\begin{figure}[H]
\centering
\includegraphics[width=.10\textwidth]{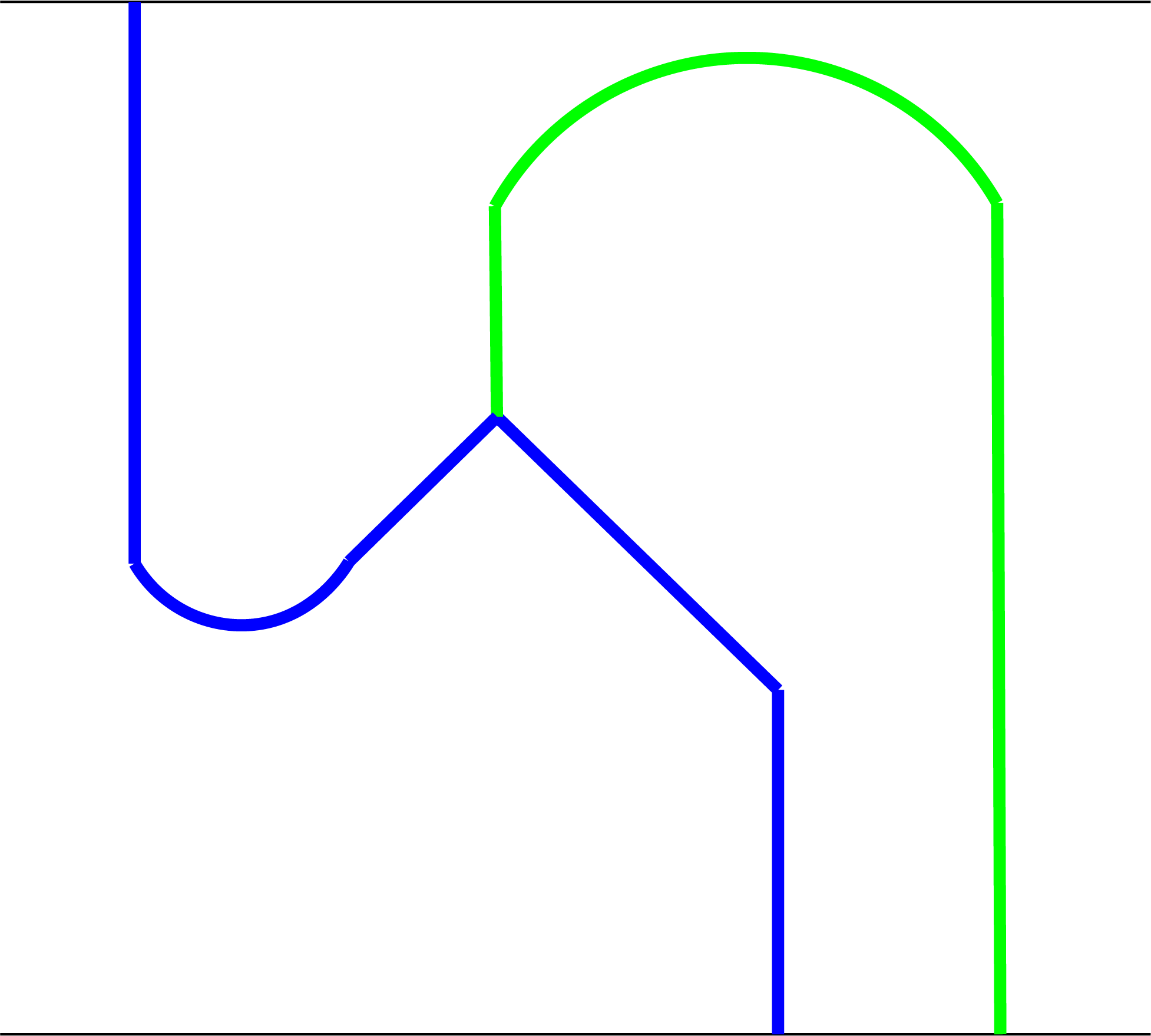}\put(8, 17){$=$} \ \ \ \ \ \ \ \ \ \includegraphics[width=.10\textwidth]{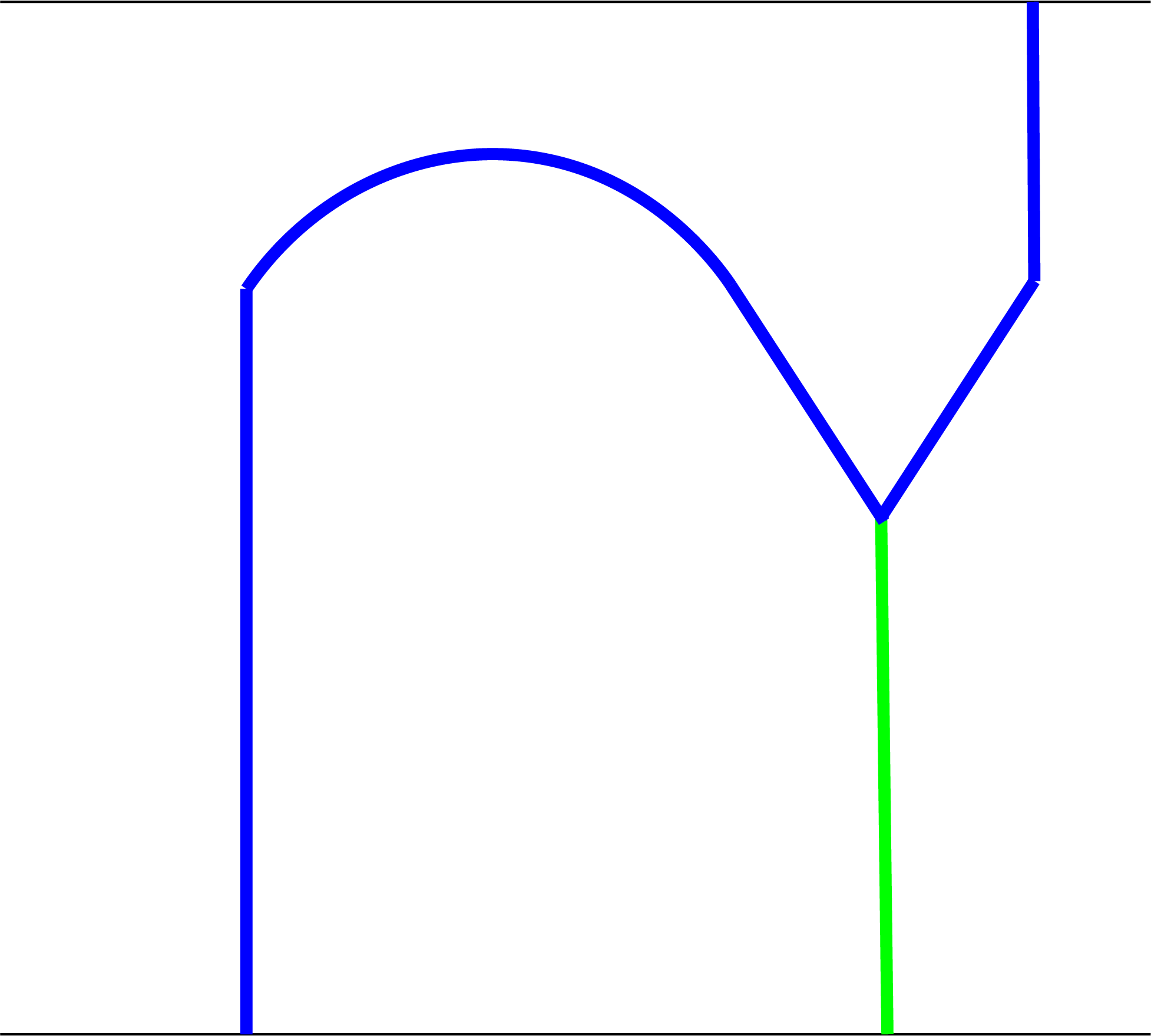}\put(8, 17){$= $} \ \ \ \ \ \ \ \ \ 
\includegraphics[width=.10\textwidth]{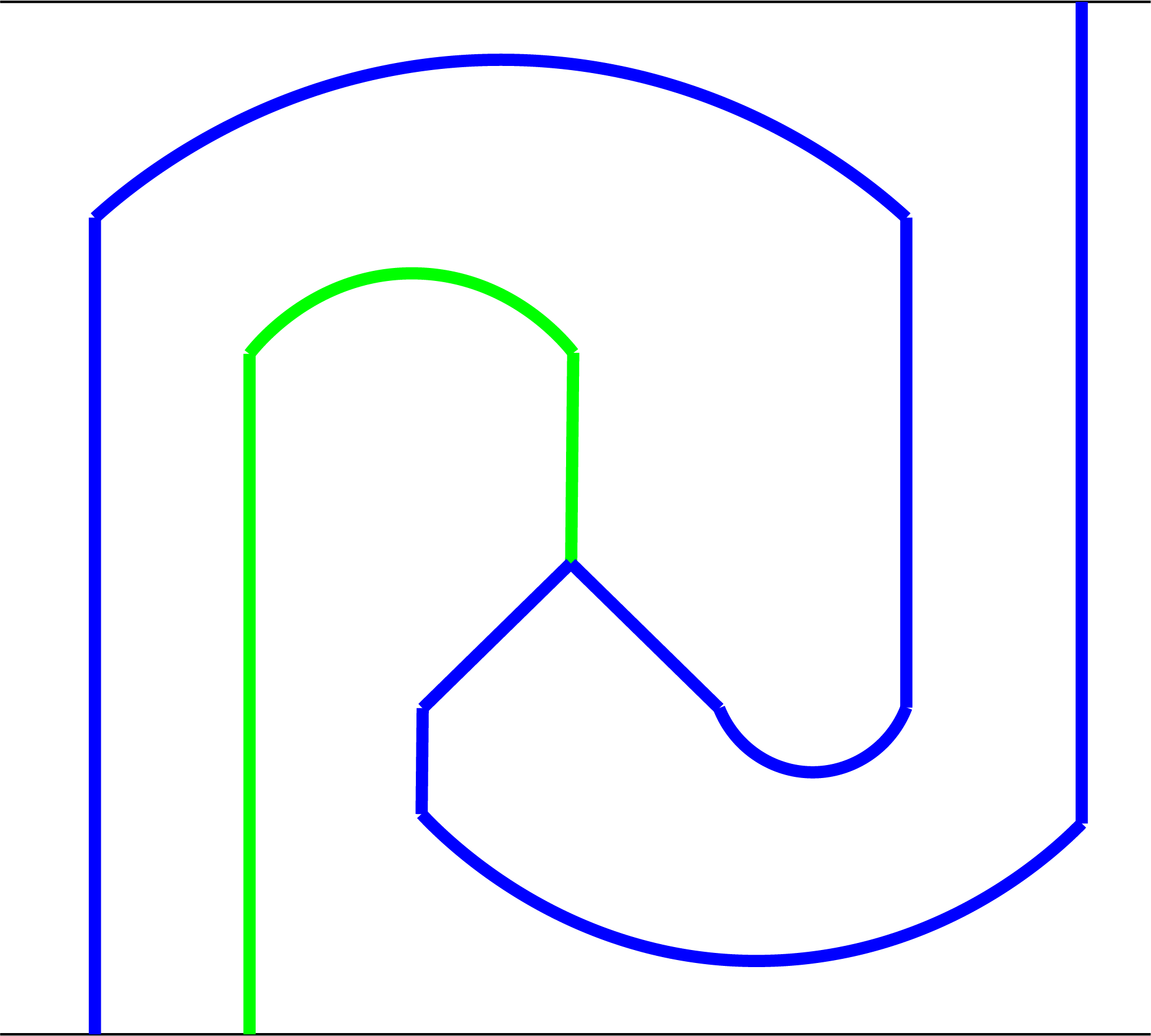}
\end{figure}
\noindent We will denote any one of these morphisms by the following trivalent vertex diagram in $\Hom_{\mathcal{D}}(\blues\greent, \blues)$.
\begin{figure}[H]
\centering
\includegraphics[width=.10\textwidth]{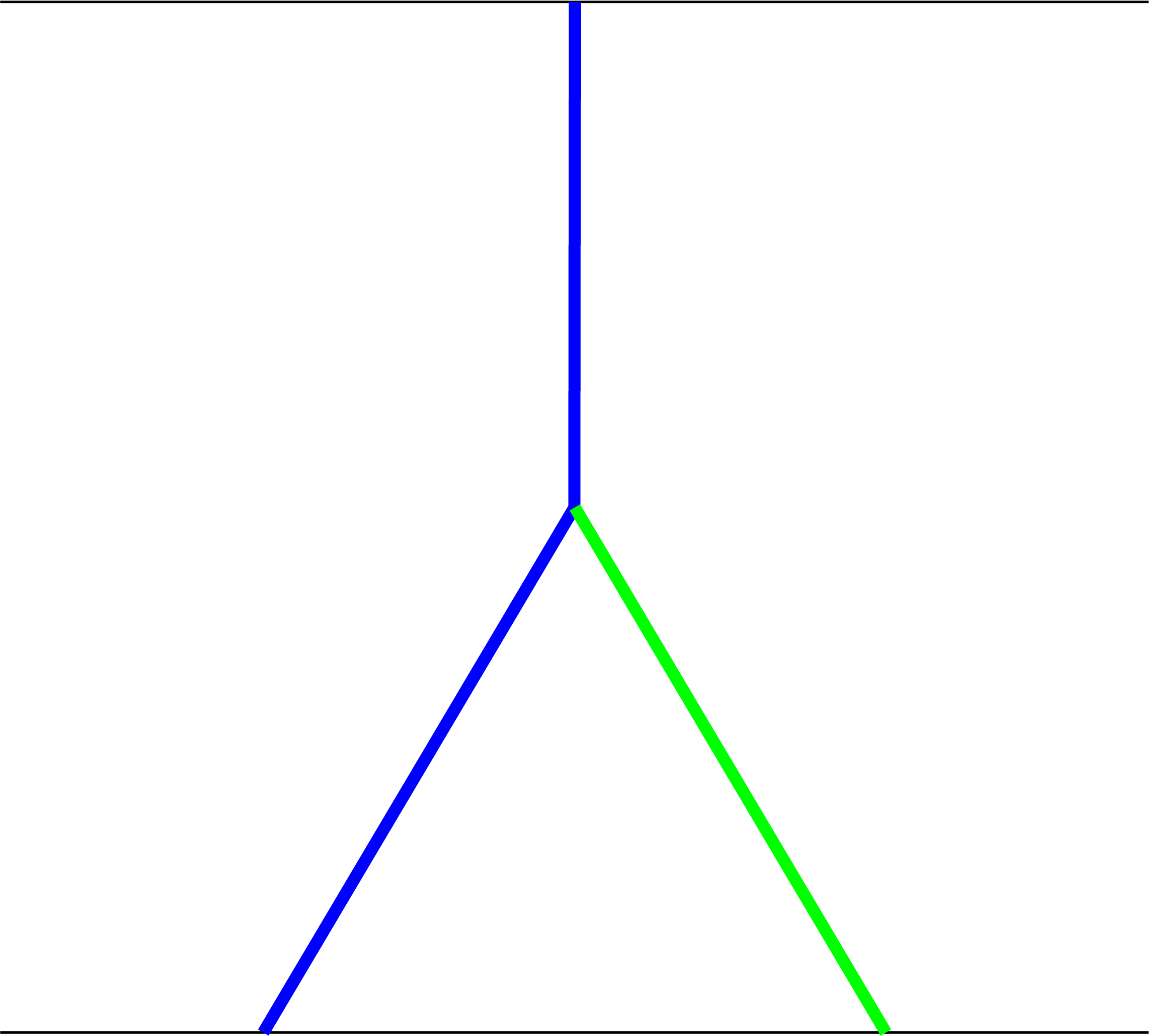}
\end{figure}

There are similar equalities for every possible vertical and horizontal reflection, and we will write the corresponding trivalent morphisms as follows. 
\begin{figure}[H]
\centering
\includegraphics[width=.10\textwidth]{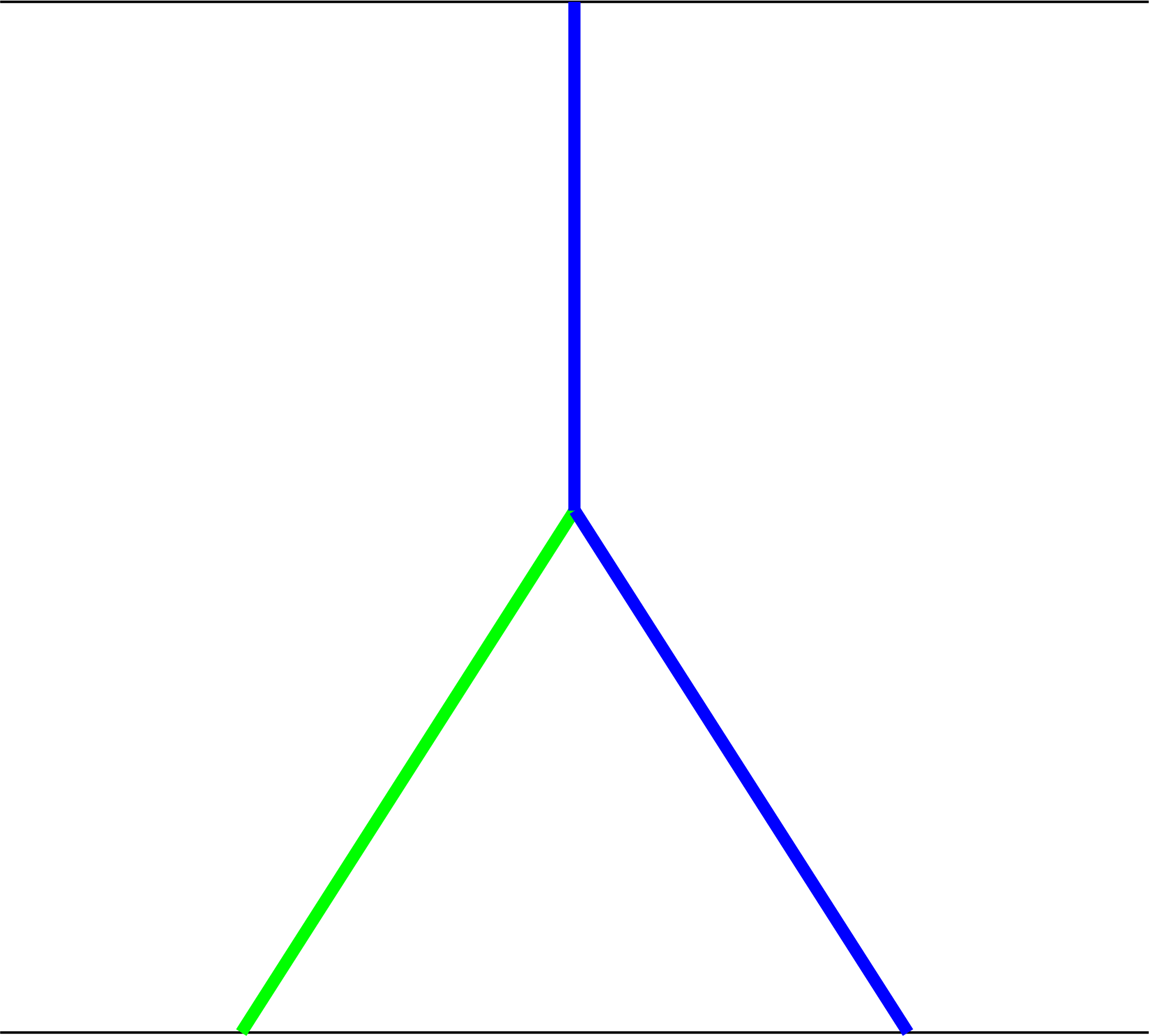}\put(8, 17){} \ \ \ \ \ \ \ \ \ \includegraphics[width=.10\textwidth]{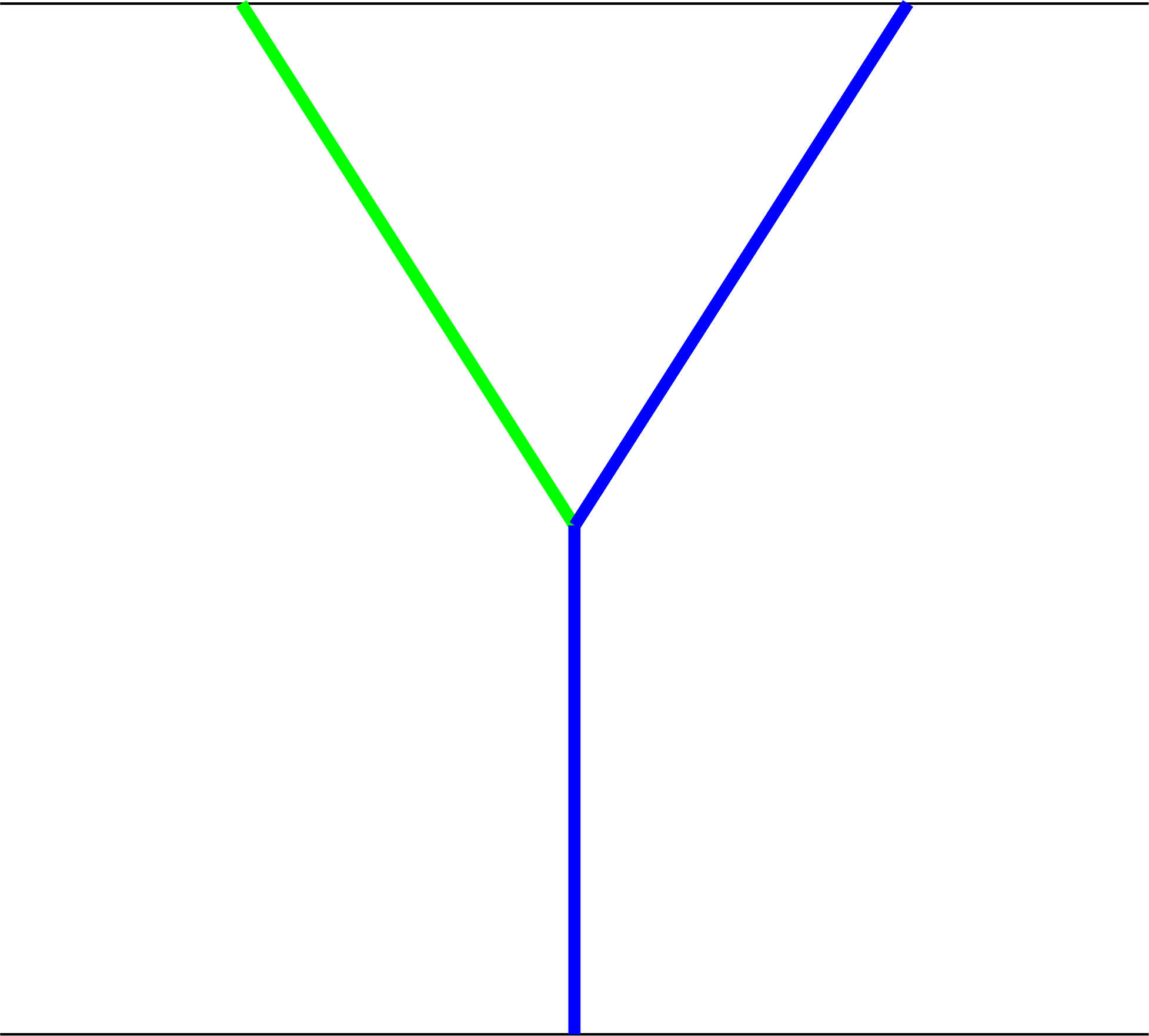}\put(8, 17){} \ \ \ \ \ \ \ \ \ 
\includegraphics[width=.10\textwidth]{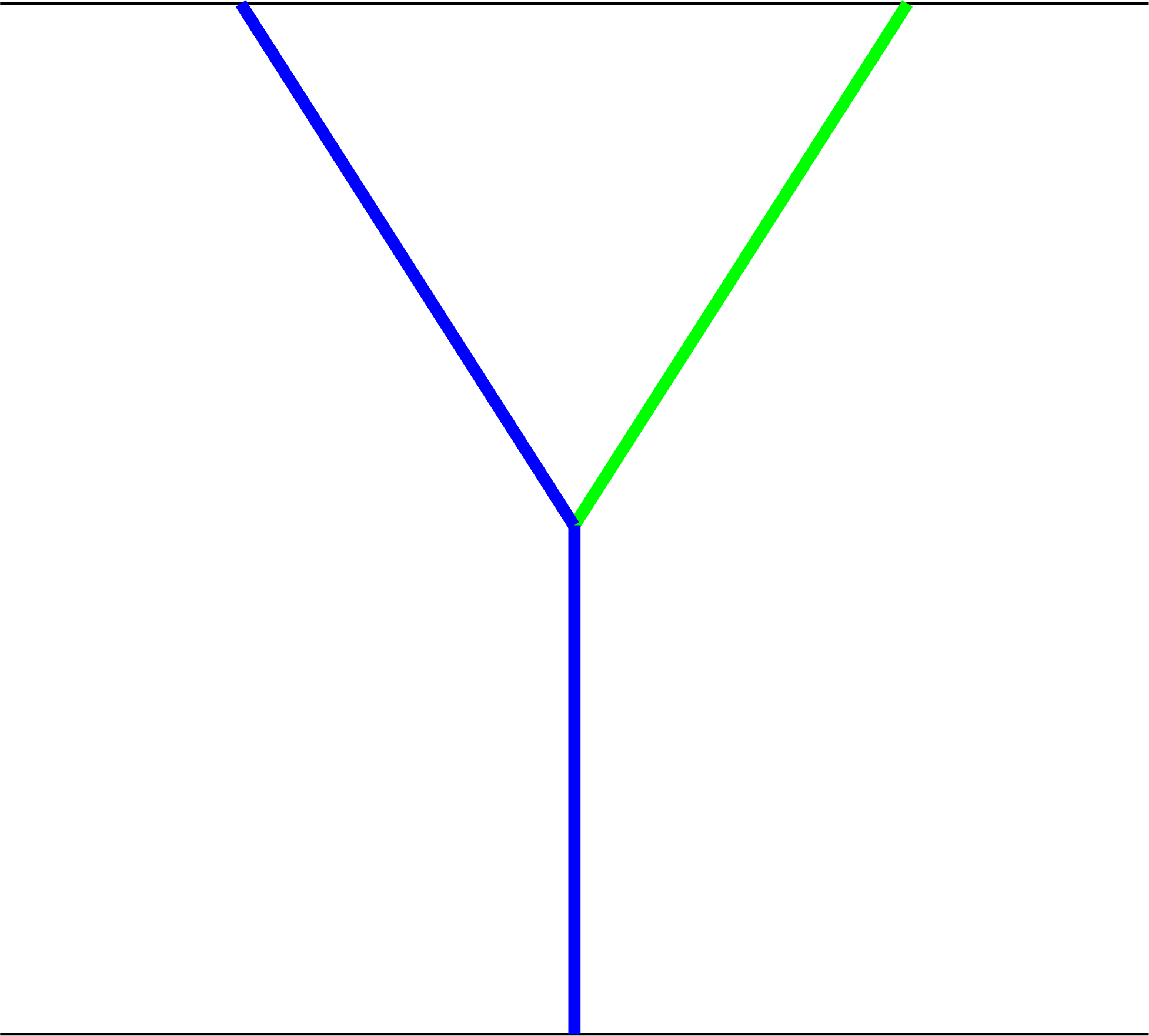}
\end{figure}
Thanks to this notation, we may now view morphisms in $\mathcal{D}$ as $\mathcal{A}$-linear combinations of isotopy classes trivalent graphs.
\end{notation}

\begin{defn}
The $\mathcal{A}$-linear monoidal category $\DD$ is the quotient of $\mathcal{D}$ by the following local relations.
\begin{figure}[H]
\centering
\includegraphics[width=0.12\textwidth]{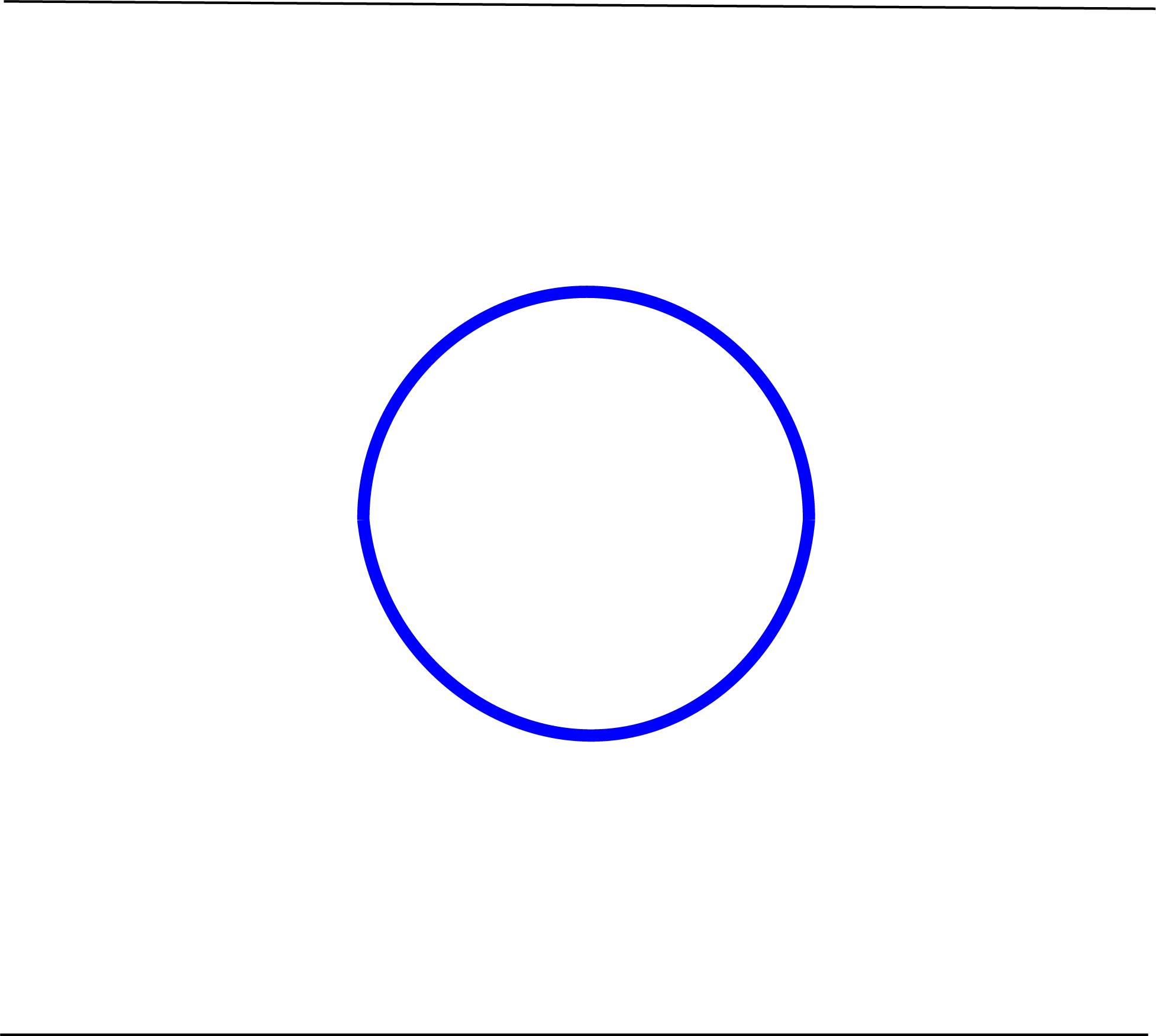}\put(5, 18){$= -\dfrac{[6]_q[2]_q}{[3]_q}$} \ \ \ \ \ \ \ \ \ \ \ \ \ \ \ \ \ \ \ \ \ \ \ \ \includegraphics[width = .12\textwidth]{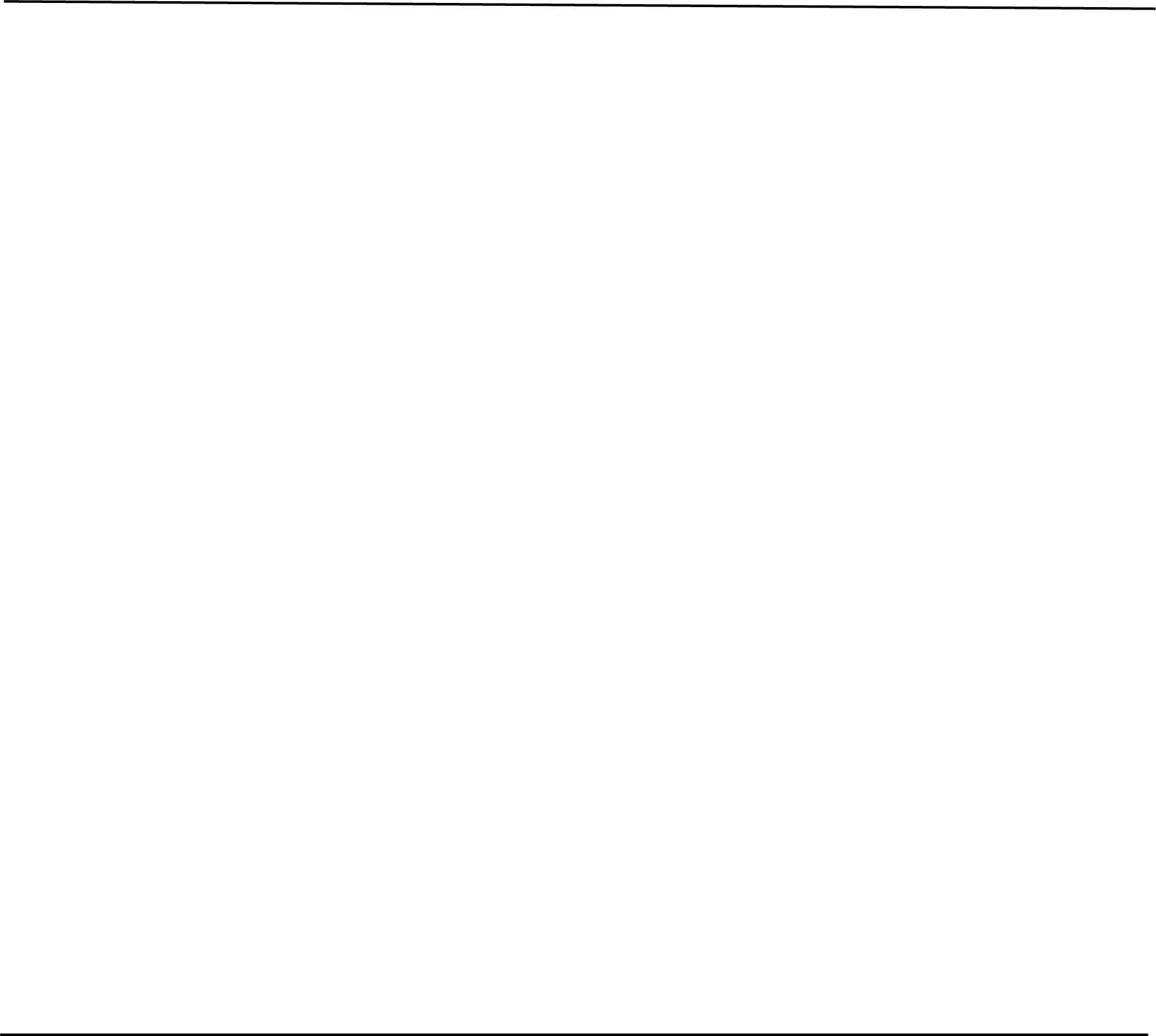}
\end{figure}
\begin{figure}[H]
\centering
\includegraphics[width=0.12\textwidth]{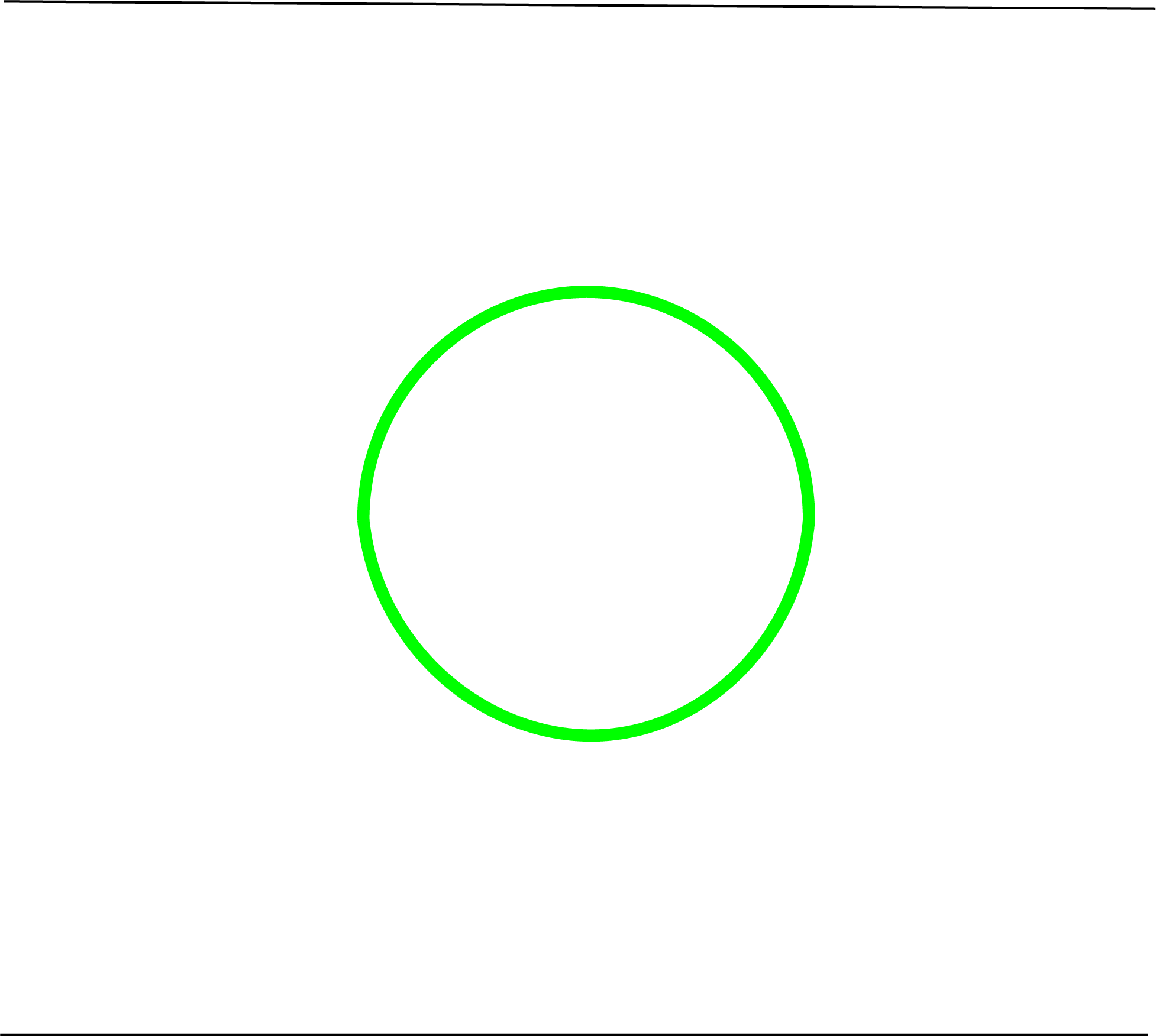}\put(5, 18){$= \dfrac{[6]_q[5]_q}{[3]_q[2]_q}$} \ \ \ \ \ \ \ \ \ \ \ \ \ \ \ \ \ \ \ \ \ \includegraphics[width = .12\textwidth]{figs1/emptyrelation}
\end{figure}
\begin{figure}[H]
\centering
\includegraphics[width=0.12\textwidth]{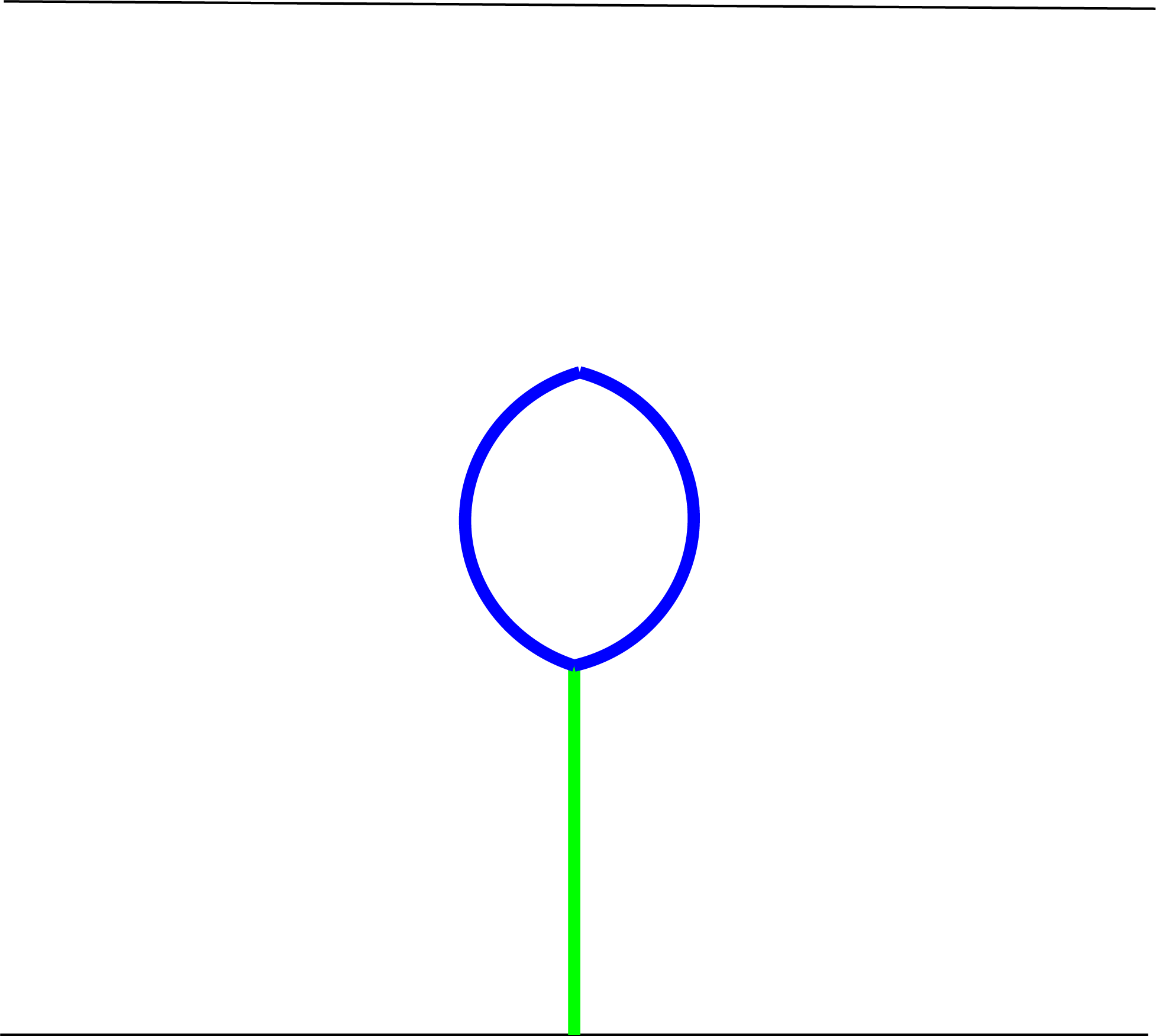}\put(5, 18){$=0$}
\end{figure}
\begin{figure}[H]
\centering
\includegraphics[width=0.12\textwidth]{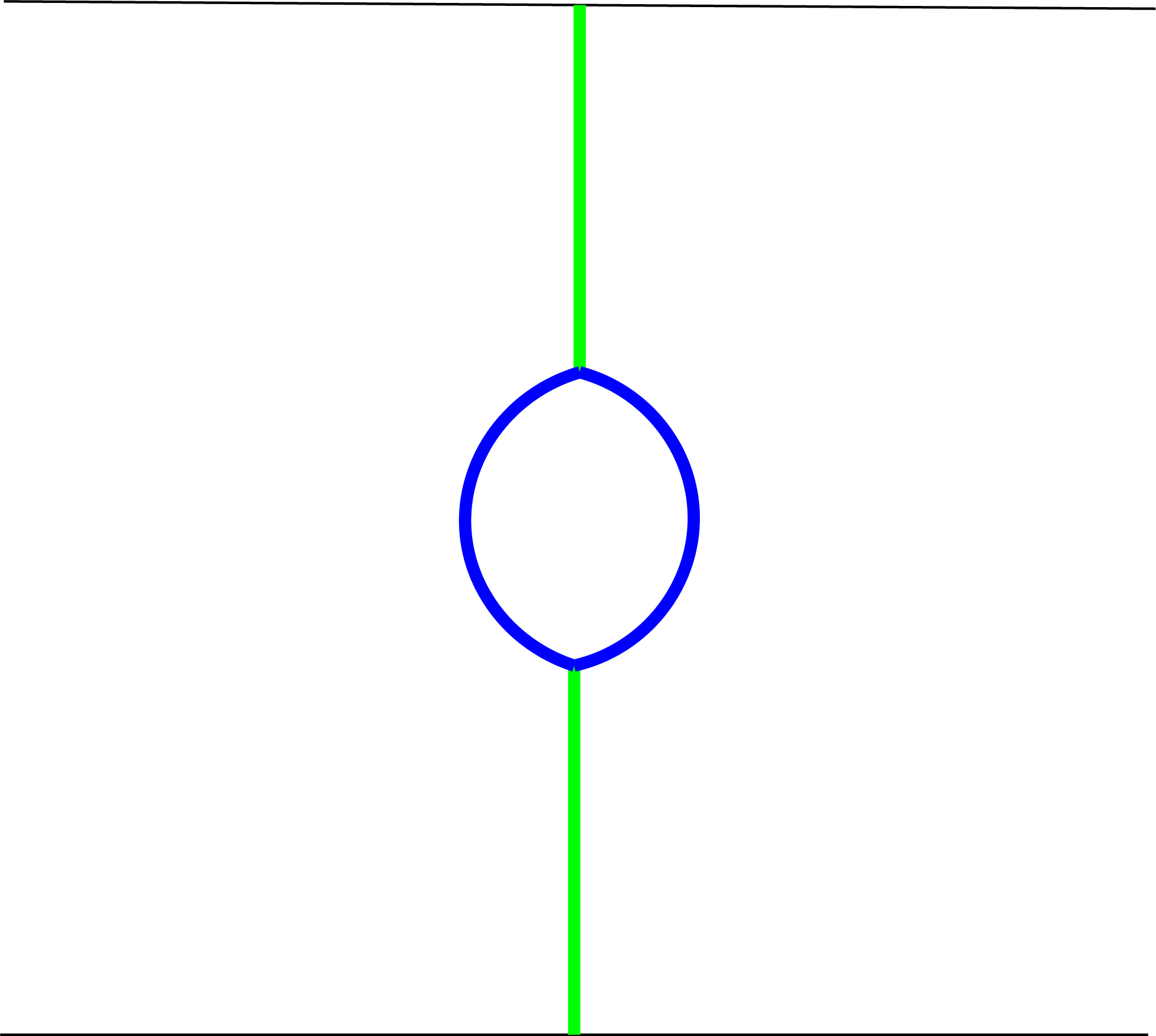}\put(5, 18){$=-[2]_q$} \ \ \ \ \ \ \ \ \ \ \ \ \ \ \ \ \ \includegraphics[width = .12\textwidth]{figs1/isotopyidt}
\end{figure}
\begin{figure}[H]
\centering
\includegraphics[width=0.12\textwidth]{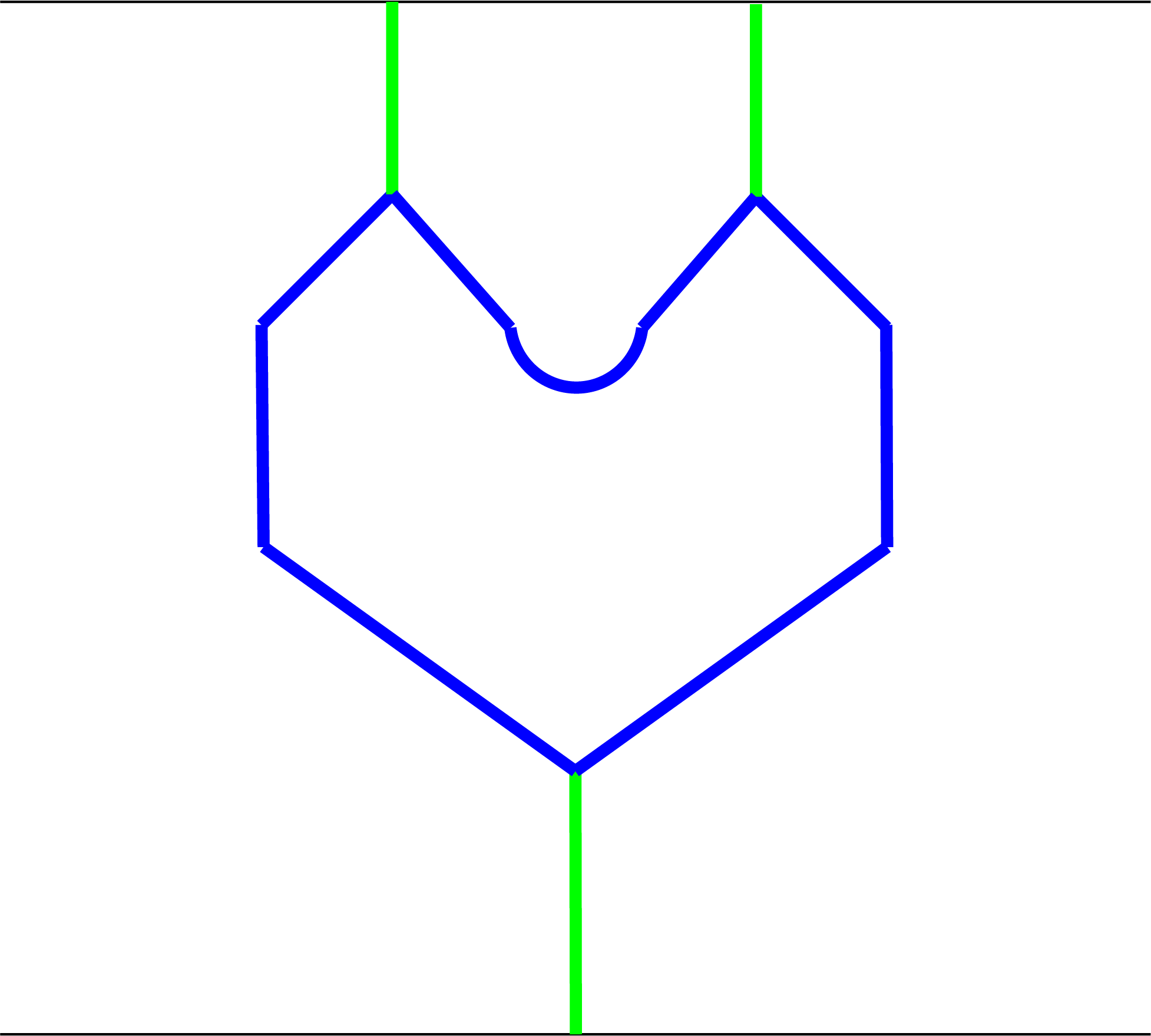}\put(5, 18){$=0$}
\end{figure}
\begin{figure}[H]
\centering
\includegraphics[width=0.12\textwidth]{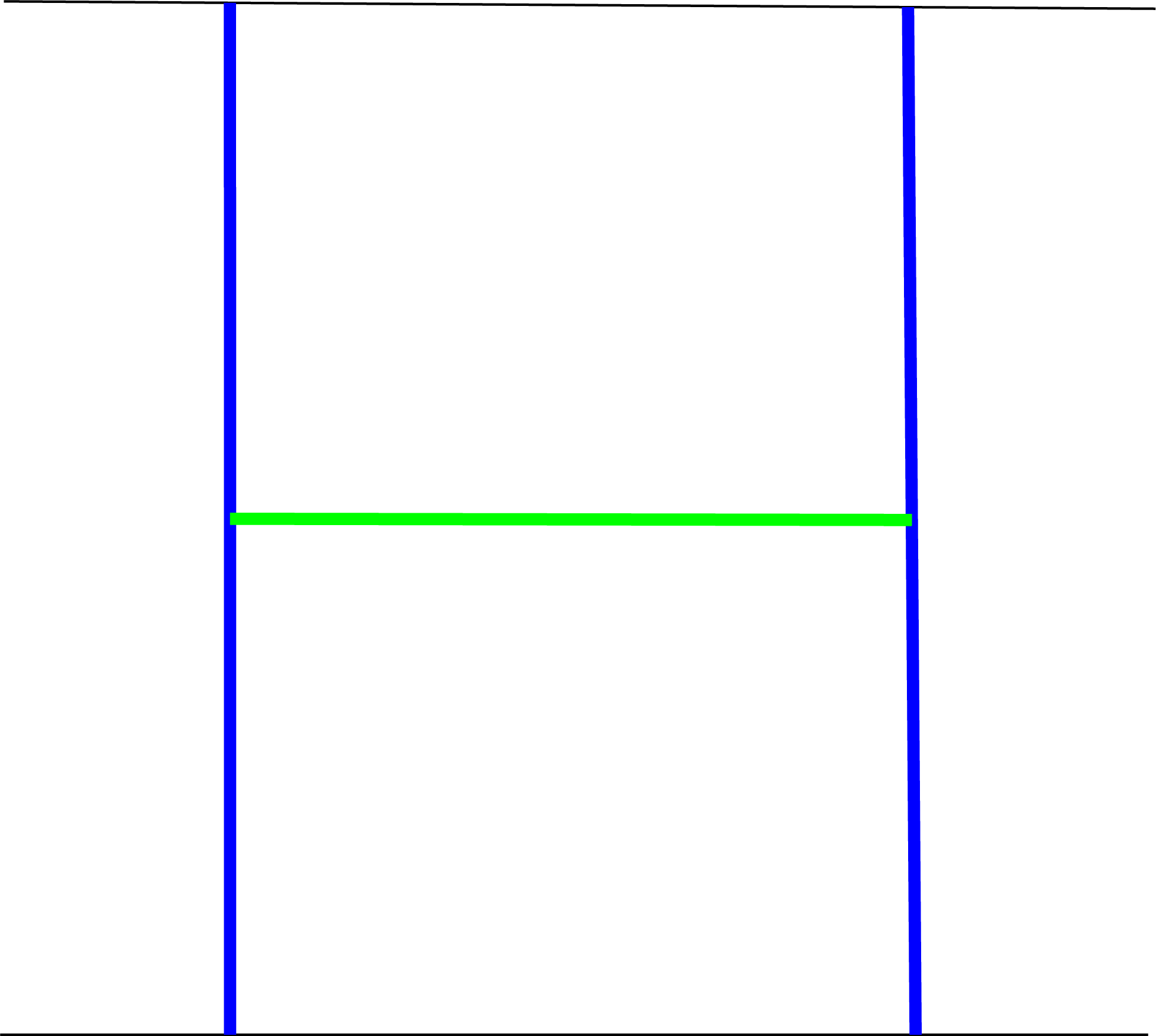}\put(5, 18){$=\dfrac{1}{[2]_q}$} \ \ \ \ \ \ \ \ \ \ \ \ \includegraphics[width=0.12\textwidth]{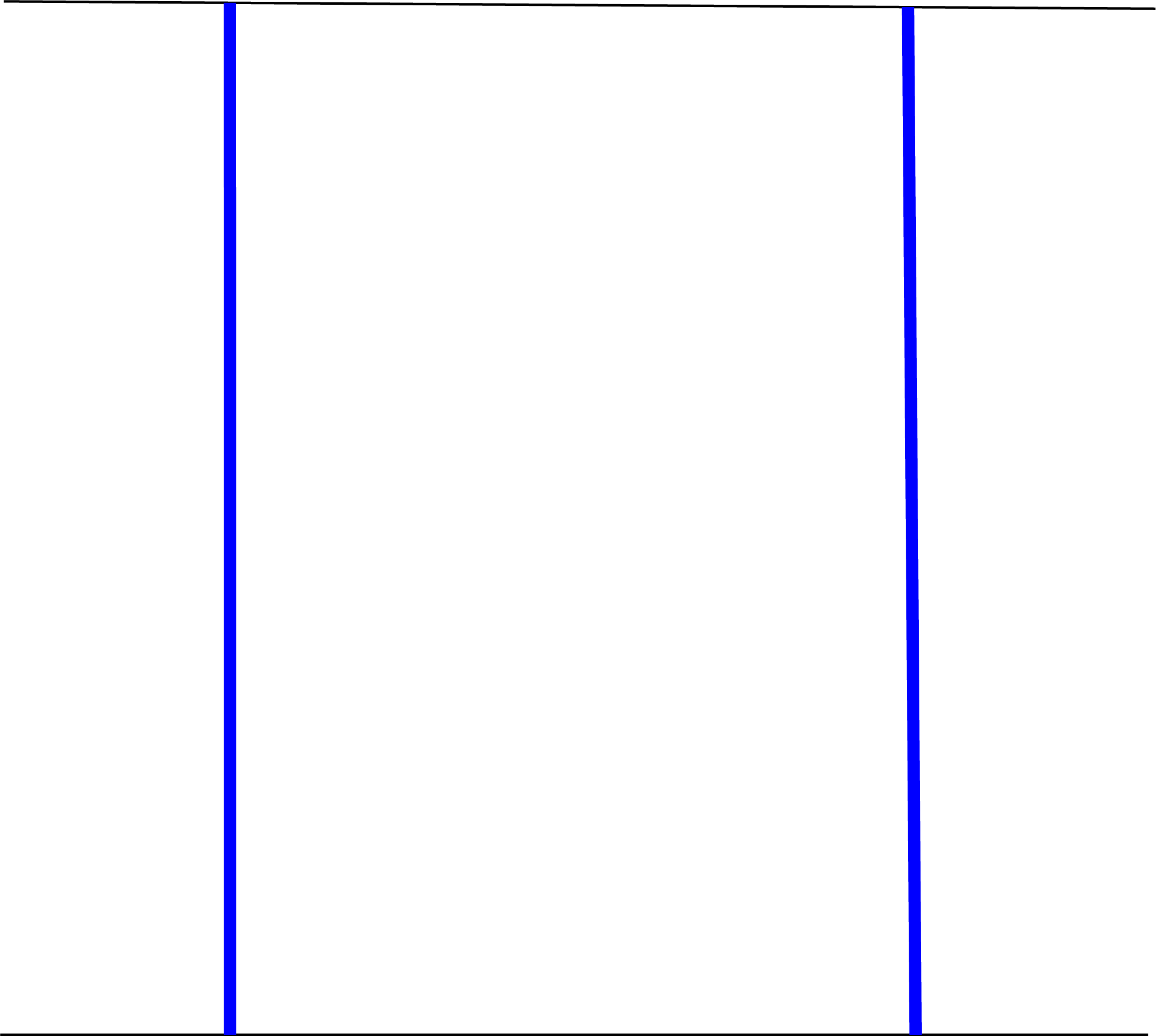}\put(5, 18){$+$}\ \ \ \ \ \ \ \ \ \ \ \ \includegraphics[width=0.12\textwidth]{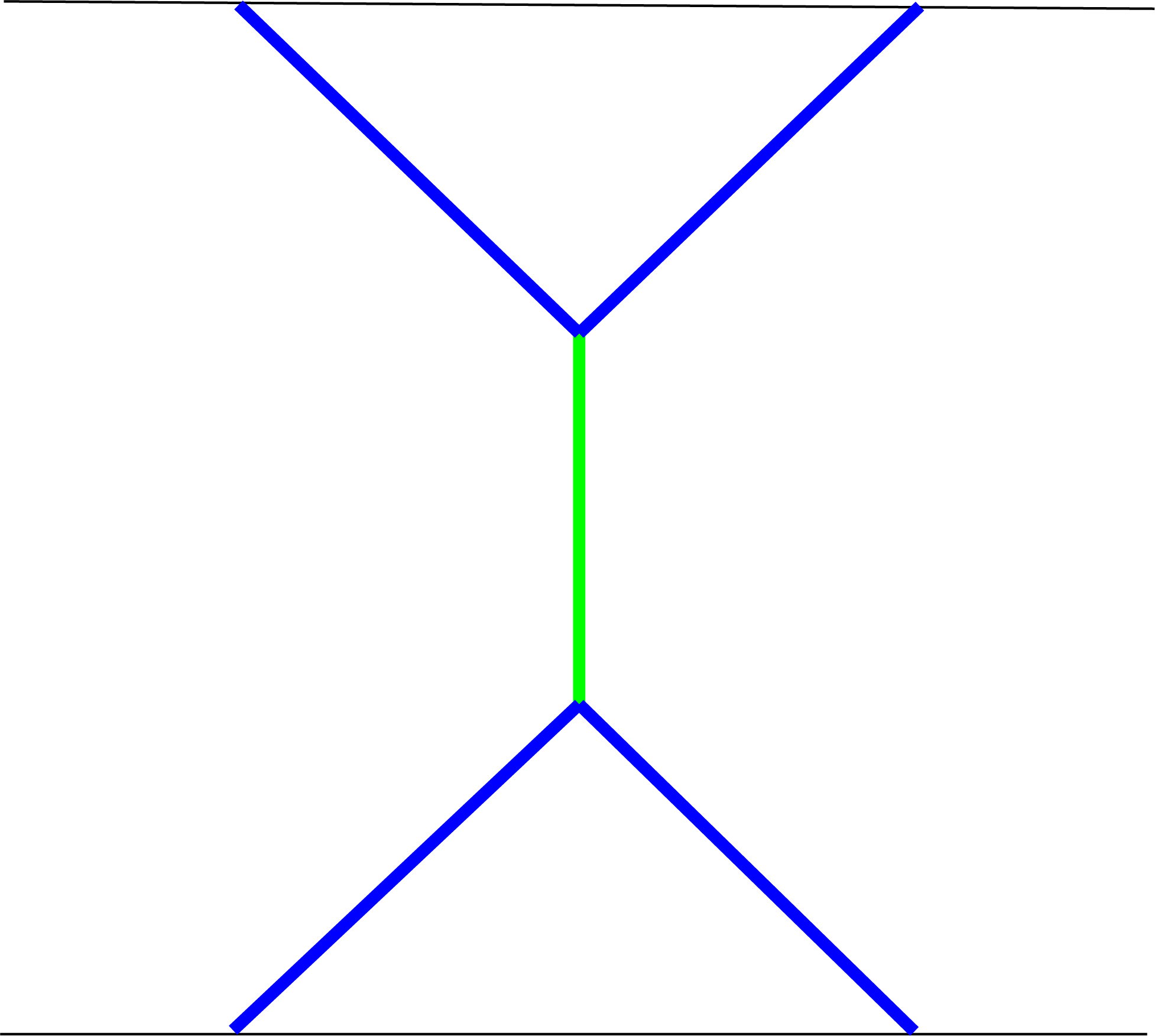}\put(5, 18){$-\dfrac{1}{[2]_q}$} \ \ \ \ \ \ \ \ \ \ \ \ \includegraphics[width=0.12\textwidth]{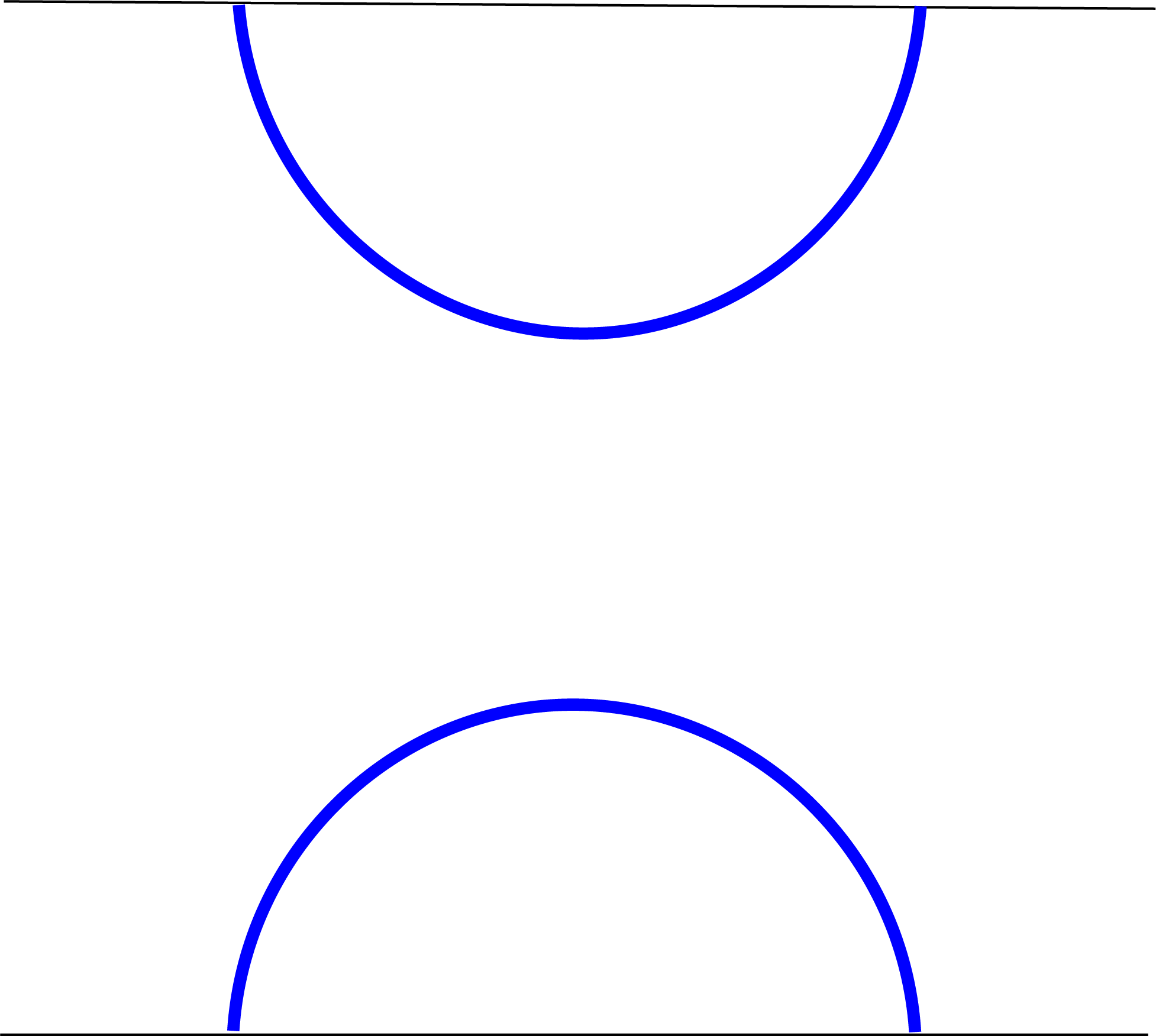}
\end{figure}
\end{defn}

\begin{notation}
When $\ak$ is an $\mathcal{A}$-algebra, we can base change the category $\DD$ to $\ak$, denoted $\ak\ot\DD$. The category $\ak\ot\DD$ has the same objects as $\DD$ and we apply $\ak \ot_{\mathcal{A}}(-)$ to homomorphism spaces. We may also write $\DDk:= \ak \ot \DD$ for short. 
\end{notation}

\begin{remark}
The coefficients in the circle relations are written as fractions but are actually elements of $\mathcal{A}$, as can be observed in the following quantum number calculations.
\begin{equation}\label{qdimcalcone}
[5]_q-[1]_q = \dfrac{\left([5]_q-[1]_q\right)[3]_q}{[3]_q} = \dfrac{[7]_q+ [5]_q+[3]_q- [3]_q}{[3]_q} = \dfrac{[6]_q[2]_q}{[3]_q}.
\end{equation} 
\begin{equation}\label{qdimcalctwo}
[7]_q-[5]_q+[3]_q = \dfrac{[8]_q +[2]_q}{[2]_q} = \dfrac{[10]_q + [8]_q+ [6]_q + [4]_q+ [2]_q}{[3]_q[2]_q}= \dfrac{[6]_q[5]_q}{[3]_q[2]_q}.
\end{equation}
\end{remark}

\begin{remark}
The category $\DD$ is almost the $B_2$ spider category in \cite{Kupe}. But we replaced $q$ with $q^2$ and rescaled the trivalent vertex by $[2]_q^{-1/2}$. The trivalent vertex in $\DD$ may seem less natural since the relations now require us to insist $[2]_q$ is invertible, but when we connect the diagrammatic category to representation theory the rescaled trivalent vertex in $\DD$ will be more natural.
\end{remark}

%\begin{remark}
%We use the letter $\mathcal{D}$ to remind ourselves we are working in a diagrammatically defined category. The definition of $\DD$ is inspired by the category $\text{Rep}(\mathfrak{sp}_4(\mathbb{C}))$ ($C_2$ is the Dynkin type of the Lie algebra $\mathfrak{sp}_4(\mathbb{C})$). 
%\end{remark}

%===========
\subsection{Decomposing Tensor Products in $\text{Rep}(\mathfrak{sp}_4(\mathbb{C}))$}
\label{subsec-plethysm}
%===========

We now recall some basic facts about $\mathfrak{sp}_4(\mathbb{C})$ and its representation theory. Some of this is worked out in detail in \cite[Lecture 16]{Fultonharris}. Then we will record some formula's describing the decomposition of certain tensor products in $\Rep(\mathfrak{sp}_4)$. 

Let $X= \mathbb{Z}\epsilon_1\oplus \mathbb{Z}\epsilon_2$ be the weight lattice for $\mathfrak{sp}_4(\mathbb{C})$. The weights $\varpi_1 = \epsilon_1$ and $\varpi_2= \epsilon_1 + \epsilon_2$ are called the fundamental weights, and $X_+= \mathbb{Z}_{\ge0}\varpi_1\oplus\mathbb{Z}_{\ge 0}\varpi_2$ is the set of dominant weights.

Let $\Fund(\symp_4(\mathbb{C}))$ be the full monoidal subcategory of $\Rep(\symp_4(\mathbb{C}))$ generated by $\VV(\varpi_1)$ and $\VV(\varpi_2)$. The decomposition
\begin{equation}
\VV(\varpi_1)\ot \VV(\varpi_1) \cong \VV(2\varpi_1) \oplus \VV(\varpi_2) \oplus \VV(0).
\end{equation}
implies there is a one-dimensional space of maps between $\VV(\varpi_1)\ot \VV(\varpi_1)$ and $\VV(\varpi_2)$. We will later prove that there is a choice for this map so that sending the trivalent vertex to the chosen map gives a well-defined monoidal functor from $\mathbb{C}\ot \DD$ to $\Fund(\symp_4(\mathbb{C}))$. We will then show that this functor is full and faithful. 

For now we will take the equivalence on faith, and use it to guide our intuition for constructing a basis for Hom spaces in $\DD$. Let $\lambda$ and $\mu$ be dominant integral weights. There is a direct sum decomposition
\begin{equation}
\VV(\lambda)\ot \VV(\mu) \cong \bigoplus_{\nu\in X(\lambda, \mu)\subset \wt(\VV(\mu))} \VV(\lambda + \nu),
\end{equation}
where $\wt(\VV(\mu))$ is the multiset of weights in $\VV(\mu)$ and $X(\lambda, \mu)$ is a submultiset. Our goal is to determine the set $X(\lambda,  \mu)$.

To simplify notation, we may write $\VV(a, b)$ in place of $\VV(a\varpi_1 + b\varpi_2)$.  The following formulas are easy to work out using classical theory. For example, one can use \cite[2.16]{PRVformula}.

\begin{subequations}\label{plethysmformula}
\begin{equation} 
\VV(a,b) \ot \VV(1,0) \cong  \begin{cases} 
      \VV(1,0), \text{if} \ a= b= 0 \\
      \VV(a+1,0)\oplus \VV(a-1, 1) \oplus \VV(a-1,0), \text{if} \ a\ge 1, b=0  \\
      \VV(1, b)\oplus \VV(1,b-1), \text{if} \ a= 0, b\ge 1 \\
      \VV(a+1, b) \oplus \VV(a-1, b+1)\oplus \VV(a-1, b) \oplus \VV(a+ 1, b-1), \\  \ \ \ \ \ \ \ \ \ \ \ \ \  \ \ \ \ \ \ \ \ \ \  \ \ \ \ \ \ \ \ \ \ \ \ \ \text{if} \ a\ge 1, b\ge 1
   \end{cases}
\end{equation}
\begin{equation}
\VV(a,b) \ot \VV(0,1) \cong  \begin{cases} 
      \VV(0,1), \text{if} \ a= b= 0 \\
      \VV(0,b+1)\oplus \VV(2,b-1) \oplus \VV(0, b-1), \text{if} \ a=0, b\ge 1  \\
      \VV(1, 1)\oplus \VV(1,0), \text{if} \ a= 1, b=0 \\
      \VV(1, b+1)\oplus \VV(1, b)\oplus \VV(3, b-1)\oplus \VV(1, b-1), \ \text{if} \ a= 1, b\ge 1 \\
      \VV(a, 1) \oplus \VV(a, 0)\oplus \VV(a-2, 1), \ \text{if} \ a\ge 2, b= 0 \\
      \VV(a, b+1)\oplus \VV(a+2, b-1)\oplus \VV(a, b-1) \oplus \VV(a,b) \oplus \VV(a-2, b+1), \\ \ \ \ \ \ \ \ \ \ \ \ \ \ \ \ \ \ \ \ \ \ \ \ \ \ \ \ \ \ \ \ \ \ \ \ \ \ \ \ \  \ \ \ \ \ \ \ \ \ \ \ \ \  \ \text{if} \ a\ge 2, b\ge 1
   \end{cases}
\end{equation}
\end{subequations}

\begin{notation}
We will write $\VV(\blues) = \VV(\varpi_1) = \VV(1, 0)$ and $\VV(\greent) = \VV(\varpi_2) = \VV(0,1)$ as well as $\wt \blues= \varpi_1$ and $\wt \greent= \varpi_2$. Also, for a sequence $\underline{w} = (w_1, \ldots, w_n)$, $w_i\in \lbrace \blues, \greent\rbrace$ we will write $\VV(\underline{w})= \VV(w_1)\ot \ldots\ot \VV(w_n)$, $\wt \underline{w} = \wt w_1 + \wt w_2 + \ldots \wt w_n$, and $\underline{w}_{\le k} = (w_1, w_2, \ldots, w_k)$. 
\end{notation}

\begin{defn} Let $\underline{w}= (w_1, \ldots, w_n)$ with $w_i\in \lbrace \blues, \greent\rbrace$. A sequence $(\mu_1, \ldots, \mu_n)$ where $\mu_i\in \wt(\VV(w_i))$ is a \textbf{dominant weight subsequence} of $\underline{w}$ if: 

\begin{enumerate}
\item $\mu_1$ is dominant;
\item $\VV(\mu_1 + \ldots + \mu_{i-1} + \mu_i)$ is a summand of $\VV(\mu_1 + \ldots + \mu_{i-1})\ot \VV(w_i)$.
\end{enumerate}
We write $E(\underline{w})$ for the set of all dominant weight subsequences of $\underline{w}$ and
\begin{equation}
E(\underline{w}, \lambda):= \lbrace (\mu_1, \ldots, \mu_n) \in E(\underline{w}) \ : \ \mu_1+ \ldots + \mu_n = \lambda\rbrace
\end{equation}
for all $\lambda \in X_+$.
\end{defn}

\begin{lemma}\label{domwtdecomp}
Let $\underline{w}= (w_1, \ldots, w_n)$, $w_i\in \lbrace \blues, \greent\rbrace$, then
\begin{equation}
\VV(\underline{w}) \cong \bigoplus_{(\mu_1, \ldots, \mu_n)\in E(\underline{w})} \VV(\mu_1 + \ldots + \mu_n).
\end{equation}
Moreover, if we denote the multiplicity of $\VV(\lambda)$ as a summand of $\VV(\underline{w})$ by $[\VV(\underline{w}):\VV(\lambda)]$, then
\begin{equation}
[\VV(\underline{w}):\VV(\lambda)] = \# E(\underline{w}, \lambda).
\end{equation}
\end{lemma}

\begin{proof}
If we begin with $\VV(\emptyset)= \mathbb{C}$ and tensor with $\VV(w_1)$, there is only one irreducible summand. This summand corresponds to the dominant weight in $\wt \VV(w_1)$, which we record as $\mu_1$. Then we tensor $\VV(w_1)$ by $\VV(w_2)$ and note that $\VV(w_1)\ot \VV(w_2)$ contains $\VV(\mu_1)\ot \VV(w_2)$ as a summand. Choose a summand of $\VV(\mu_1)\ot \VV(w_2)$ and record this choice by the weight $\mu_2\in \wt \VV(w_2)$ so that the chosen summand is isomorphic to $\VV(\mu_1 + \mu_2)$. Next, we tensor $\VV(w_1)\ot \VV(w_2)$ by $\VV(w_3)$, observe that $\VV(w_1)\ot \VV(w_2)\ot \VV(w_3)$ contains a summand isomorphic to $\VV(\mu_1+ \mu_2)\ot \VV(w_3)$, and choose a weight $\mu_3\in \wt \VV(w_3)$ so that $\VV(\mu_1+ \mu_2+ \mu_3)$ is a summand of $\VV(\mu_1+ \mu_2)\ot \VV(w_3)$. Iterating this procedure, we end up with a sequence of weights $(\mu_1, \ldots, \mu_n)$, which is a dominant weight subsequence of $\underline{w}$, and a summand in $\VV(\underline{w})$ isomorphic to $\VV(\mu_1 + \ldots + \mu_n)$. Furthermore, all summands of $\VV(\underline{w})$ can be realized uniquely as the end result of the process we just described.
\end{proof}

\begin{lemma}
Let $\underline{u}= (u_1 ,\ldots, u_n)$ be a sequence with $u_i\in \lbrace \blues, \greent\rbrace$, then 
\begin{equation}\label{LLsdimhom}
\dim \Hom_{\mathfrak{sp}_4(\mathbb{C})}(\VV(\underline{w}), \VV(\underline{u})) =\sum_{\lambda\in X_+}[\VV(\underline{w}): \VV(\lambda)][\VV(\underline{u}): \VV(\lambda)]. 
\end{equation}
\end{lemma}

\begin{proof}
Thanks to Lemma \eqref{domwtdecomp}, this is consequence of Schur's lemma. 
\end{proof}

%===========
\subsection{Motivating the Light Ladder Algorithm}
\label{subsec-motivatell}
%===========

We outline a construction of a basis of homomorphism spaces in the category $\Fund(\symp_4(\mathbb{C}))$. This is a special case of a much more general construction \cite{tiltcellular}. 

Suppose that $(\mu_1, \ldots, \mu_m)\in E(\underline{w}, \lambda)$. For $i= 1, \ldots, m$ there is a projection map $P_{(\mu_1, \ldots, \mu_i)}: \VV(w_1)\ot \ldots\ot \VV(w_i)\longrightarrow \VV(\mu_1 + \ldots + \mu_i)$. The map $P_{(\mu_1, \ldots, \mu_i)}$ is the projection $P_{(\mu_1, \ldots, \mu_{i-1})}: \VV(w_1)\ot \ldots \ot \VV(w_{i-1})\longrightarrow \VV(\mu_1 + \ldots + \mu_{i-1})$ postcomposed with the projection $p_{\mu_i}: \VV(\mu_1 + \ldots + \mu_{i-1})\ot \VV(w_i)\longrightarrow \VV(\mu_1 + \ldots + \mu_i)$. 

Let $(\nu_1, \ldots, \nu_n)\in E(\underline{v}, \lambda)$. Now, for $i= 1,\ldots,n$ there are inclusion maps $I^{(\nu_1, \ldots, \nu_i)}: \VV(\nu_1 + \ldots + \nu_i)\longrightarrow \VV(u_1)\ot \ldots \ot \VV(u_i)$. Composing the projection with the inclusion we get a map $I^{(\nu_1, \ldots, \nu_n)}\circ P_{(\mu_1, \ldots, \mu_m)}: \VV(\underline{w})\longrightarrow \VV(\underline{u})$, factoring through $\VV(\lambda)$. 

Since $[\VV(\lambda): \VV(\underline{w})]= E(\underline{w}, \lambda)$ and $[\VV(\lambda): \VV(\underline{u})]= E(\underline{u}, \lambda)$, the maps 
\begin{equation}\label{PIbasis}
\bigcup_{\substack{ \lambda\in X_+ \\ (\mu_1, \ldots, \mu_m)\in E(\underline{w}, \lambda) \\ (\nu_1, \ldots, \nu_n)\in E(\underline{u}, \lambda)}} \lbrace I^{(\nu_1, \ldots, \nu_n)}\circ P_{(\mu_1, \ldots, \mu_m)}\rbrace
\end{equation}
form a basis in $\Hom_{\mathfrak{sp}_4(\mathbb{C})}(\VV(\underline{w}), \VV(\underline{u}))$. 

The maps $P_{(\mu_1, \ldots, \mu_n)}$ are built inductively out of the $p_{\mu_i}$'s in a way that is analogous to how we will define light ladder diagrams in terms of elementary light ladder diagrams. The inclusion map $I^{(\nu_1, \ldots, \nu_n)}: \VV(\lambda)\longrightarrow \VV(\underline{u})$ is analogous to what we will call upside down light ladder diagrams. We will define double ladder diagrams as the composition of a light ladder diagram and an upside down light ladder diagram, in analogy with the $I\circ P$'s. Then our work will be to argue that double ladder diagrams are a basis. 

\begin{remark}
The projection and inclusion maps we discuss here are not the image of the light ladder diagrams under a functor $\DD\longrightarrow \Fund(\symp_4(\mathbb{C}))$. There are at least two reasons for this. The first being that the object $\VV(\lambda)$ is not in the category $\Fund(\symp_4(\mathbb{C}))$, so we have to construct light ladder maps not from $\VV(\underline{w})$ to $\VV(\lambda)$, but from $\VV(\underline{w})$ to $\VV(\underline{x})$ where $\wt \underline{x}= \lambda$. 

The second reason is that we want to construct a basis for the diagrammatic category which descends to a basis in $\Fund$ for fields other than $\mathbb{C}$. Over other fields the representation theory is no longer semisimple so $\VV(\lambda)$ may not be a summand of $\VV(\underline{w})$. There will still be the same number of maps from $\VV(\underline{w})$ to a suitable version of $\VV(\lambda)$ but they may not be inclusions and projections. 
\end{remark}

%===========
\subsection{Light Ladder Algorithm}
\label{subsec-algorithm}
%===========

Now we define some morphisms in the diagrammatic category. 

\begin{defn}\label{elementaryll}
An \textbf{elementary light ladder diagram} is one of the following diagrams in $\DD$. We will say that $L_{\mu}$ is the elementary light ladder diagram of weight $\mu$.

\begin{figure}[H]
\centering
\includegraphics[width=0.15\textwidth]{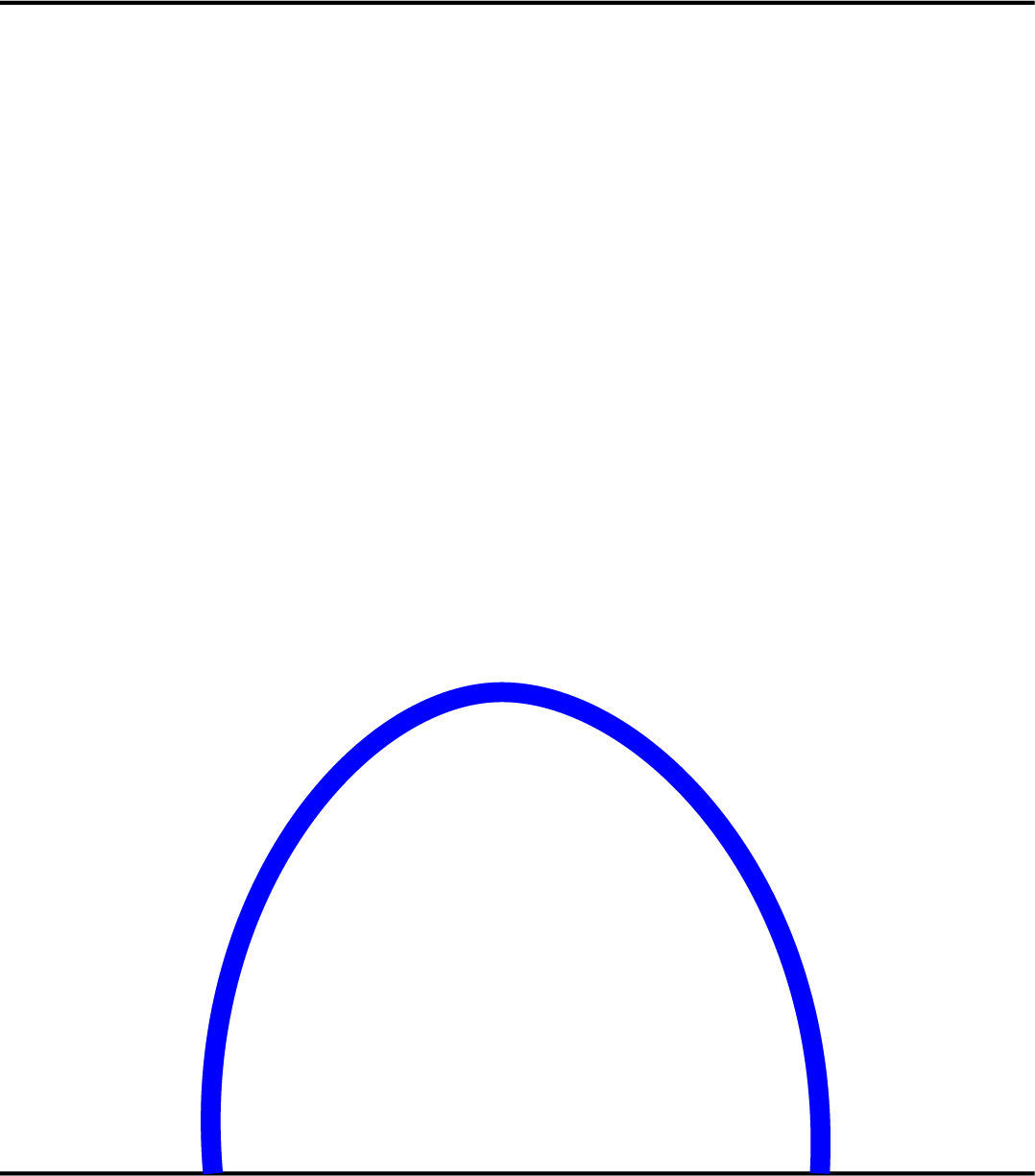}\put(-50,-15){$L_{(-1,0)}$} \ \ \ \ \ \includegraphics[width=0.15\textwidth]{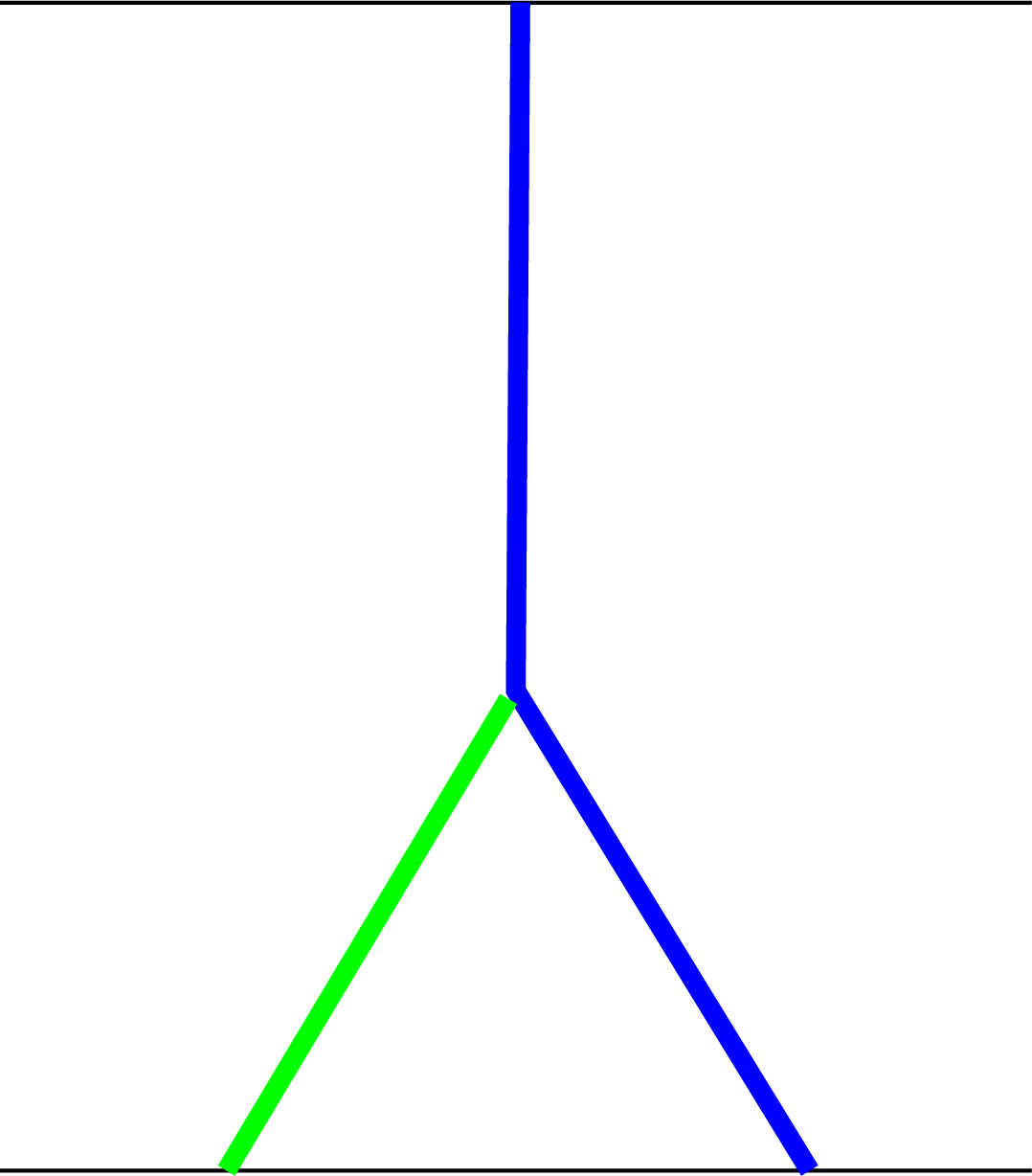}\put(-50,-15){$L_{(1,-1)}$} \ \ \ \ \ \includegraphics[width=0.15\textwidth]{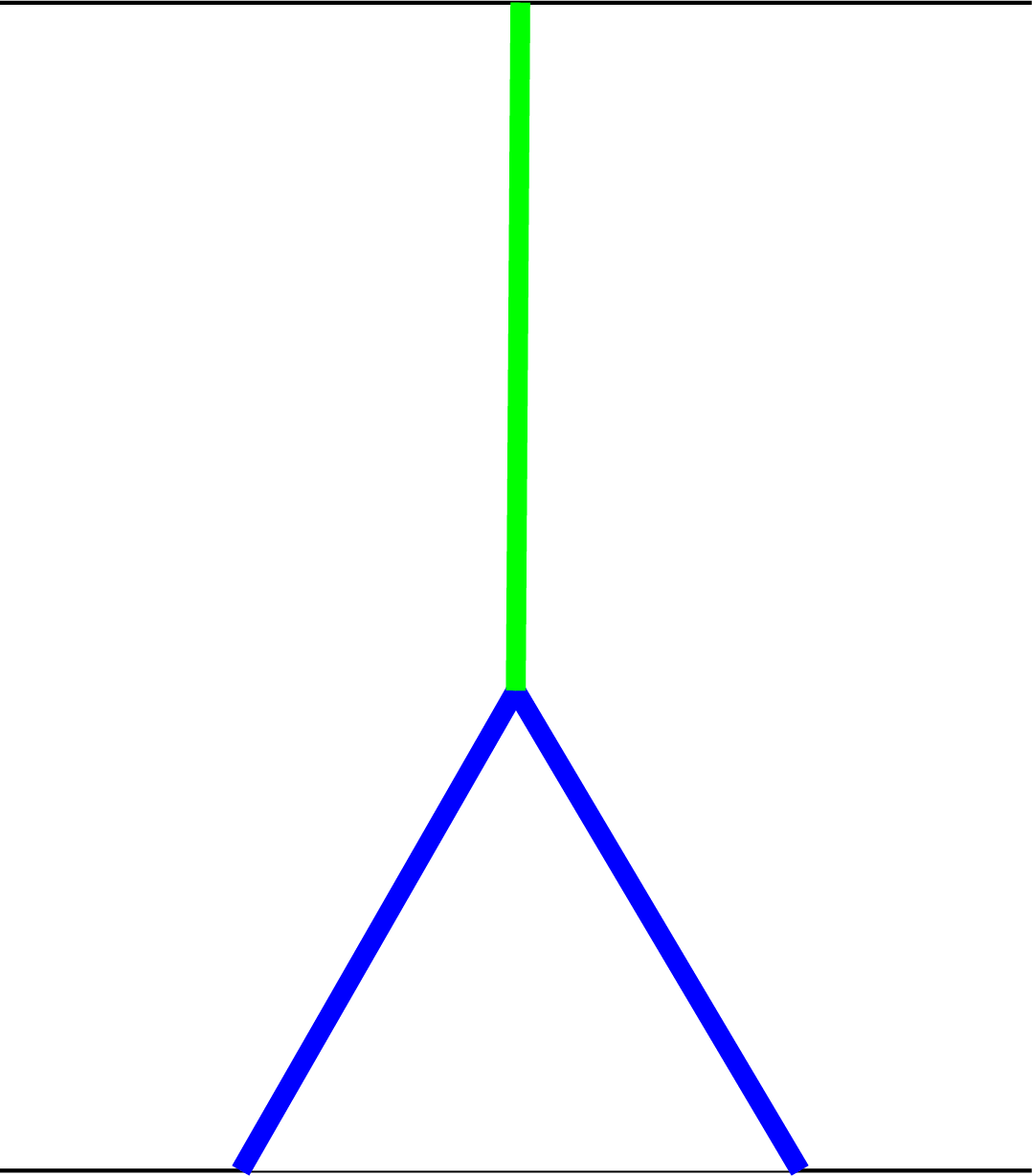}\put(-50,-15){$L_{(-1,1)}$} \ \ \ \ \ \includegraphics[width=0.15\textwidth]{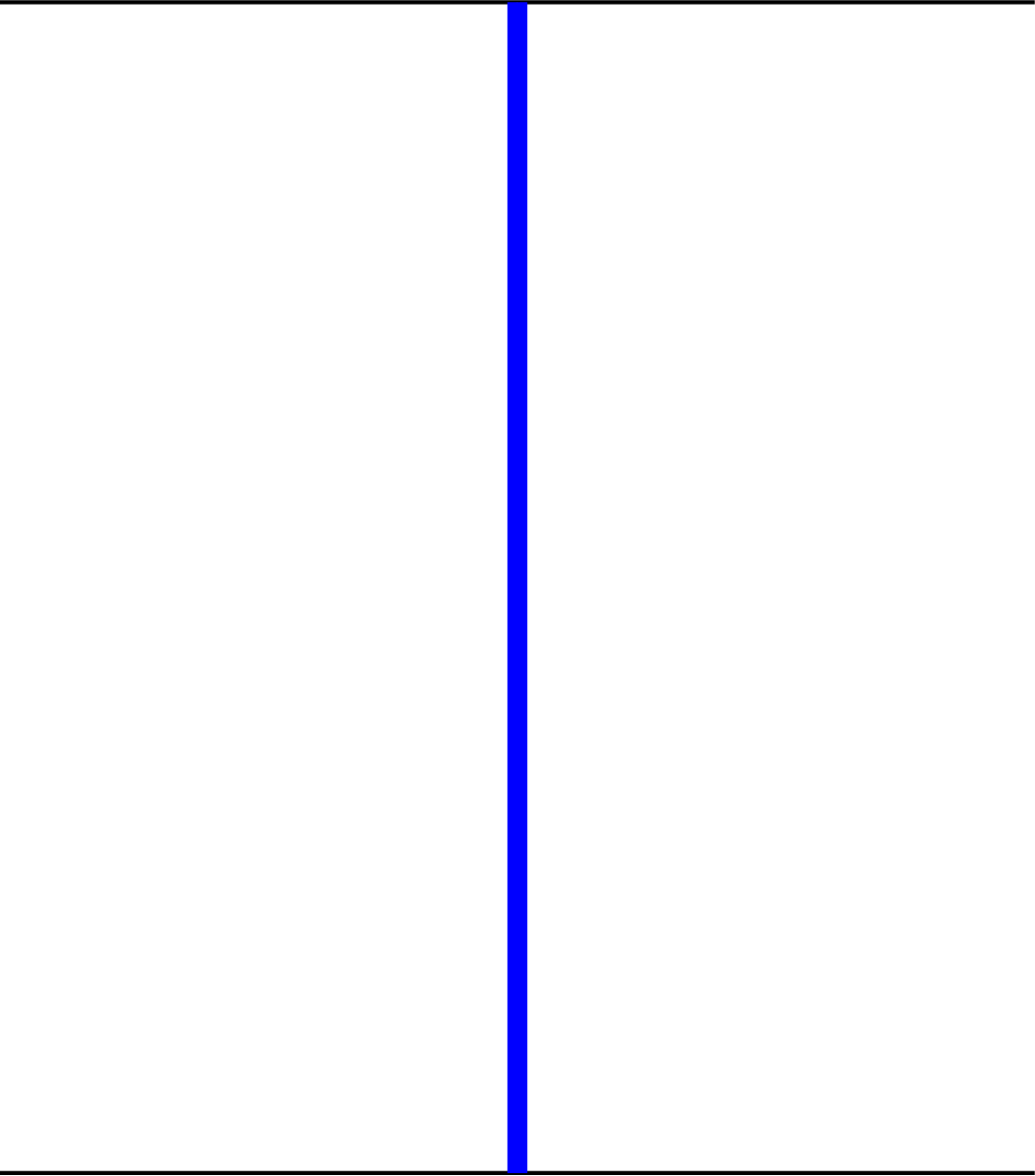}\put(-50,-15){$L_{(1, 0)}$}
\end{figure}

\begin{figure}[H]
\centering
\includegraphics[width=0.15\textwidth]{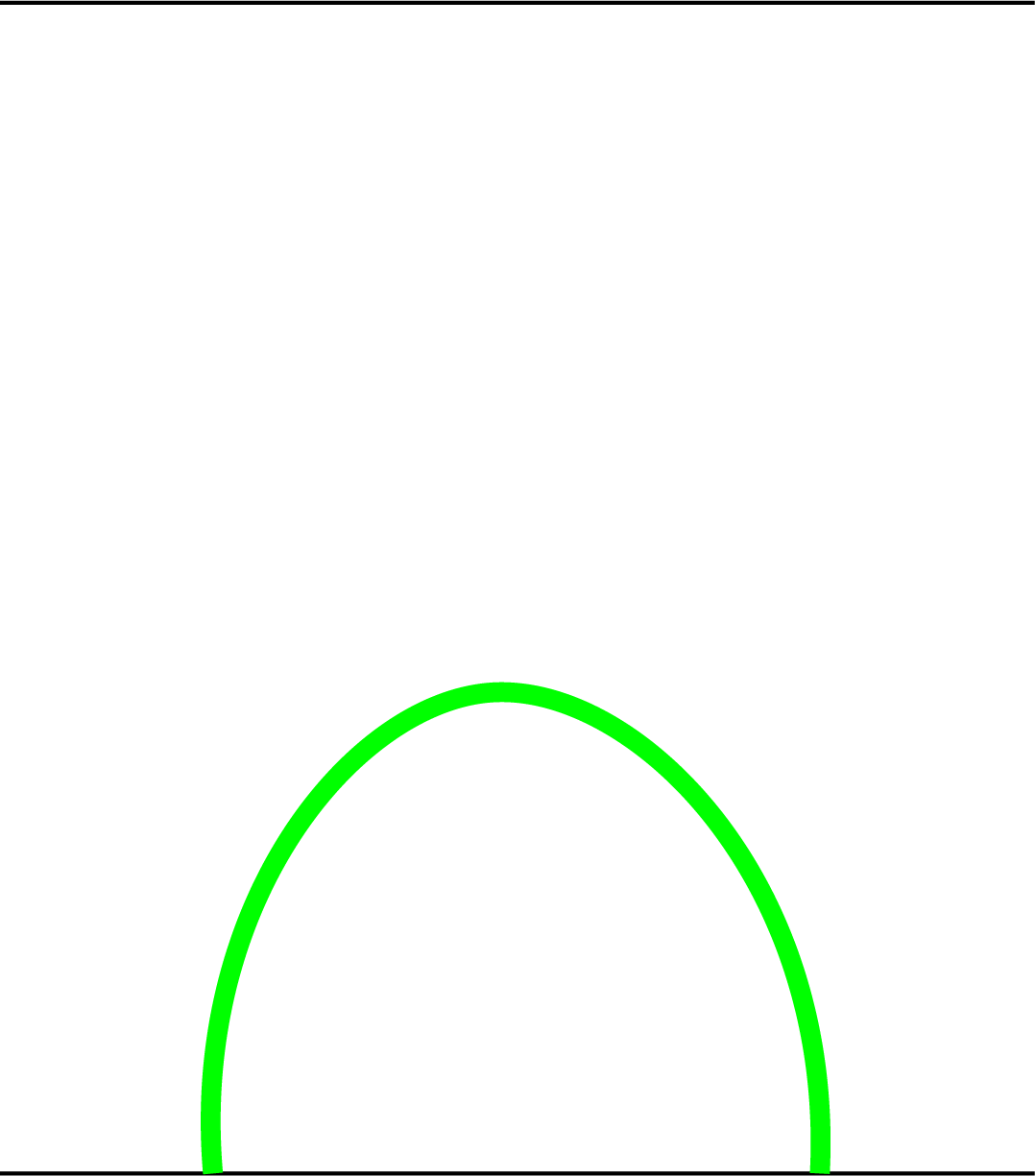}\put(-50,-15){$L_{(0,-1)}$}\ \ \ \ \ \includegraphics[width=0.15\textwidth]{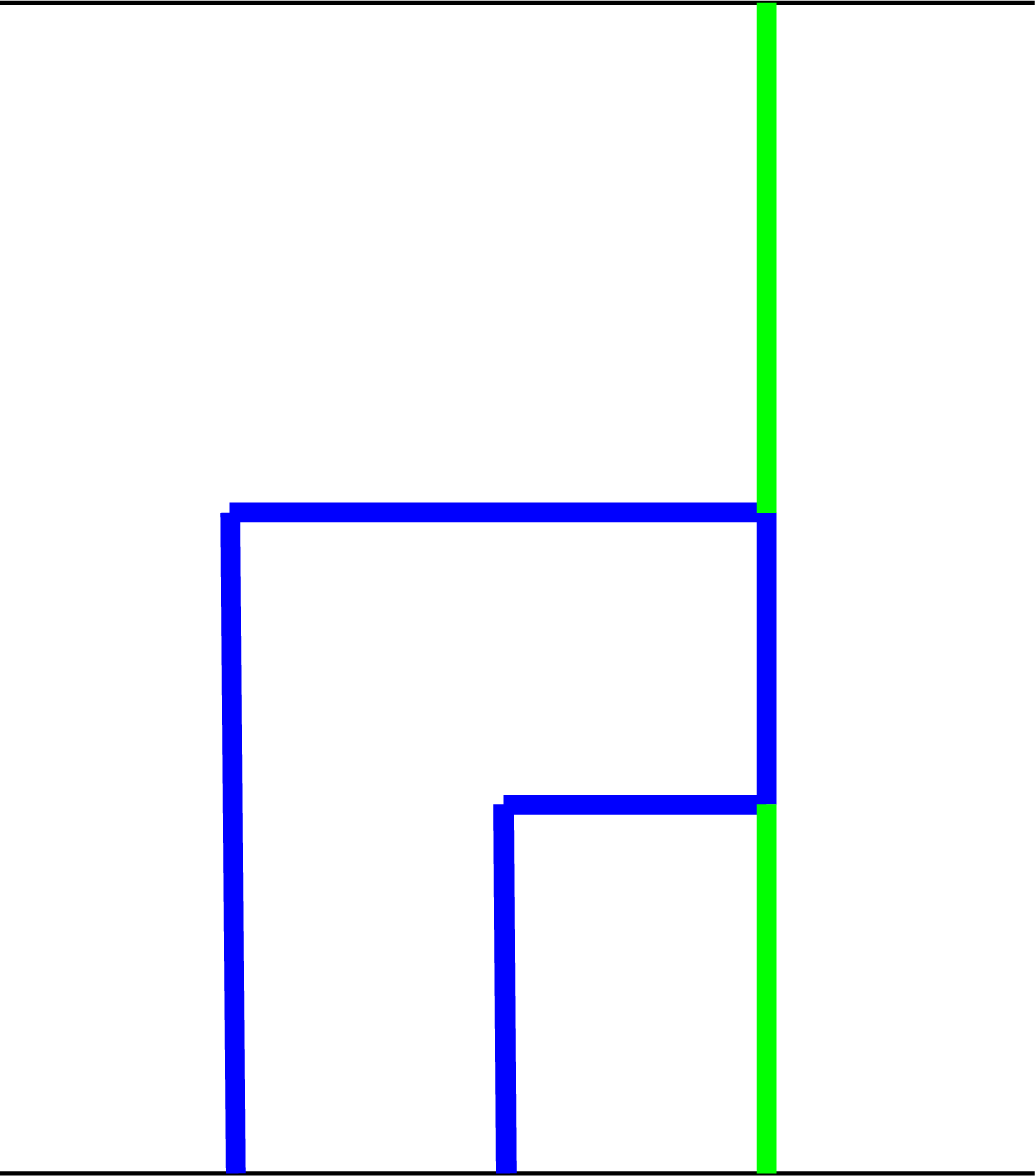}\put(-50,-15){$L_{(-2,1)}$} \ \ \ \ \ \includegraphics[width=0.15\textwidth]{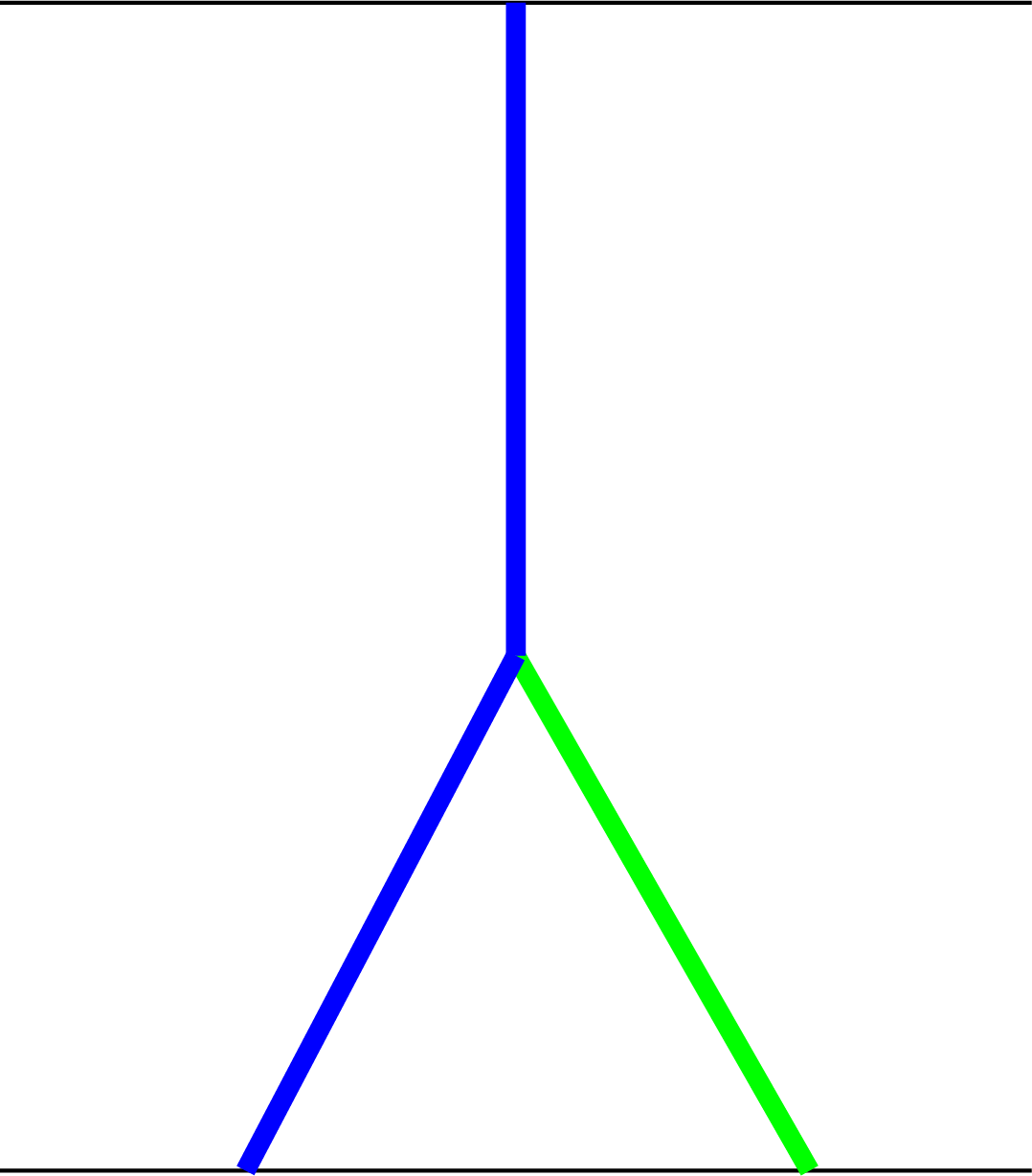}\put(-50,-15){$L_{(0,0)}$} \ \ \ \ \ \includegraphics[width=0.15\textwidth]{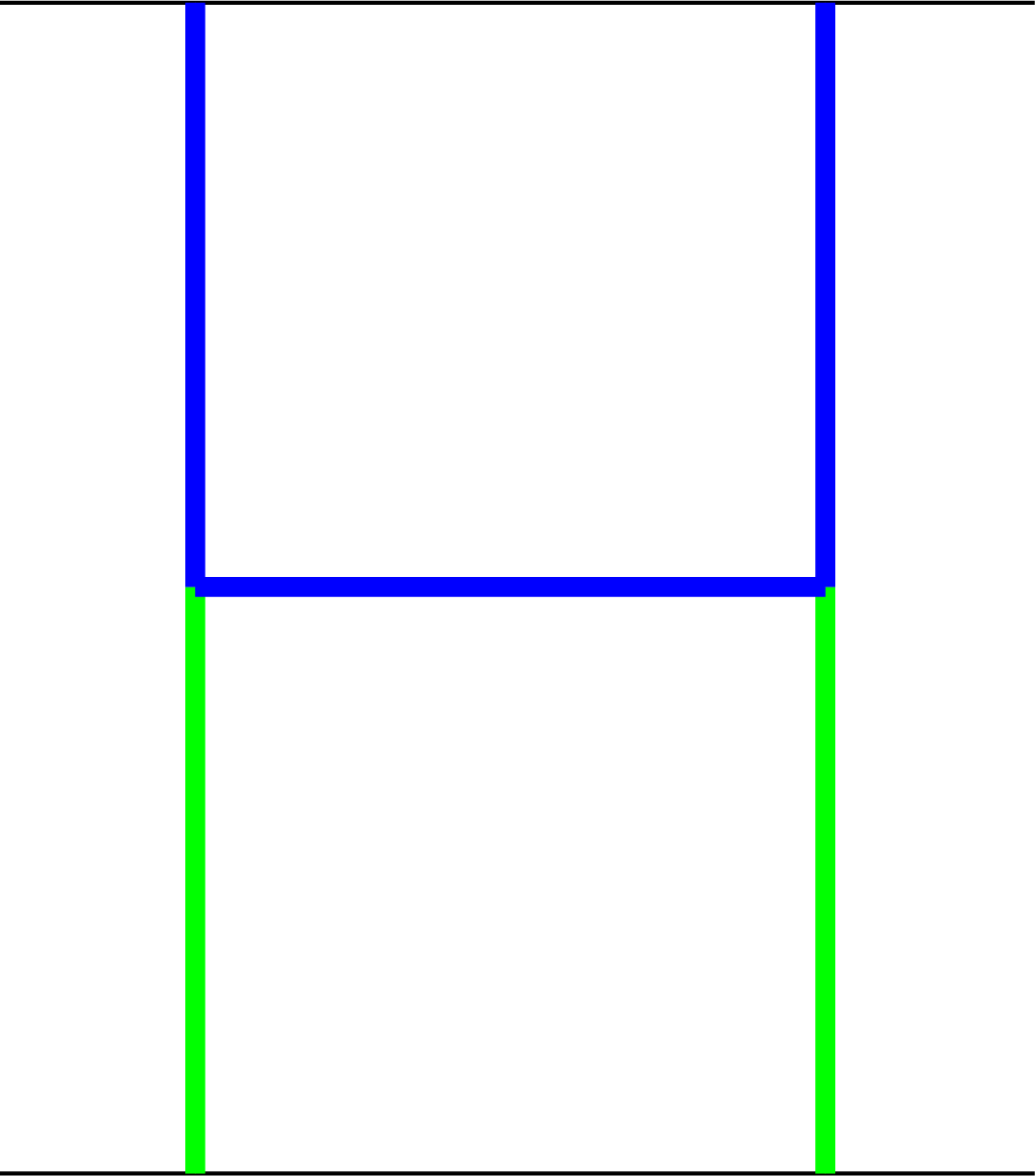}\put(-50,-15){$L_{(2,-1)}$} \ \ \ \ \ 
\includegraphics[width=0.15\textwidth]{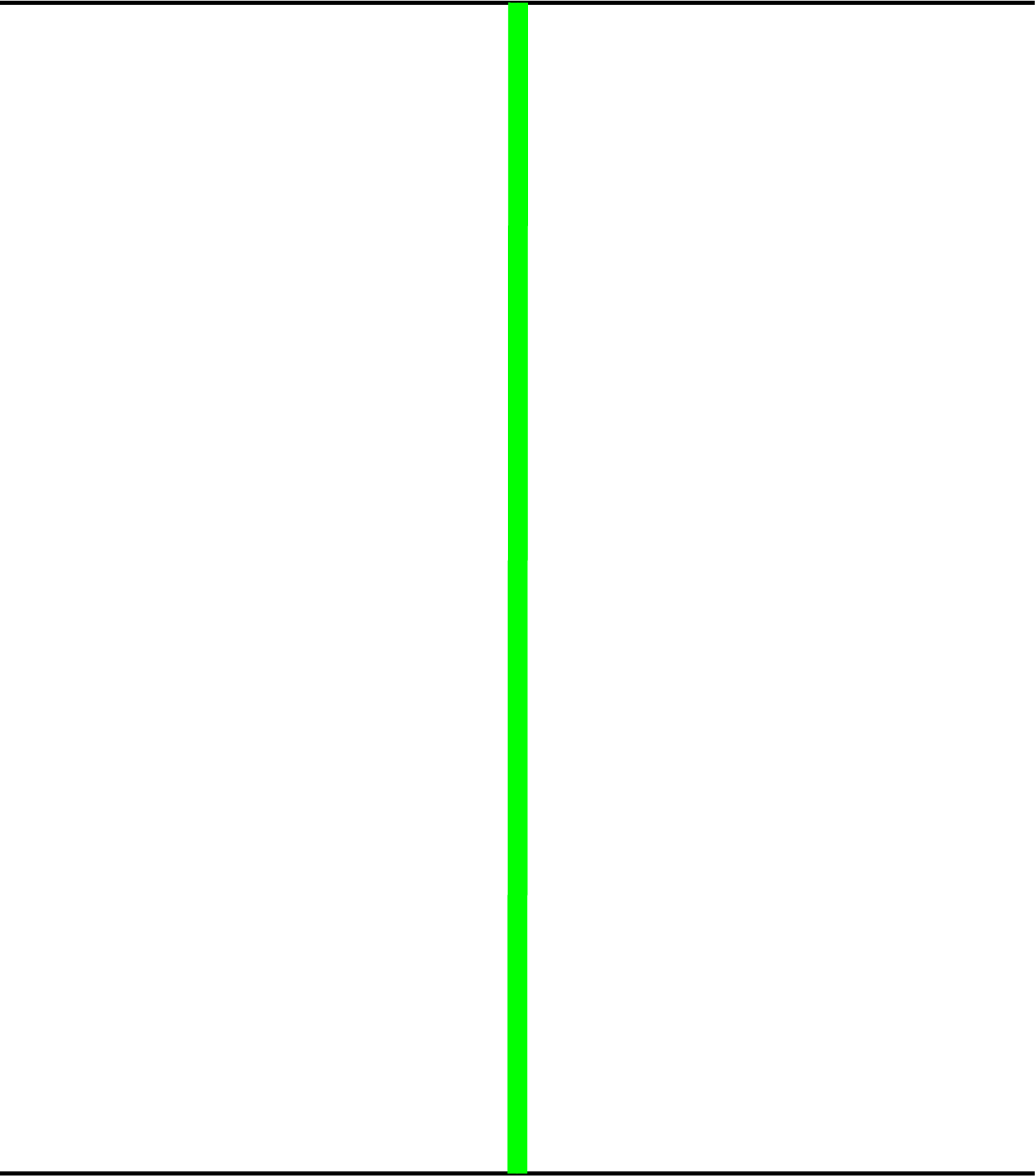}\put(-50,-15){$L_{(0,1)}$}
\end{figure}
\end{defn}

\begin{remark}
If $L_{\mu}:\underline{u}\ast \rightarrow \underline{w}$, for $\ast \in \lbrace \blues, \greent\rbrace$, then $\mu\in \wt \VV(\ast)$ and $\wt \underline{w} = \wt \underline{u} + \mu$. 
\end{remark}

\begin{defn}A \textbf{neutral diagram} is any diagram which is the horizontal and/or vertical composition of identity maps and the following \textbf{basic neutral diagrams}.

\begin{figure}[H]
\centering
\includegraphics[width=0.15\textwidth]{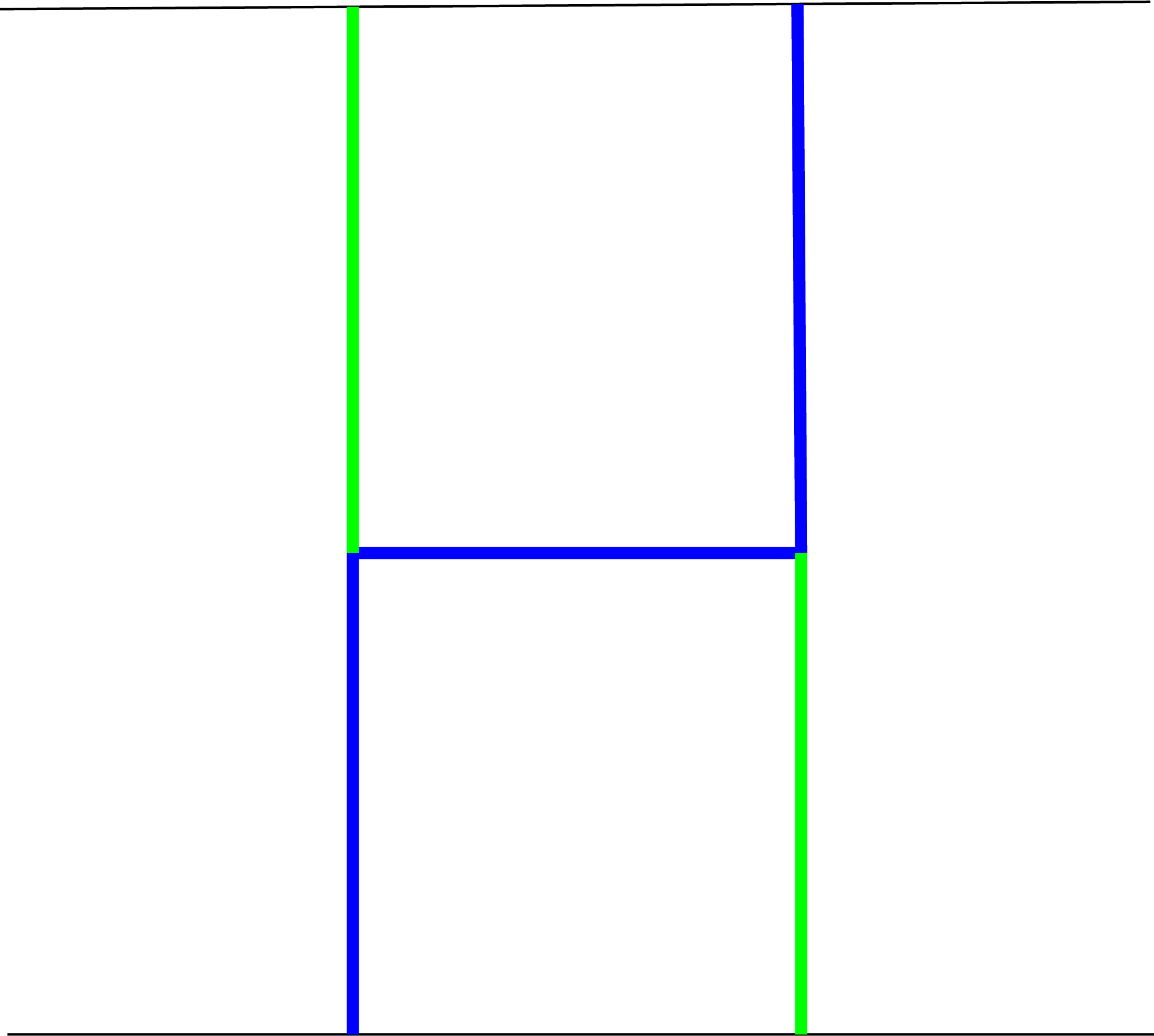}\put(-45,-15){$N_{\blues\greent}^{\greent\blues}$} \ \ \ \ \ \includegraphics[width=0.15\textwidth]{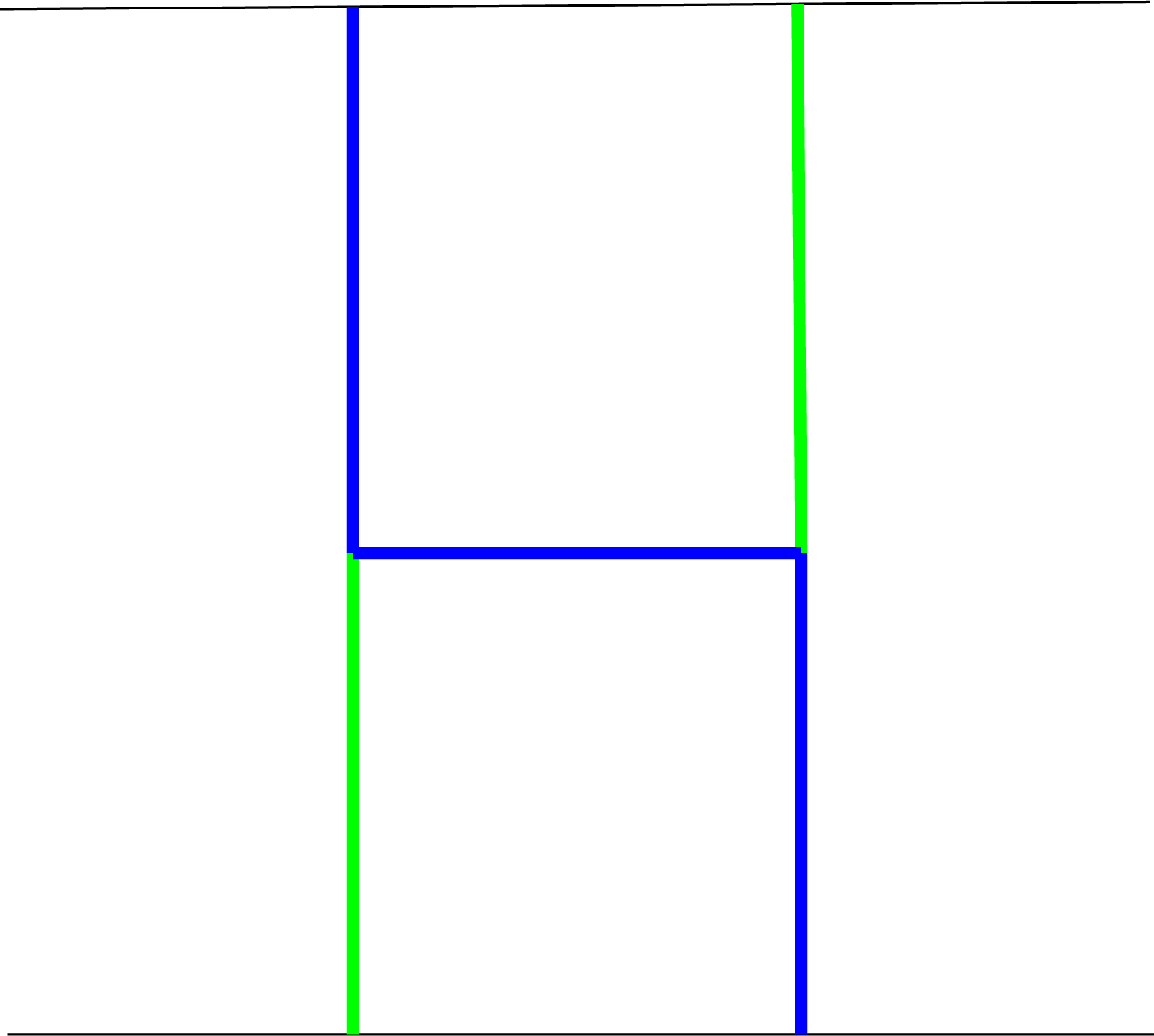}\put(-45,-15){$N_{\greent\blues}^{\blues\greent}$} \ \ \ \ \ 
\end{figure}
\end{defn}

\begin{figure}[H]
\caption{A neutral diagram from $\blues\blues\greent\greent\greent\blues$ to $\greent\greent\blues\greent\blues\blues$.}
\centering
\includegraphics[width=0.2\textwidth]{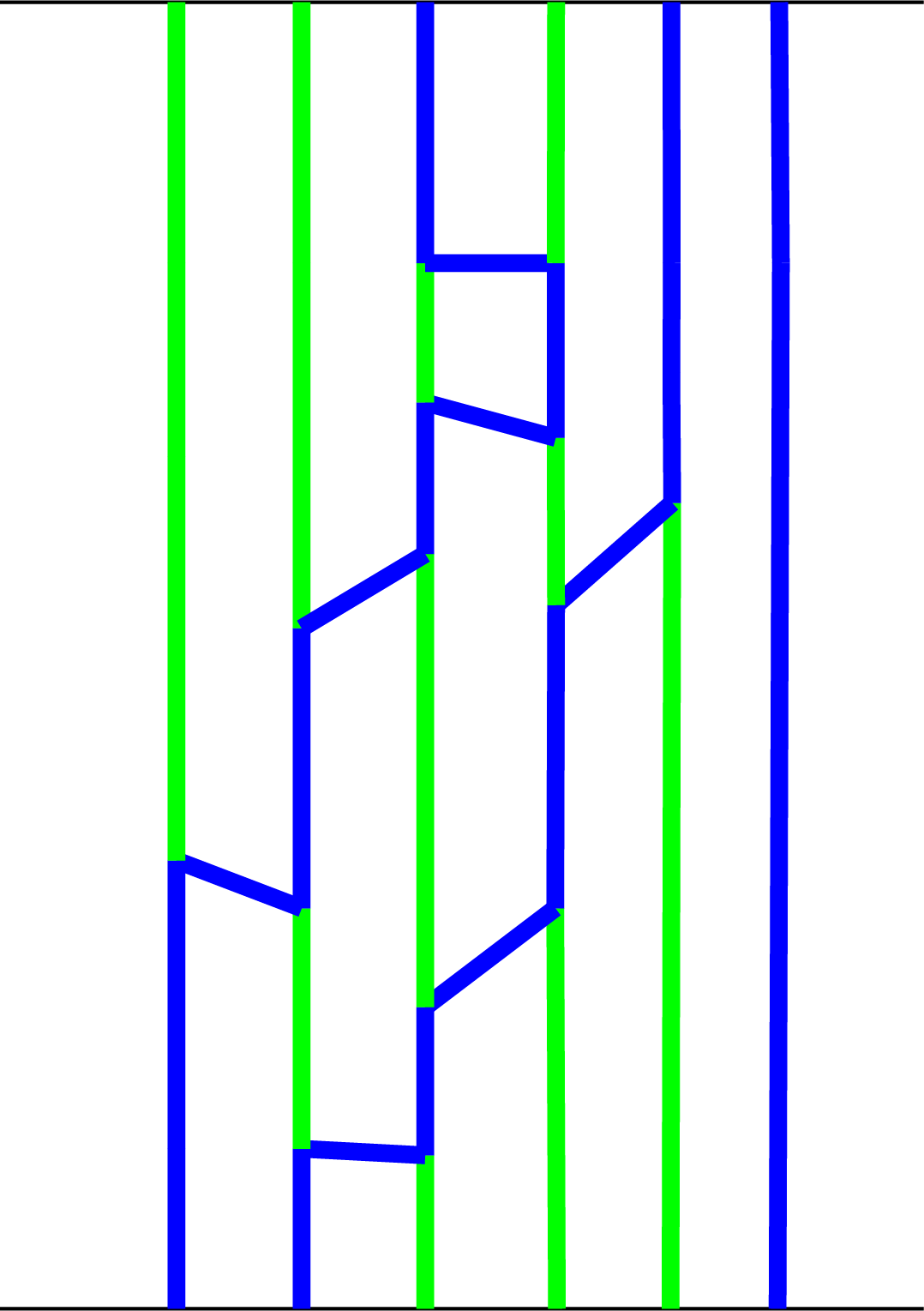} 
\end{figure}

\begin{defn}\label{lightladder}
Fix an object $\underline{w}  = (w_1, \ldots, w_n)$ in $\DD$, a dominant weight subsequence $\vec{\mu} = (\mu_1, \ldots, \mu_n)\in E(\underline{w})$, and an object $\underline{v}= (v_1, \ldots v_m)$ in $\DD$ such that $\wt \underline{v}= \mu_1 + \ldots + \mu_n$. We will describe an algorithm, which we will refer to as the \textbf{light ladder algorithm}, to construct a diagram in $\DD$ with source $\underline{w}$ and target $\underline{v}$. This diagram will be denoted $LL_{\underline{w}, \vec{\mu}}^{\underline{v}}$ and we will call it a \textbf{light ladder diagram}.

We define the diagrams inductively, starting by defining $LL_{\emptyset, (\emptyset)}^{\emptyset}$ to be the empty diagram. Suppose we have constructed $LL_{\underline{w}_{\le n-1}, (\mu_1, \ldots, \mu_{n-1})}^{\underline{u}}$, where $\wt (\underline{u}) = \mu_1 + \ldots + \mu_{n-1}$. Then we define 
\begin{equation}
LL_{\underline{w}, (\mu_1, \ldots, \mu_n)}^{\underline{v}} = N_{?}^{\underline{v}}\circ \left(\id \ot L_{\mu_n}\right)\circ \left(N_{\underline{u}}^{?}\ot \id\right) \circ \left(LL_{\underline{w}_{\le n-1}, (\mu_1, \ldots, \mu_{n-1})}^{\underline{u}}\ot \id_{w_n}\right)
\end{equation}
where $N_{?}^{?}$ is a neutral diagram with appropriate source (subscript) and target (superscript). 
\end{defn}

\begin{figure}[H]
\caption{A schematic for the inductive definition of a light ladder diagram $LL_{\underline{w}, (\mu_1, \ldots, \mu_n)}^{\underline{v}}$.}
\centering
\includegraphics[width=.6\textwidth]{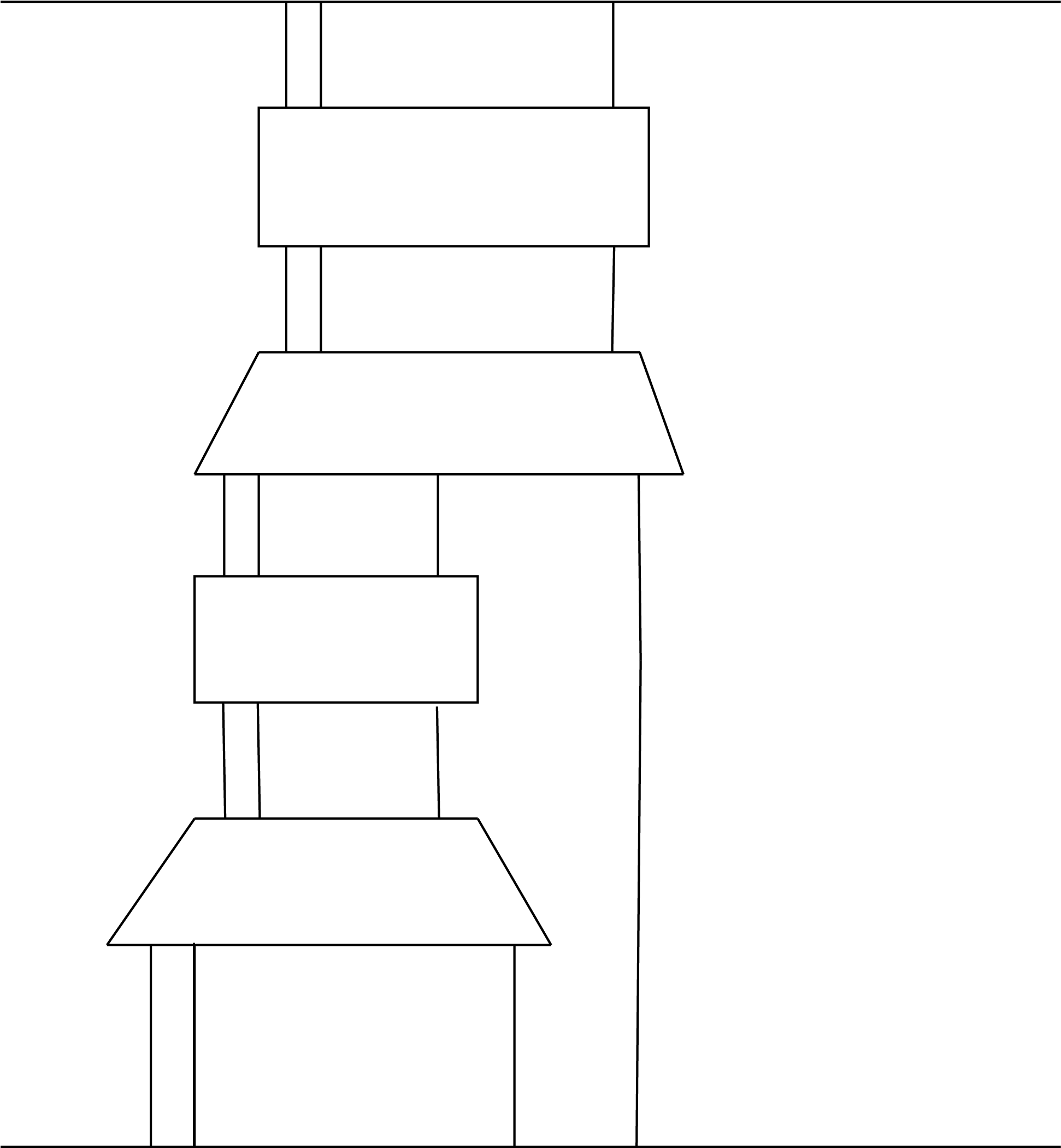}\put(-228,60){$LL_{\underline{w}_{\le n-1}, (\mu_1, \ldots, \mu_{n-1})}^{\underline{u}}$}\put(-190, 125){$N_{\underline{u}}^{?}$}\put(-170, 180){$\id\ot L_{\mu_n}$}\put(-165, 240){$N_{?}^{\underline{v}}$}\put(-265, 20){$\underline{w}_{\le n-1}$}\put(-125, 20){$w_n$}\put(-225, 95){$\underline{u}$}\put(-225, 152){$?$}\put(-205, 210){$?$}\put(-210, 270){$\underline{v}$}
\end{figure}

To further aid the readers understanding of the light ladder construction we give an example and some clarifying comments.  

\begin{example}
The light ladder diagram $LL_{\greent\blues\greent\blues\greent, ((0,1), (1, -1), (-1,0), (2, -1))}^{\blues\blues}$

\begin{figure}[H]
\centering
\includegraphics[width=.35\textwidth]{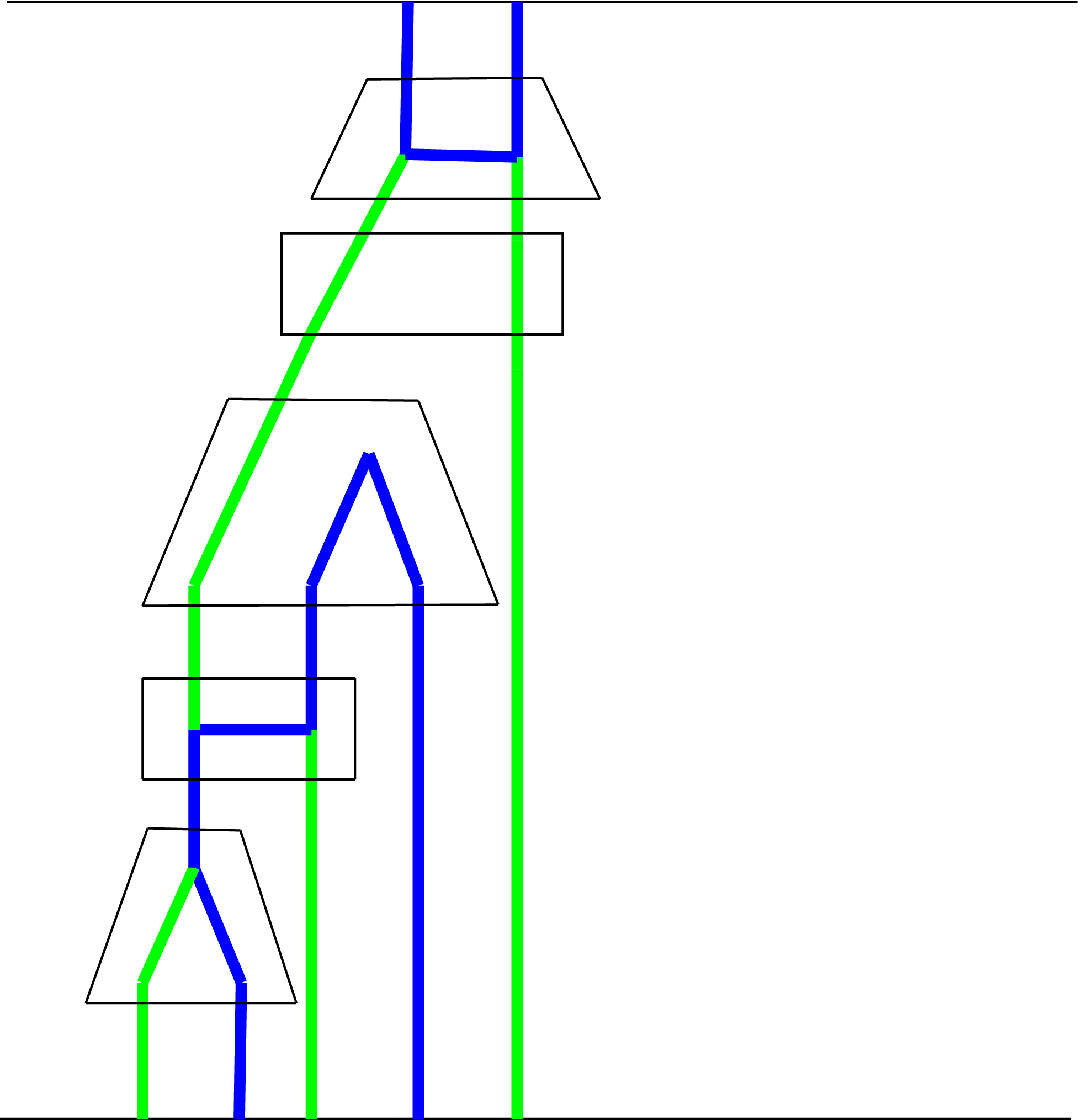}
\end{figure}

Our convention of rectangles and trapezoids is to indicate whether a diagram is a neutral diagram or a diagram of the form $\id\ot \ \text{elementary light ladder diagram}$. We omitted the first step corresponding to $\mu= (0,1)$.
\end{example}

The elementary light ladder diagrams have fixed source and target. As a result one can construct $LL_{\underline{w}_{\le n-1}, (\mu_1, \ldots, \mu_{n-1})}^{\underline{u}}$, then see that $\VV(\mu_n)$ is a summand of $\VV(\mu_1 + \ldots + \mu_{n-1})\ot \VV(w_n)$, but still not guarantee there is an object $\underline{y}$ in $\DD$ so that $\underline{y}w_n$ is the source of $\id \ot L_{\mu_n}$. 

\begin{figure}[H]
\caption{An example of what can go wrong without neutral diagrams.}
\centering
\includegraphics[width=.3\textwidth]{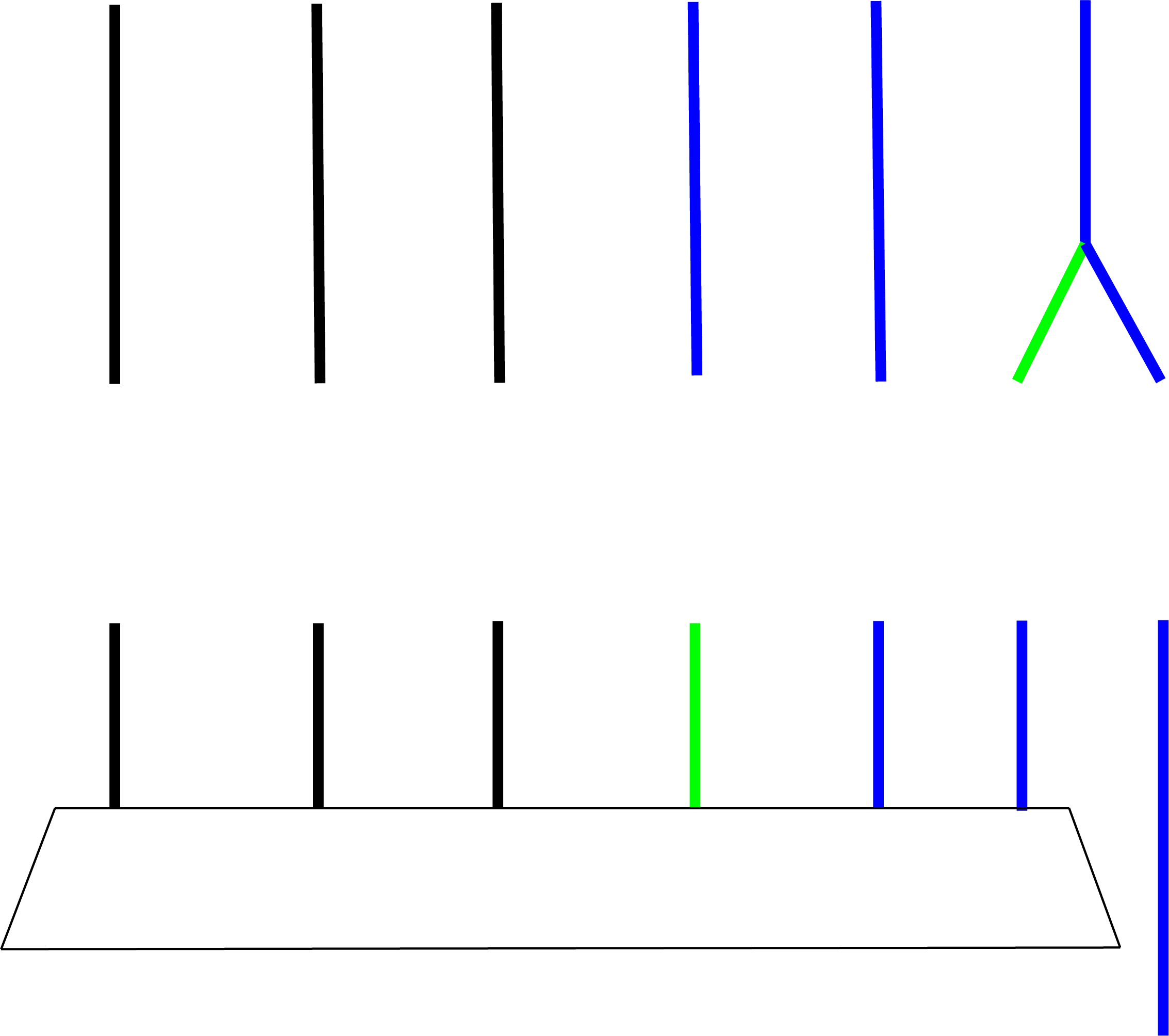}
\end{figure}

Basic neutral diagrams encode isomorphisms $\blues\greent \rightarrow \greent\blues$ and $\greent\blues\rightarrow \blues \greent$, while arbitrary neutral diagrams encode isomorphisms $\underline{w}\rightarrow \underline{w}'$. Intuitively, one could think that neutral diagrams are built out of colored crossings which interchange $\blues\greent$ and $\greent\blues$. 

\begin{figure}[H]
\caption{Colored crossing diagrams.}
\centering
\includegraphics[width=.2\textwidth]{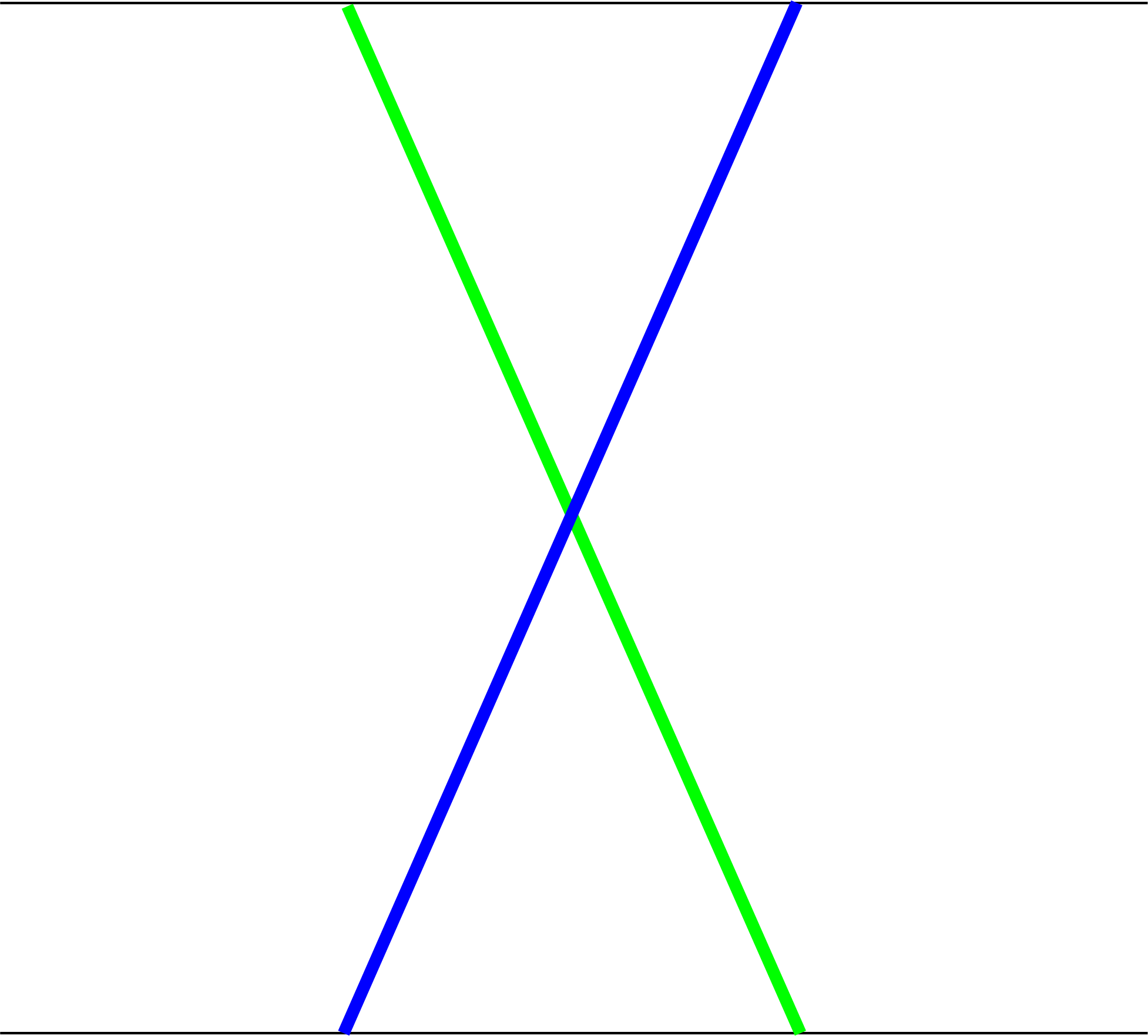} \ \ \ \ \ \ \ \ \ \ \ \ \ \ \ \ \ \ \ \ 
\includegraphics[width=.2\textwidth]{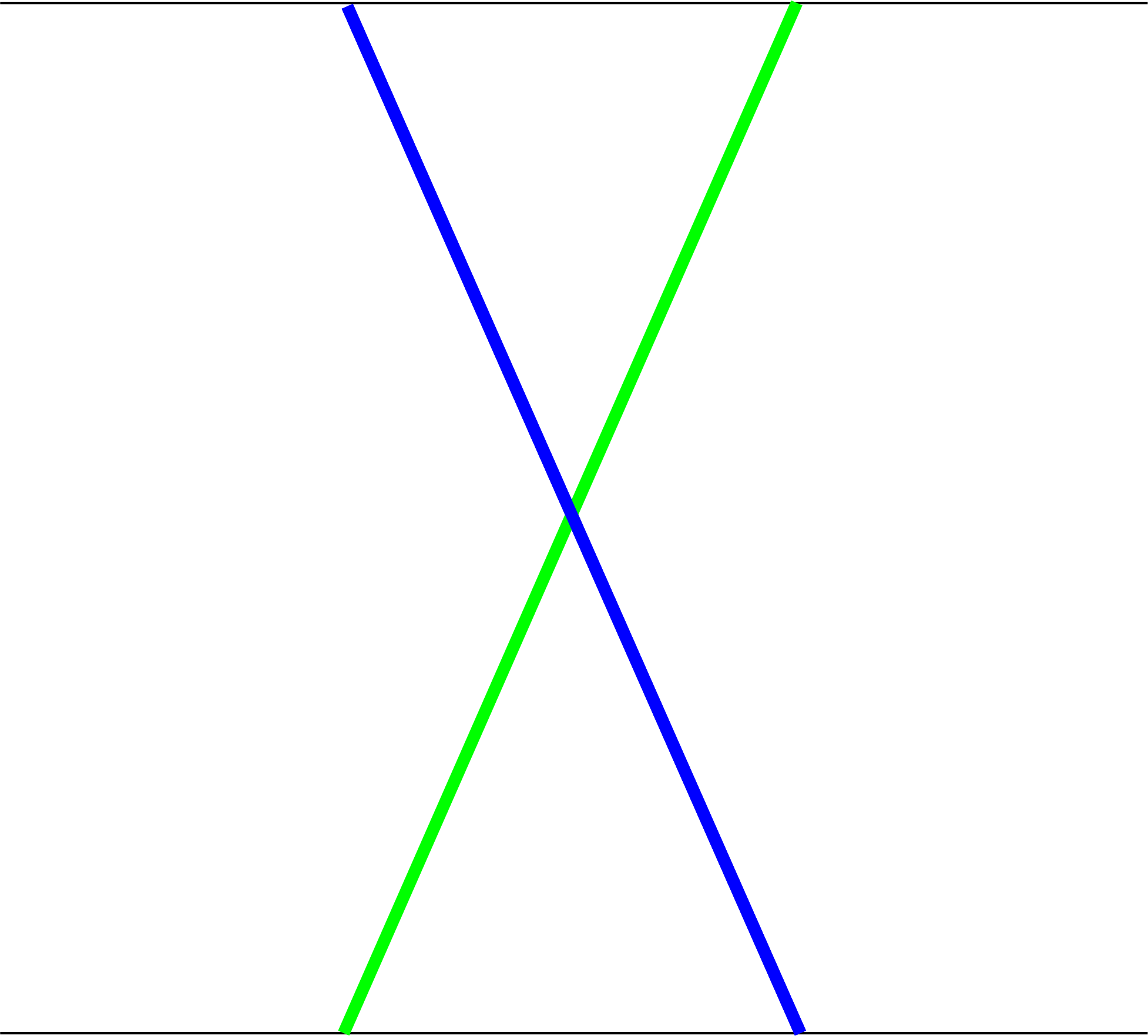}
\end{figure}

\begin{rmk}
The reason we use basic neutral diagrams instead of colored crossing diagrams is that the latter are not allowed as graphs in our definition of morphisms in $\DD$. But there is an obvious way to convert a colored crossing diagram to neutral diagram, and vice versa. 
\end{rmk}

\begin{lemma}
Given two sequences $\underline{w}$ and $\underline{w}'$ such that $\wt \underline{w} = \wt \underline{w}'$, there is a diagram built from colored crossings connecting $\underline{w}$ to $\underline{w}'$. 
\end{lemma}
\begin{proof}
Suppose that $\wt \underline{w} = a\varpi_1 + b \varpi_2 = \wt\underline{w}'$. Connect both $\underline{w}$ and $\underline{w}'$ via colored crossing diagrams to the standard sequence $\blues^{\ot a}\ot \greent^{\ot b}$ and then compose the diagram from $\underline{w}$ to the standard sequence with the vertical flip of the diagram from the standard diagram to $\underline{w}'$. 
\end{proof}
\begin{lemma}\label{shuffleneutrals}
Given two sequences $\underline{w}$ and $\underline{w}'$ such that $\wt \underline{w} = \wt \underline{w}'$, there is a neutral diagram connecting $\underline{w}$ to $\underline{w}'$. 
\end{lemma}
\begin{proof}
Replacing the colored crossings with the associated basic neutral diagrams, we obtain a neutral diagram $\underline{w}\rightarrow \underline{w}'$. 
\end{proof}

The following lemma uses this observation to fix the problem, in the light ladder algorithm, of elementary diagrams having fixed source and target. 

\begin{lemma}
Let $(\mu_1, \ldots, \mu_n)\in E(\underline{w})$ (in particular, $\VV(\mu_n)$ is a summand of $\VV(\mu_1 + \ldots + \mu_{n-1})\ot \VV(w_n)$). Suppose we have constructed $LL_{\underline{w}_{\le n-1}, (\mu_1, \ldots, \mu_{n-1})}^{\underline{u}}$. There is an object $\underline{y}$ in $\DD$ and a neutral map $N_{\underline{w}_{\le n-1}}^{\underline{y}}$ such that $\underline{y}\ot w_n$ is the source of $\id \ot L_{\mu_n}$. 
\end{lemma}

\begin{proof}
We will argue this for the elementary diagram $L_{(1, -1)}$, so $\mu_n = (1, -1)$ and $w_n= \blues$. The arguments for the rest of the cases follow the same pattern. From the tensor product decomposition formulas \eqref{plethysmformula} we see that $\VV(1, -1)$ being a summand of $\VV(\mu_1 + \ldots + \mu_{n-1})$ implies that, if $\mu_1 + \ldots + \mu_{n-1}= a\varpi_1 + b\varpi_2$, then $b\ge 1$. Thus, in the sequence $\underline{u}= (u_1, \ldots, u_k)$  there is some $k$ so that $u_k= \greent$. By Lemma \eqref{shuffleneutrals} there is a neutral diagram from the sequence $\underline{u}$ to a sequence which ends in $\greent$. The target of this neutral diagram will be an object $\underline{y}$ such that $\underline{y} \ot \blues$ is the source of $\id\ot L_{(1, -1)}$. 
\end{proof}

\begin{figure}[H]
\caption{Using a neutral map to fix the problem.}
\centering
\includegraphics[width=.3\textwidth]{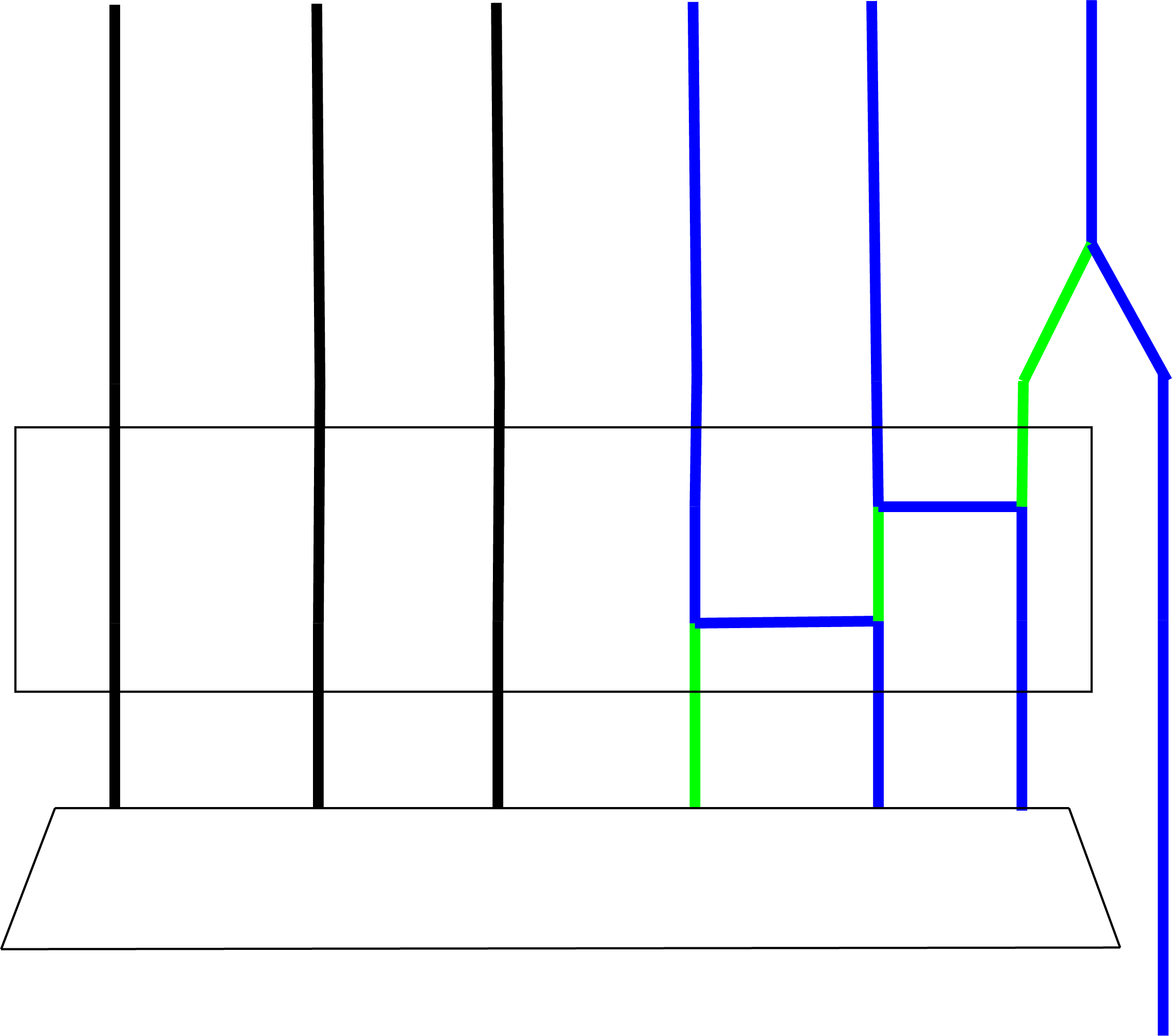}
\end{figure}

Comparing the tensor product decompositions in \eqref{plethysmformula} with the elementary light ladder diagrams it is apparent that dominant weight subsequences always produce a light ladder diagram. Neutral diagrams from one word to another are not unique. The choice of neutral diagram could result in several different light ladder diagrams for a given dominant weight subsequence. 

\begin{remark}
For any $\underline{w}$ and $\underline{u}$ so that $\wt \underline{w} = \wt \underline{u}$, there is a distinguished choice of neutral diagram corresponding to the minimal coset representative in the symmetric group realizing the shuffle from one sequence to the other. However, we do not require that we choose particular elements as our neutral diagrams in the light ladder algorithm.
\end{remark}

%===========
\subsection{Double Ladders}\label{doubleladslabel}
\label{subsec-doublealgorithm}
%===========

We define a contravariant endofunctor $\mathbb{D}$ on the category $\DD$ by requiring that $\mathbb{D}$ fixes objects and turns diagrams upside down. Note that $\mathbb{D}^2= \id_{\DD}$, so $\mathbb{D}$ is a duality on the category.

\begin{defn}\label{upsidedownlightladder}
Let $LL_{\underline{w}, \vec{\mu}}^{\underline{v}}$ be a light ladder diagram. The associated \textbf{upside down light ladder diagram} is defined to be
\begin{equation}
\mathbb{D}(LL_{\underline{w}, \vec{\mu}}^{\underline{v}}). 
\end{equation}
\end{defn}

\begin{figure}[H]
\caption{An upside down light ladder diagram $\mathbb{D}(LL_{\blues\blues\greent\blues\greent\blues, ((1, 0), (1, 0), (-2, 1), (1, 0), (2, -1), (-1, 1))}^{\blues\greent\blues})$.}
\centering
\includegraphics[width=.20\textwidth]{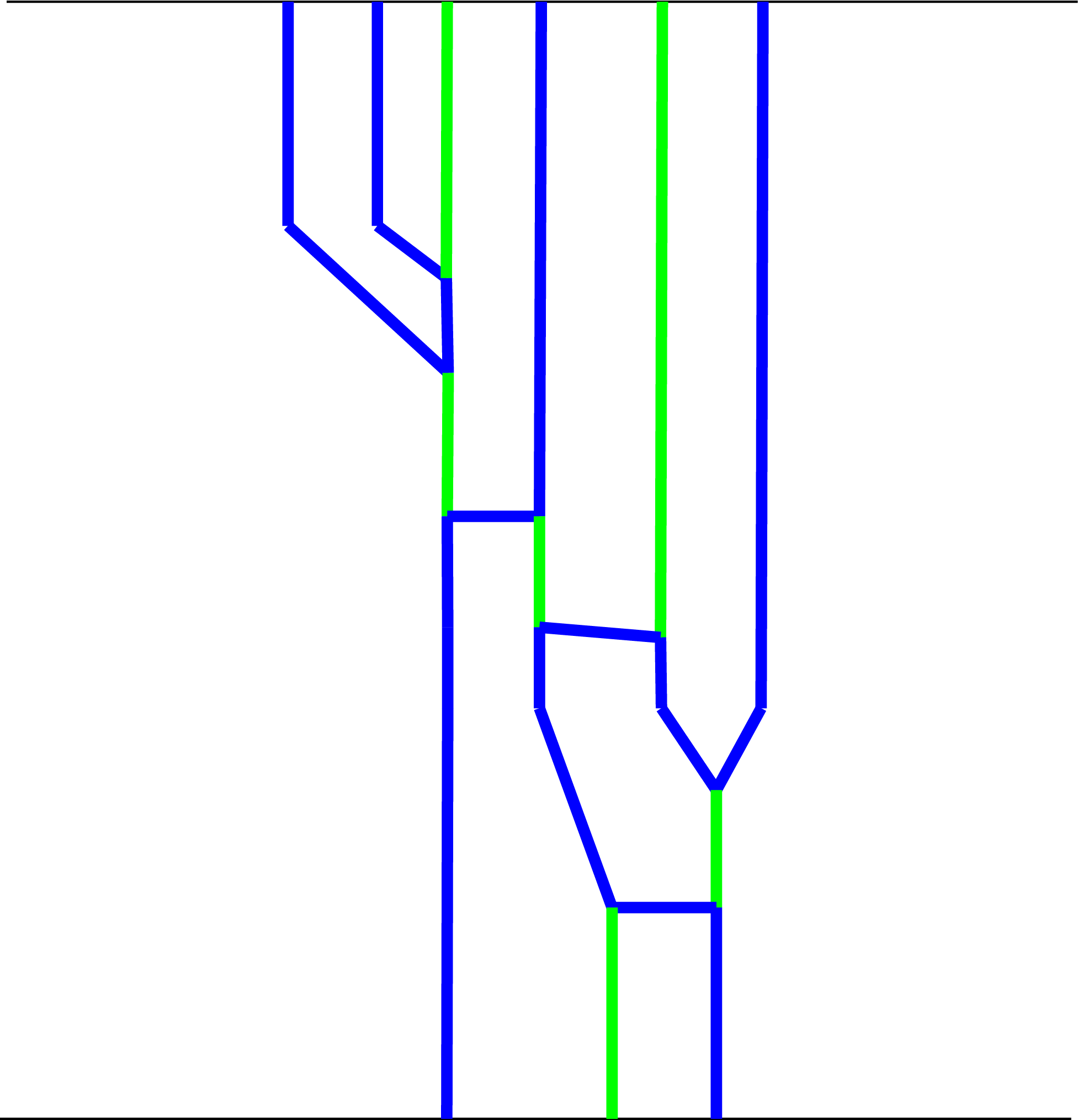}
\end{figure}

For each dominant weight $\lambda$ fix a word $\underline{x}_{\lambda}$ in the alphabet $\lbrace\blues, \greent\rbrace$ corresponding to a sequence of fundamental weights which sum to $\lambda$. For all words $\underline{w}$ and for each dominant weight subsequence $\vec{\mu}  \in E(\underline{w}, \lambda)$, we choose one light ladder diagram from $\underline{w}$ to $\underline{x}_{\lambda}$. If $\underline{w}= \underline{x}_{\lambda}$ and each $\mu_i$ is dominant, then we choose the identity diagram. From now on we denote this chosen light ladder diagram by $L_{\underline{w}, \vec{\mu}}$.  

\begin{remark}
The choice of $LL_{\underline{x}_{\lambda}, \vec{\lambda}} = \id_{\underline{x}_{\lambda}}$ when the $\lambda_i$ are all dominant is not essential for our arguments, but does ensure our construction is aligned with other conventions. For example this is required in the definition of an object adapted cellular category in \eqref{ELauda}. 
\end{remark}

\begin{defn} 
If $\underline{w}$ and $\underline{u}$ are fixed words in $\lbrace \blues, \greent\rbrace$ and $\lambda$ is a dominant weight, then for $\vec{\mu} \in E(\underline{w}, \lambda)$ and $\vec{\nu}\in E(\underline{u}, \lambda)$ we obtain a \textbf{double ladder diagram} (associated to our choices of $\underline{x}_{\lambda}$'s and our choices of light ladder diagrams)
\begin{equation}\label{LLdefn}
\mathbb{LL}_{\underline{w}, \vec{\mu}}^{\underline{u}, \vec{\nu}}= \mathbb{D}(LL_{\underline{u}, \vec{\nu}})\circ LL_{\underline{w}, \vec{\mu}}.
\end{equation}
\end{defn}

\begin{remark}
One reason for fixing an $\underline{x}_{\lambda}$ for all $\lambda$ is so the composition on the right hand side of \eqref{LLdefn} is well defined. 
\end{remark}

\begin{figure}[H]
\caption{A schematic for the double ladder diagram}
\centering
\includegraphics[width=.40\textwidth]{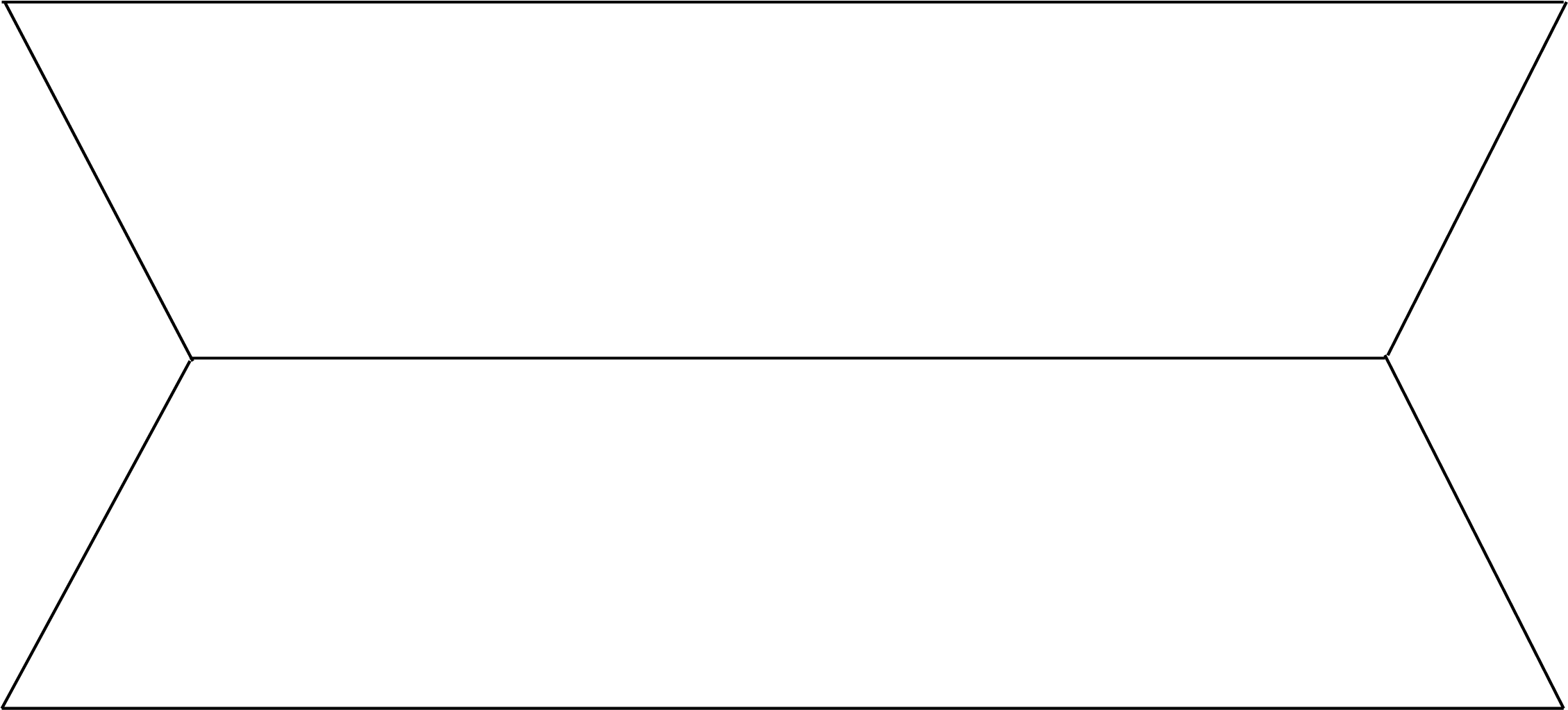}\put(-105,20){$LL_{\underline{w}, \vec{\mu}}$}\put(-105, 60){$\mathbb{D}(LL_{\underline{u}, \vec{\nu}})$}
\end{figure}

\begin{remark}
Note that light ladder diagrams ending in $\underline{x}_{\lambda}$ are double ladder diagrams, where the upside-down light ladder happens to be the identity diagram.
\end{remark}

\begin{defn}
We define the set of all double ladder diagrams from $\underline{w}$ to $\underline{u}$ factoring through $\lambda$ (associated to our choice of $\underline{x}_{\lambda}$'s and light ladder diagrams) to be
\begin{equation}
\mathbb{LL}_{\underline{w}}^{\underline{u}}(\lambda)= \big\lbrace \mathbb{LL}_{\underline{w}, \vec{\mu}}^{\underline{u}, \vec{\nu}} \ : \ \vec{\mu}\in E(\underline{w}, \lambda), \ \vec{\nu}\in E(\underline{u}, \lambda)\big\rbrace,
\end{equation}
and define the set of all double ladder diagrams from $\underline{w}$ to $\underline{u}$ (associated to our choice of $\underline{x}_{\lambda}$'s and light ladder diagrams) to be
\begin{equation}
\mathbb{LL}_{\underline{w}}^{\underline{u}} = \bigcup_{\lambda\in X_+}\mathbb{LL}_{\underline{w}}^{\underline{u}}(\lambda)
\end{equation}
\end{defn}

\begin{remark}
Anytime we write $\mathbb{LL}_{\underline{w}, \vec{\mu}}^{\underline{u}, \vec{\nu}}$ or $\mathbb{LL}_{\underline{w}}^{\underline{u}}$, we have already fixed choices of $\underline{x}_{\lambda}$'s and choices of light ladder diagrams. The notation does not account for these choices, but we will not be comparing double ladders for different choices so the notation should not lead to confusion.
\end{remark}

%===========
\subsection{Relating Non-Elliptic Webs to Double Ladders}
\label{subsec-reduction}
%===========

Our next goal is to define an evaluation functor from $\DD$ to the category $\Fund(\mathfrak{sp}_4)$, and then to prove that the functor is an equivalence. That the functor is an equivalence will follow from showing that double ladder diagrams span the category $\DD$, and map to a set of linearly independent morphisms in $\Fund(\mathfrak{sp}_4)$. This approach is modeled on the work on type $A$ webs in \cite{elias2015light}, where most of the work goes into showing that double ladder diagrams span the diagrammatic category. Checking linear independence is comparatively easy once you know the functor explicitly. But for $\DD$, the extra work to show double ladders span can be circumvented by bootstrapping known results about $B_2$ webs which we recall below. 

Kuperberg's paper \cite[pp. 14-15]{Kupe} introduces a tetravalent vertex in the $B_2$ web category which can be used to remove all internal double edges. Let $\textbf{B}$ be the set of $B_2$ diagrams with no internal double edges and with no faces having one, two, or three adjacent edges. These diagrams are called non-elliptic in \cite{Kupe}. There are local relations in the $B_2$ category (now including the tetravalent vertex) which can be used to reduce triangular faces, bigons, monogons, and circles to sums of diagrams with fewer crossings (i.e. $\textbf{B}$ is the set of irreducible webs with respect to the relations). It follows that the set $\textbf{B}$ spans the $B_2$ category over $\mathbb{Z}[q, q^{-1}]$. Let $\textbf{B}_{\underline{w}}$ be the set of diagrams in $\textbf{B}$ with $\underline{w}$ on the boundary. One of the main results of \cite{Kupe} is that 
\begin{equation}\label{equinumeration}
\# \textbf{B}_{\underline{w}} = \dim \VV(\underline{w})^{\mathfrak{sp}_4(\mathbb{C})}.
\end{equation}

If we work in the $\mathcal{A}$-linear category $\DD$, there is an analogous $90$ degree rotation invariant morphism, which we will call the tetravalent vertex, in $\End_{\DD}(\blues\blues)$.
\begin{figure}[H]
\centering
\includegraphics[width=0.1\textwidth]{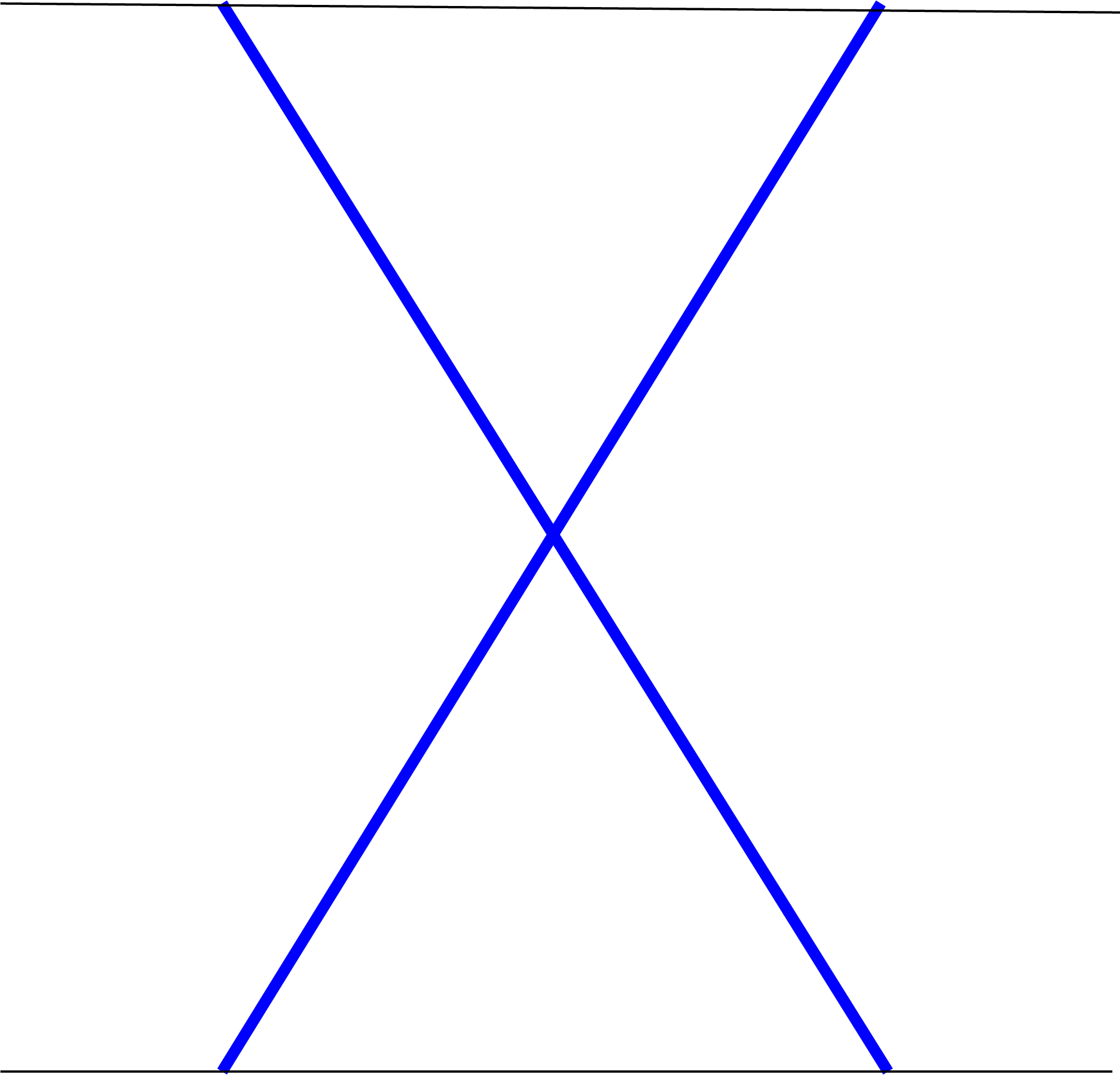}\put(5, 16){$:=$} \ \ \ \ \ \ \ \ \ \includegraphics[width=0.1\textwidth]{figs1/I=HI}\put(5, 16){$-\dfrac{1}{[2]_q}$} \ \ \ \ \ \ \ \ \ \ \ \ \includegraphics[width=0.1\textwidth]{figs1/I=Hcupcap}
\end{figure}
\noindent Since $[2]_q$ is invertible in our ground ring, we can use this tetravalent vertex to remove all internal green label edges in any diagram in $\DD$. The tetravalent vertex satisfies the following relations in $\DD$ 
\begin{figure}[H]
\centering
\includegraphics[width=0.1\textwidth]{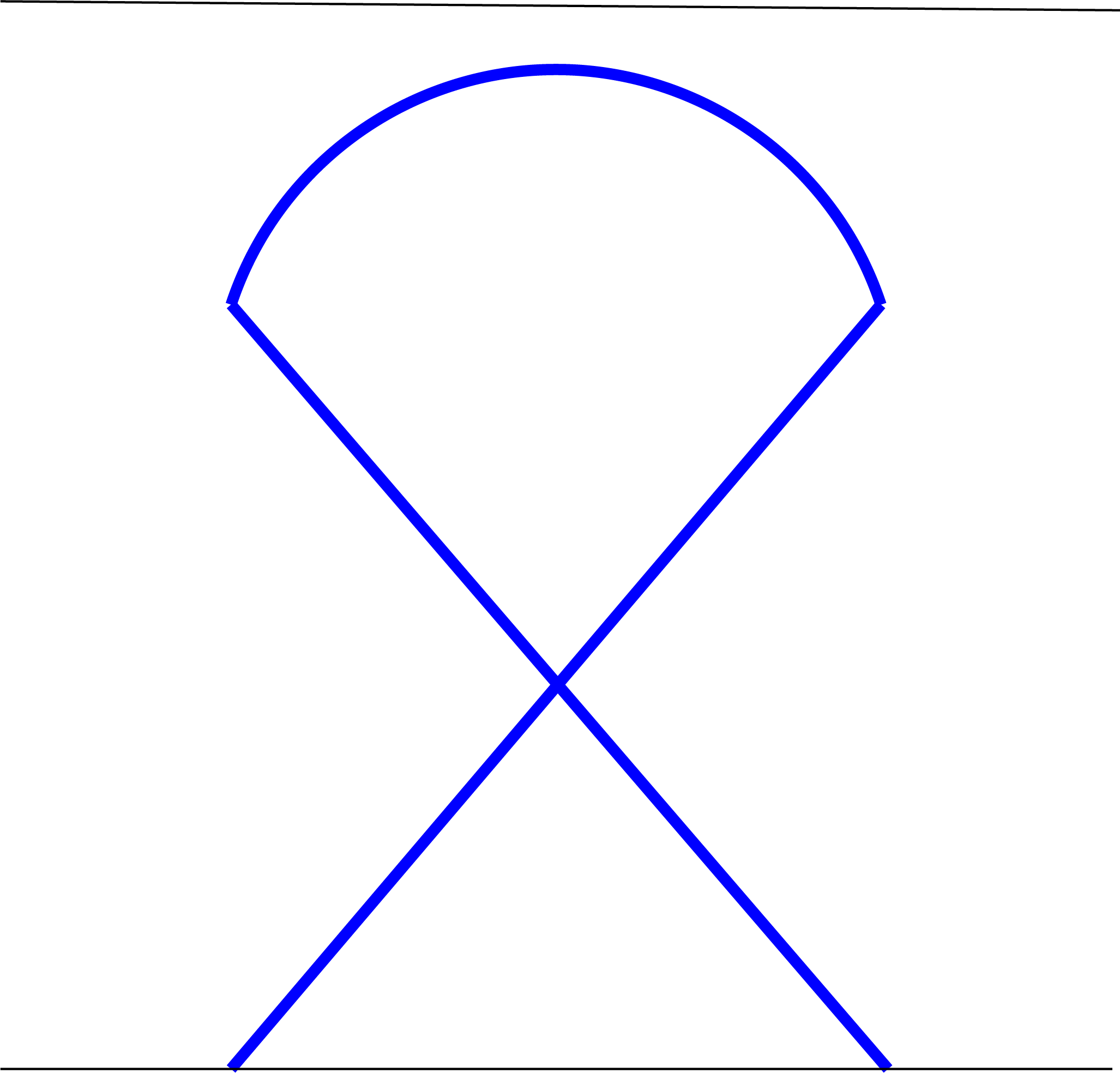}\put(5, 16){$=\dfrac{[6]_q}{[3]_q}$} \ \ \ \ \ \ \ \ \ \ \ \ \ \ \includegraphics[width=0.1\textwidth]{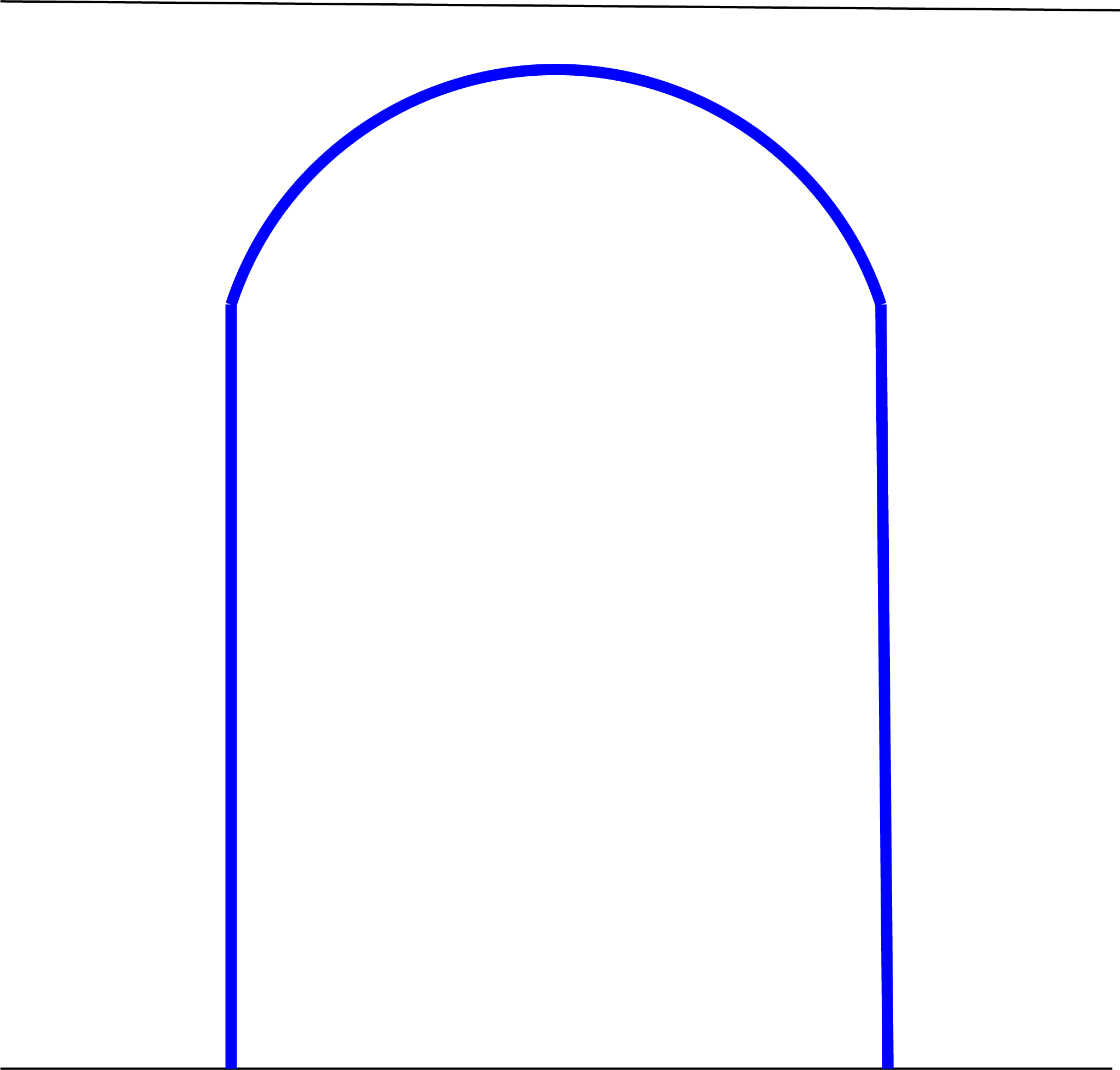}
\end{figure}
\begin{figure}[H]
\centering
\includegraphics[width=0.1\textwidth]{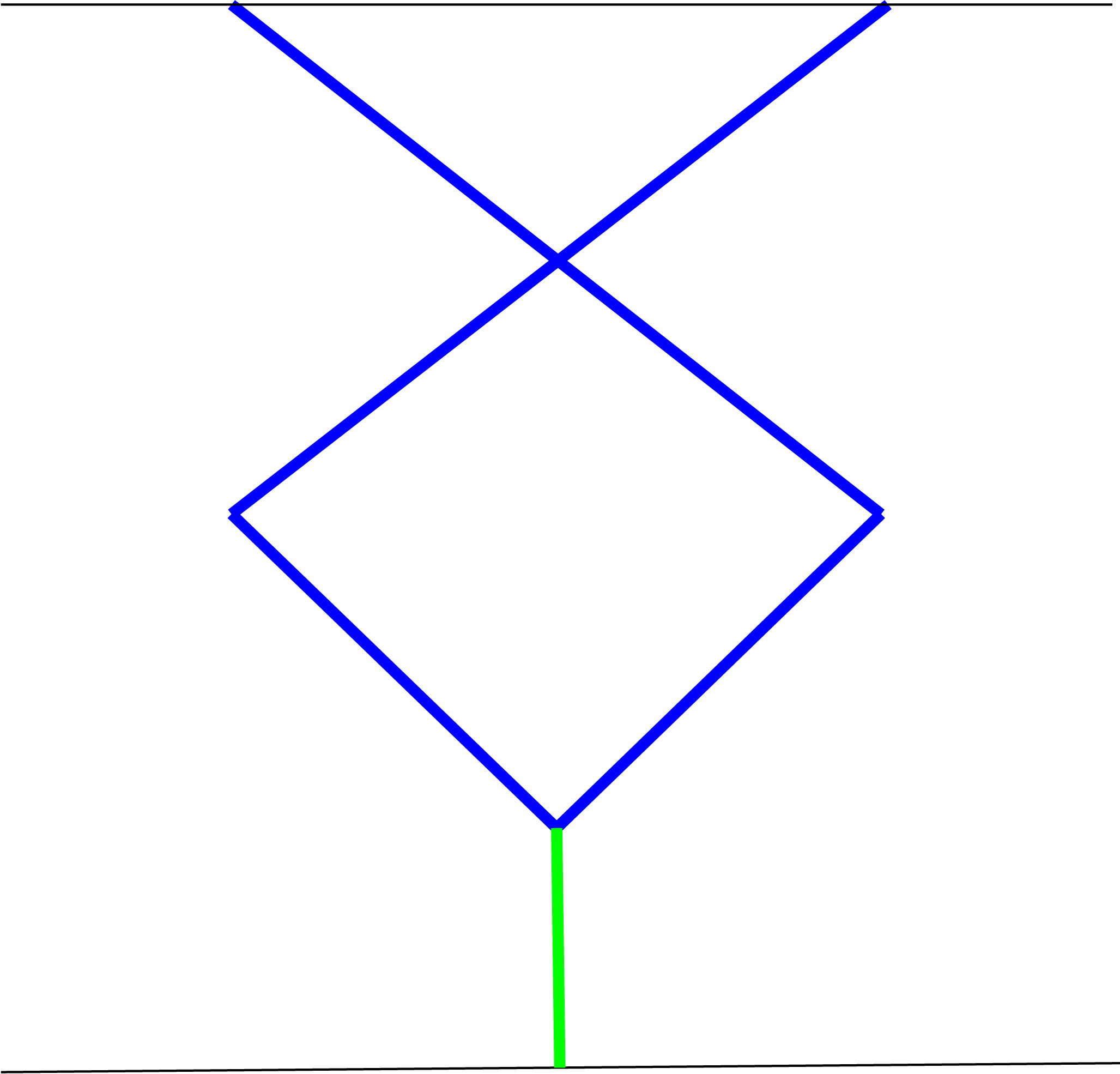}\put(5, 16){$=-[2]_q$} \ \ \ \ \ \ \ \ \ \ \  \ \ \ \includegraphics[width=0.1\textwidth]{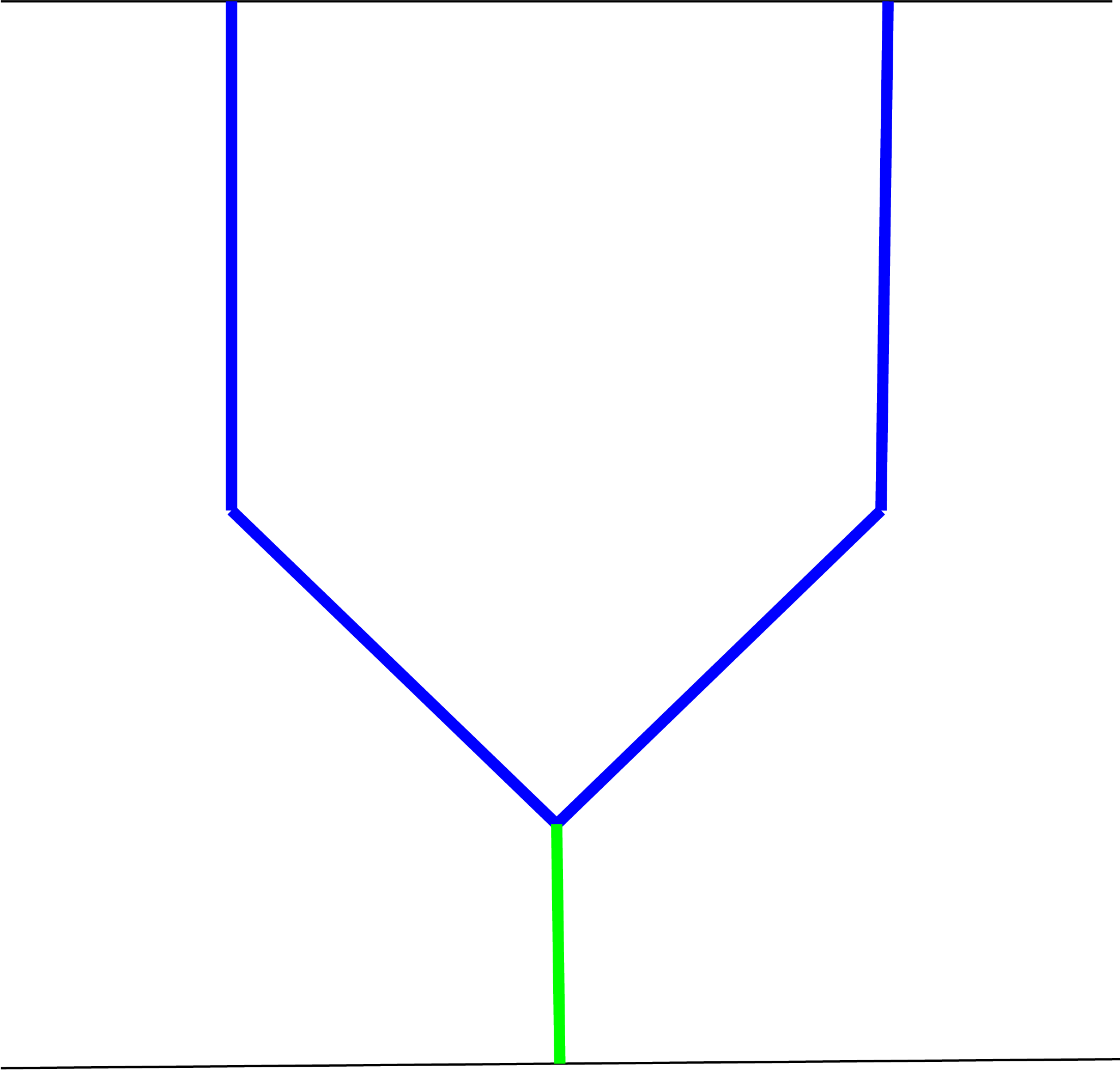}
\end{figure}
\begin{figure}[H]
\centering
\includegraphics[width=0.1\textwidth]{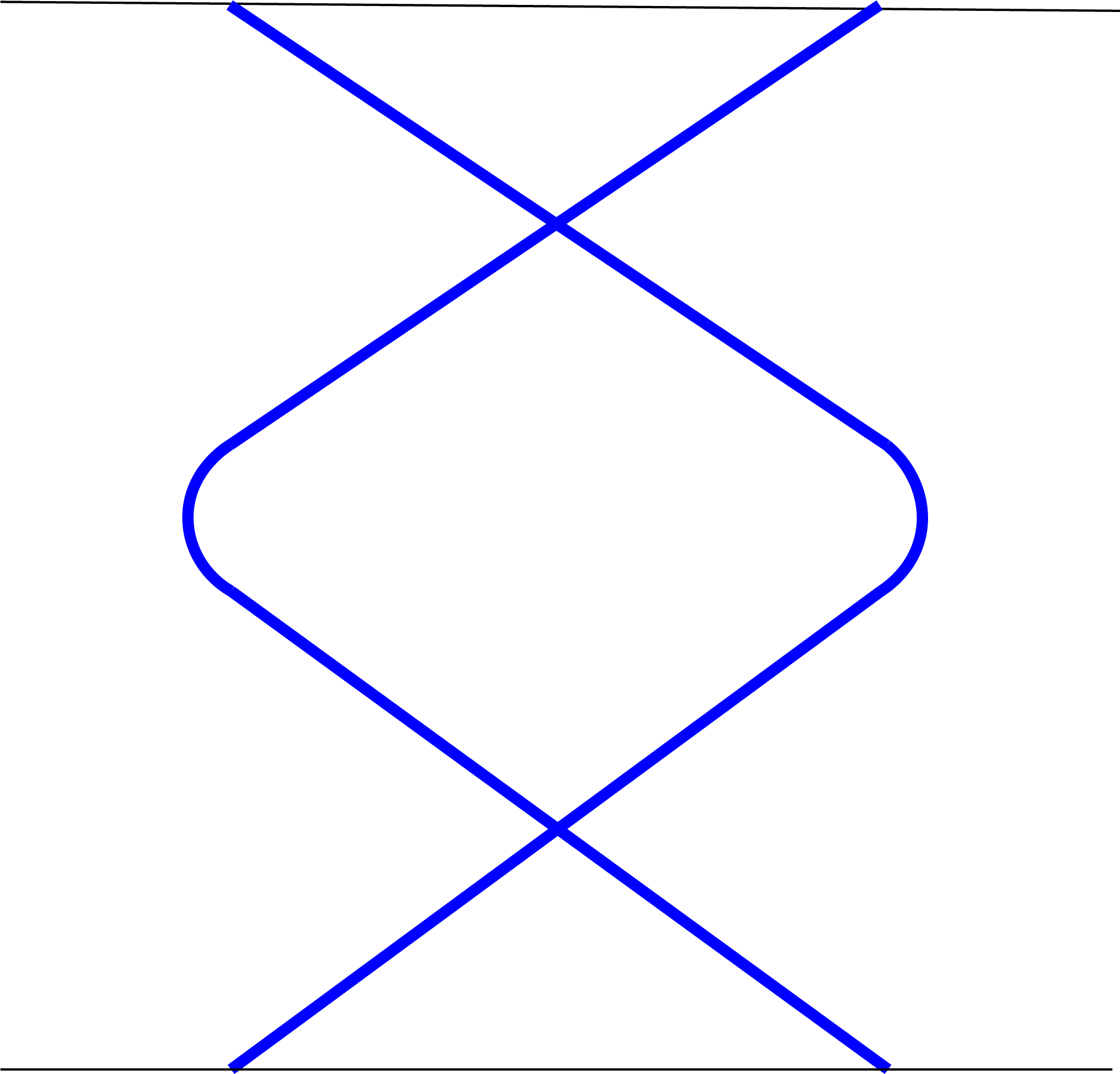}\put(5, 16){$=-[2]_q$} \ \ \ \ \ \ \ \ \ \ \ \ \ \ \includegraphics[width=0.1\textwidth]{figs1/X}\put(5, 16){$-\dfrac{[4]_q}{[2]_q}$} \ \ \ \ \ \ \ \ \ \ \ \ \includegraphics[width=0.1\textwidth]{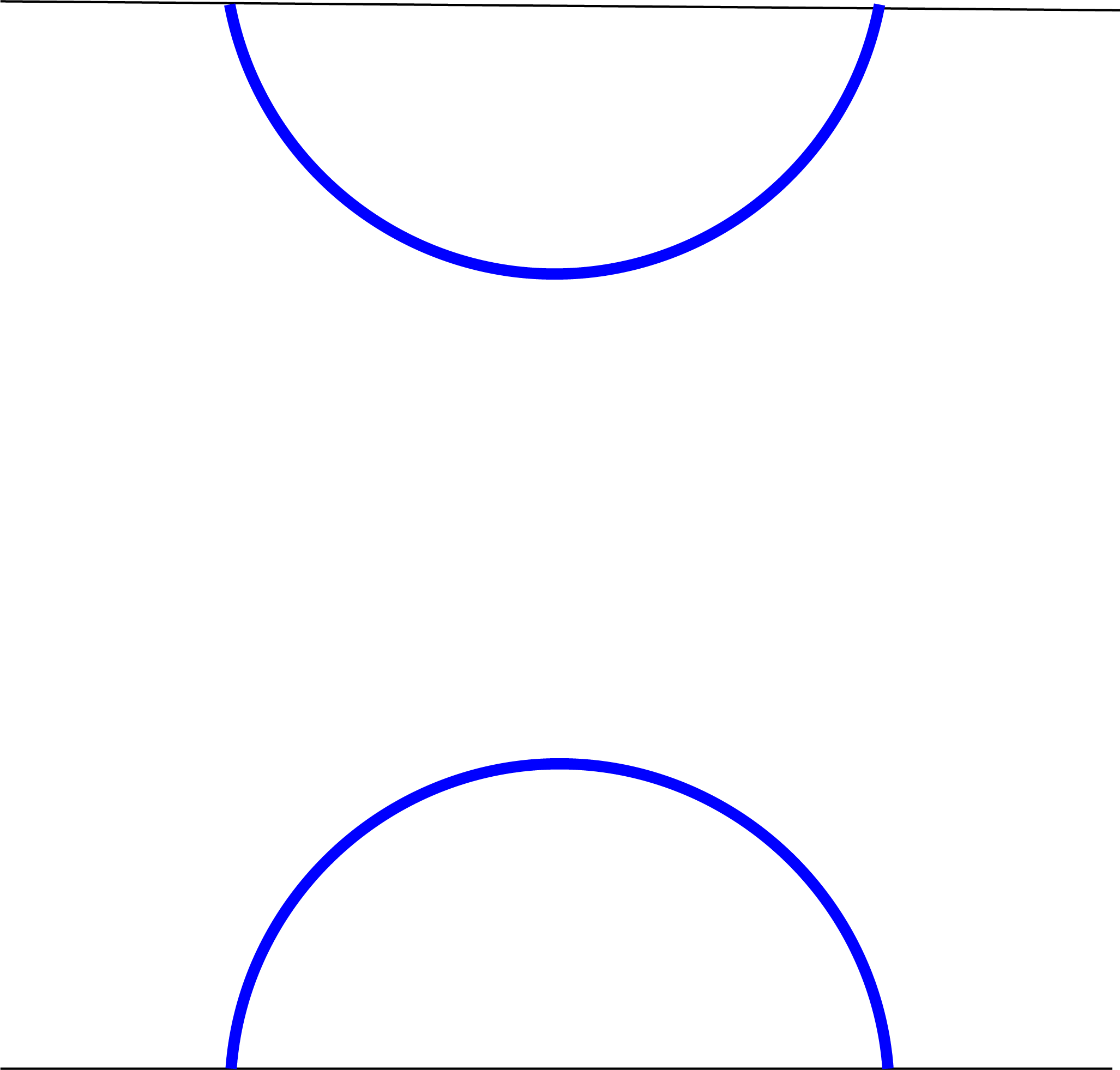}
\end{figure}
\begin{figure}[H]
\centering
\includegraphics[width=0.1\textwidth]{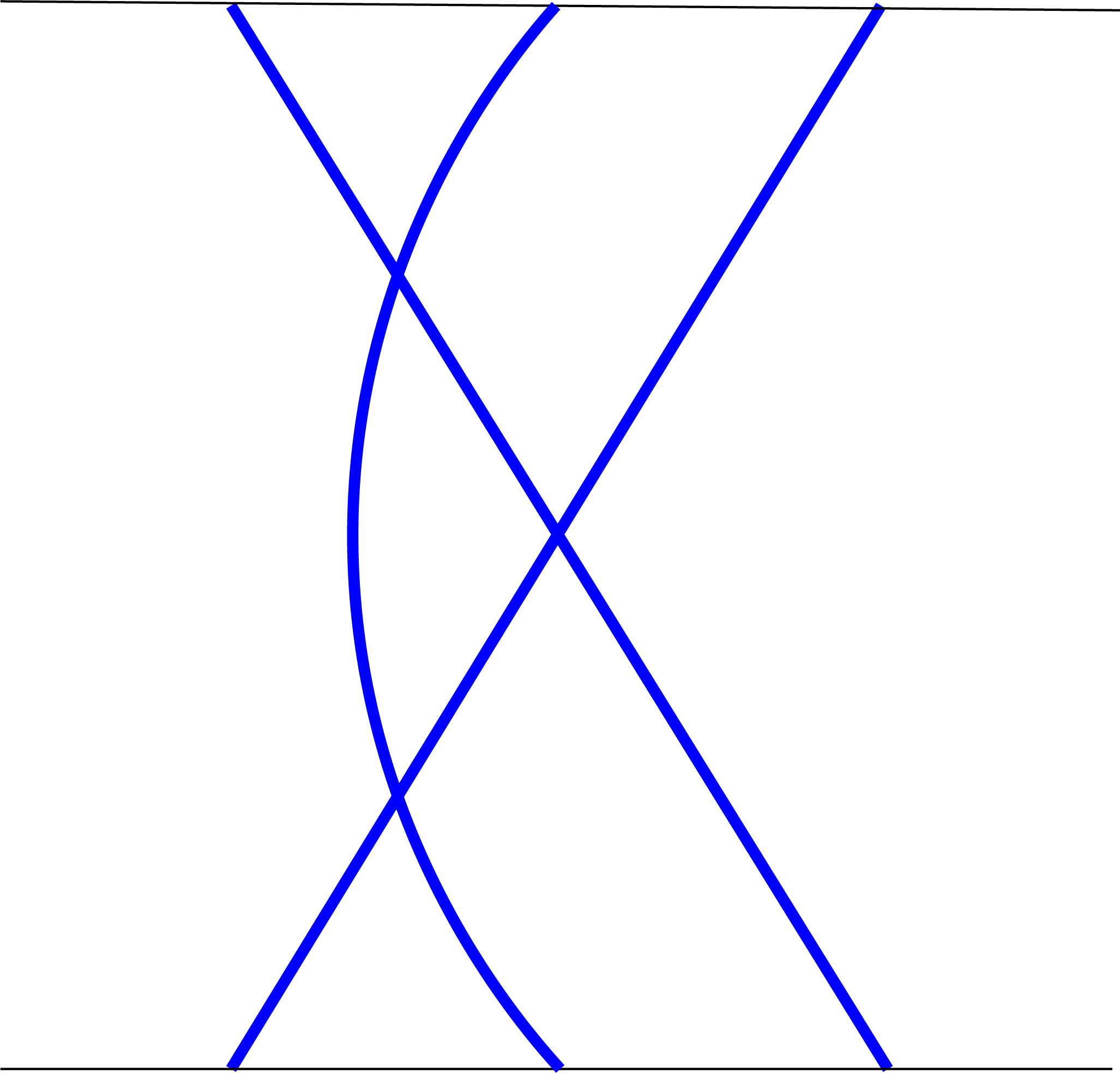}\put(5, 16){$=$} \ \ \ \ \ \ \ \ \ \ \ \ \includegraphics[width=0.1\textwidth]{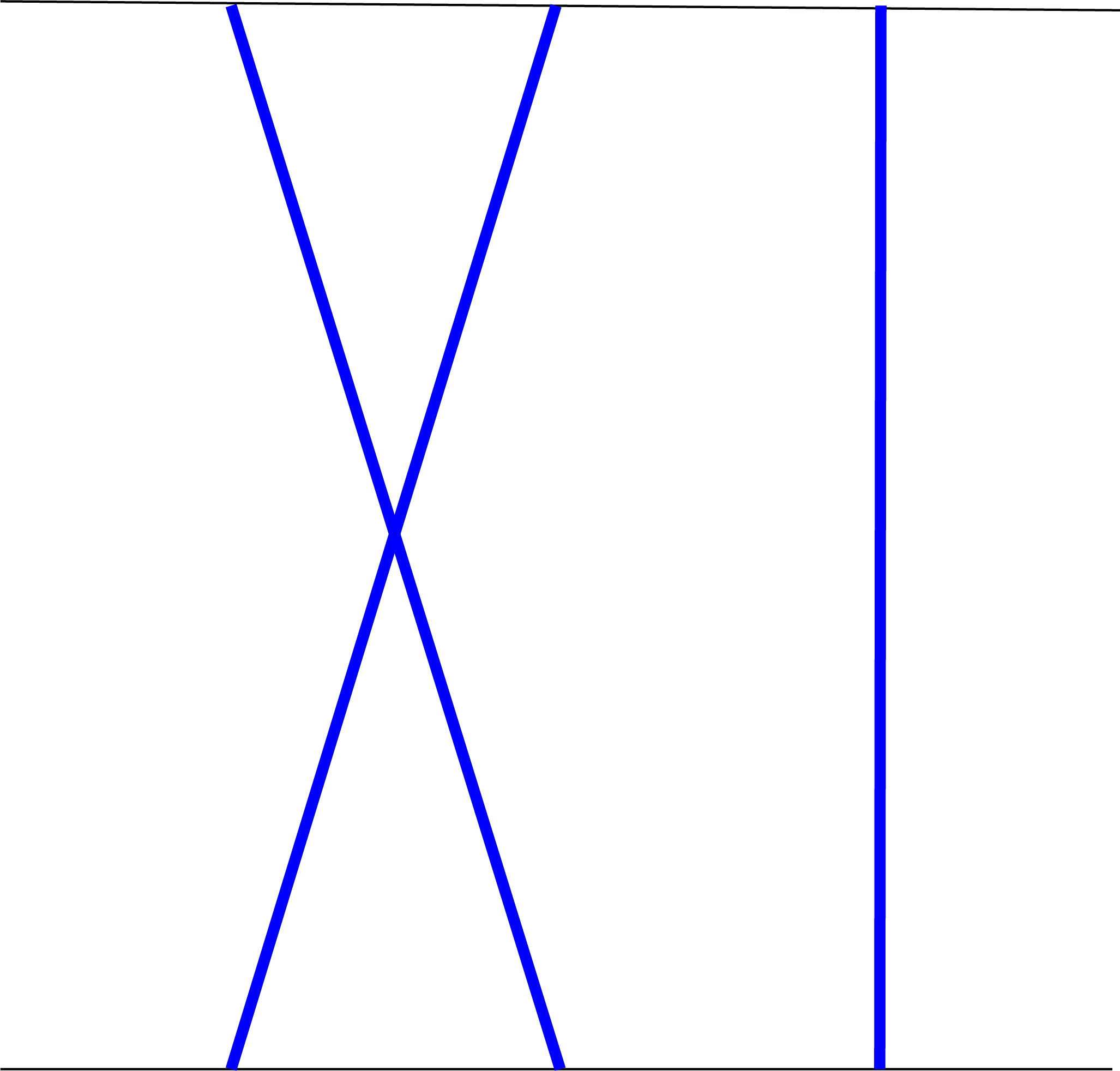}\put(5, 16){$+$}\ \ \ \ \ \ \ \ \ \includegraphics[width=0.1\textwidth]{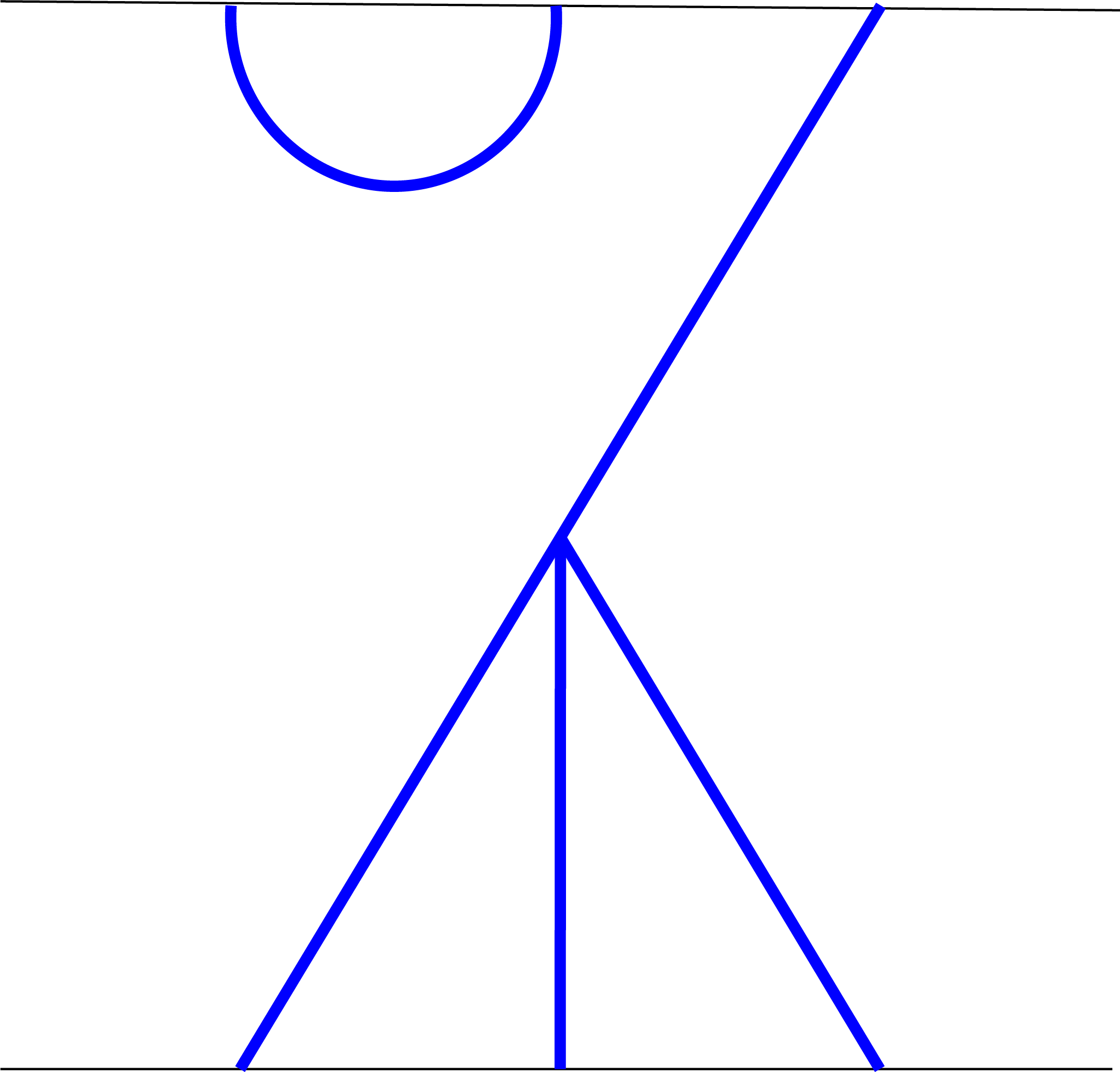}\put(5, 16){$+$} \ \ \ \ \ \ \ \ \ \ \ \ \includegraphics[width=0.1\textwidth]{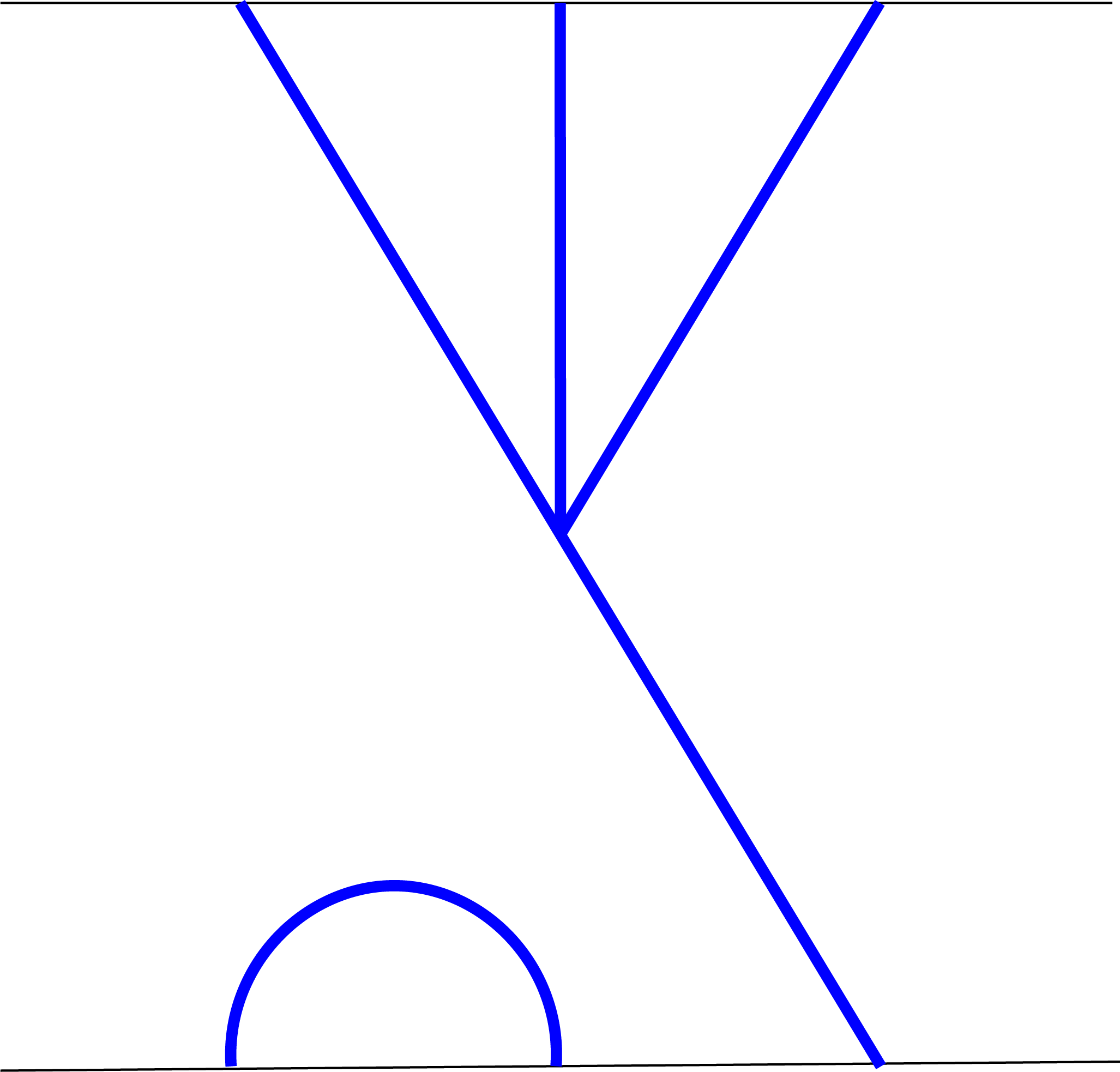}\put(5, 16){$+[2]_q$} \ \ \ \ \ \ \ \ \ \ \ \ 
\includegraphics[width=0.1\textwidth]{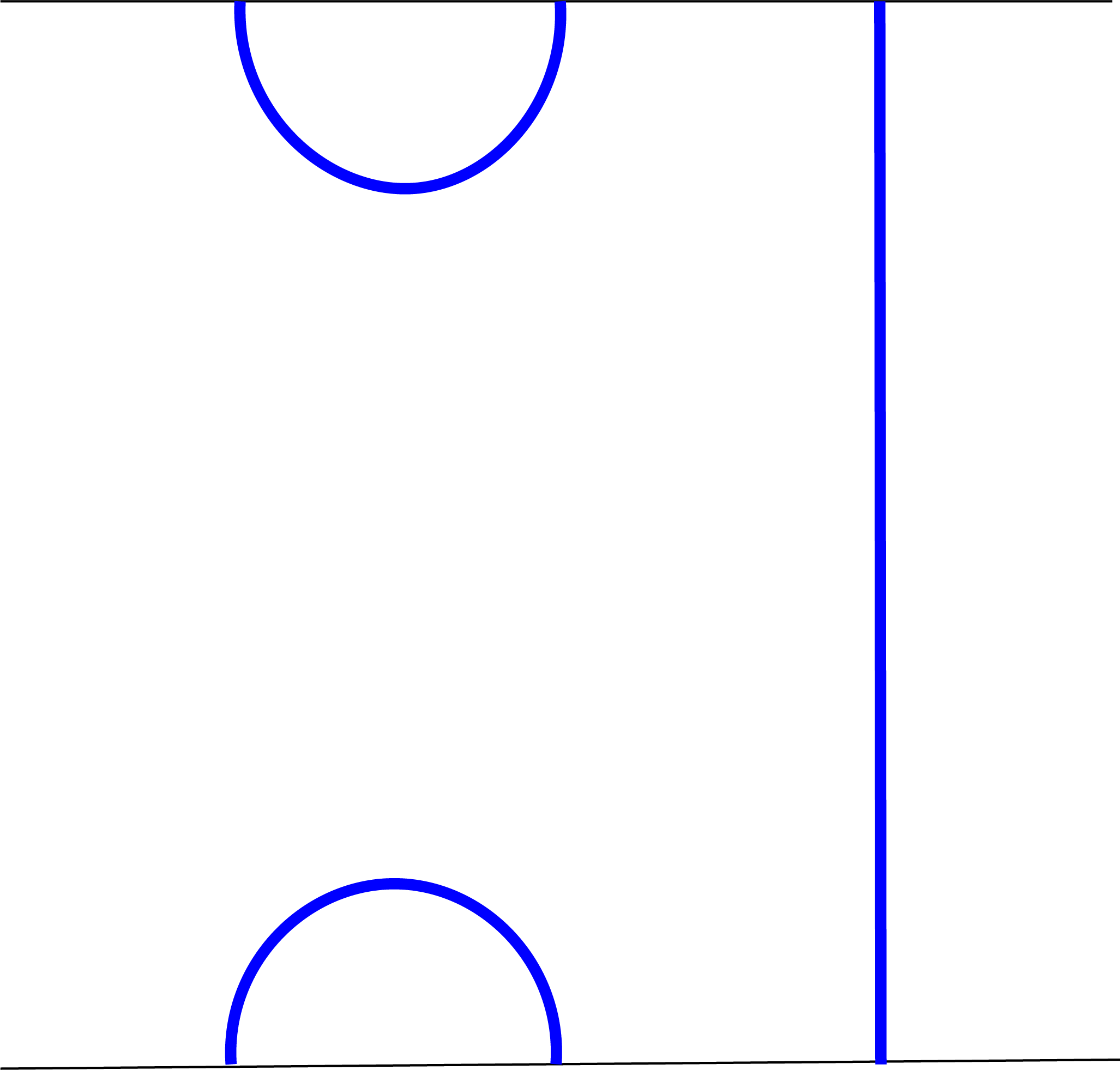}
\end{figure}
\begin{figure}[H]
\centering
\includegraphics[width=0.1\textwidth]{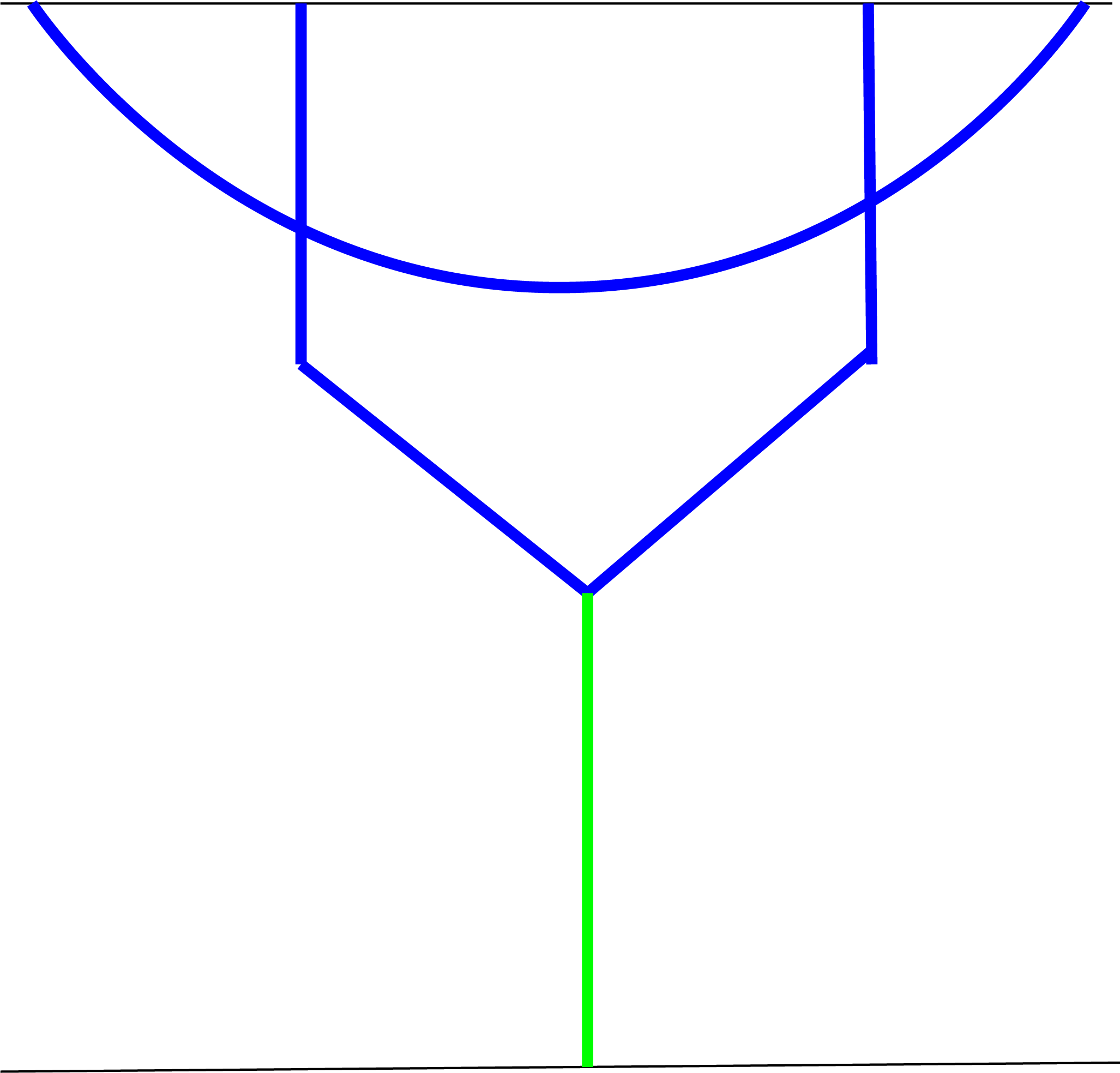}\put(5, 16){$=$} \ \ \ \ \ \ \ \ \ \ \ \ \ \ \includegraphics[width=0.1\textwidth]{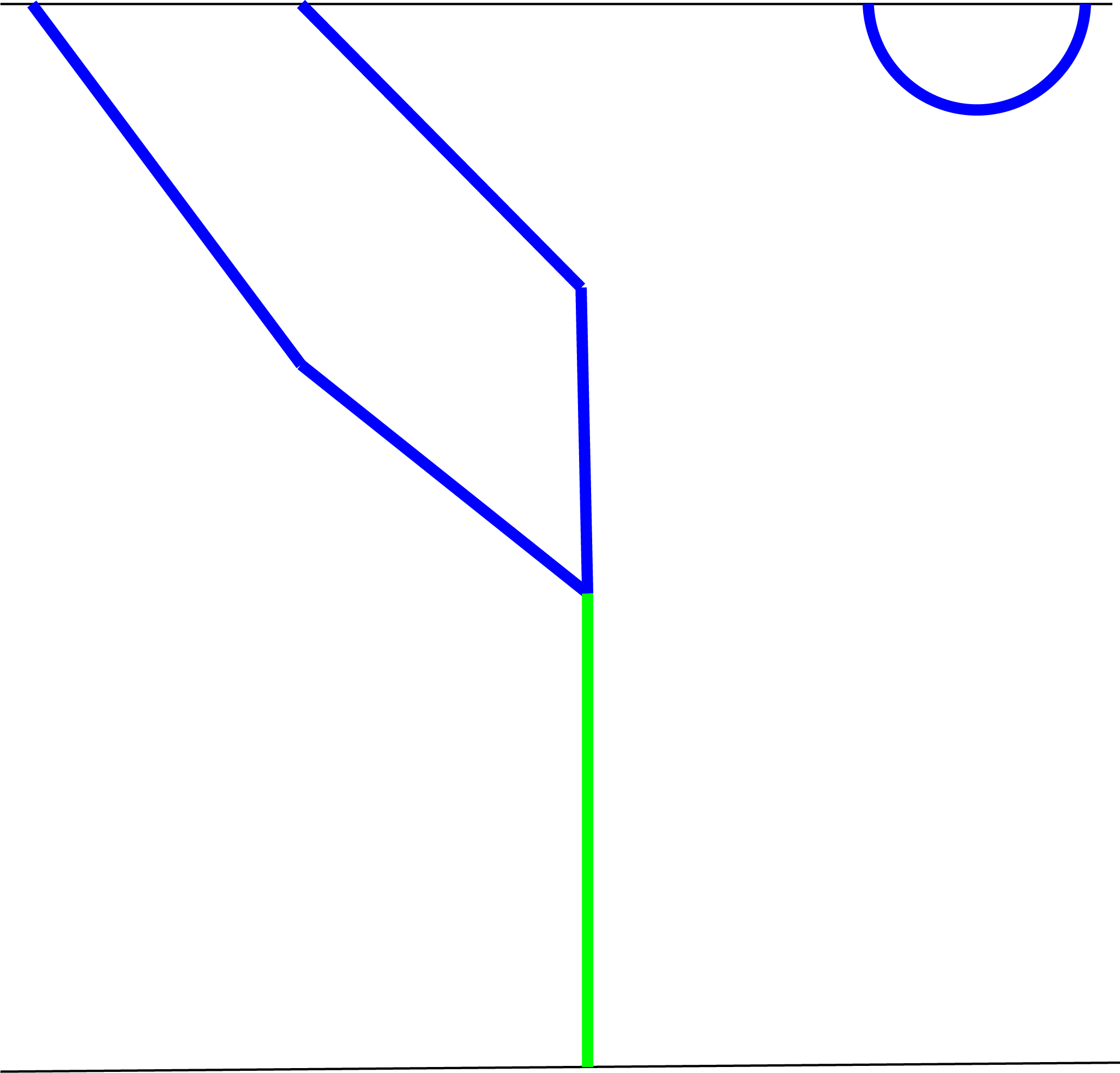}\put(5, 16){$+$} \ \ \ \ \ \ \ \ \ \ \ \ \includegraphics[width=0.1\textwidth]{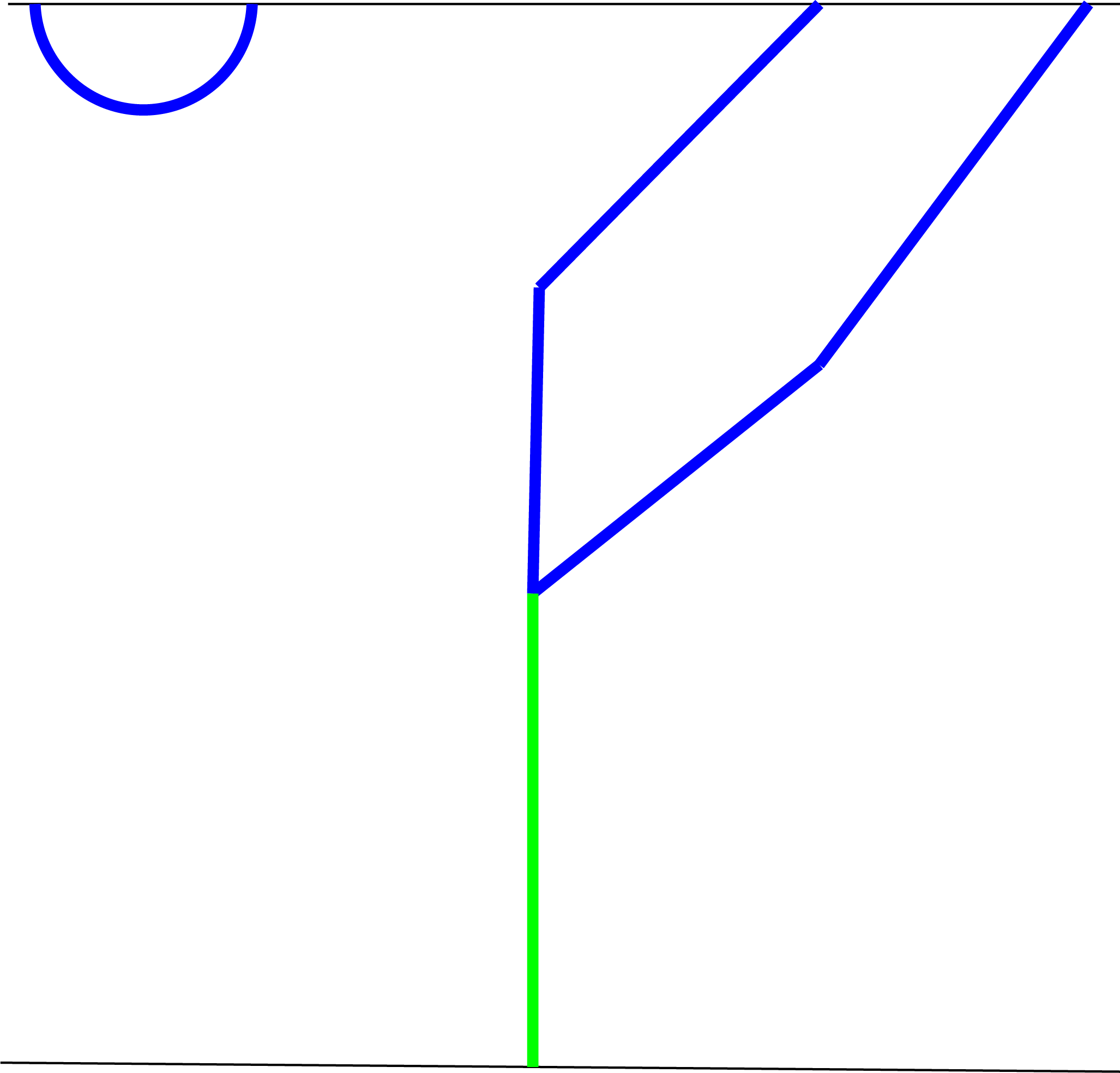}
\end{figure}
\begin{figure}[H]
\centering
\includegraphics[width=0.1\textwidth]{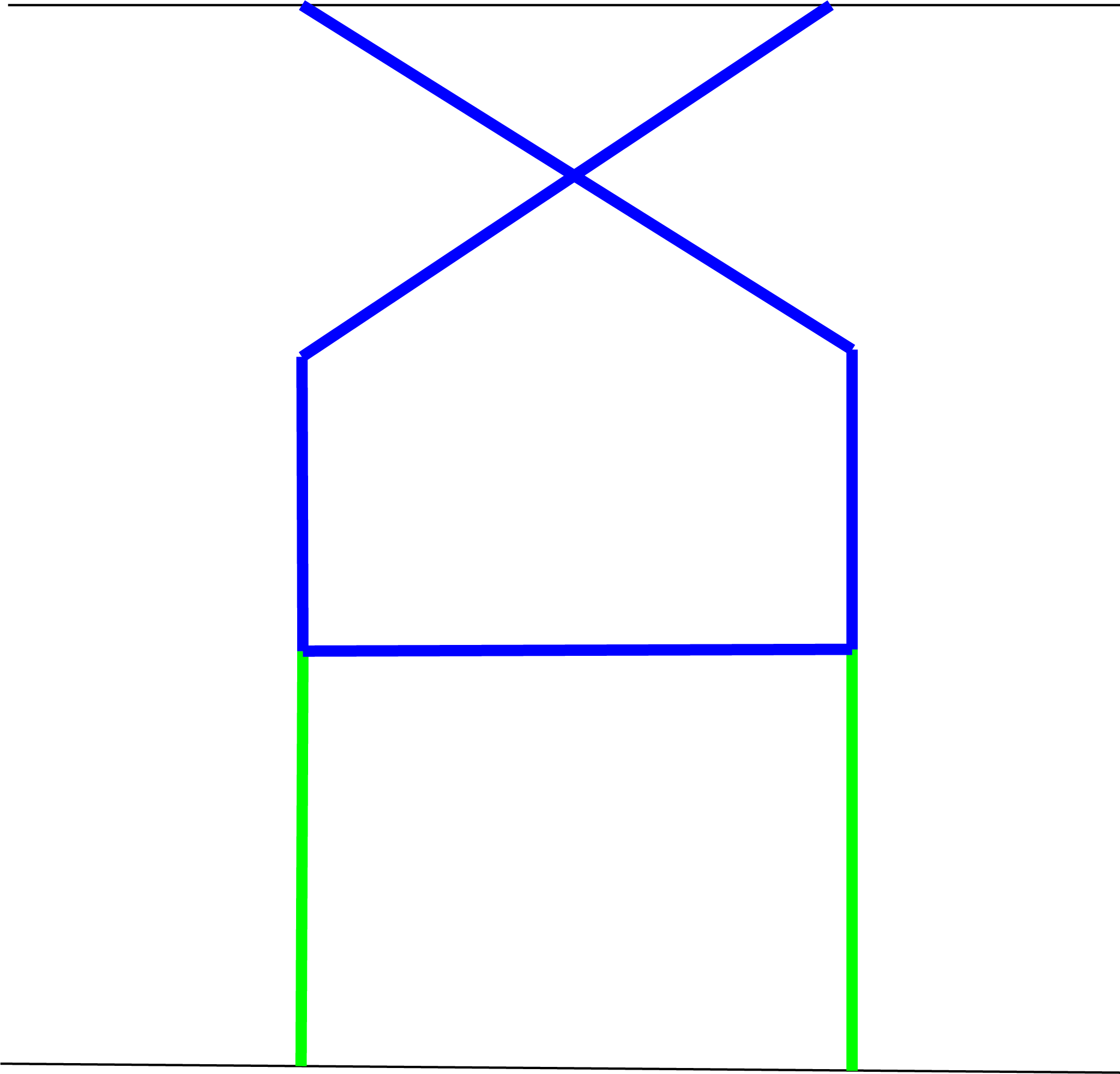}\put(5, 16){$=$} \ \ \ \ \ \ \ \ \ \ \ \ \ \ \includegraphics[width=0.1\textwidth]{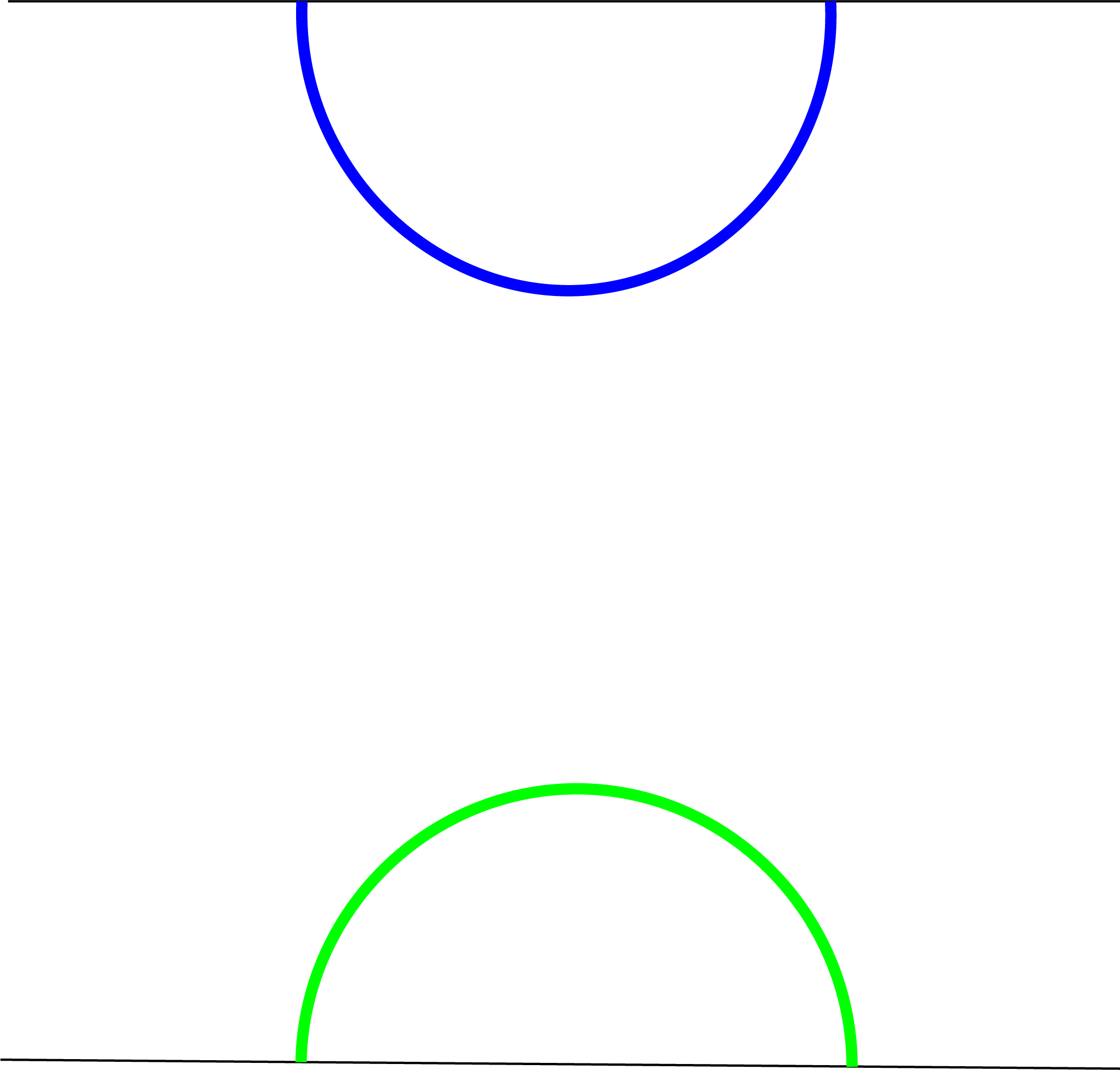}
\end{figure}

\begin{remark}
Due to the identity $[2n]_q/[n]_q = [n+1]_q- [n-1]_q$, the coefficients in these relations all lie in the ring $\mathbb{Z}[q, q^{-1}]$. 
\end{remark}

\begin{defn}
A \textbf{face} of a diagram in $\DD$ is a simply connected component of the complement of the diagram, which does not touch the boundary. 
\end{defn}

\begin{defn}
A \textbf{non-elliptic diagram} in $\DD$ is a diagram such that all faces have more than three sides (i.e a diagram with no triangular faces, bigons, monogons, or circles).
\end{defn}

\begin{defn}
An \textbf{internal $\greent$ edge} of a diagram in $\DD$ is a $\greent$ edge in the diagram which does not connect to the boundary. 
\end{defn}

\begin{figure}[H]
\caption{An example of an elliptic web with internal $\greent$ edges in $\DD$. The only face is the interior of the $\greent$ circle.}
\centering
\includegraphics[width=.1\textwidth]{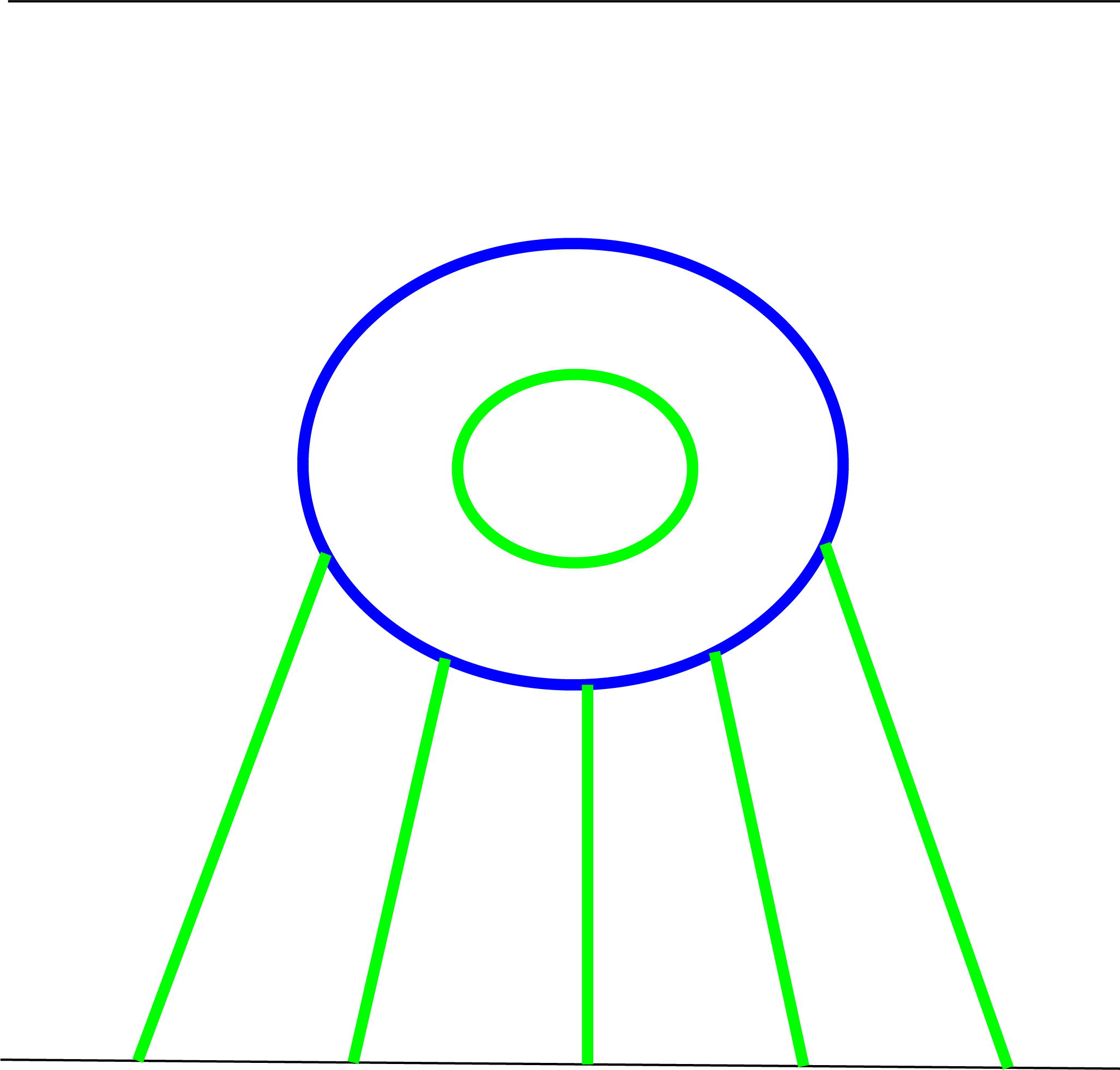}
\end{figure}

\begin{figure}[H]\label{nonellipticweb}
\caption{An example of a non-elliptic web with no internal $\greent$ edges in $\DD$. There is only one face and it has five sides.}
\centering
\includegraphics[width=.1\textwidth]{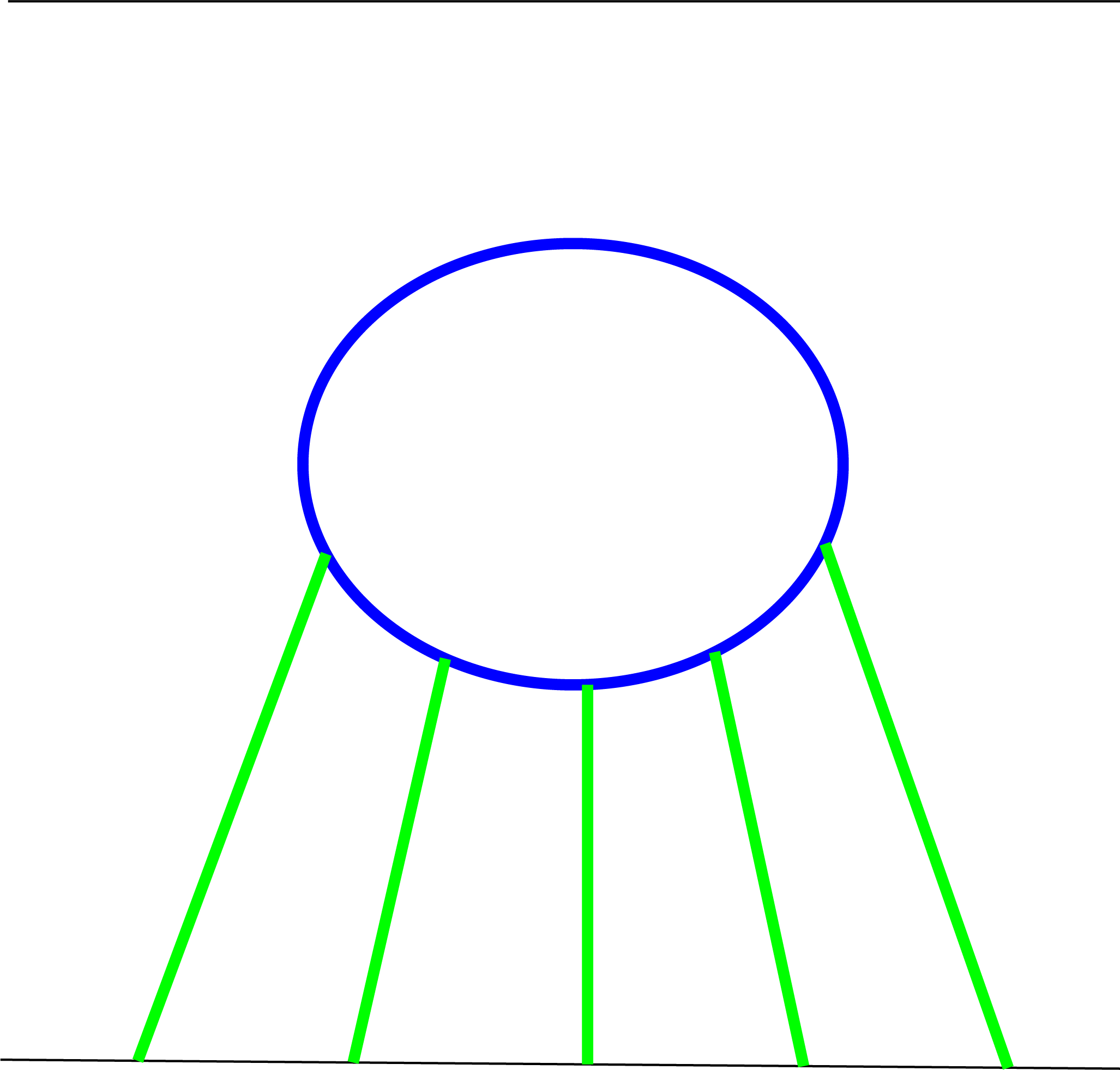}
\end{figure}

\begin{defn}
The set $\textbf{D}$ is the collection of all non-elliptic diagrams in $\mathcal{D}_{C_2}$ with no internal $\greent$ edges, and the set $\textbf{D}_{\underline{w}}^{\underline{u}}$ is the set of diagrams in $\textbf{D}\cap \Hom_{\mathcal{D}_{C_2}}(\underline{w}, \underline{u})$. 
\end{defn}

\begin{lemma}
The set $\textbf{D}$ spans $\DD$ over $\mathcal{A}$.
\end{lemma}
\begin{proof}
Let $D$ be an arbitrary diagram in $\DD$. We will argue that $D$ is a linear combination of non-elliptic webs with no internal $\greent$ edges. If a $\greent$ edge does not connect to a trivalent vertex, then you can use the bigon relation to introduce one. Thus, every $\greent$ edge either connects to the boundary of $D$, or connects two trivalent vertices. Using the tetravalent vertex to remove all pairs of trivalent vertices, we can rewrite $D$ as a linear combination of diagrams with no internal $\greent$ edges. Thus, we may assume that $D$ is a diagram with no internal $\greent$ edges. Using the defining relations in $\DD$ along with the tetravalent relations, we can remove all faces with less than four edges. 
\end{proof}

\begin{remark}
In order to introduce a trivalent vertex, we used the bigon relation backwards, which required $[2]_q^{-1}\in \mathcal{A}$. 
\end{remark}

\begin{lemma}\label{hominequality}
Let $\ak$ be a field and let $q\in \ak^{\times}$ be such that $q+ q^{-1} \ne 0$. Then
\begin{equation}\label{kupinequality}
\dim \Hom_{\mathbb{\ak}\ot\DD}(\underline{w}, \underline{u}) \le \dim \Hom_{\mathfrak{sp}_4(\mathbb{C})}(\VV(\underline{w}), \VV(\underline{u})).
\end{equation}. 
\end{lemma}

\begin{proof}
There is an obvious bijection between the set $\textbf{B}$ and the set $\textbf{D}$. The result then follows from \eqref{equinumeration}.
\end{proof}

\begin{remark}
We sketch a more direct argument to deduce the inequality \eqref{kupinequality}. The dimension of the $\mathfrak{sp}_4(\mathbb{C})$ invariants in $\VV(\blues)^{\ot 2n}$ is known to be equal to the number of matchings of $2n$ points on the boundary of a disc so that there is no $6$-point star in the matching \cite{1990IMA....19..191S}\cite[8.4]{Kupe}. One can argue that the local condition of being non-elliptic implies the global condition of having no six point star. Then, noting that non-elliptic diagrams have a unique representative up to isotopy (there are no potential Reidemeister moves), it follows that there is a bijection between non-elliptic diagrams and matchings without a $6$-point star. This proves that the inequality \eqref{kupinequality} holds when $\underline{w} = \blues ^{\ot a}$ and $\underline{u}= \blues^{\ot b}$ for some $a, b\in \mathbb{Z}_{\ge 0}$. Since $\greent$ is a direct summand of $\blues\blues$ it follows that \eqref{kupinequality} holds for any words $\underline{w}$ and $\underline{u}$ in the alphabet $\lbrace \blues, \greent \rbrace$. 
\end{remark}

We have defined a set $\mathbb{LL}_{\underline{w}}^{\underline{u}}$ of double ladders in $\DD$. It follows from the construction of $\mathbb{LL}_{\underline{w}}^{\underline{u}}$ and \eqref{LLsdimhom} that
\begin{equation}
\#\mathbb{LL}_{\underline{w}}^{\underline{u}} = \sum_{\lambda\in X_+}\#E(\underline{w}, \lambda)\#E(\underline{u}, \lambda)= \dim \Hom_{\mathfrak{sp}_4(\mathbb{C})}(\VV(\underline{w}), \VV(\underline{u})).
\end{equation}
We want to show linear independence of the set of double ladders, or equivalently that the inequality of dimensions in \eqref{kupinequality} is in fact an equality, for a general choice of base ring $\ak$. To this end we will define an evaluation functor from the diagrammatic category $\DDk$ to the representation theoretic category $\Fund(\ak\ot U_q^{\mathcal{A}}(\mathfrak{sp}_4))$, and interpret the image of the evaluation functor in terms of tilting modules. If we can show that the image of the double ladder diagrams under the evaluation functor is a linearly independent set, then \eqref{hominequality} will imply that the double ladder diagrams must be linearly independent in $\DD$. This implies that the inequality in \eqref{kupinequality} is an equality, and it follows that the evaluation functor maps bases to bases, so is fully faithful. 

\begin{remark}
Since $\textbf{D}$ spans $\DDk$ and is in a non-canonical bijection with the set of double ladder diagrams (for fixed choices of $\underline{x}_{\lambda}$ and fixed choices of light ladders), linear independence of the double ladder diagrams over $\ak$ implies that both sets are bases. 

Note that double ladders have many internal $\greent$ label edges while the diagrams in $\textbf{D}$ will have none. On the other hand, sometimes the double ladder diagrams will be non-elliptic webs with no internal $\greent$ edges. A good exercise for the reader is to rewrite the diagram in figure \eqref{nonellipticweb} as a double ladder diagram. A hint is that a double ladder diagram in $\Hom_{\DD}(\greent^{\ot 5}, \emptyset)$ will just be a light ladder diagram $LL_{\greent^{\ot 5}, ?}^{\emptyset}$. 
\end{remark}

%%%%%%%%%%%%%%%%%%%%%%%%%
\section{The Evaluation Functor and Tilting Modules}
\label{sec-evaltilt}
%%%%%%%%%%%%%%%%%%%%%%%%%

%===========
\subsection{Defining the Evaluation Functor on Objects}
\label{subsec-objects}
%===========

We are now going to be more precise about what representation category associated to $\mathfrak{sp}_4$ we are considering. The discussion below is well-known, but we reproduce it here to help the reader follow certain calculations which come later. 

Our main reference for quantum groups is Jantzen's book \cite{JantzenQgps}. Recall that $\mathfrak{sp}_4(\mathbb{C})$ gives rise to a root system $\Phi$ and a Weyl group $W$. We choose simple roots $\Delta = \lbrace \alpha_s = \epsilon_1- \epsilon_2, \alpha_t = 2\epsilon_2\rbrace$. There is a unique $W$ invariant symmetric form $(-,-)$ on the root lattice $\mathbb{Z}\Phi$ such that the short roots pair with themselves to be $2$. This is the form $(\epsilon_i, \epsilon_j)= \delta_{ij}$, restricted to the root lattice. For $\alpha \in \Phi$ we define the coroot $\alpha^{\vee} = 2\alpha/(\alpha, \alpha)$, in particular $\alpha_s^{\vee} = \alpha_s$ and $\alpha_t^{\vee} = \alpha_t/2$ and the Cartan matrix $((\alpha_i^{\vee}, \alpha_j))$ is 
\[
\begin{pmatrix}
\alpha_s^{\vee}(\alpha_s) & \alpha_s^{\vee}(\alpha_t)\\
\alpha_t^{\vee}(\alpha_s) & \alpha_t^{\vee}(\alpha_t)
\end{pmatrix} = \begin{pmatrix}
2 & -2\\
-1 & 2
\end{pmatrix}.
\]

Define the algebra $U_q(\mathfrak{sp}_4)$ as the $\mathbb{Q}(q)$ algebra given by generators 
\[
F_{s},F_t,  K_{s}^{\pm 1},K_t^{\pm1} E_{s}, E_t
\]
and relations 
\begin{itemize}
\item $K_sK_s^{-1} = 1= K_sK_s^{-1}, K_tK_t^{-1} = 1= K_t^{-1}K_t, K_sK_t = K_t K_s$
\item $K_tE_t = q^4 E_tK_t,  K_tE_s= q^{-2}E_sK_t$
\item $K_sE_t = q^{-2} E_tK_s,  K_sE_s = q^2E_sK_s$
\item $K_tF_t= q^{-4} F_tK_t, K_tF_s = q^2F_sK_t$
\item $K_sF_t= q^2F_tK_s, K_sF_s = q^{-2} F_sK_s$
\item $E_tF_s = F_sE_t, E_sF_t = F_tE_s$
\item $E_tF_t= F_tE_t + \dfrac{K_t - K_t{^-1}}{q^2- q^{-2}}$
\item $E_sF_s = F_sE_s + \dfrac{K_s- K_s^{-1}}{q- q^{-1}}$
\item $E_t^2E_s - \dfrac{[4]_q}{[2]_q}E_tE_sE_t + E_sE_t^2 = 0$
\item $E_s^3E_t - [3]_qE_s^2E_tE_s + [3]_qE_sE_tE_s^2 - E_sE_t^3$
\end{itemize}
Our convention is $[n]_q := \dfrac{q^n- q^{-n}}{q- q^{-1}}$ and $[n]_q!= [n]_q[n-1]_q\ldots[2]_q[1]_q$. 

Recall that $\mathcal{A} = \mathbb{Z}[q, q^{-1}, [2]_q^{-1}]$. Let $U^{\mathcal{A}}_q(\mathfrak{sp}_4)$ be the unital $\mathcal{A}$-subalgebra of $U_q(\mathfrak{sp}_4)$ spanned by $K_{s}^{\pm 1}, K_t^{\pm 1}$, and the divided powers  
\[
E_s^{(n)} = \dfrac{E_s^n}{[n]_q!}, F_s^{(n)}= \dfrac{F_s^n}{[n]_q!}, E_t^{(n)} = \dfrac{E_t}{[n]_{{q^2}}!}, F_t^{(n)} = \dfrac{F_t}{[n]_{{q^2}}!}
\]
for all $n\in \mathbb{Z}_{\ge 1}$. So $U_q^{\mathcal{A}}(\mathfrak{sp}_4)$ is Lusztig's divided powers quantum group \cite{Andersen:1991wl}. 

Let $V^{\mathcal{A}}(\varpi_1)$ denote the free $\mathcal{A}$ module with basis 
\begin{equation}
v_{(1,0)}, v_{(-1,1)}, v_{(1,-1)}, v_{(-1,0)},
\end{equation}
and action of $U_q^{\mathcal{A}}(\mathfrak{sp}_4)$ given by:
\begin{equation}\label{fund1}
v_{(-1,0)}\substack{\xrightarrow{E_s= 1}\\ \xleftarrow[F_s = 1]{}} v_{(1, -1)} \substack{\xrightarrow{E_t= 1}\\ \xleftarrow[F_t= 1]{}} v_{(-1, 1)}  \substack{\xrightarrow{E_s= 1}\\ \xleftarrow[F_s = 1]{}} v_{(1, 0)}. 
\end{equation}
Also, let $V^{\mathcal{A}}(\varpi_2)$ denote the free $\mathcal{A}$ module with basis
\begin{equation}
v_{(0,1)} , v_{(2, -1)}, v_{(0,0)}, v_{(-2, 1)}, v_{(0,-1)},
\end{equation}
and action of $U_q^{\mathcal{A}}(\mathfrak{sp}_4)$ given by: 
\begin{equation}\label{fund2}
v_{(0, -1)}\substack{\xrightarrow{E_t= 1} \\ \xleftarrow[F_t = 1]{}} v_{(-2, 1)}  \substack{\xrightarrow{E_s= 1} \\ \xleftarrow[\text{$F_s = [2]_q$}]{}} v_{(0, 0)} \substack{\xrightarrow{E_s= [2]_q}\\ \xleftarrow[F_s = 1]{}} v_{(2, -1)} \substack{\xrightarrow{E_t= 1}\\ \xleftarrow[F_t = 1]{}} v_{(0,1)}.
\end{equation}
The elements $K_{\alpha}$ act on the basis vectors by
\begin{equation}
K_s\cdot v_{(i, j)} = q^iv_{(i, j)} \ \ \ \ \text{and} \ \ \ \ \ K_t \cdot v_{(i, j)} = q^{2j}v_{(i, j)}.
\end{equation}

Our convention is that whenever we do not indicate the action of $E_{\alpha}$ or $F_{\alpha}$ they act by zero. The action of higher divided powers on these modules can be extrapolated from the given data. For example, $F_s^{(2)}v_{(-2, 1)} = v_{(-2, 1)}$. 

\begin{remark}
Why are we using $\mathcal{A}$ instead of $\mathbb{Z}[q, q^{-1}]$? When $[2]_q= 0$, the Weyl module $\ak \ot V^{\mathcal{A}}(\varpi_2)$ is not irreducible and the correct choice of combinatorial category seems to be the $\mathbb{Z}[q^{1/2}, q^{-{1/2}}]$-linear monoidal category generated by $V_q$ and $\Lambda^2(V_q)$. The module $\Lambda^2(V_q)$ is the $\mathbb{Z}[q^{1/2}, q^{-{1/2}}]$-basis 
\begin{equation}
v_{(0,1)} , v_{(2, -1)}, X_0, Y_0, v_{(-2, 1)}, v_{(0,-1)},
\end{equation}
and action of $U_q^{\mathcal{A}}(\mathfrak{sp}_4)$ given by: 
\begin{equation}\label{fund2}
v_{(0, -1)}\substack{\xrightarrow{E_t= 1} \\ \xleftarrow[F_t = 1]{}} v_{(-2, 1)}  \substack{\xrightarrow{E_s} \\ \xleftarrow[\text{$F_s$}]{}} X_0 \oplus Y_0 \substack{\xrightarrow{E_s}\\ \xleftarrow[\text{$F_s$}]{}} v_{(2, -1)} \substack{\xrightarrow{E_t= 1}\\ \xleftarrow[F_t = 1]{}} v_{(0,1)}.
\end{equation}
where
\begin{equation}
\begin{split}
 E_s\cdot Y_0 = q^{-1/2}v_{(2, -1)}& \ \ \ \ \ \  E_s\cdot X_0 = q^{1/2} v_{(2, -1)} \\
 E_s\cdot v_{(-2, 1)}& = q^{1/2}X_0 + q^{-1/2} Y_0\\
 F_s\cdot Y_0 = q^{-1/2}v_{(-2, 1)}& \ \ \ \ \ \  F_s \cdot X_0 = q^{1/2}v_{(-2, 1)}\\
F_s\cdot v_{(2, -1)}& = q^{1/2}X_0 + q^{-1/2}Y_0. 
\end{split}
\end{equation}

The module $V^{\mathcal{A}}(\varpi_2)$ can be defined over $\mathbb{Z}[q^{1/2}, q^{-1/2}]$. There is a map from $V^{\mathcal{A}}(\varpi_2)$ into $\Lambda^2(V_q)$, such that $v_{(0, 0)}\mapsto q^{1/2}X_0 + q^{-1/2}Y_0$. Moreover, the cokernel of this inclusion map will be isomorphic to the trivial module. Thus, $\Lambda^2(V_q)$ is filtered by Weyl modules, and the filtration splits when $[2]_q\ne 0$. If $[2]_q= 0$, then $\Lambda^2(V_q)$ is indecomposable with socle and head isomorphic to the trivial module, and middle subquotient isomorphic to the irreducible module of highest weight $\varpi_2$. 
\end{remark}

The algebra $U_q(\mathfrak{sp}_4)$ is a Hopf algebra with structure maps $(\Delta, S, \epsilon)$ defined on generators by 
\begin{itemize}
\item $\Delta(E_{\alpha}) = E_{\alpha} \ot 1 + K_{\alpha} \ot E_{\alpha}, \Delta(F_{\alpha}) = F_{\alpha} \ot K_{\alpha}^{-1} + 1\ot F_{\alpha}, \Delta(K_{\alpha}) = K_{\alpha} \ot K_{\alpha}$
\item $S(E_{\alpha}) = - K_{\alpha}^{-1} E_{\alpha}, S(F_{\alpha})= - F_{\alpha}K_{\alpha}, S(K_{\alpha}) = K_{\alpha}^{-1}$ 
\item $\epsilon(E_{\alpha}) = 0, \epsilon(F_{\alpha}) = 0, \epsilon(K_{\alpha}) = 1$.
\end{itemize}
Furthermore, the algebra $U_q^{\mathcal{A}}(\mathfrak{sp}_4)$ is a sub-Hopf-algebra of $U_q(\mathfrak{sp}_4)$ \cite{Andersen:1991wl}. Therefore, $U^{\mathcal{A}}_q(\mathfrak{sp}_4)$ will act on the tensor product of representations through the coproduct $\Delta$. 

Using the antipode $S$, we can define an action of $U_q^{\mathcal{A}}(\mathfrak{sp}_4)$ on 

\begin{equation}
V^{\mathcal{A}}(\varpi_1)^* = \Hom_{\mathcal{A}}(V^{\mathcal{A}}(\varpi_1), \mathcal{A})
\end{equation}
by
\begin{equation}\label{dualfund1}
-q^4v_{(1,0)}^*\substack{\xrightarrow{E_s= 1}\\ \xleftarrow[F_s = 1]{}} q^3v_{(-1, 1)}^* \substack{\xrightarrow{E_t= 1}\\ \xleftarrow[F_t= 1]{}} -q^2 v_{(1, -1)}^*  \substack{\xrightarrow{E_s= 1}\\ \xleftarrow[F_s = 1]{}} v_{(-1, 0)}^*, 
\end{equation}
and on
\begin{equation}
V^{\mathcal{A}}(\varpi_2)^* = \Hom_{\mathcal{A}}(V^{\mathcal{A}}(\varpi_2), \mathcal{A})
\end{equation}
by
\begin{equation}\label{dualfund2}
q^6 v_{(0, -1)}^*\substack{\xrightarrow{E_t= 1}\\ \xleftarrow[F_t = 1]{}} -q^4v_{(-2, 1)}^* \substack{\xrightarrow{E_s= 1}\\ \xleftarrow[\text{$F_s= [2]_q$}]{}} q^2[2]_qv_{(0, 0)}^* \substack{\xrightarrow{E_s= [2]_q}\\ \xleftarrow[F_s = 1]{}} -q^2v_{(2, -1)}^*\substack{\xrightarrow{E_t= 1}\\ \xleftarrow[F_t = 1]{}} v_{(0,1)}^*. 
\end{equation}

Comparing \eqref{fund1} and \eqref{dualfund1} we see there is an isomorphism of $U_q^{\mathcal{A}}(\mathfrak{sp}_4)$ modules
\begin{equation}
\varphi_1: V^{\mathcal{A}}(\varpi_1)\rightarrow V^{\mathcal{A}}(\varpi_1)^*
\end{equation}
such that basis elements in \eqref{fund1} are sent to the basis elements in \eqref{dualfund1}. By comparing \eqref{fund2} and \eqref{dualfund2} we similarly obtain an isomorphism
\begin{equation}
\varphi_2: V^{\mathcal{A}}(\varpi_2)\rightarrow V^{\mathcal{A}}(\varpi_2)^*
\end{equation}
sending basis elements in \eqref{fund2} to the basis elements in \eqref{dualfund2}.

In \eqref{subsec-defeval} we will define a monoidal functor from $\DD$ to  $U_q^{\mathcal{A}}(\mathfrak{sp}_4)-\text{mod}$. The functor will send $\blues$ to $V^{\mathcal{A}}(\varpi_1)$ and $\greent$ to $V^{\mathcal{A}}(\varpi_2)$. The dual modules $V^{\mathcal{A}}(\varpi_1)^*$ and $V^{\mathcal{A}}(\varpi_2)^*$ will not be in the image of the functor $\eval$. However, the maps $\varphi_1$ and $\varphi_2$ are fixed isomorphisms of these dual modules with modules which are in the image of the functor. 

%===========
\subsection{Caps and Cups}
\label{subsec-caps and cups}
%===========

\begin{lemma}
If $V$ is any finite rank $\mathcal{A}$ lattice with basis $e_i$, define maps:
\begin{equation}
\mathcal{A}\xrightarrow{u} V\ot \Hom_{\mathcal{A}}(V, \mathcal{A})\xrightarrow{c} \mathcal{A}
\end{equation}
\begin{equation}
\mathcal{A}\xrightarrow{u'}  \Hom_{\mathcal{A}}(V, \mathcal{A})\ot V\xrightarrow{c'} \mathcal{A}
\end{equation}
where $u(1) = \sum e_i\ot e_i^*$, $u'(1)= \sum e_i^*\ot e_i$, $c(v\ot f)= f(v)$, and $c'(f\ot v)= f(v)$. Then
\begin{equation}\label{interchangeV}
(\id_V\ot c')\circ(u\ot\id_V) = \id_V = (c\ot \id_V)\circ (\id_V\ot u')
\end{equation}
and
\begin{equation}\label{interchangeVdual}
(\id_{V^*}\ot c)\circ (u'\ot \id_{V^*})= \id_{V^*} = (c'\ot \id_{V^*})\circ (\id_{V^*}\ot u).
\end{equation}
\end{lemma}

\begin{proof}
We will show that 
\[
(\id_V\ot c')\circ(u\ot\id_V) = \id_V 
\]
the arguments to establish the other three equalities in \eqref{interchangeV} and \eqref{interchangeVdual} are similar. 

Let $v\in V$. Since $e_i$ is a basis for $V$ we can write $v= \sum v_ie_i$ for some $v_i \in \mathcal{A}$. Thus,
\[
(\id_V\ot c')\circ(u\ot \id_V)(v) =(\id_V\ot c')( \sum e_i\ot e_i^*\ot v) = \sum e_i \cdot e_i^*(v) =\sum v_i e_i =  v.
\]
\end{proof}

\begin{lemma}
Fix an isomorphism $\varphi: V\rightarrow V^{*}$ and write $\textbf{cap}= c'\circ (\varphi\ot \id)$ and $\textbf{cup}= (\id \ot \varphi^{-1})\circ u$. Then  
\begin{equation}\label{cupcapisotopy}
(\id_V\ot \textbf{cap})\circ (\textbf{cup}\ot \id_V) = \id_V = (\textbf{cap}\ot \id_V)\circ(\id_V\ot \textbf{cup}).
\end{equation}

\end{lemma}

\begin{proof}
Using $\varphi \circ\varphi^{-1} = \id= \varphi^{-1} \circ \varphi$, \eqref{cupcapisotopy} follows easily from \eqref{interchangeV} and \eqref{interchangeVdual}.
\end{proof}

The $\mathcal{A}$-linear maps
\begin{equation}\label{cupcapdef}
\mathcal{A}\xrightarrow{\textbf{cup}_i:= (\id \ot \varphi_i^{-1}) \circ u_i} V^{\mathcal{A}}(\varpi_i)\ot V^{\mathcal{A}}(\varpi_i) \xrightarrow{\textbf{cap}_i :=  c_i'\circ (\varphi_i\ot \id)} \mathcal{A}, \ \ \ \ \  \text{for} \ i=1, 2,
\end{equation}
are actually maps of $U_q^{\mathcal{A}}(\mathfrak{sp}_4)$ modules, where $\mathcal{A}$ is the trivial module. The functor $\eval$ will send the cups and caps from the diagrammatic category to the maps $\textbf{cup}_i$ and $\textbf{cap}_i$

The module $V^{\mathcal{A}}(\varpi_1)$ has basis
\begin{equation}\label{fund1basis}
\lbrace v_{(1,0)}, v_{(-1,1)} = F_sv_{(1, 0)}, v_{(1,-1)} = F_tF_sv_{(1, 0)}, v_{(-1,0)} = F_sF_tF_sv_{(1,0)} \rbrace,
\end{equation}
and the module $V^{\mathcal{A}}(\varpi_2)$ has basis
\begin{equation}\label{fund2basis}
\lbrace v_{(0,1)} , v_{(2, -1)}= F_tv_{(0,1)}, v_{(0,0)} = F_sF_tv_{(0,1)}, v_{(-2, 1)} = F_s^{(2)}F_tv_{(0,1)}, v_{(0,-1)} = F_tF_s^{(2)}F_tv_{(0,1)}\rbrace.
\end{equation}
With respect to these bases, we can write $\cupa: \mathcal{A}\rightarrow V^{\mathcal{A}}(\varpi_1)\ot V^{\mathcal{A}}(\varpi_1)$ as
\begin{equation}\label{cup1formula}
1\mapsto -q^{-4}v_{(1,0)}\ot v_{(-1,0)} + q^{-3} v_{(-1,1)} \ot v_{(1, -1)} - q^{-1} v_{(1, -1)} \ot v_{(-1, 1)} + v_{(-1, 0)}\ot v_{(1, 0)},
\end{equation}
and $\cupb  : \mathcal{A} \rightarrow V^{\mathcal{A}}(\varpi_2)\ot V^{\mathcal{A}}(\varpi_2)$ as
\begin{equation}\label{cup2formula}
1\mapsto q^{-6} v_{(0,1)}\ot v_{(0,-1)} - q^{-4}v_{(2, -1)}\ot v_{(-2, 1)} +\dfrac{q^{-2}}{[2]_q}v_{(0,0)}\ot v_{(0,0)} - q^{-2} v_{(-2,1)}\ot v_{(2, -1)} + v_{(0,-1)}\ot v_{(0,1)}. 
\end{equation}

To record the maps $\textbf{cap}_i$ in our basis we use the matrices
\begin{equation}\label{cap1formula}
\capa   (v_i \ot v_j) = \begin{blockarray}{ccccc}
 & v_{(-1, 0)} & v_{(1,-1)} & v_{(-1, 1)} & v_{(1, 0)}\\
\begin{block}{c(cccc)}
v_{(-1, 0)}& 0 & 0 & 0 &-q^4\\
v_{(1, -1)}& 0 & 0 & q^3 & 0\\
v_{(-1, 1)}& 0 & -q & 0 & 0\\
v_{(1, 0)}& 1 & 0 &0 &0 \\
\end{block}
\end{blockarray}
 \end{equation}
and
\begin{equation}\label{cap2formula}
\capb  (v_i \ot v_j) = \begin{blockarray}{cccccc}
 & v_{(0, -1)} & v_{(-2, 1)} & v_{(0,0)} & v_{(2, -1)} & v_{(0, 1)}\\
\begin{block}{c(ccccc)}
v_{(0, -1)}& 0 & 0 & 0 & 0 & q^6\\
v_{(-2, 1)}& 0 & 0 & 0 & -q^4 & 0\\
v_{(0,0)}& 0 & 0 & q^2[2]_q & 0 & 0\\
v_{(2, -1)}& 0 & -q^2 & 0 & 0 & 0\\
v_{(0, 1)}& 1 & 0 & 0 & 0 & 0 \\
\end{block}.
\end{blockarray}
\end{equation}

\begin{example}
We give two calculations to clarify how we arrived at these formulas:
\[
\capb  (v_{(0, 0)} \ot v_{(0, 0)}) = c'_2\circ (\varphi_2\ot \id)(v_{(0, 0)} \ot v_{(0, 0)}) = c'_2(q^2[2]_qv_{(0, 0)}^* \ot v_{(0, 0)}) = q^2[2]_q
\]
and
\begin{align*}
\cupa  (1) &= (\id \ot \varphi_1^{-1})\circ u_1(1) \\
&= v_{(-1, 0)}\ot \varphi_1^{-1}(v_{(-1, 0)}^*) + v_{(1, -1)}\ot \varphi_1^{-1}(v_{(1, -1)}^*)  \\
& \ \ \ \ \ \ \ \ \ \ \ \ \ \ \ \ \ \  \ \ \ \ \ \ \ \ \ \ \ \ \ \ \ \ \ \ \ \ \ \ \ \ +v_{(-1, 1)}\ot \varphi_1^{-1}(v_{(-1, 1)}^*) + v_{(1, 0)}\ot \varphi_1^{-1}(v_{(1, 0)}) \\
&= -q^{-4}v_{(1,0)}\ot v_{(-1,0)} + q^{-3} v_{(-1,1)} \ot v_{(1, -1)} - q^{-1} v_{(1, -1)} \ot v_{(-1, 1)} + v_{(-1, 0)}\ot v_{(1, 0)}.
\end{align*}
\end{example}

The maps $\textbf{cup}_i$ and $\textbf{cap}_i$ in $U_q^{\mathcal{A}}(\mathfrak{sp}_4)-\text{mod}$ are going to correspond to the colored cap and cup maps in $\DD$. In which case, the equation \eqref{cupcapisotopy} corresponds to the isotopy relations
\begin{figure}[H]
\centering
\includegraphics[width=.10\textwidth]{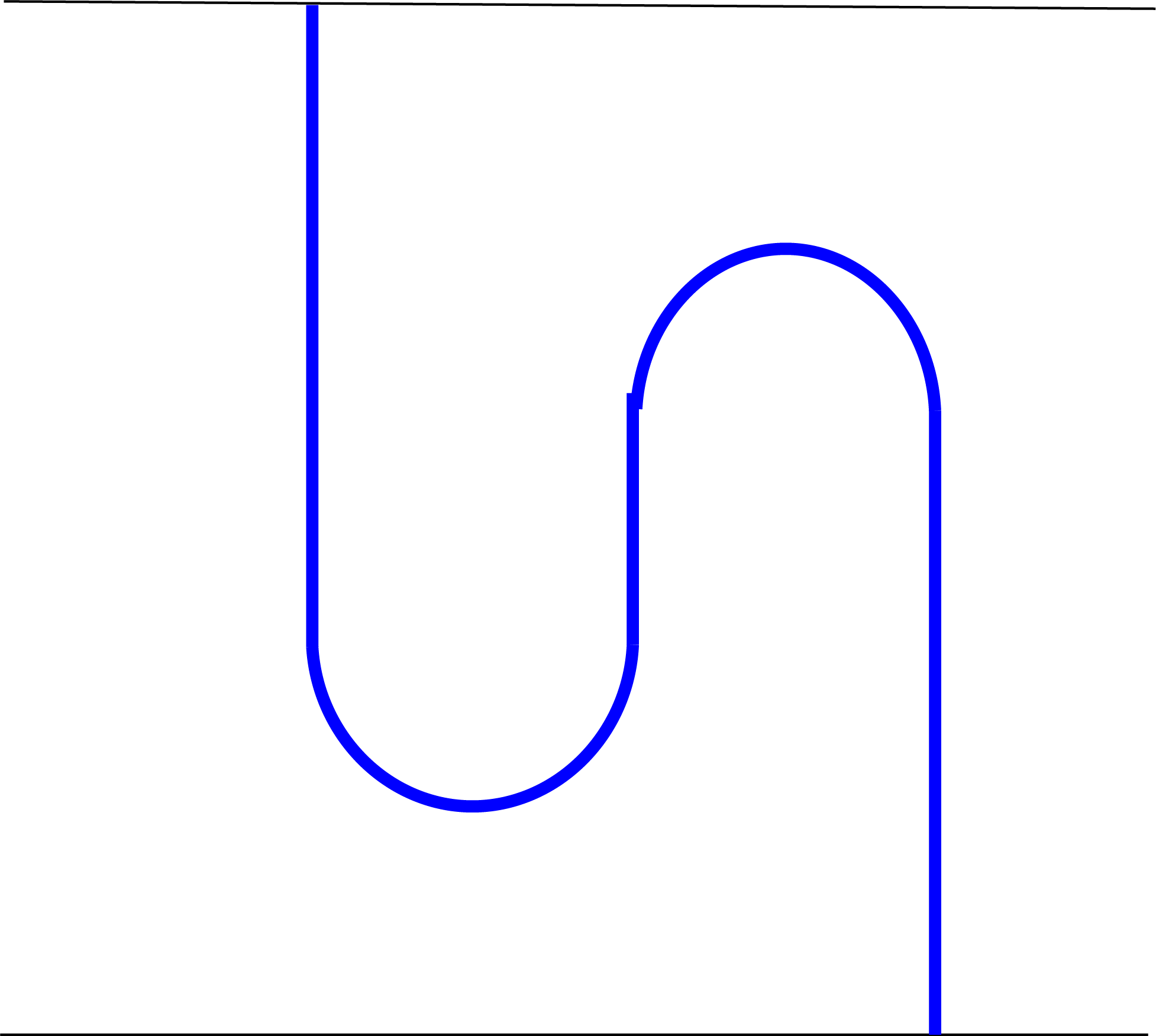}\put(8, 17){$=$} \ \ \ \ \ \ \ \ \ \includegraphics[width=.10\textwidth]{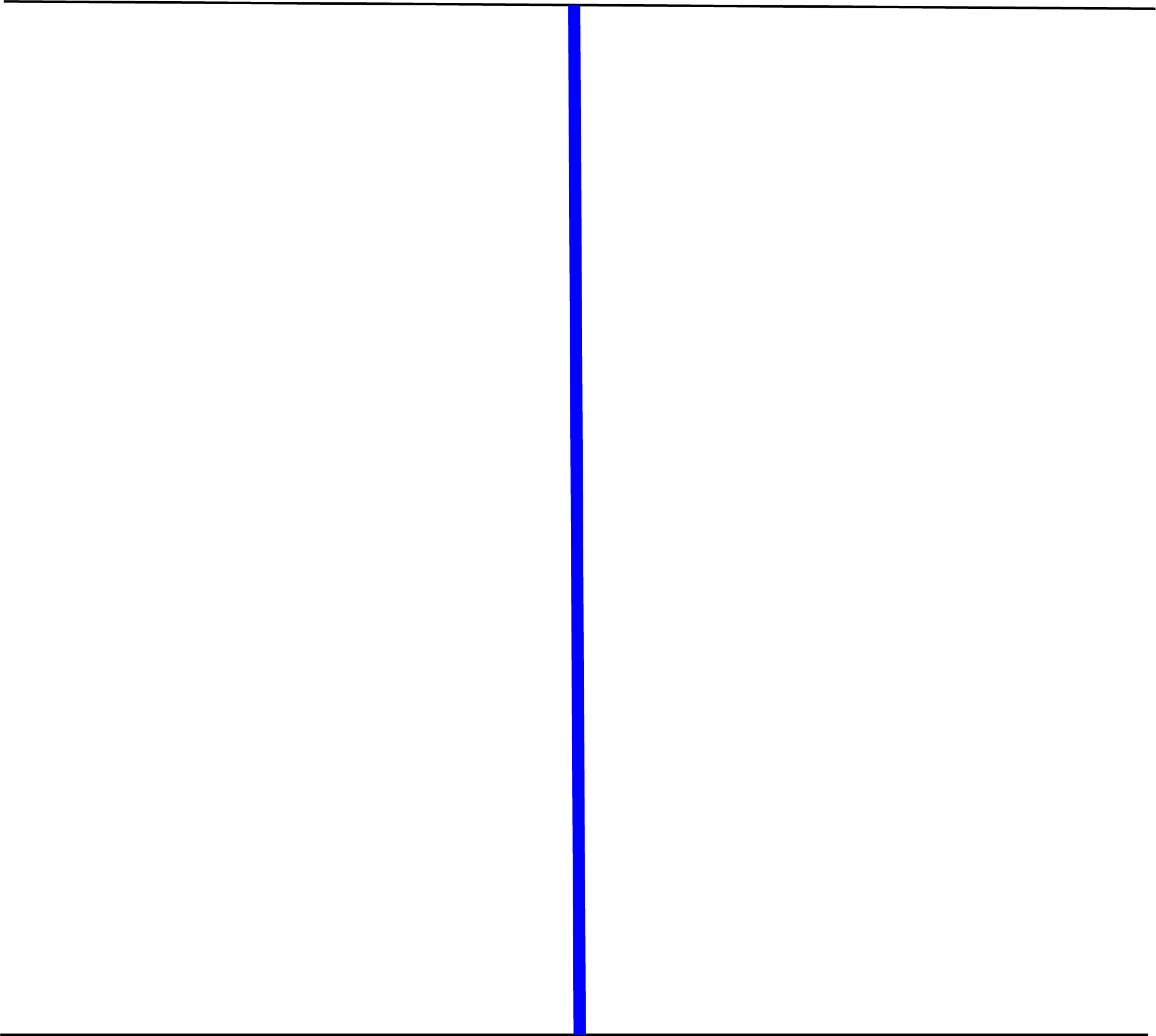}\put(8, 17){$= $} \ \ \ \ \ \ \ \ \ 
\includegraphics[width=.10\textwidth]{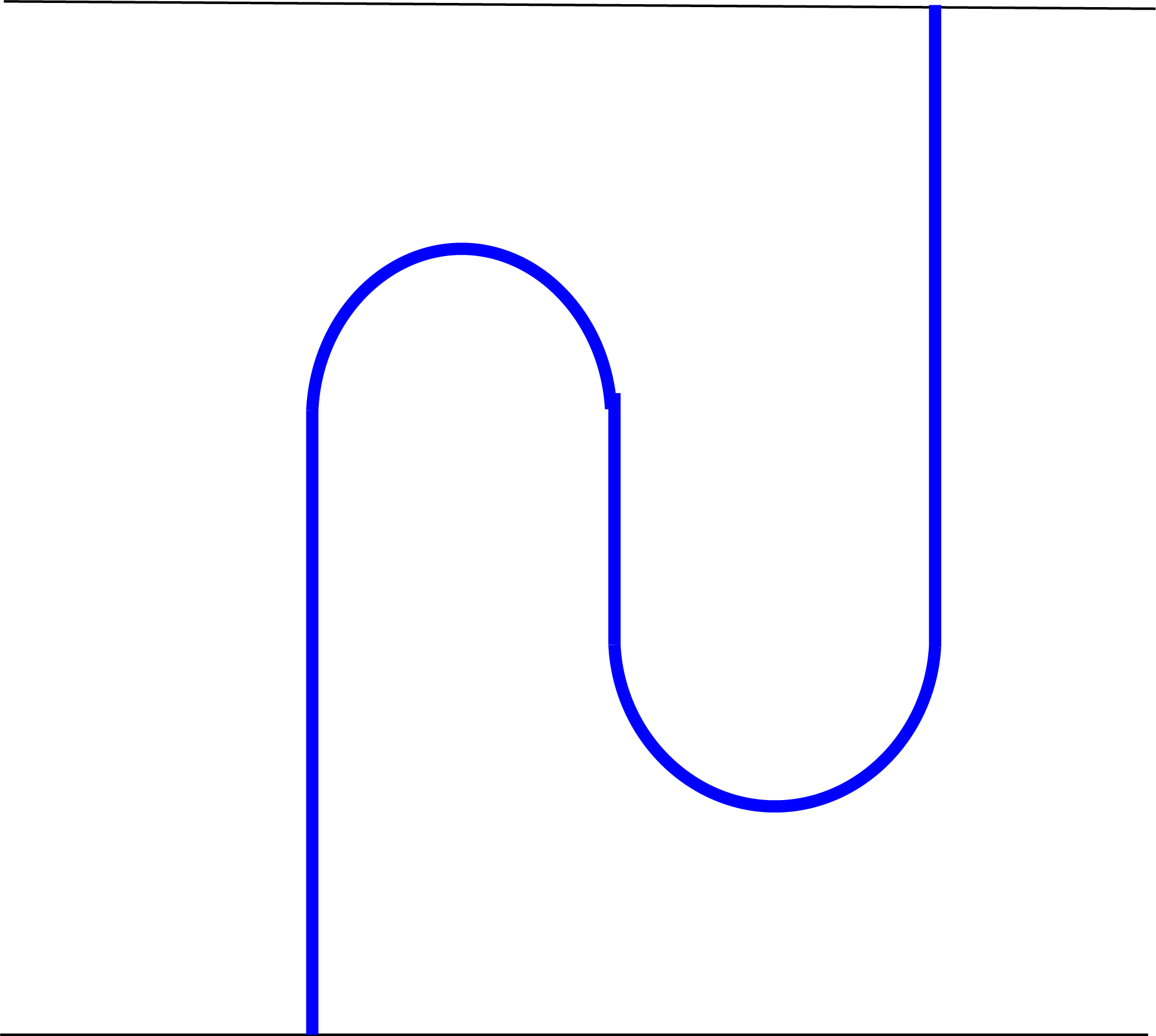}
\end{figure}
\begin{figure}[H]
\centering
\includegraphics[width=.10\textwidth]{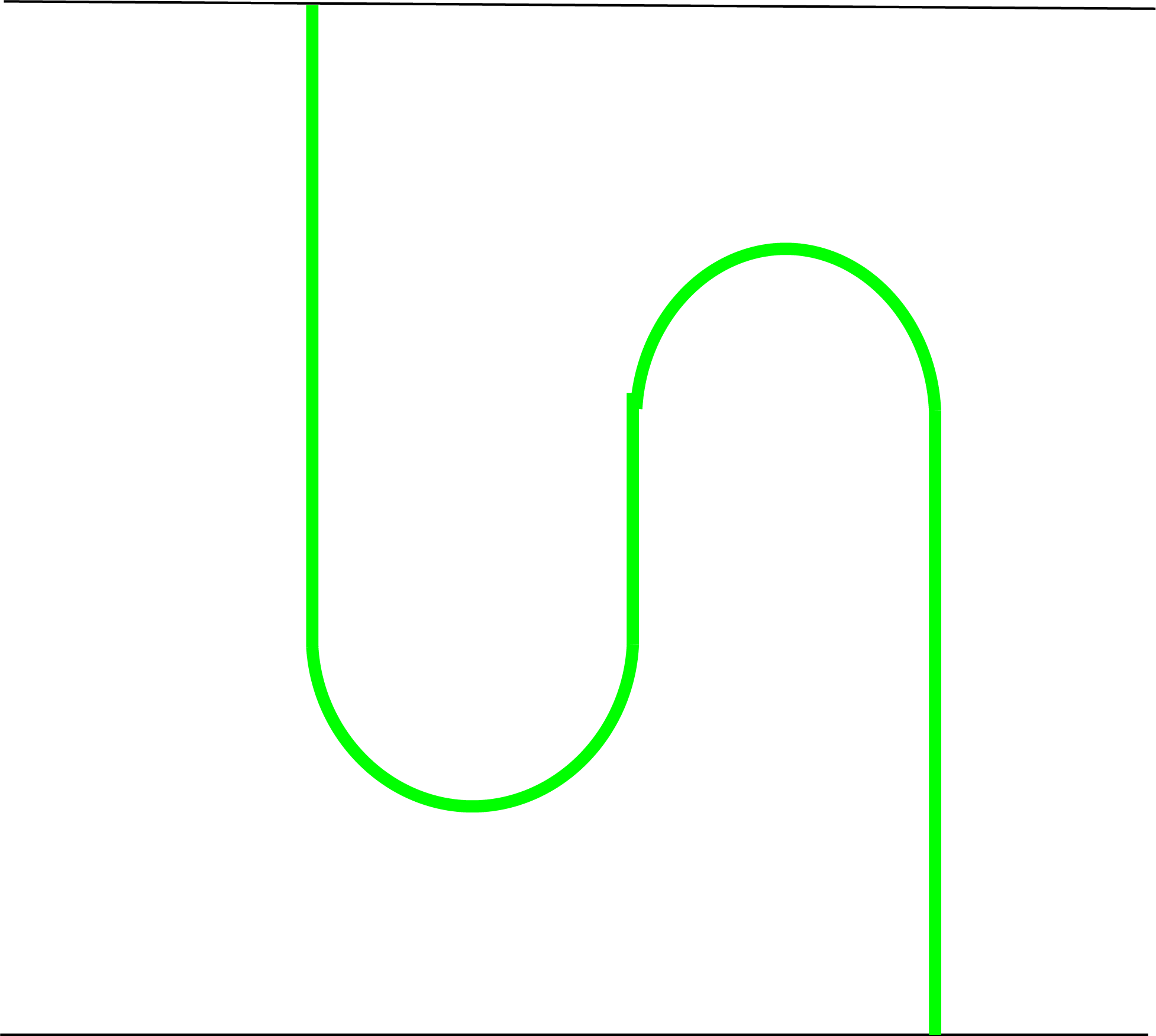}\put(8, 17){$=$} \ \ \ \ \ \ \ \ \ \includegraphics[width=.10\textwidth]{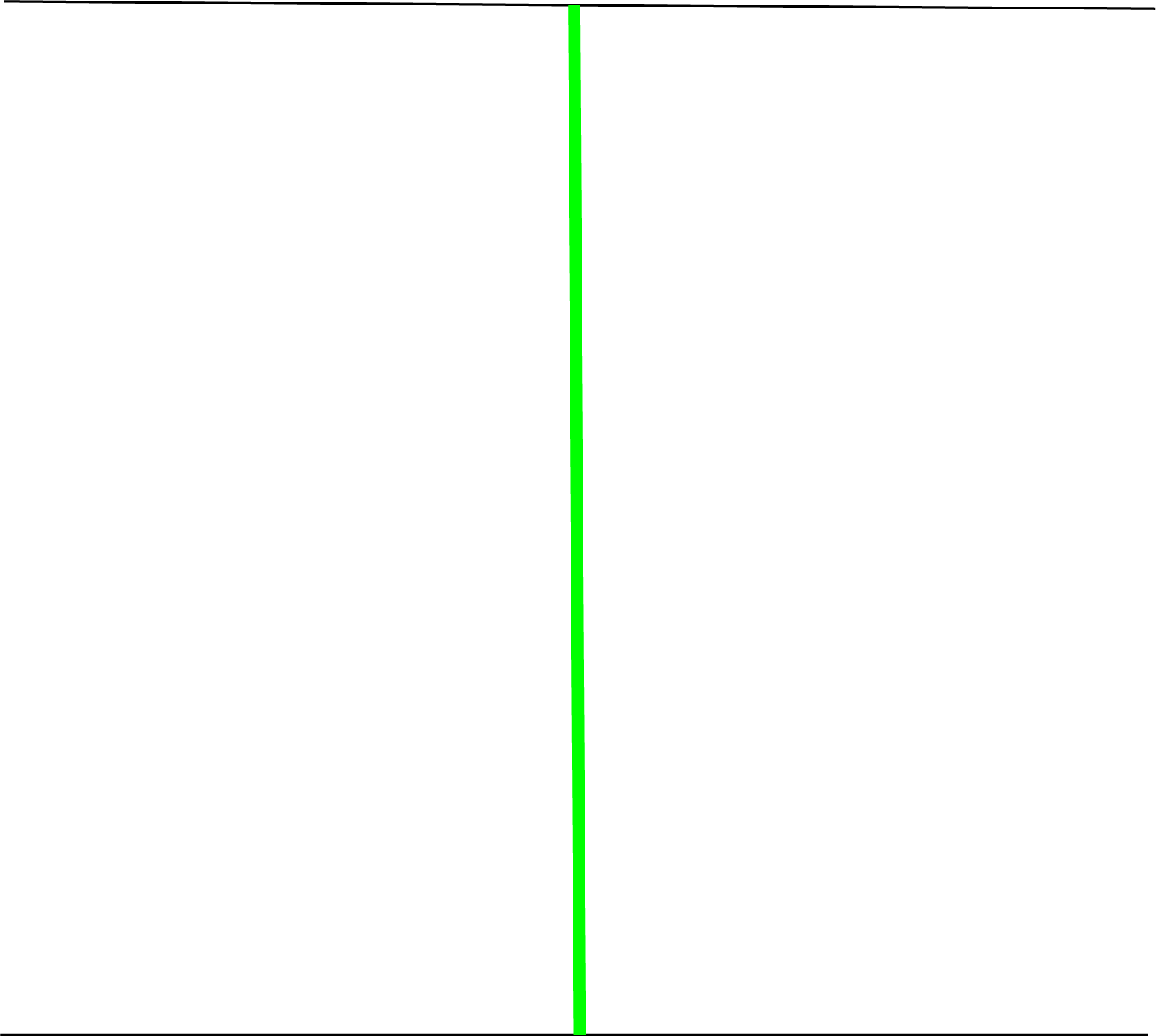}\put(8, 17){$= $} \ \ \ \ \ \ \ \ \ 
\includegraphics[width=.10\textwidth]{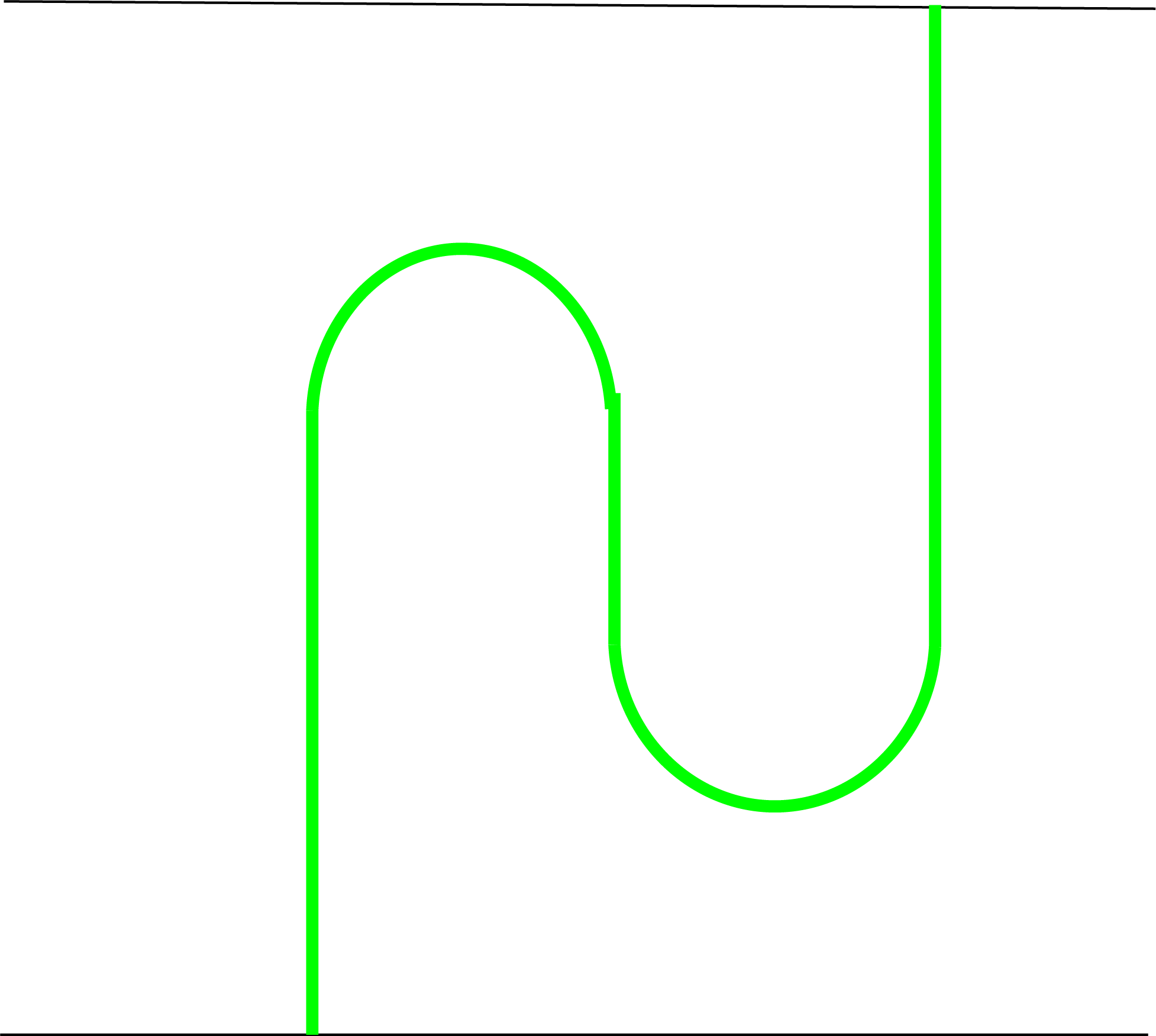}
\end{figure}

%===========
\subsection{Trivalent Vertices}
\label{subsec-trivalents}
%===========

Consider the module $V^{\mathcal{A}}(\varpi_1)\ot V^{\mathcal{A}}(\varpi_1)$. We observe that the vector $q^{-1}v_{(1,0)}\ot v_{(0,1)} - v_{(0,1)}\ot v_{(1, 0)}$ is annihilated by $E_s$ and $E_t$. The action of $K_s$ scales this vector by $1$ and the action of $K_t$ scales the vector by $q^2$. There is an $\mathcal{A}$-linear map
\begin{equation}\label{iformula}
\begin{split}
\ii: V^{\mathcal{A}}(\varpi_2)&\rightarrow V^{\mathcal{A}}(\varpi_1)\ot V^{\mathcal{A}}(\varpi_1) \\
v_{(0,1)} &\mapsto q^{-1}v_{(1,0)} \ot v_{(-1,1)} - v_{(-1,1)}\ot v_{(1,0)} \\
v_{(2,-1)} &\mapsto q^{-1}v_{(1, 0)}\ot v_{(1, -1)} - v_{(1, -1)} \ot v_{(1,0)} \\
v_{(0,0)} &\mapsto q^{-1} v_{(1,0)} \ot v_{(-1,0)} + q^{-2}v_{(-1,1)} \ot v_{(1,-1)} \\
& \ \ \ \ \ \ \ - v_{(1,-1)}\ot v_{(-1,1)}  -q^{-1} v_{(-1,0)} \ot v_{(1,0)} \\
v_{(-2, 1)} &\mapsto q^{-1} v_{(-1,1)} \ot v_{(-1,0)} - v_{(-1,0)} \ot v_{(-1,1)} \\
v_{(0,-1)} &\mapsto q^{-1} v_{(1, -1)} \ot v_{(-1,0)} - v_{(-1,0)} \ot v_{(1, -1)}.
\end{split}
\end{equation}
Using the explicit description of $V^{\mathcal{A}}(\varpi_2)$ in \eqref{fund2}, one checks that $\ii$ is a map of $U_q^{\mathcal{A}}(\mathfrak{sp}_4)$-modules by computing the action of the generators of $U^{\mathcal{A}}_q(\mathfrak{sp}_4)$ on the vectors appearing on the right hand side of \eqref{iformula}. The morphism $\ii$ will correspond to the following diagram.
\begin{figure}[H]
\centering
\includegraphics[width=.10\textwidth]{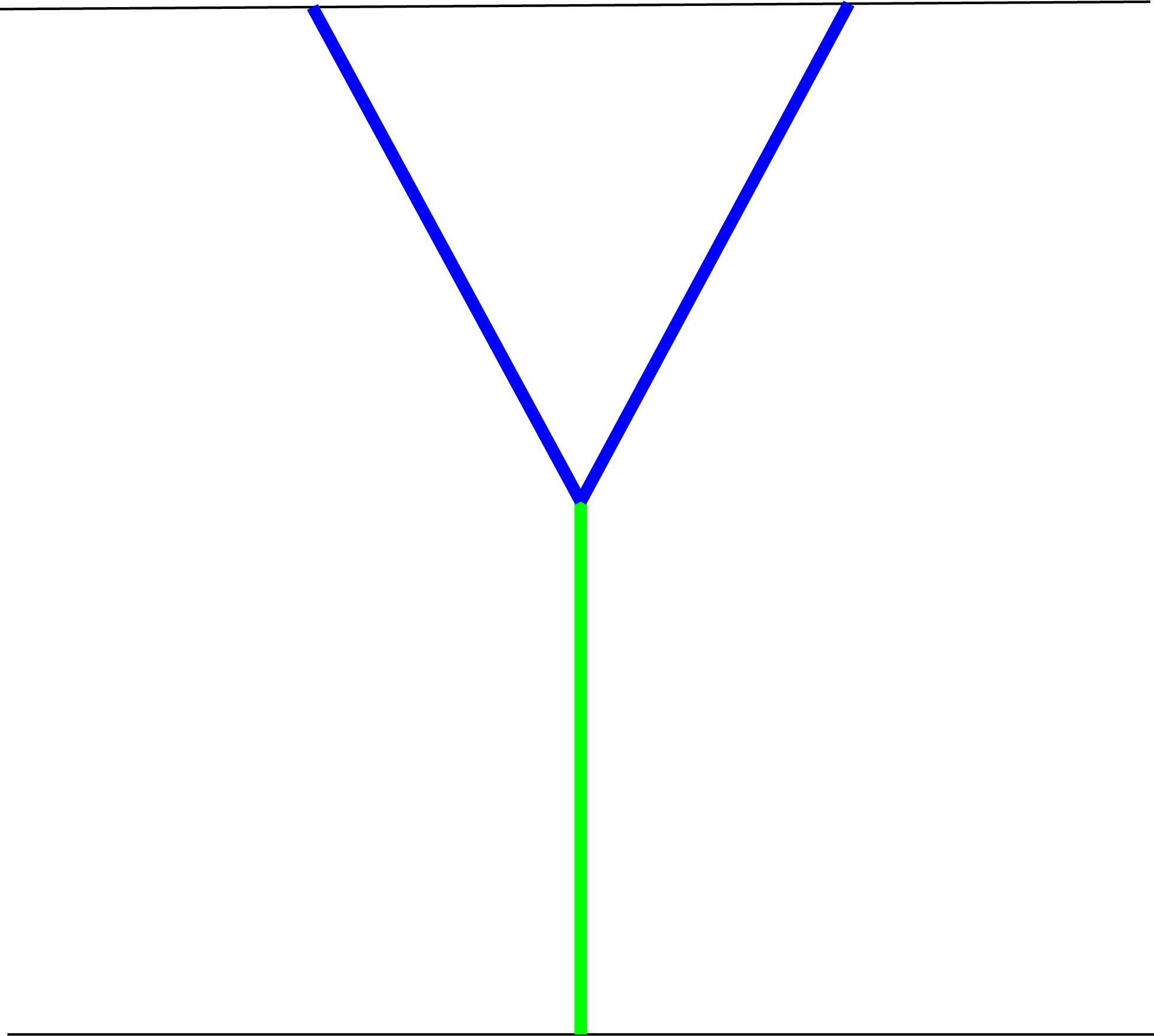}
\end{figure}

One can also check the equality of the following two elements of $\Hom_{U_q^{\mathcal{A}}(\mathfrak{sp}_4)}(V^{\mathcal{A}}(\varpi_1)\ot V^{\mathcal{A}}(\varpi_1), V^{\mathcal{A}}(\varpi_2))$:
\begin{equation}\label{LHSp}
(\id  \ot \capa   )\circ (\id  \ot \id  \ot \capa    \ot \id ) \circ (\id  \ot \ii\ot\id \ot \id ) \circ (\cupb   \ot \id  \ot\id ) 
\end{equation}
and
\begin{equation}\label{RHSp}
(\capa    \ot \id )\circ(\id \ot \capa   \ot \id  \ot \id )\circ (\id \ot \id \ot \ii \ot \id )\circ (\id  \ot \id  \ot \cupb  ).
\end{equation}
Then we will unambiguously denote both maps by $\pp$. In the graphical calculus this corresponds to the following.
\begin{figure}[H]
\centering
\includegraphics[width=.10\textwidth]{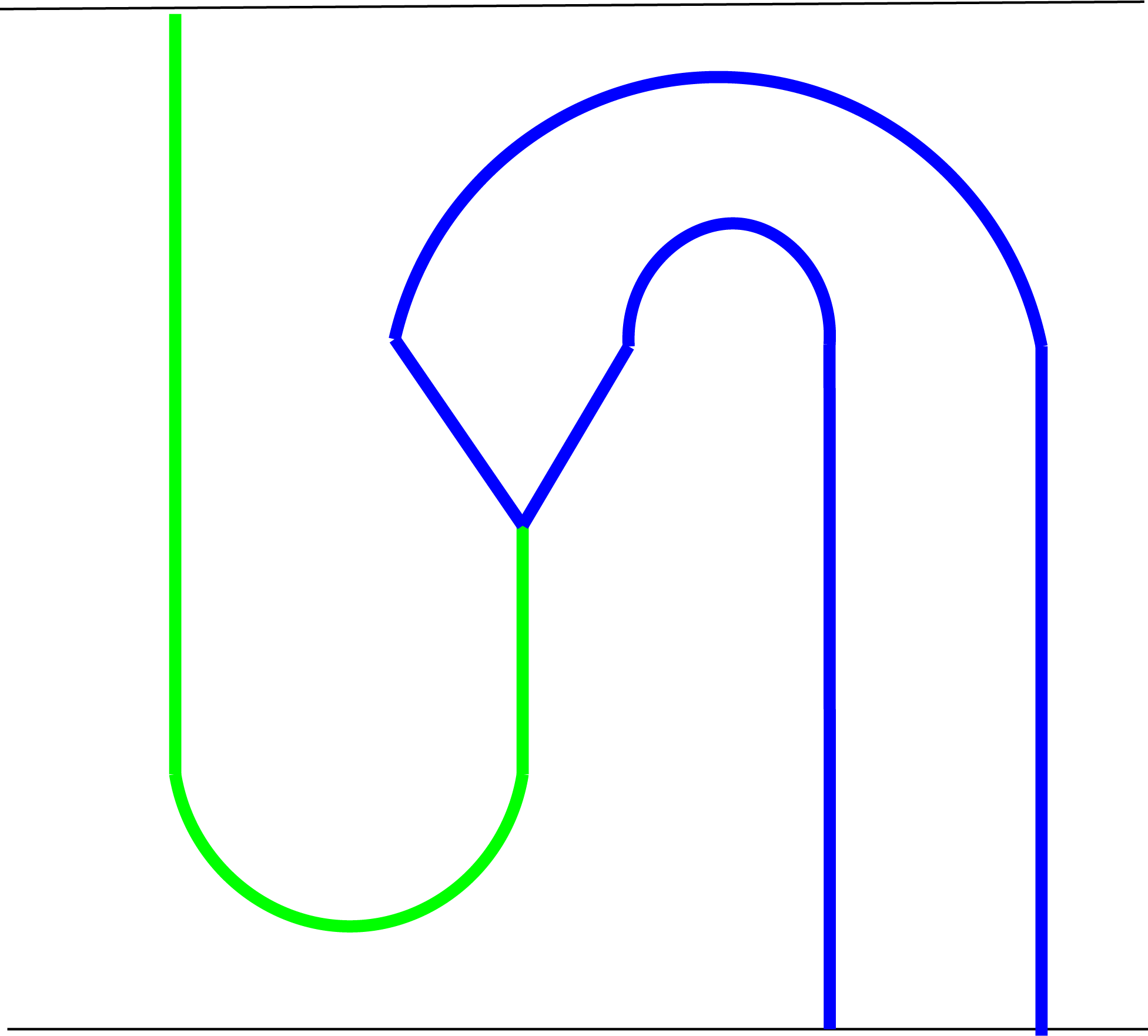}\put(8, 17){$=$}\put(-32, -12){\eqref{LHSp}} \ \ \ \ \ \ \ \ \ \includegraphics[width=.10\textwidth]{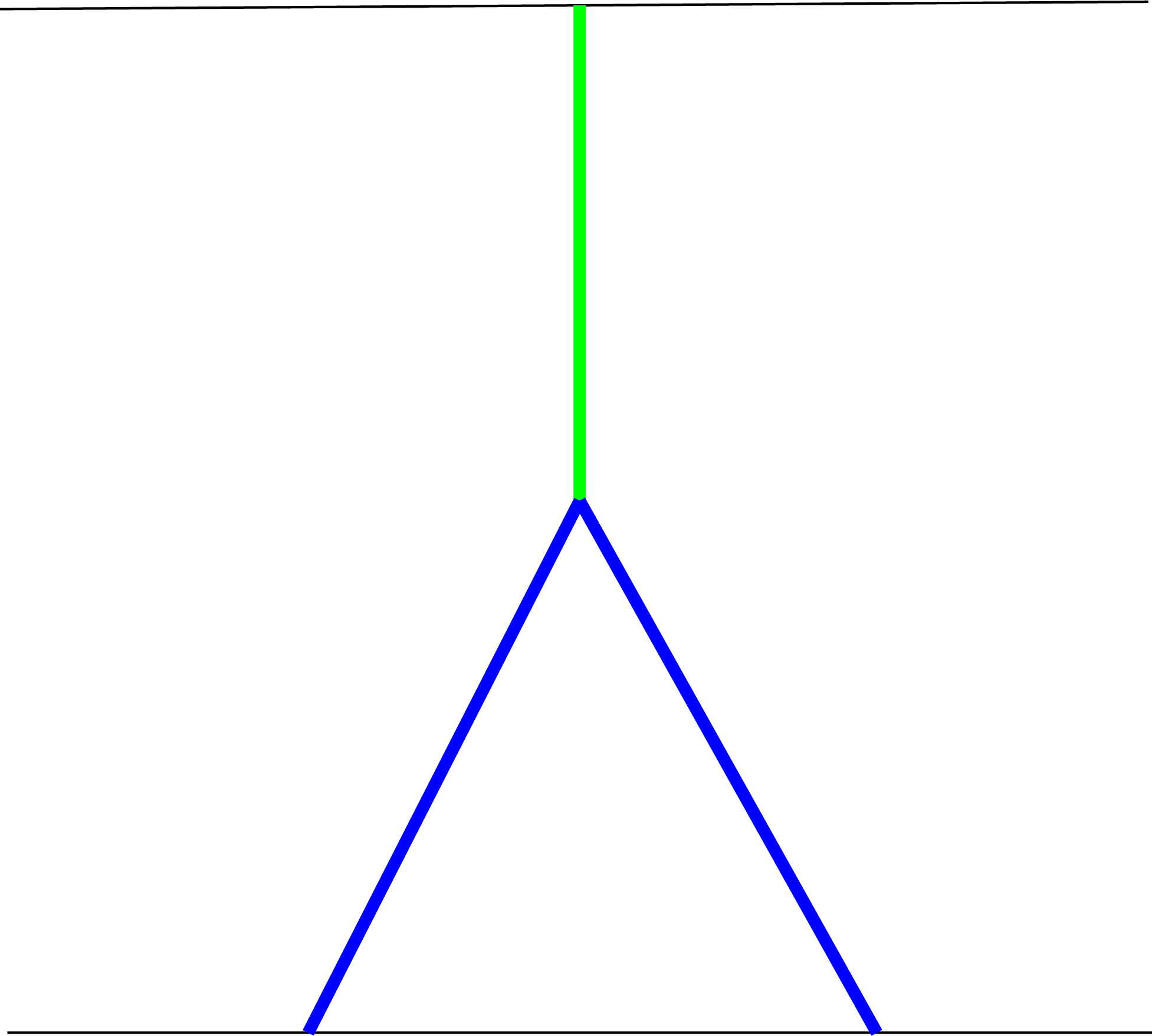}\put(8, 17){$= $}\put(-25, -10){$\pp$} \ \ \ \ \ \ \ \ \ 
\includegraphics[width=.10\textwidth]{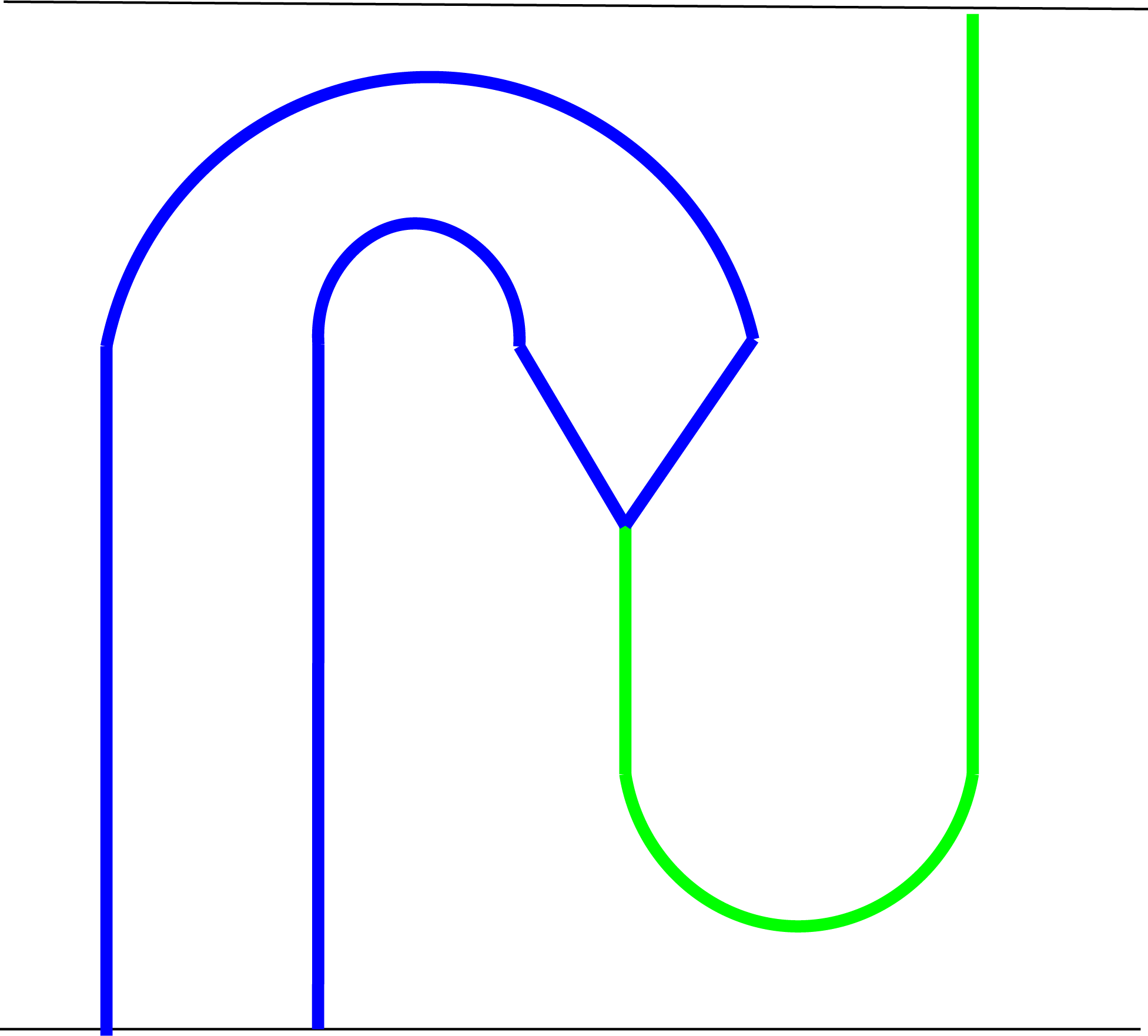}\put(-32, -12){\eqref{RHSp}}
\end{figure}
\noindent The equality of \eqref{LHSp} and \eqref{RHSp} follows from verifying that both maps act on a basis as follows. 
\begin{equation}\label{formulaforpi}
\pp: V^{\mathcal{A}}(\varpi_1)\ot V^{\mathcal{A}}(\varpi_1) \rightarrow V^{\mathcal{A}}(\varpi_2).
\end{equation}
\begin{align*}
v_{(1, 0)}\ot v_{(1, 0)}&\mapsto 0           &  v_{(-1, 1)}\ot v_{(1, 0)} &\mapsto qv_{(0, 1)}            \\
v_{(1, 0)}\ot v_{(-1, 1)}&\mapsto -v_{(0, 1)}         &  v_{(-1, 1)}\ot v_{(-1, 1)}&\mapsto 0  \\
v_{(1, 0)}\ot v_{(1, -1)}&\mapsto -v_{(2, -1)}   &  v_{(-1, 1)}\ot v_{(1, -1)}&\mapsto \dfrac{-1}{[2]_q}v_{(0, 0)}    \\
v_{(1, 0)}\ot v_{(-1, 0)}&\mapsto \dfrac{-q}{[2]_q}v_{(0,0)}   &  v_{(-1, 1)}\ot v_{(-1, 0)}&\mapsto -v_{(-2, 1)}  \\       
\end{align*}

\begin{align*}
v_{(1, -1)}\ot v_{(1, 0)}&\mapsto qv_{(2, -1)}           &  v_{(-1, 0)}\ot v_{(1, 0)} &\mapsto \dfrac{q}{[2]_q}v_{(0,0)}           \\
v_{(1, -1)}\ot v_{(-1, 1)}&\mapsto \dfrac{q^2}{[2]_q}v_{(0,0)}        &  v_{(-1, 0)}\ot v_{(-1, 1)}&\mapsto qv_{(-2, 1)} \\
v_{(1, -1)}\ot v_{(1, -1)}&\mapsto 0   &  v_{(-1, 0)}\ot v_{(1, -1)}&\mapsto qv_{(0, -1)}   \\
v_{(1, -1)}\ot v_{(-1, 0)}&\mapsto -v_{(0,-1)}   &  v_{(-1, 0)}\ot v_{(-1, 0)}&\mapsto 0        
\end{align*}

\begin{remark}
We sketch a method to compute \eqref{LHSp} evaluated on $v_{(-1, 1)}\ot v_{(1, -1)}$, the other calculations follow the same pattern. The $\capa$'s in the definition of \eqref{LHSp} are only non-zero on basis vectors of the form $v_{\mu}\ot v_{-\mu}$. Also, in the formula for $\ii$ \eqref{iformula} the only basis vector with a tensor of the form $v_{(-1, 1)}\ot v_{(1,-1)}$ is $v_{(0, 0)}$. Therefore, \eqref{LHSp} acts as
\begin{equation}
\begin{split}
v_{(-1, 1)}\ot v_{(1, -1)} &\mapsto (\id \ot \capa)\circ (\id \ot\id \ot \capa\ot \id)\left(\dfrac{q^{-2}}{[2]_q} v_{(0,0)}\ot \ii(v_{(0,0)})\ot v_{(-1, 1)}\ot v_{(1, -1)}\right) \\
&= q^{-2}\capa(v_{(1, -1)}\ot v_{(-1, 1)})\capa(v_{(-1, 1)} \ot v_{(1, -1)})\dfrac{q^{-2}}{[2]_q} v_{(0,0)}\\
&= q^{-2}q^3(-q)\dfrac{q^{-2}}{[2]_q} v_{(0,0)} \\
&= \dfrac{-1}{[2]_q}v_{(0, 0)}.
\end{split}
\end{equation}
\end{remark}

%One can use either \eqref{RHSp} or \eqref{LHSp} to compute that $p$ acts as
%\begin{align*}
%v_{(1, 0)}\ot v_{(1, 0)} &\mapsto \\
%v_{(1, 0)}\ot v_{(-1, 1)}&\mapsto  \\
%v_{(1,0)}\ot v_{(1, -1)}&\mapsto \\
%v_{(1, 0)}\ot v_{(-1, 0)}&\mapsto \\
%v_{(-1,1)}\ot v_{(1, 0)}&\mapsto \\
%v_{(-1, 1)}\ot v_{(-1, 1)} &\mapsto \\
%v_{(-1, 1)}\ot v_{(1, -1)}&\mapsto  \\
%v_{(-1,1)}\ot v_{(-1, 0)}&\mapsto \\
%v_{(1, -1)}\ot v_{(1, 0)}&\mapsto \\
%v_{(1,-1)}\ot v_{(-1, 1)}&\mapsto \\
%v_{(1, -1)} \ot v_{(1, -1)}&\mapsto \\
%v_{(1, -1)}\ot v_{(-1, 0)}&\mapsto  \\
%v_{(-1,0)}\ot v_{(1, 0)}&\mapsto \\
%v_{(-1, 0)}\ot v_{(-1, 1)}&\mapsto \\
%v_{(-1,0)}\ot v_{(1, -1)}&\mapsto \\
%v_{(-1, 0)} \ot v_{(-1, 0)}&\mapsto \\
%\end{align*}

%===========
\subsection{The Definition of the Evaluation Functor}
\label{subsec-defeval}
%===========

\begin{thm}\label{evalthm}
There is a monoidal functor
\[
\eval:\DD \rightarrow U_q^{\mathcal{A}}(\mathfrak{sp}_4)-\text{mod}.
\]
defined on objects by defining 
$\eval(\blues)= V^{\mathcal{A}}(\varpi_1)$ and $\eval(\greent) = V^{\mathcal{A}}(\varpi_2)$ and then extending monoidally. The functor $\eval$ is defined on morphisms by first defining
\begin{figure}[H]
\centering
\includegraphics[width=.10\textwidth]{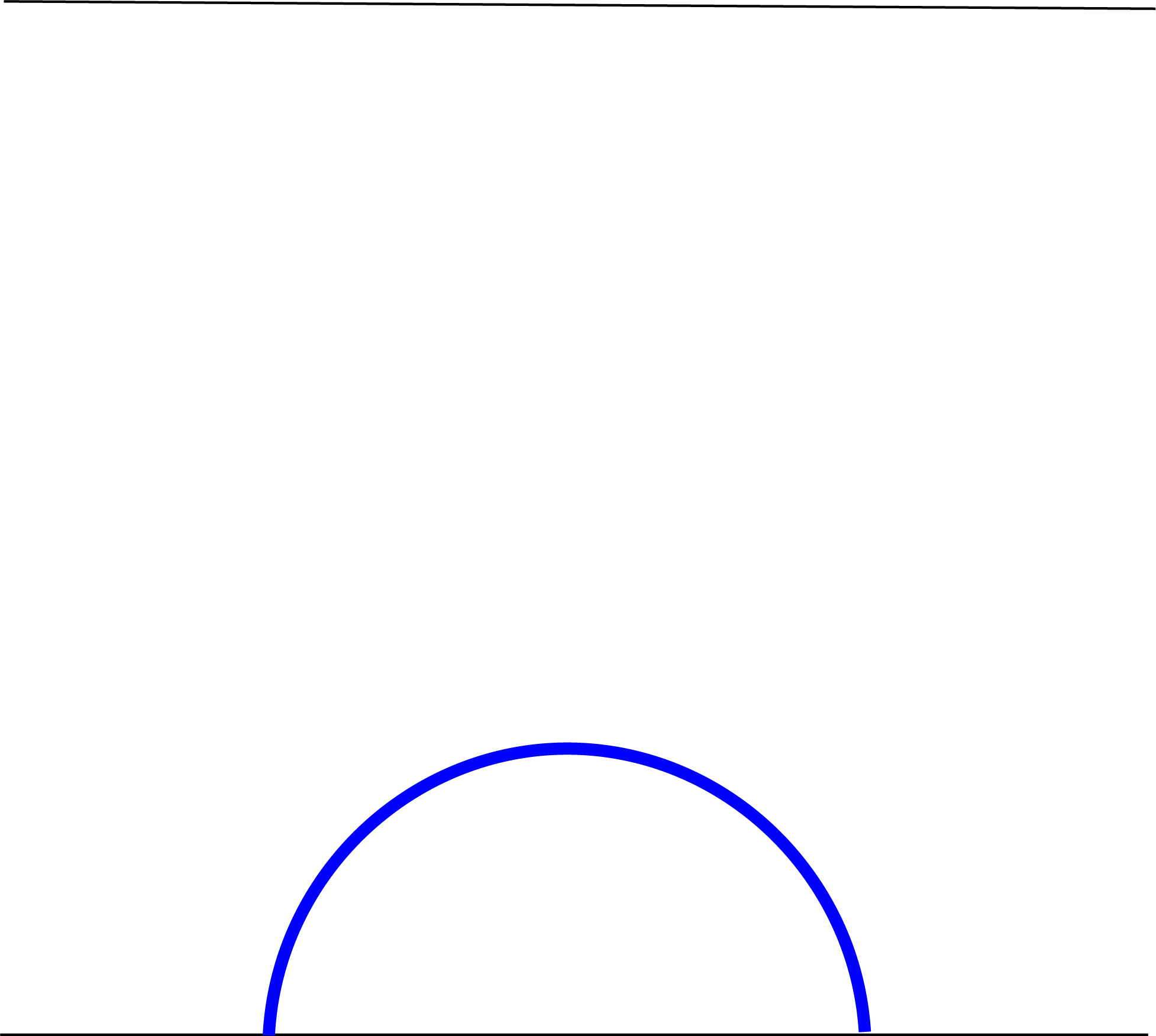}\put(10, 20){$\mapsto \capa  $} \ \ \ \ \ \ \ \ \ \ \ \ \ \ \ \ \ \ \ \  \ \ \ \ \ \ \ \ \ \ \ \ \ \ \ \ \ \ \ \ 
\includegraphics[width=.10\textwidth]{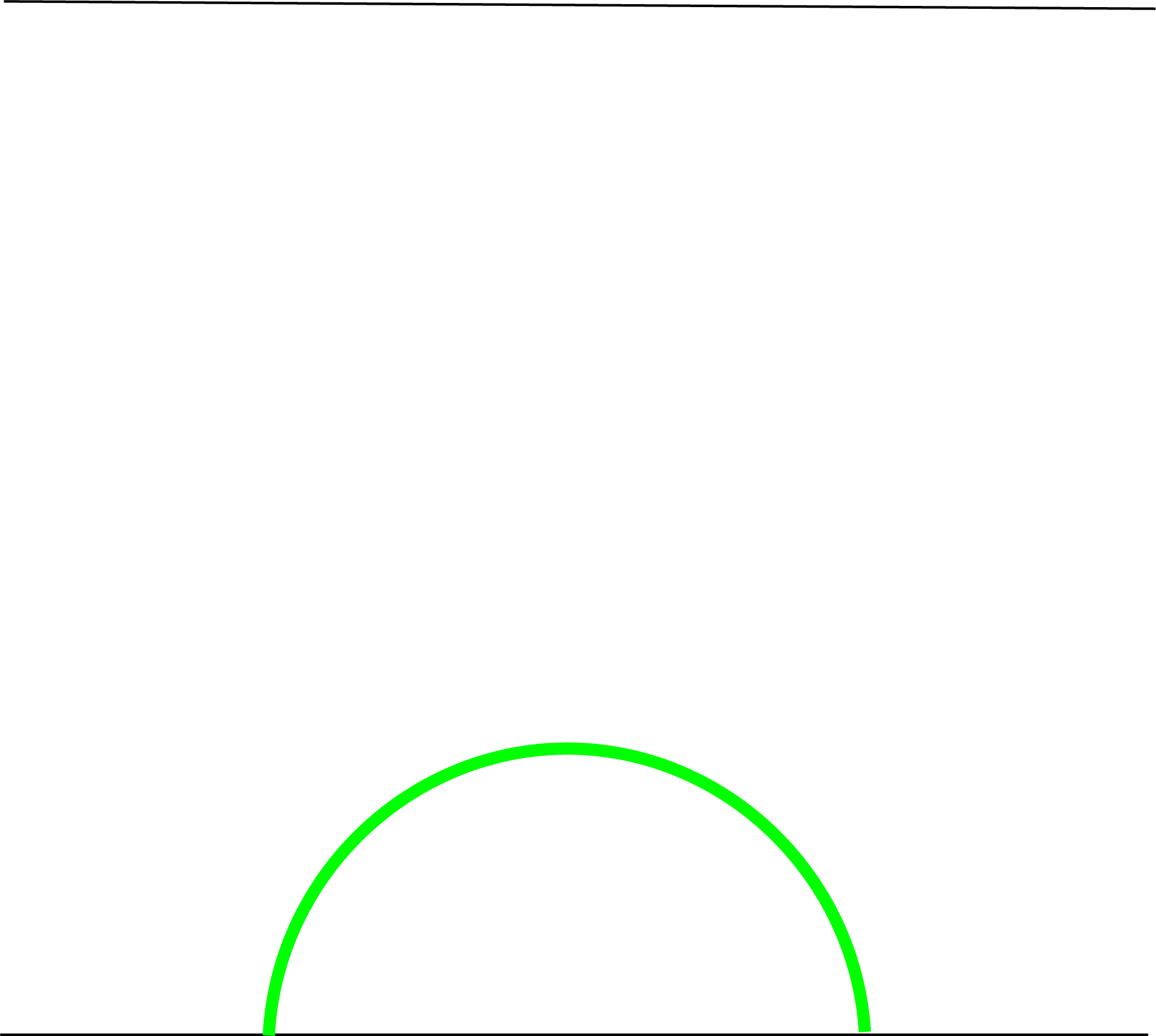}\put(10, 20){$\mapsto \capb  $}  \ \ \ \ \ \ \ \ \ \ \ \ \ \ \ \ \ \ \ \  \ \ \ \ \ \ \ \ \ \ \ \ \ \ \ \ \ \ \ \ 
\includegraphics[width=.10\textwidth]{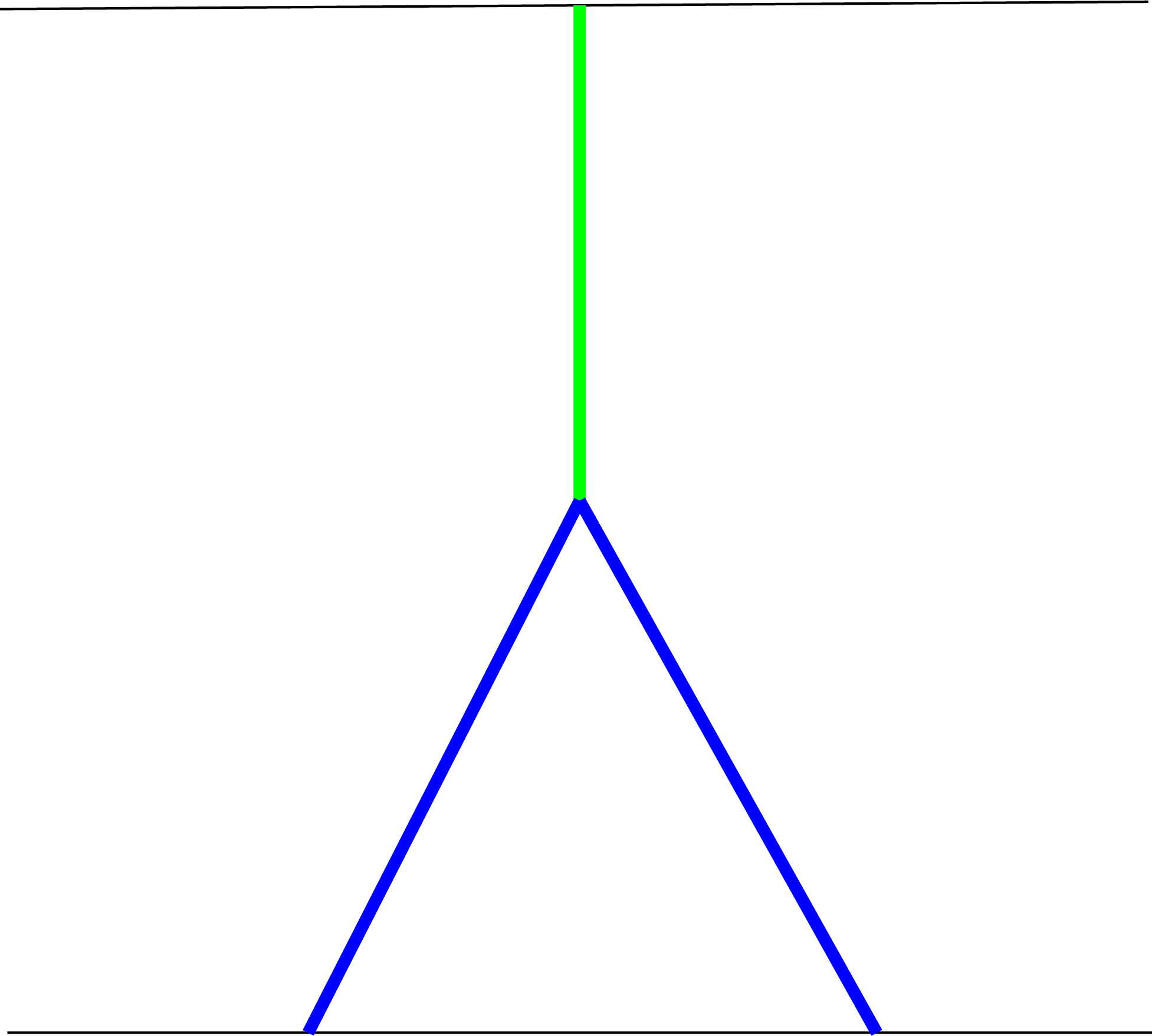}\put(10, 20){$\mapsto \pp$}
\end{figure}
\begin{figure}[H]
\centering
\includegraphics[width=.10\textwidth]{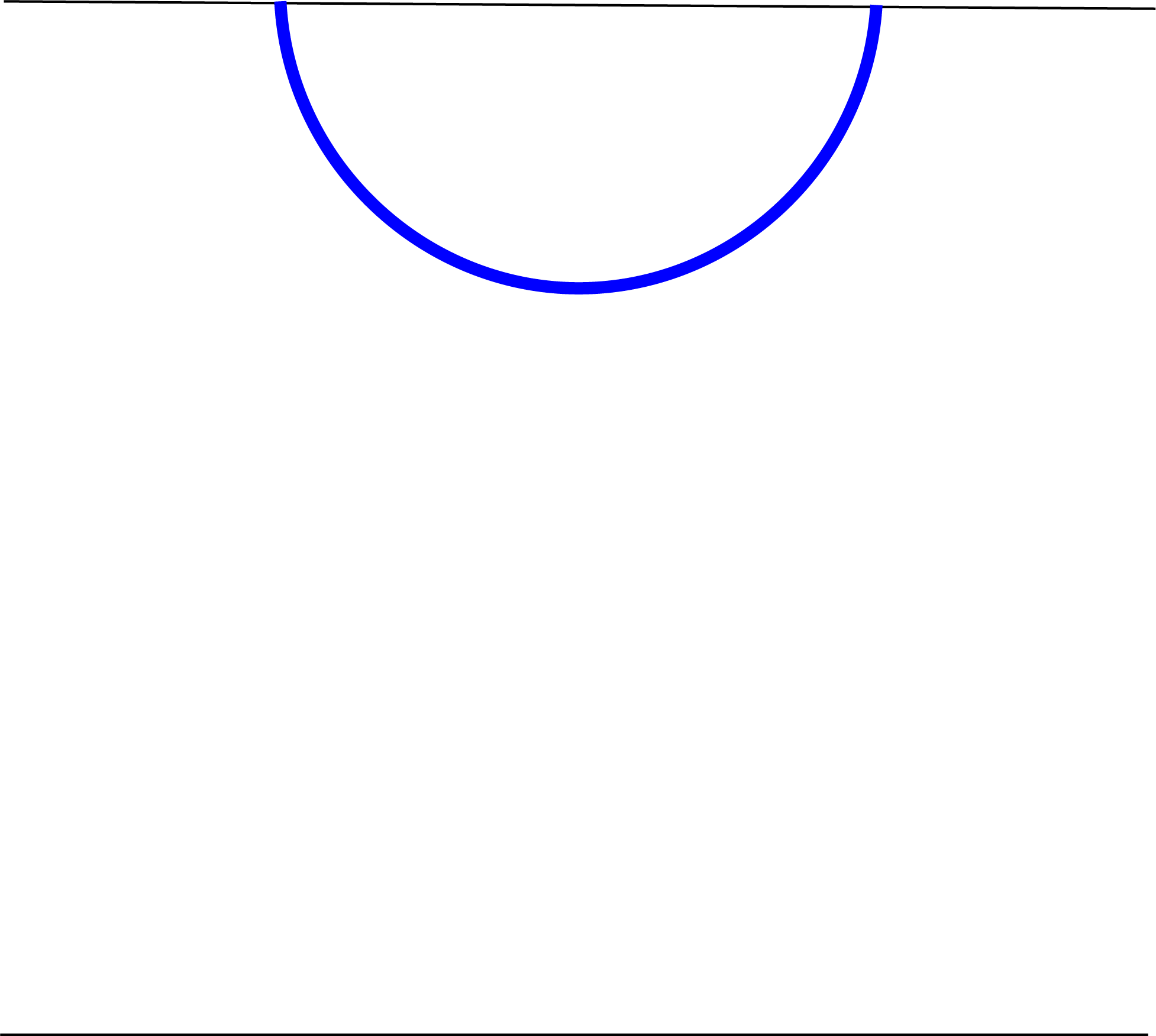}\put(10, 20){$\mapsto \cupa   $} \ \ \ \ \ \ \ \ \ \ \ \ \ \ \ \ \ \ \ \  \ \ \ \ \ \ \ \ \ \ \ \ \ \ \ \ \ \ \ \ 
\includegraphics[width=.10\textwidth]{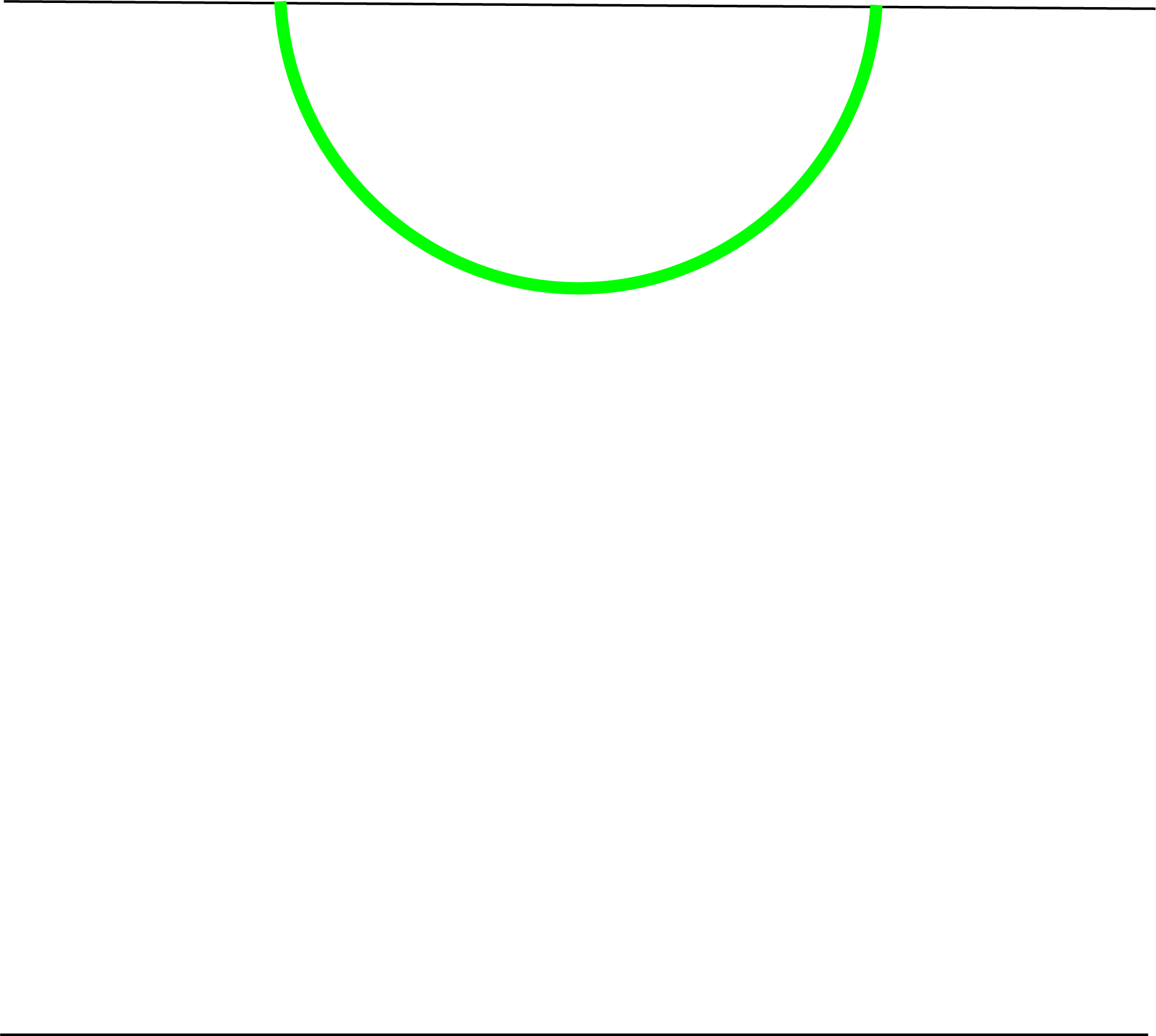}\put(10, 20){$\mapsto \cupb  $}  \ \ \ \ \ \ \ \ \ \ \ \ \ \ \ \ \ \ \ \  \ \ \ \ \ \ \ \ \ \ \ \ \ \ \ \ \ \ \ \ 
\includegraphics[width=.10\textwidth]{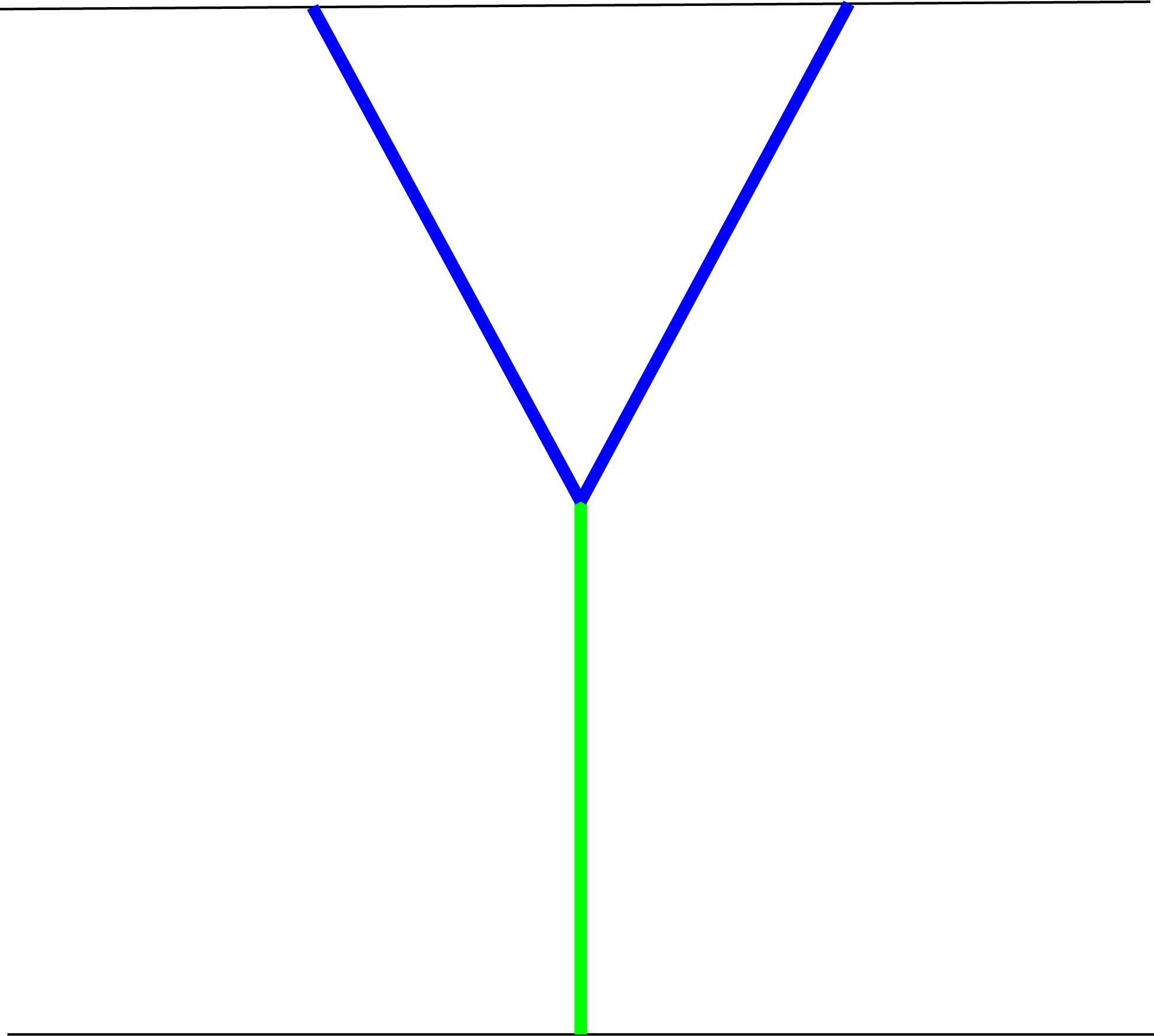}\put(10, 20){$\mapsto \ii$}
\end{figure}
\noindent and then extending $\mathcal{A}$-linearly so that horizontal concatenation of diagrams corresponds to tensor product of morphisms in $U_q^{\mathcal{A}}(\mathfrak{sp}_4)-\text{mod}$ and vertical composition of diagrams corresponds to composition of morphisms in $U_q^{\mathcal{A}}(\mathfrak{sp}_4)-\text{mod}$. 
\end{thm}

\begin{example}
We illustrate how $\eval$ is defined on objects and on morphisms:
\[
\eval(\blues \greent\greent) = V^{\mathcal{A}}(\blues\greent\greent) = V^{\mathcal{A}}(\varpi_1)\ot V^{\mathcal{A}}(\varpi_1) \ot V^{\mathcal{A}}(\varpi_2)
\]
\begin{figure}[H]
\centering
\hspace*{-9cm}\includegraphics[width=0.1\textwidth]{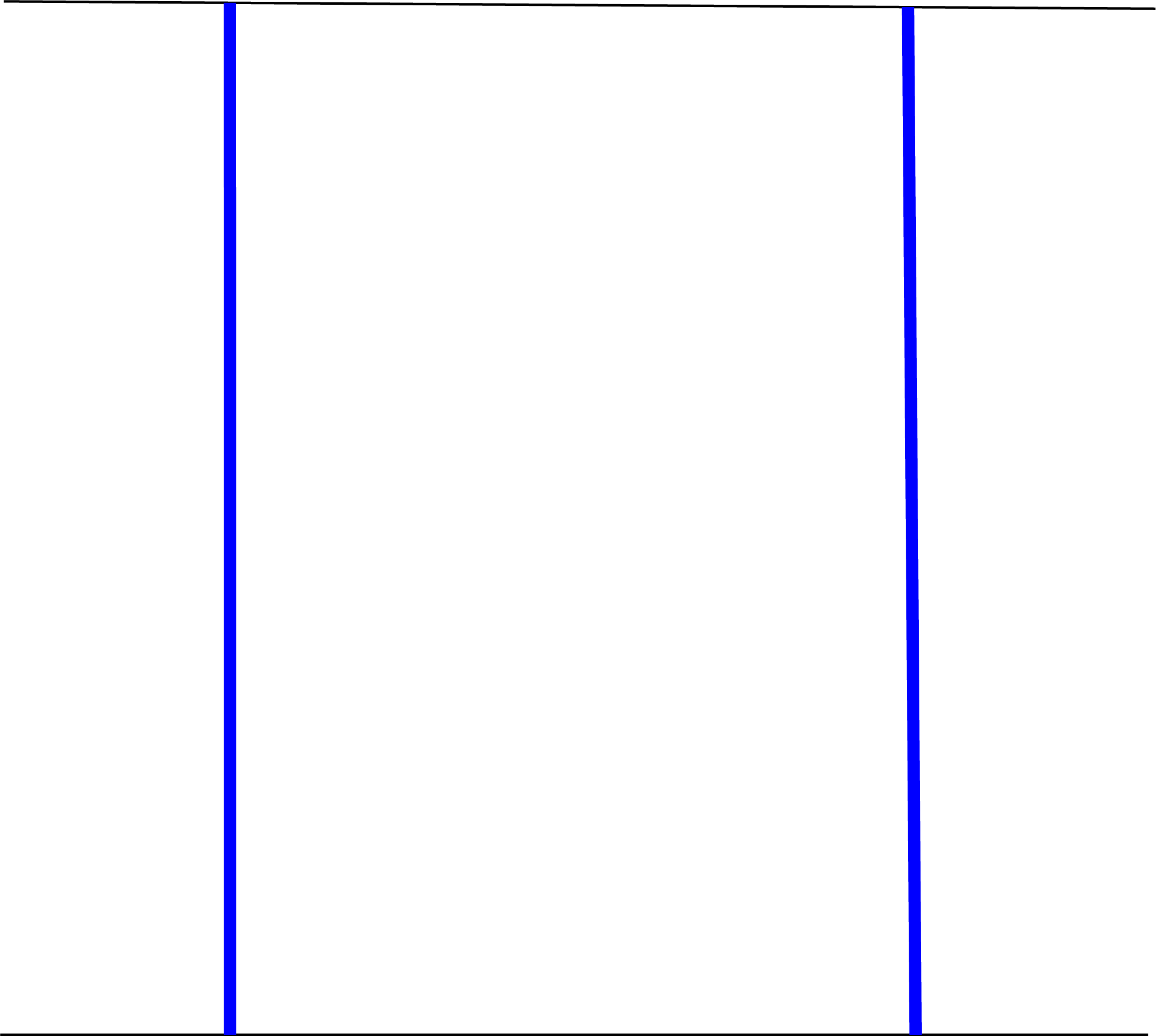}\put(5, 16){$-$}\ \ \ \ \ \ \ \ \includegraphics[width=0.1\textwidth]{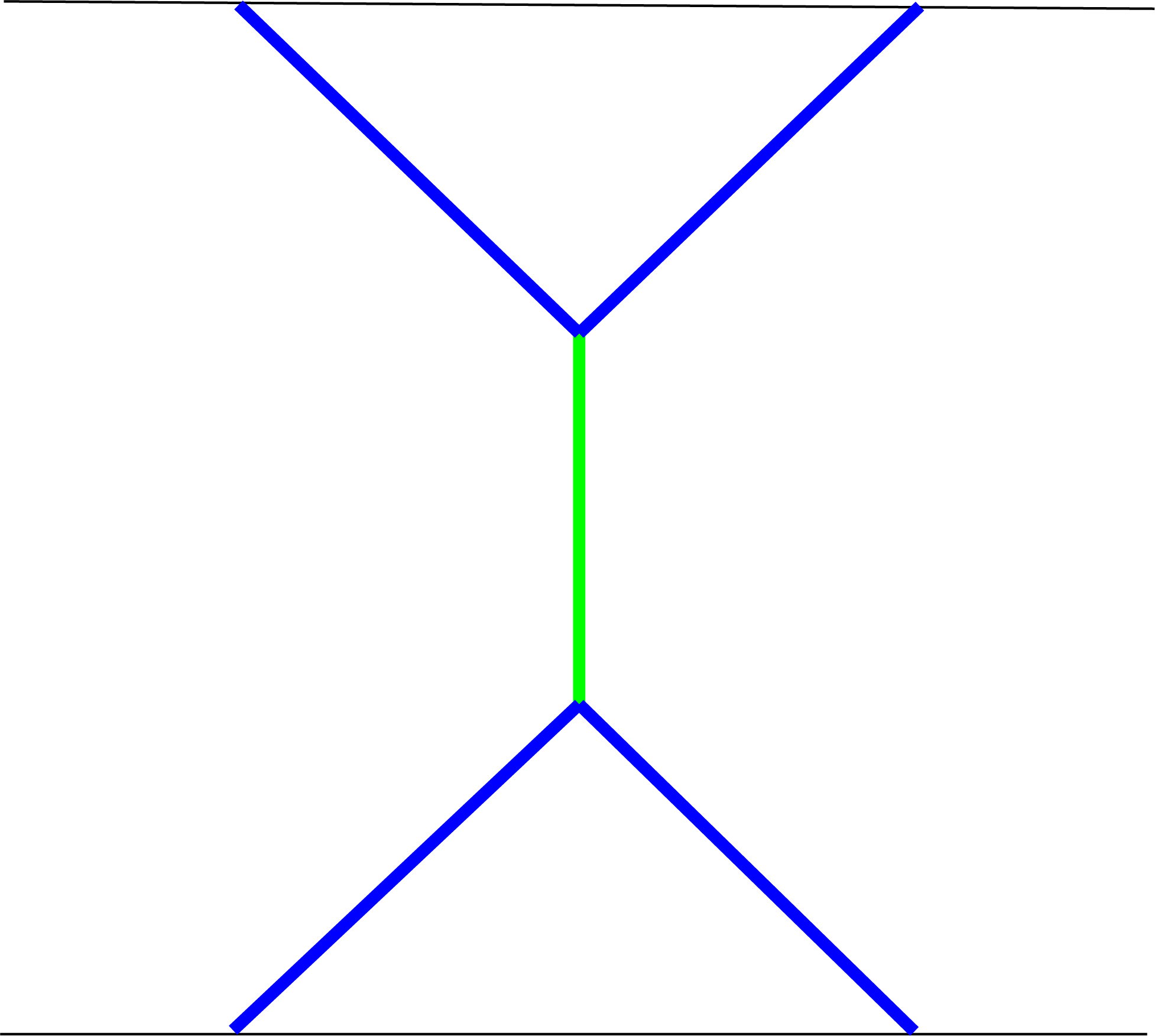}\put(5, 16){$+\dfrac{1}{[2]_q}$} \ \ \ \ \ \ \ \ \ \ \ \ \includegraphics[width=0.1\textwidth]{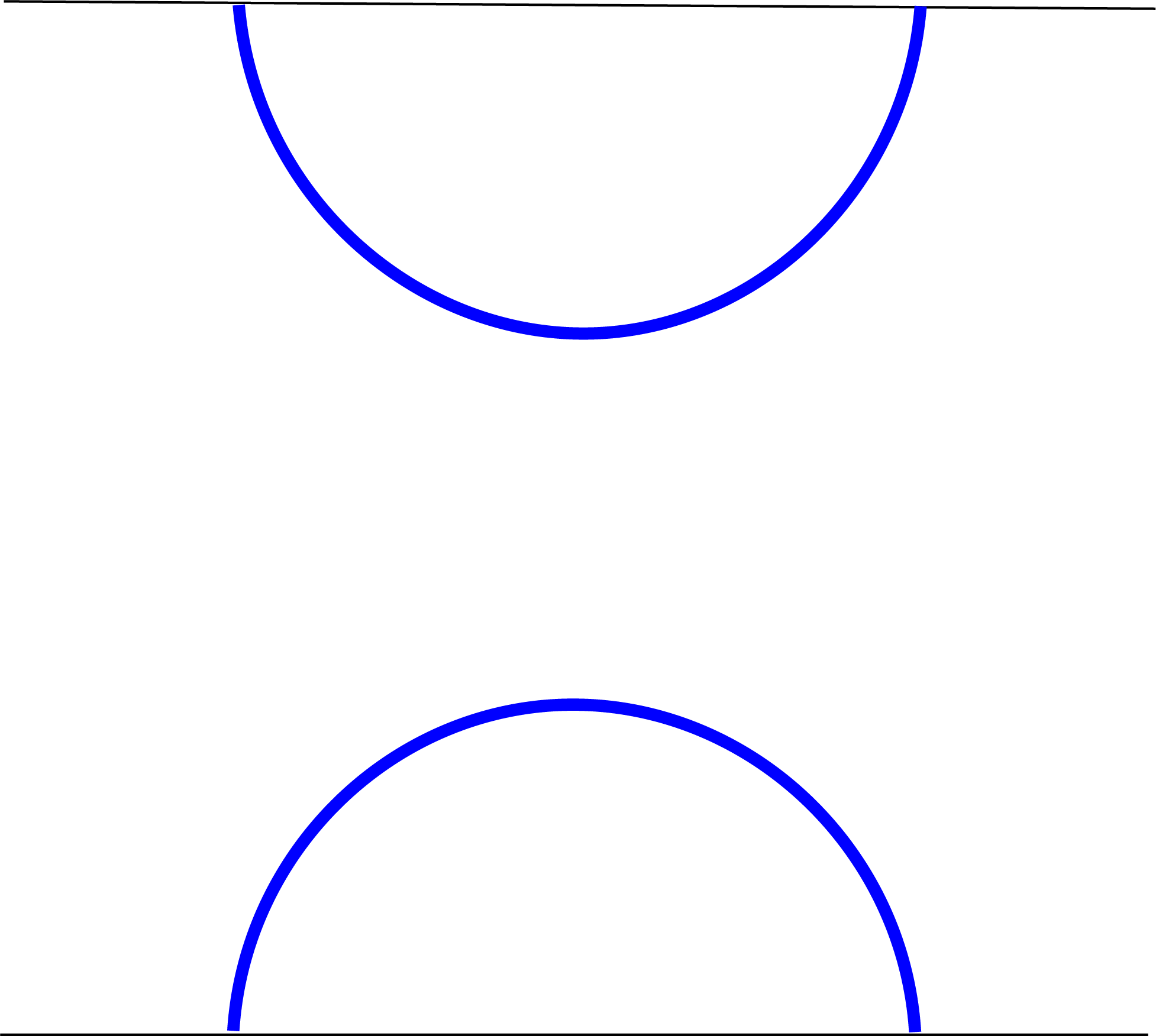} \put(12, 16){$\xmapsto{\eval}\id \ot \id- \ii\circ \pp + \dfrac{1}{[2]_q}\cupa  \circ \capa $.}
\end{figure}
\end{example}

%===========
\subsection{Checking Relations}
\label{subsec-relationcheck}
%===========
Since $\DD$ is defined by generators and relations, in order to verify the theorem we must check that the diagrammatic relations hold in $U_q^{\mathcal{A}}(\mathfrak{sp}_4)-\text{mod}$.

\begin{proof}[Proof of Theorem \ref{evalthm}]

The isotopy relations follow from \eqref{cupcapisotopy} and the equality of \eqref{LHSp} and \eqref{RHSp}. 

To verify the relation 
\begin{figure}[H]
\centering
\includegraphics[width=0.10\textwidth]{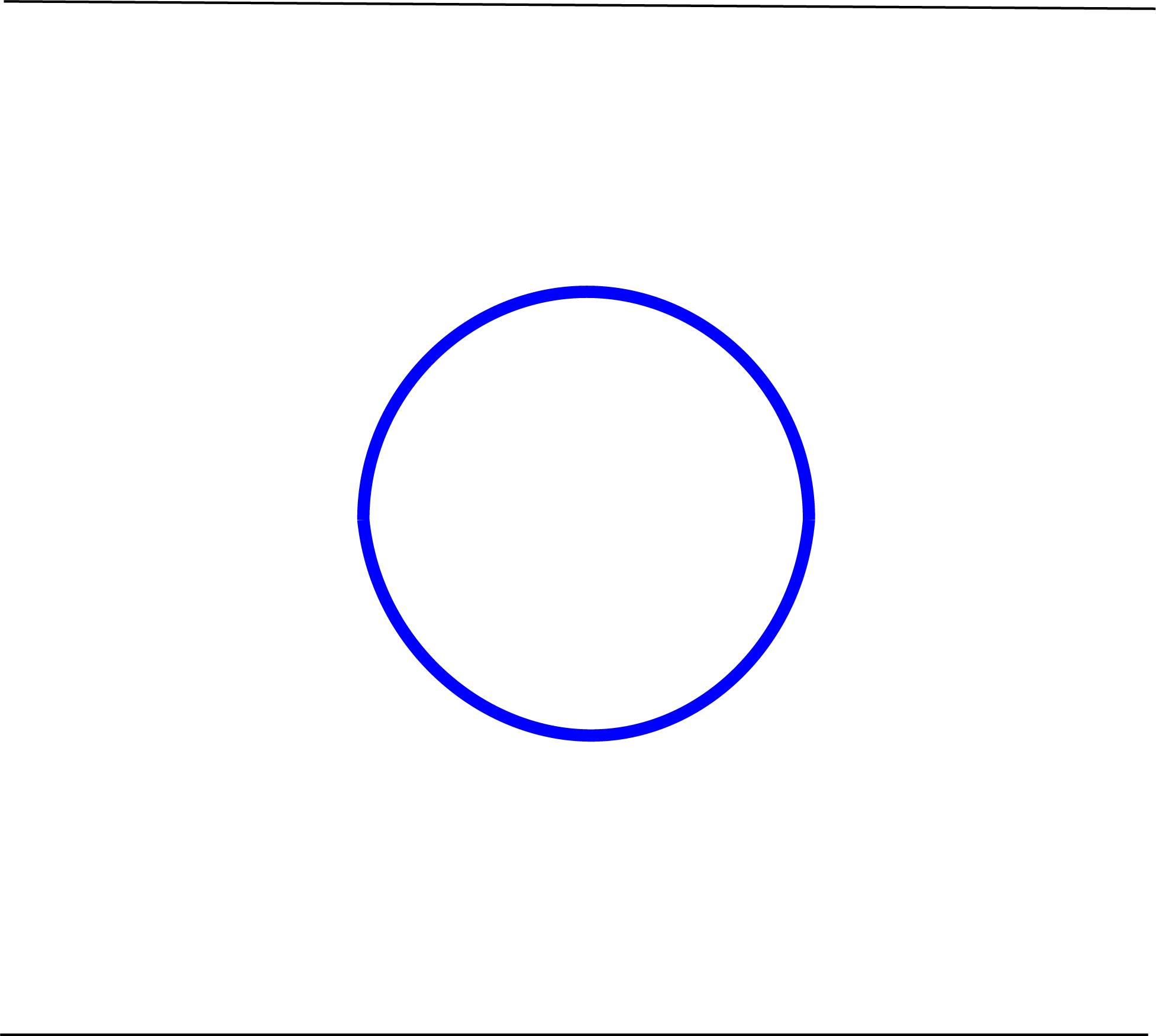}\put(5, 16){$= -\dfrac{[6]_q[2]_q}{[3]_q}$} \ \ \ \ \ \ \ \ \ \ \ \ \ \ \ \ \ \ \ \ \ \ \ \ \includegraphics[width = .10\textwidth]{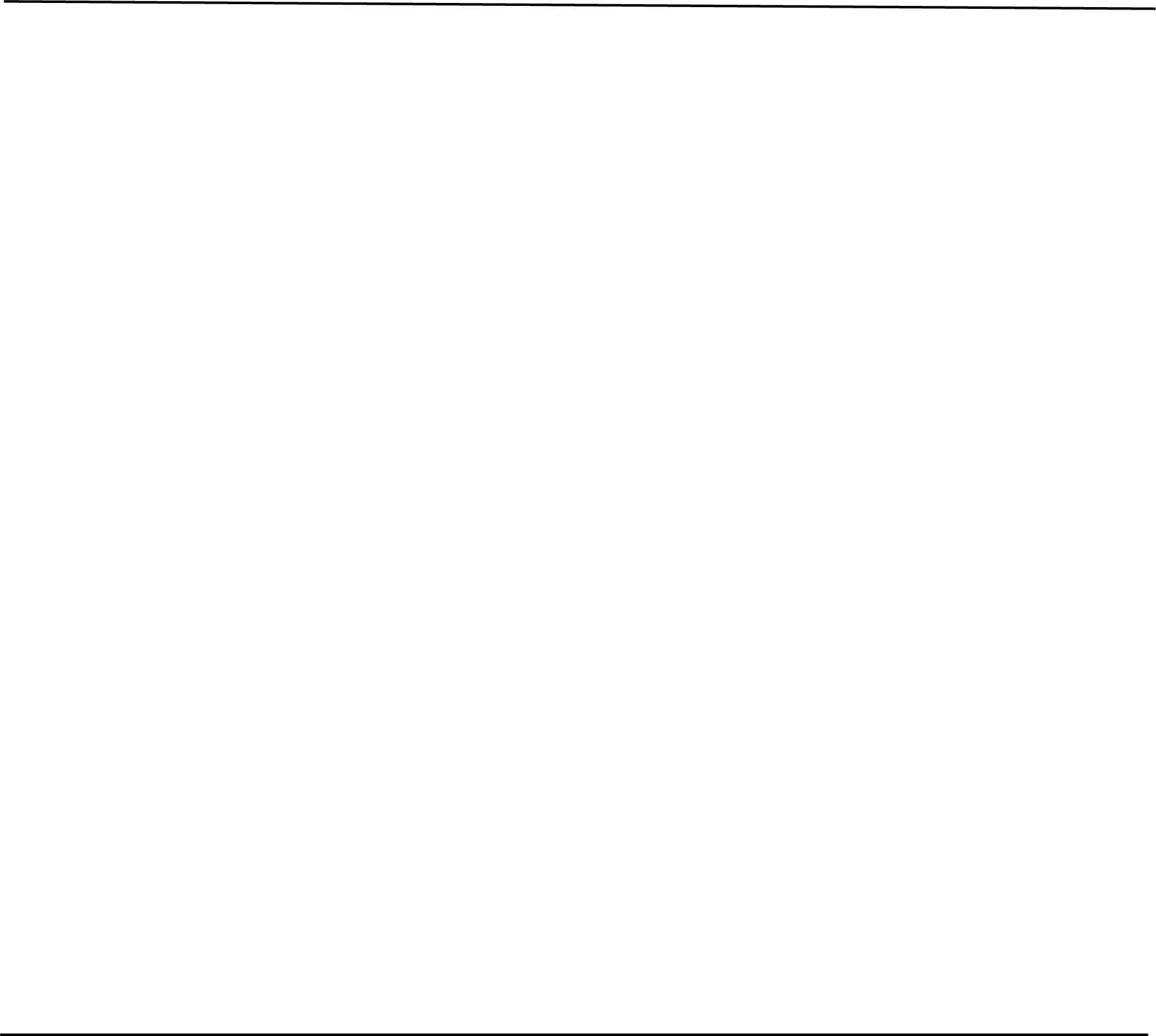}
\end{figure}
\noindent it suffices to show that
\begin{equation}\label{deseqfromqnums}
\capa   \circ \cupa  (1) = -\dfrac{[6]_q[2]_q}{[3]_q}.
\end{equation}
\noindent Using \eqref{cup1formula} and \eqref{cap1formula} we find
\begin{equation}
\capa   \circ \cupa  (1) = -q^{-4}\cdot 1 + q^{-3} \cdot (-q) - q^{-1} \cdot q^3 + 1\cdot (-q^4) = -\left([5]_q- [1]_q\right).
\end{equation}
The desired equality \eqref{deseqfromqnums} comes from the quantum number calculation in \eqref{qdimcalcone}.

One can similarly argue that the relation
\begin{figure}[H]
\centering
\includegraphics[width=0.10\textwidth]{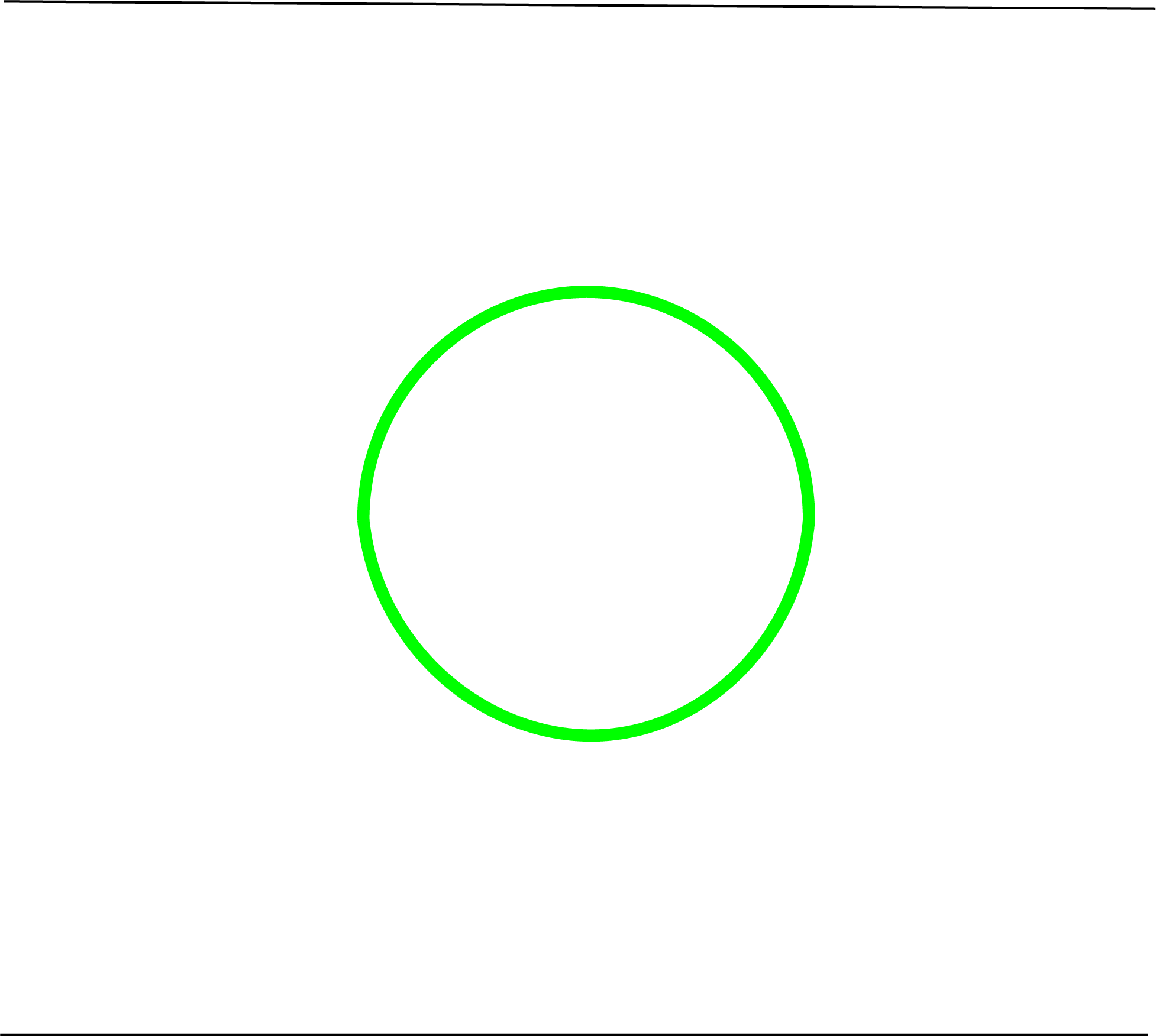}\put(5, 16){$= \dfrac{[6]_q[5]_q}{[3]_q[2]_q}$} \ \ \ \ \ \ \ \ \ \ \ \ \ \ \ \ \ \ \ \ \ \includegraphics[width = .10\textwidth]{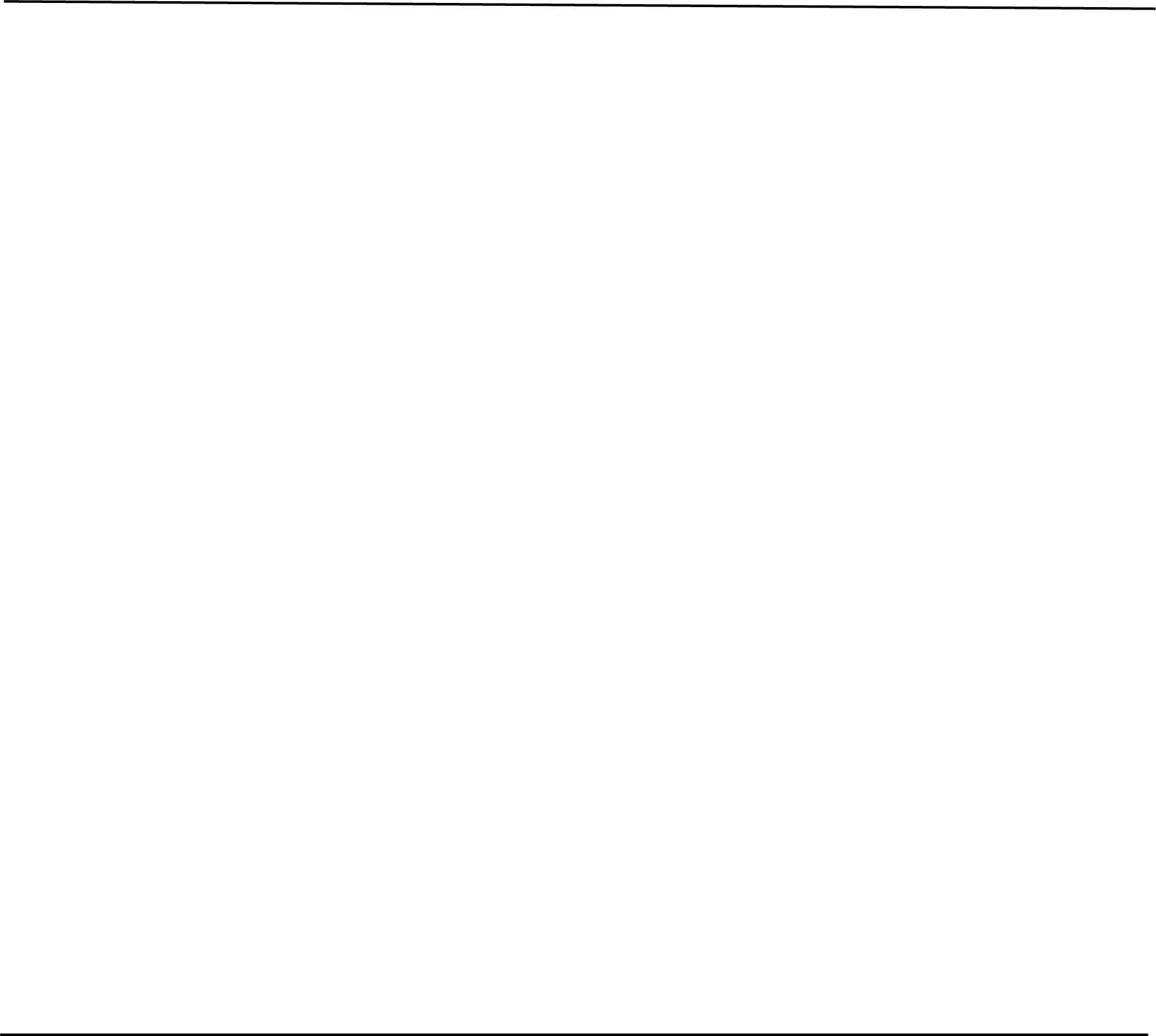}
\end{figure}
\noindent is satisfied. Use \eqref{cup2formula} and \eqref{cap2formula} to compute
\begin{equation}
\capb  \circ \cupb  (1) = q^{-6} + q^{-2} + 1+ q^2+ q^6 = [7]_q- [5]_q+[3]_q,
\end{equation}
then use \eqref{qdimcalctwo} to deduce 
\begin{equation}
\capb  \circ \cupb  (1) = \dfrac{[6]_q[5]_q}{[3]_q[2]_q}.
\end{equation}

To check the monogon relation
\begin{figure}[H]
\centering
\includegraphics[width=0.10\textwidth]{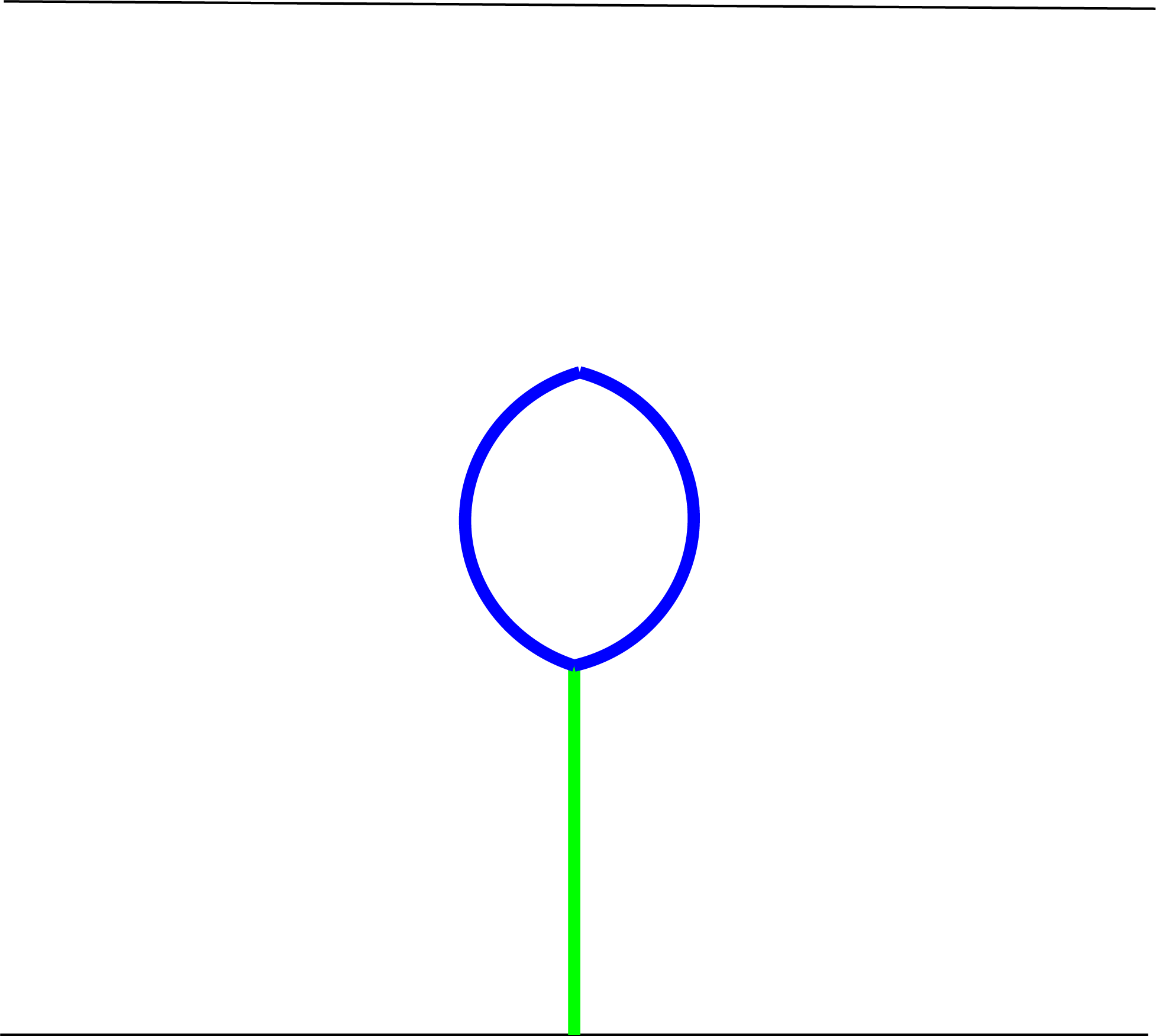}\put(5, 15){$=0$}
\end{figure}
\noindent and the bigon relation
\begin{figure}[H]
\centering
\includegraphics[width=0.10\textwidth]{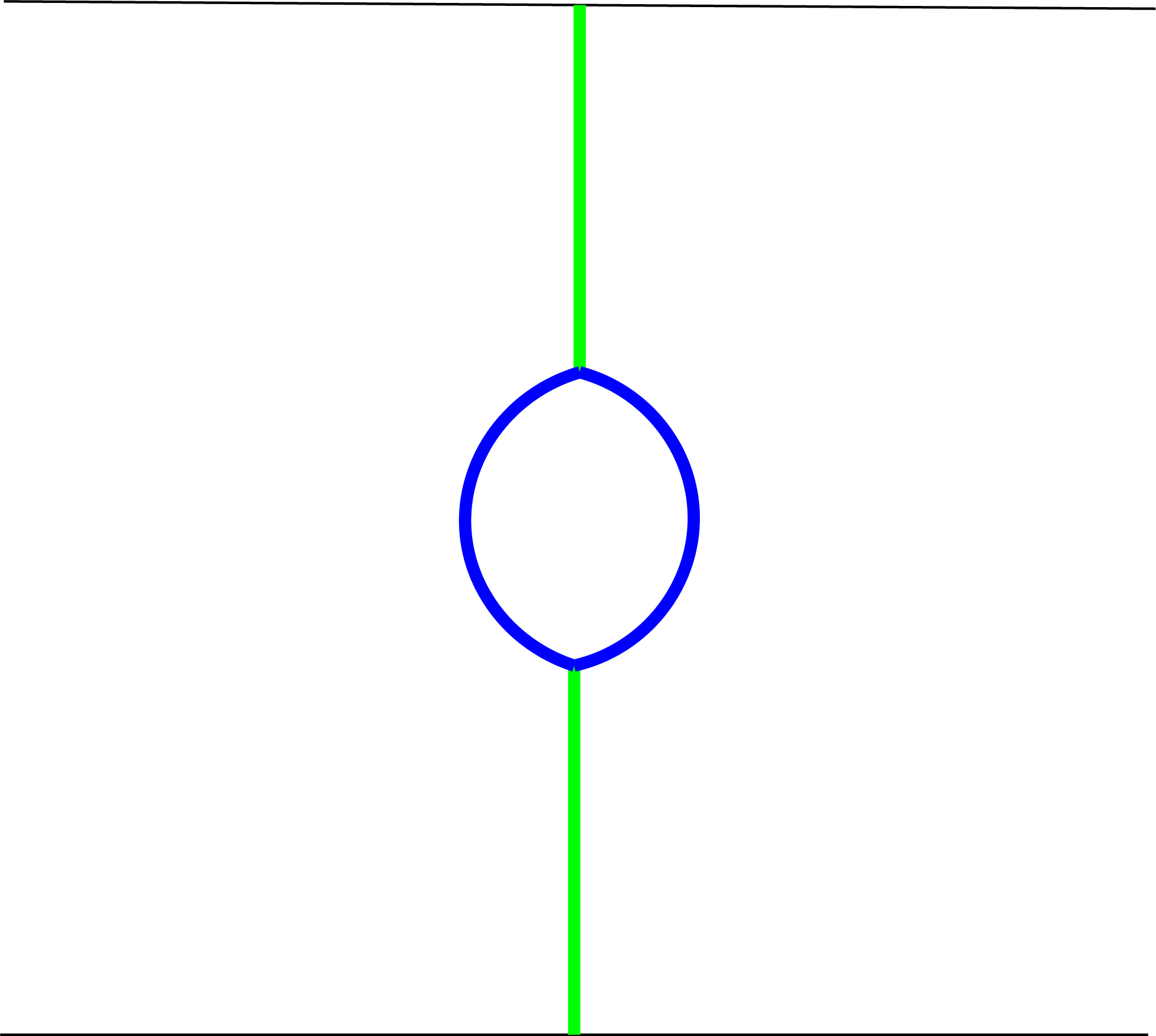}\put(5, 15){$=-[2]_q$} \ \ \ \ \ \ \ \ \ \ \ \ \ \ \ \ \ \includegraphics[width = .10\textwidth]{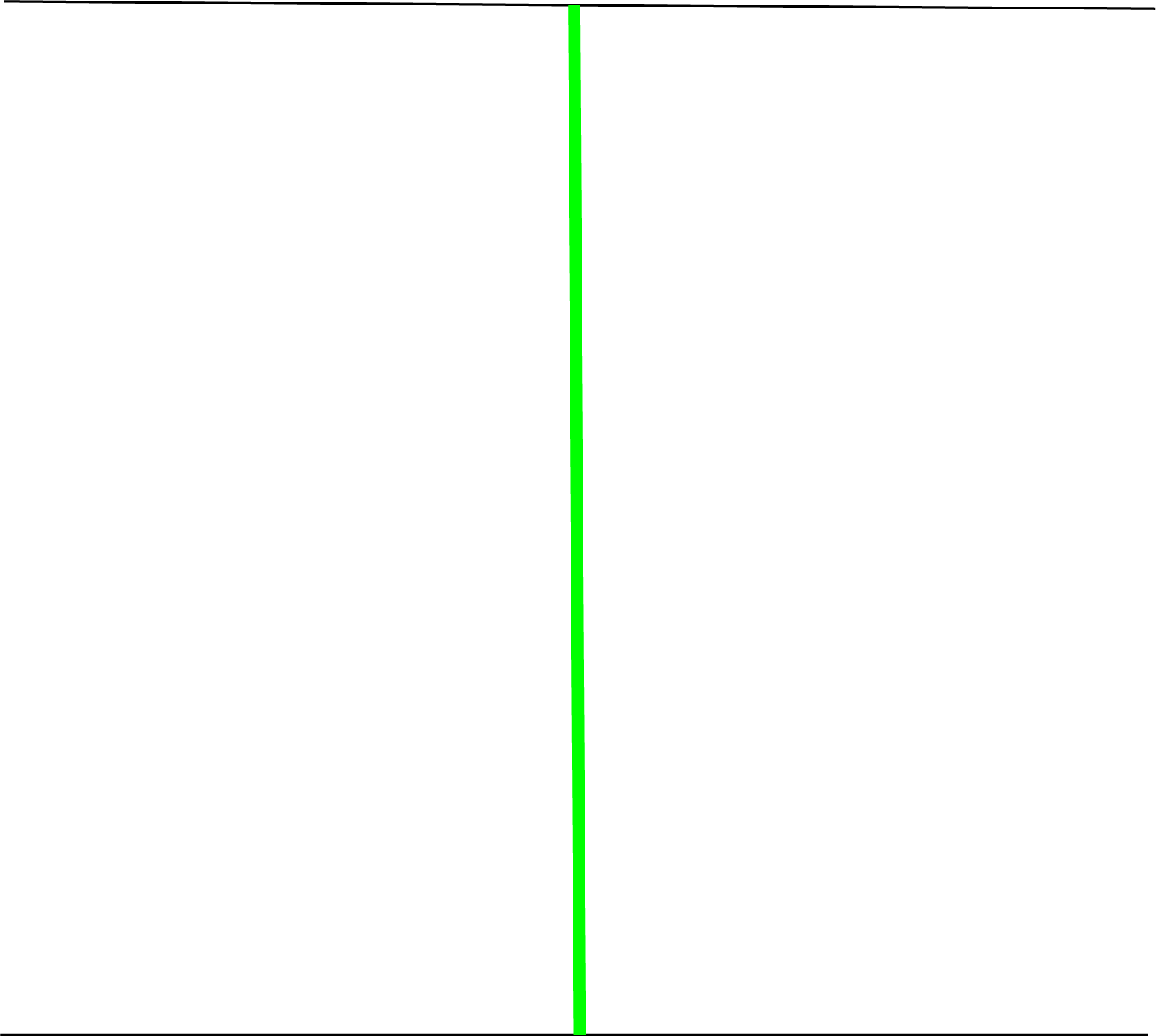}
\end{figure}

\noindent we need to show $\capa   \circ \ii = 0$ and $\pp\circ \ii= -[2]_q\id$ respectively. Since the module $V^{\mathcal{A}}(\varpi_2)$ is generated by the highest weight vector $v_{(0, 1)}$ it suffices to show that $\capa   \circ \ii(v_{(0, 1)})= 0$ and $\pp\circ \ii(v_{(0, 1)}) = -[2]_q v_{(0, 1)}$. The calculations go as follows:
\begin{equation}
\begin{split}
\capa   \circ \ii(v_{(0, 1)})&\stackrel{\eqref{iformula}}{=}  \capa (q^{-1}v_{(1, 0)}\ot v_{(-1, 1)} - v_{(-1, 1)}\ot v_{(1, 0)}) \\
&\stackrel{\eqref{cup1formula}}{=} 0. 
\end{split}
\end{equation}
and
\begin{equation}
\begin{split}
\pp   \circ \ii(v_{(0, 1)})&\stackrel{\eqref{iformula}}{=}  \pp(q^{-1}v_{(1, 0)}\ot v_{(-1, 1)} - v_{(-1, 1)}\ot v_{(1, 0)}) \\
&\stackrel{\eqref{formulaforpi}}{=} -q^{-1}v_{(0,1)} + qv_{(0, 1)} \\
&\stackrel{}{=} -[2]_qv_{(0, 1)}. 
\end{split}
\end{equation}

Verifying the trigon relation
\begin{figure}[H]
\centering
\includegraphics[width=0.10\textwidth]{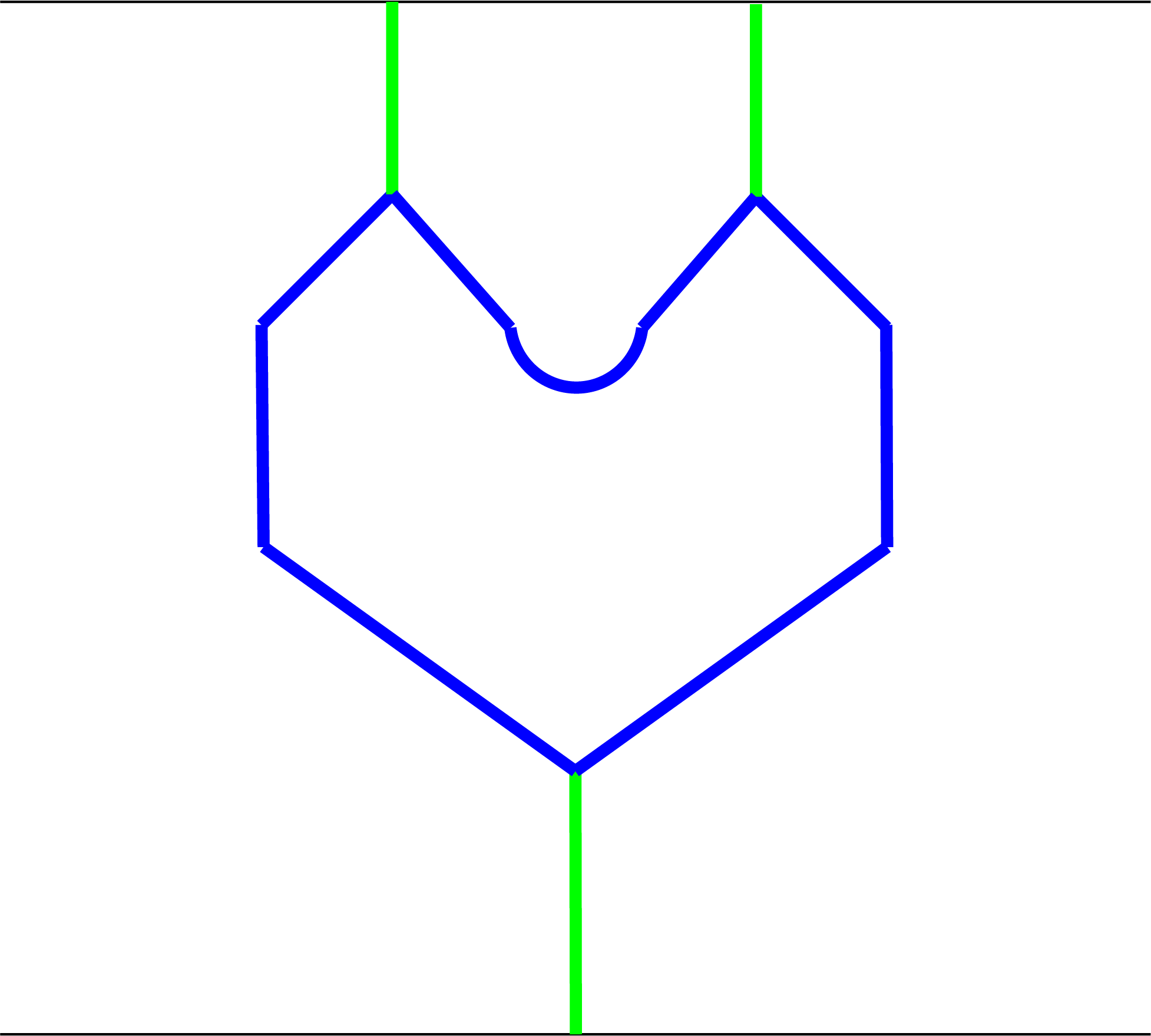}\put(5, 15){$=0$}
\end{figure}
\noindent
\noindent is left as an exercise (Hint: apply $(\pp\ot \pp )\circ(\id \ot \cupa \ot \id)\circ \ii$ to the vector $v_{(0,1)}$ and use \eqref{iformula} and \eqref{cup1formula} and \eqref{formulaforpi}).

Now we endeavor to check 
\begin{figure}[H]
\centering
\includegraphics[width=0.1\textwidth]{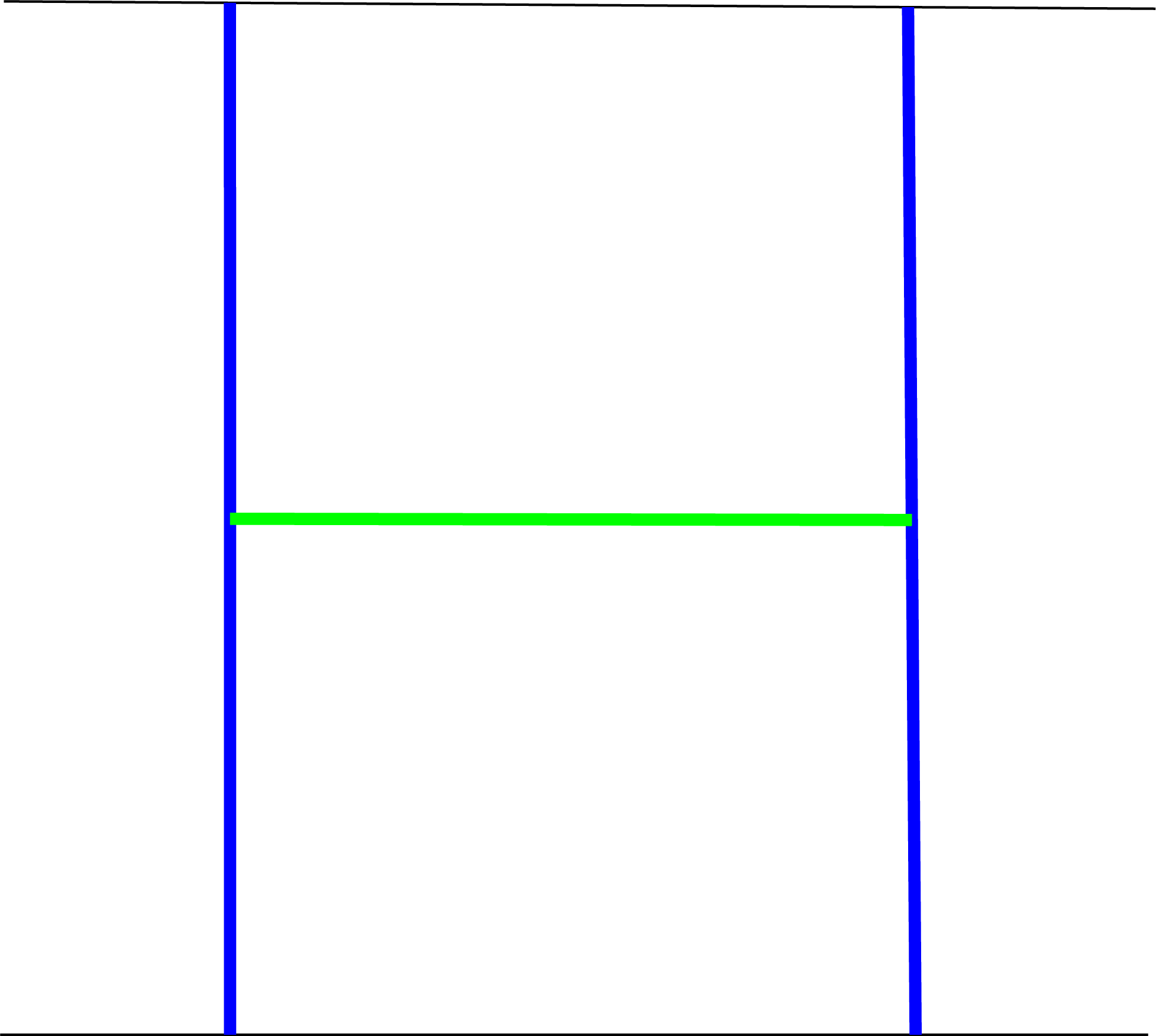}\put(5, 16){$=\dfrac{1}{[2]_q}$} \ \ \ \ \ \ \ \ \ \ \ \ \includegraphics[width=0.1\textwidth]{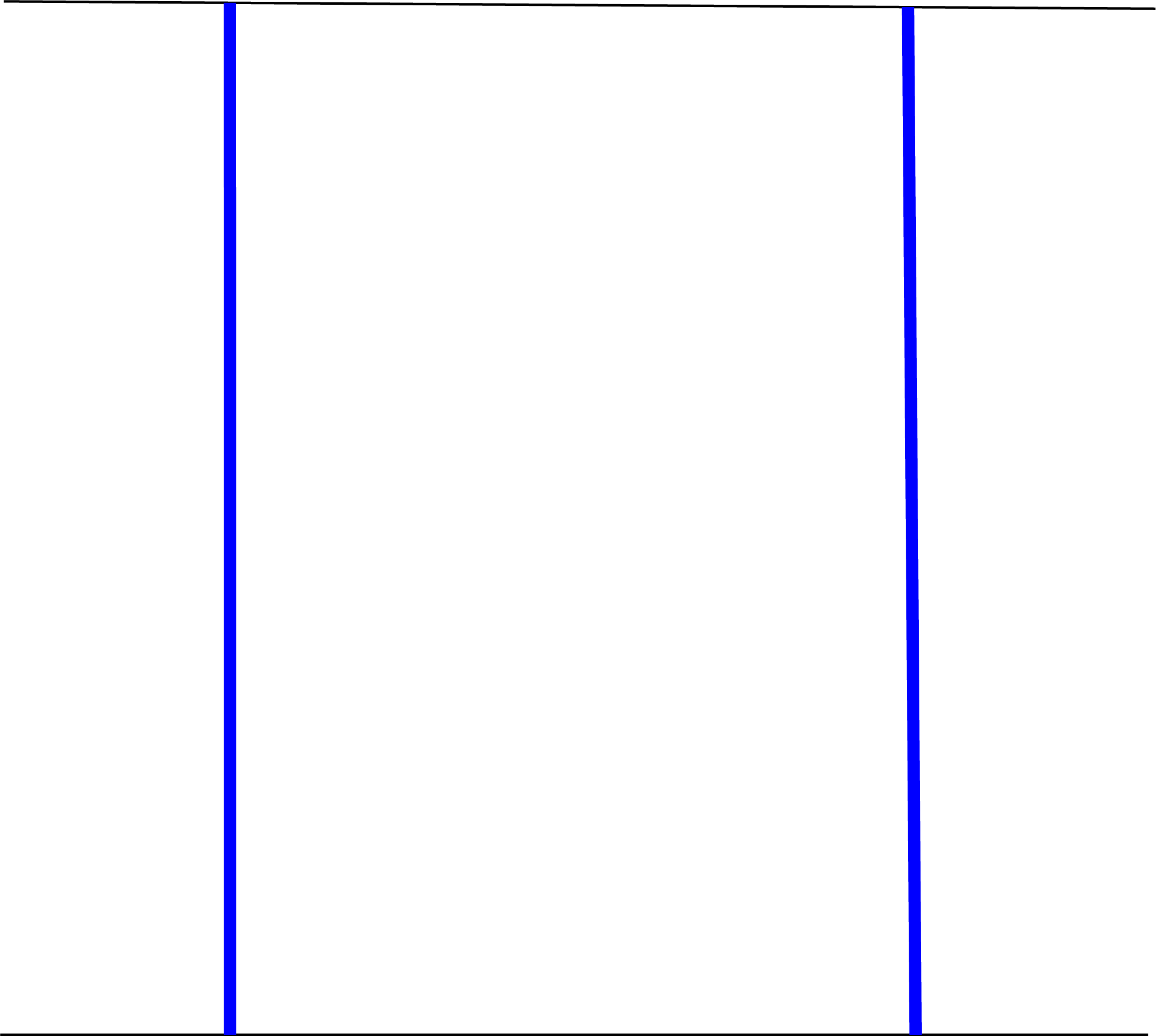}\put(5, 16){$+$}\ \ \ \ \ \ \ \ \ \ \ \ \includegraphics[width=0.1\textwidth]{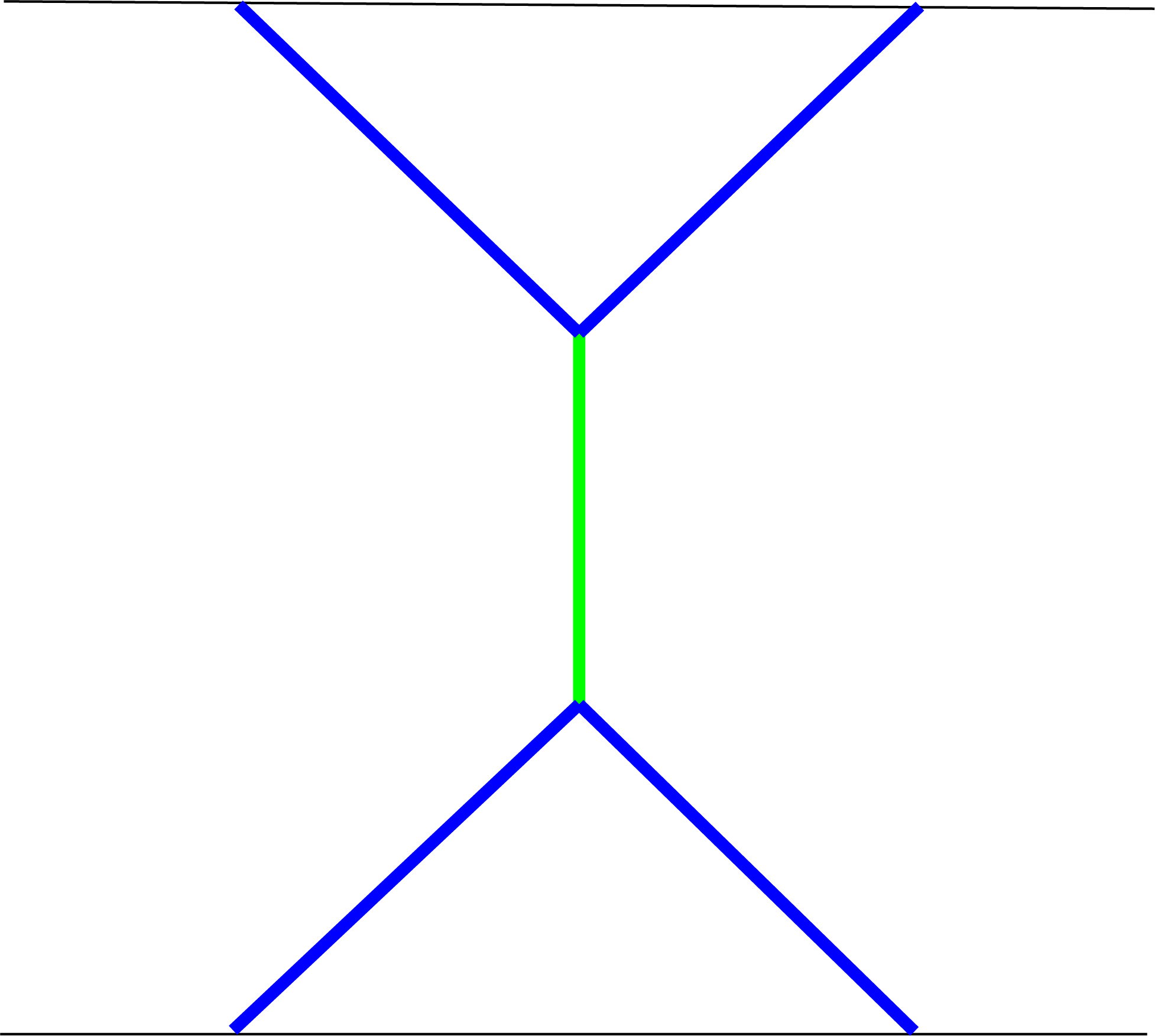}\put(5, 16){$-\dfrac{1}{[2]_q}$} \ \ \ \ \ \ \ \ \ \ \ \ \includegraphics[width=0.1\textwidth]{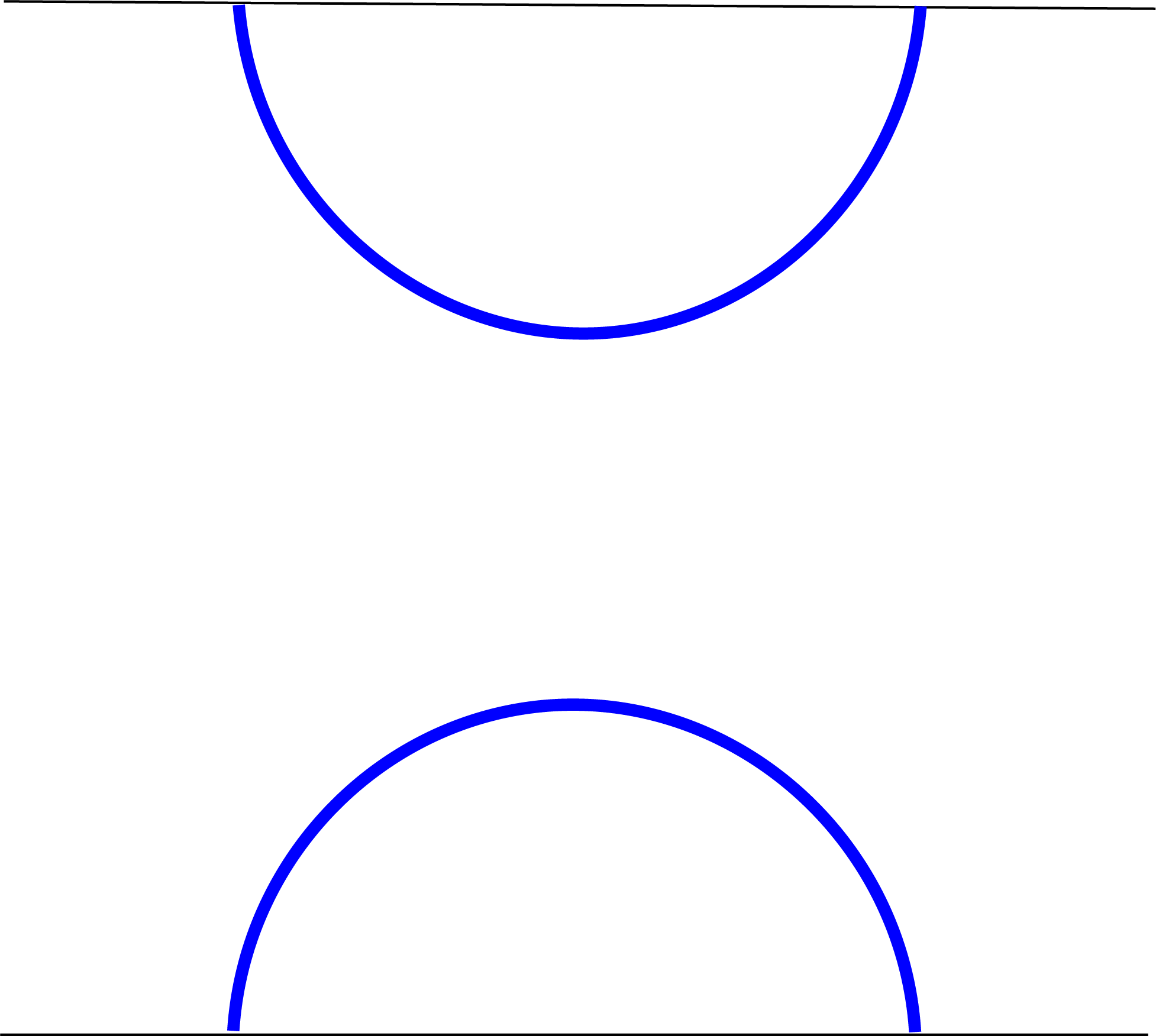}
\end{figure}
\noindent Precomposing with $\id\ot \cupa  $ is an $\mathcal{A}$-linear map
\begin{equation}
\Hom_{U_q^{\mathcal{A}}(\mathfrak{sp}_4)}(V^{\mathcal{A}}(\varpi_1)^{\ot 2}, V^{\mathcal{A}}(\varpi_1)^{\ot 2}) \longrightarrow \Hom_{U_q^{\mathcal{A}}(\mathfrak{sp}_4)}(V^{\mathcal{A}}(\varpi_1), V^{\mathcal{A}}(\varpi_1)^{\ot 3}),
\end{equation}
while postcomposing with $\id\ot\id \ot \capa   $ is an $\mathcal{A}$-linear map in the other direction. From \eqref{cupcapisotopy} it follows that the two maps are mutually inverse isomorphisms of $\mathcal{A}$-modules, so we can instead check the following relation.
\begin{figure}[H]
\centering
\includegraphics[width=0.1\textwidth]{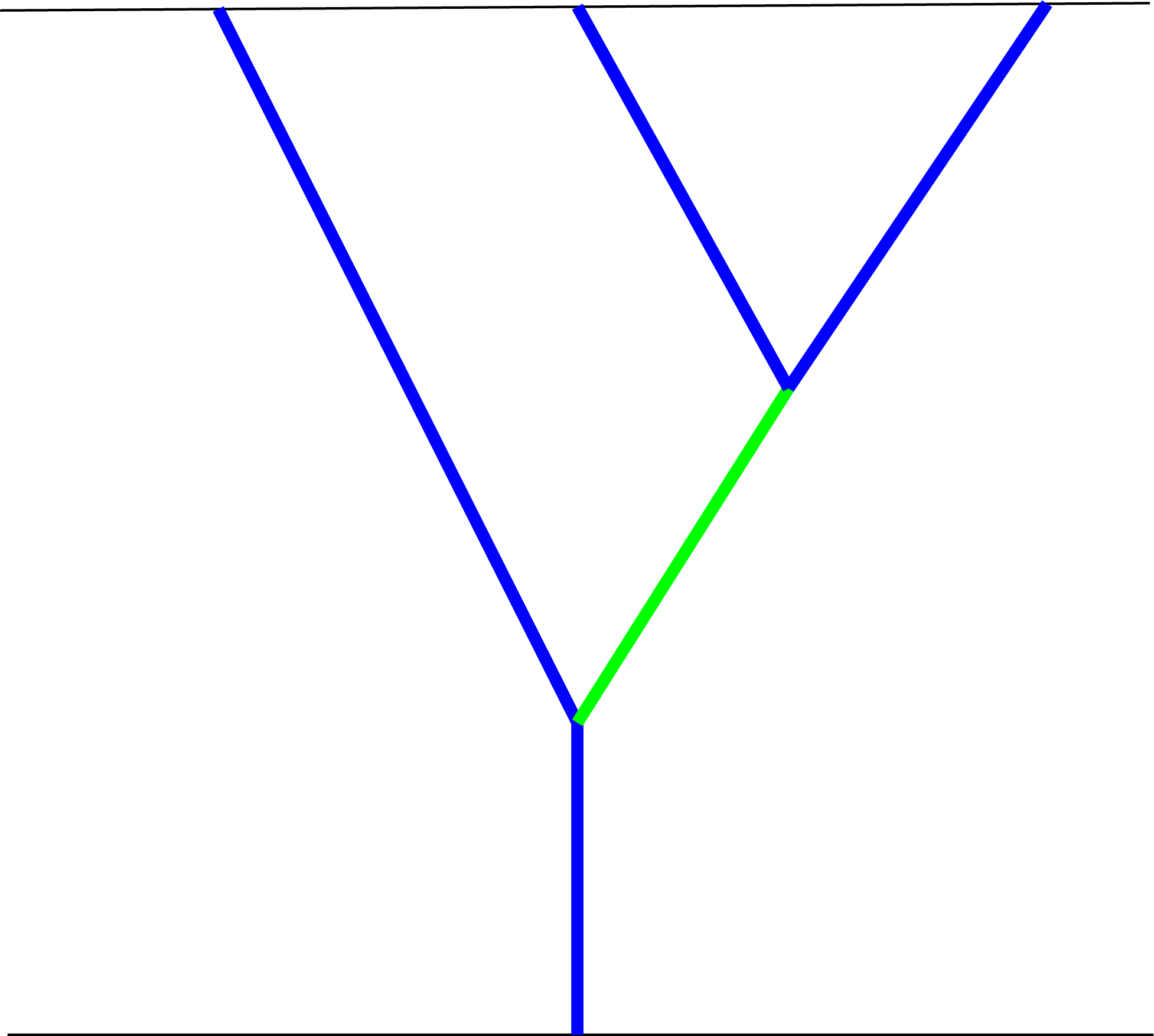}\put(-60, 16){$[2]_q$}\put(5, 16){$-[2]_q$} \ \ \ \ \ \ \ \ \ \includegraphics[width=0.1\textwidth]{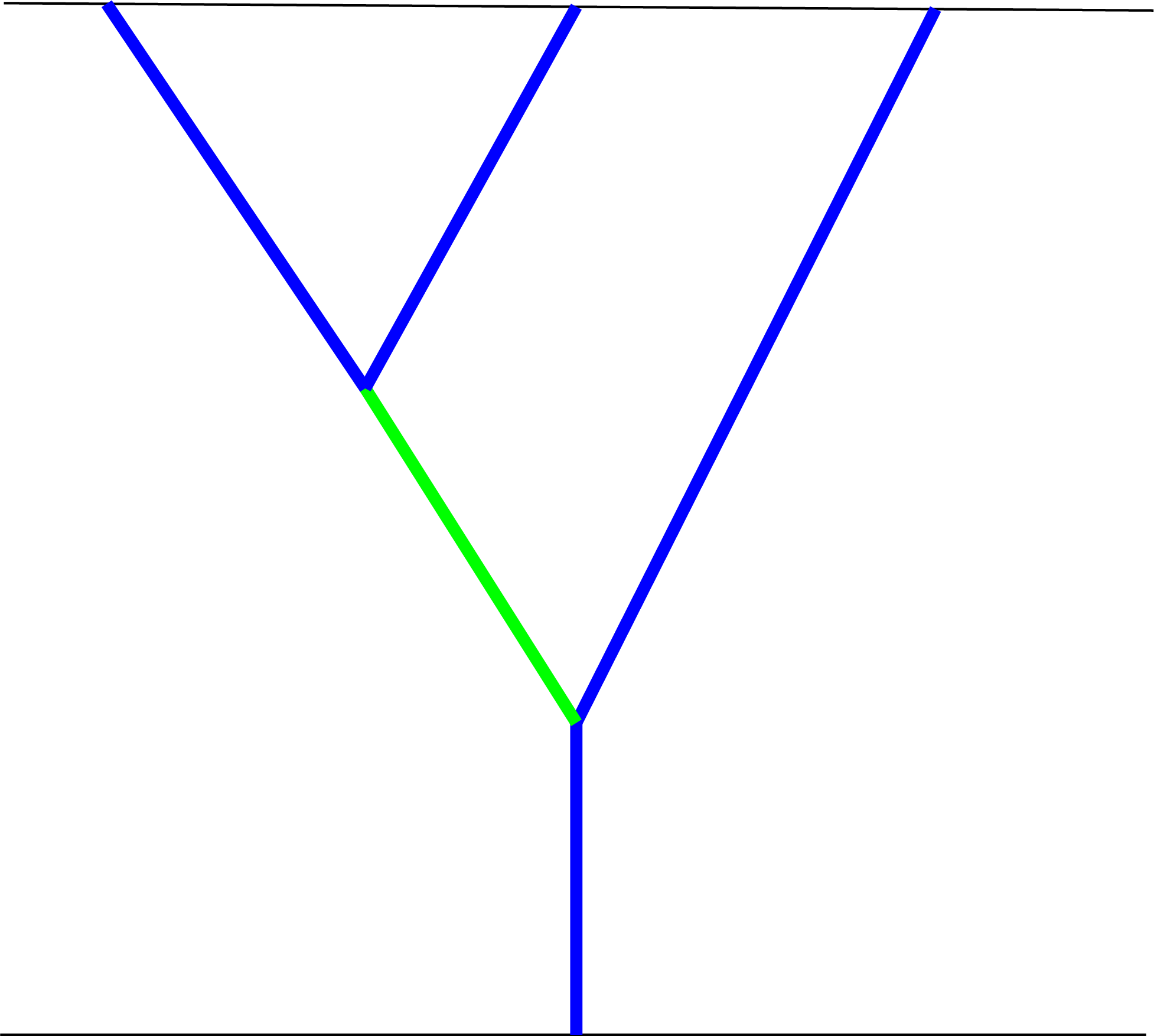}\put(5, 16){$=$}\ \ \ \ \ \ \ \ \ \ \ \ \ \ \ \ \ \includegraphics[width=0.1\textwidth]{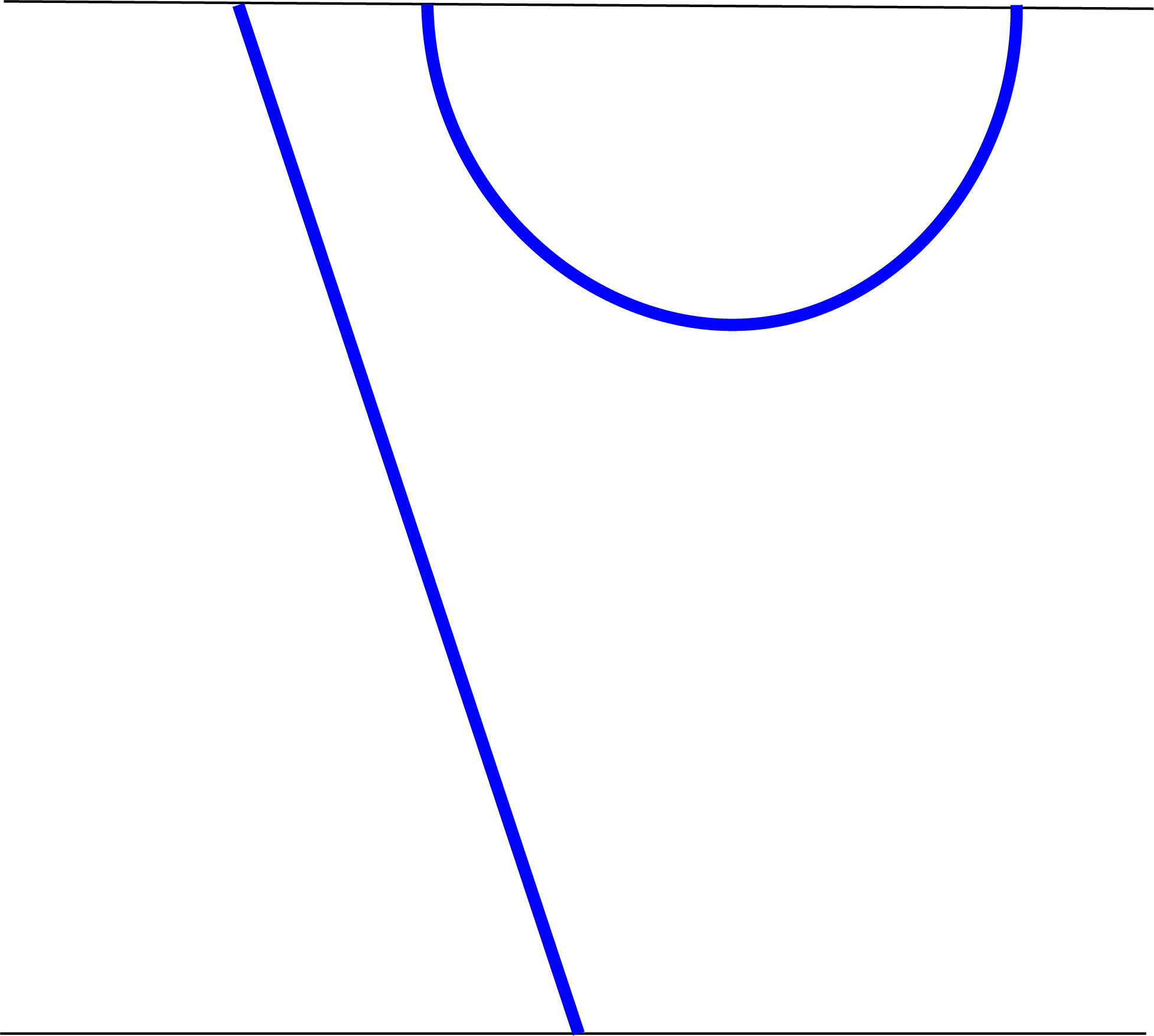}\put(5, 16){$-$} \ \ \ \ \ \ \ \ \ \ \includegraphics[width=0.1\textwidth]{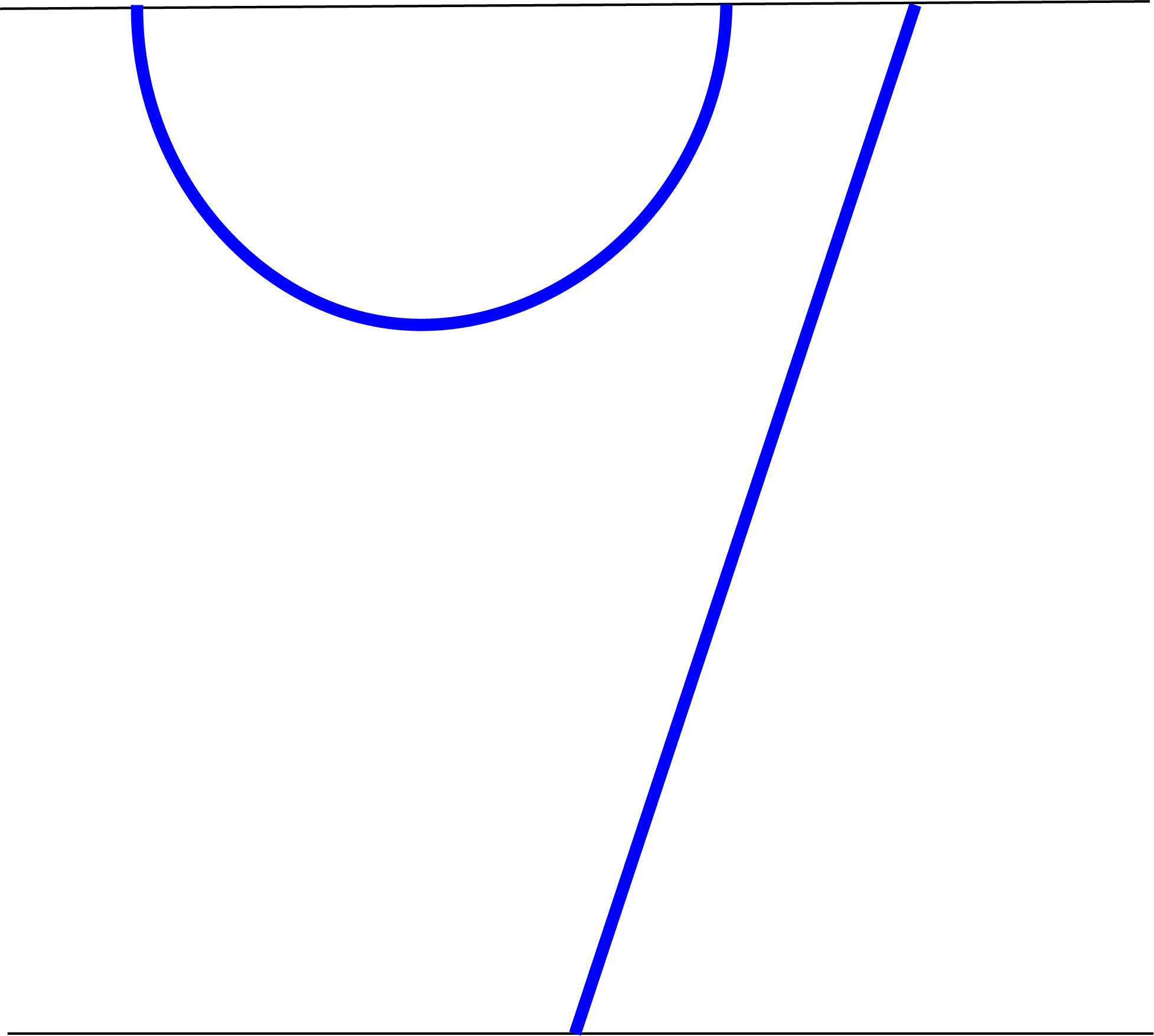}
\end{figure}
From the discussion in remark \eqref{isotopyrmk} it follows that we need to show 
\begin{equation}\label{IHLHS}
[2]_q(\id \ot \ii)\circ(\id \ot \pp)\circ (\cupa \ot \id  ) -[2]_q(\ii\ot \id) (\pp\ot \id)\circ(\id\ot \cupa)
\end{equation}
is equal to 
\begin{equation}\label{IHRHS}
\id \ot \cupa  - \cupa  \ot \id.
\end{equation}
Since $V^{\mathcal{A}}(\varpi_1)$ is generated by the vector $v_{(1, 0)}$ it suffices to check that \eqref{IHLHS} and \eqref{IHRHS} send $v_{(1,0)}$ to the same vector in $V^{\mathcal{A}}(\varpi_1)\ot V^{\mathcal{A}}(\varpi_1)\ot V^{\mathcal{A}}(\varpi_1)$.

From \eqref{cup1formula}, \eqref{formulaforpi}, and \eqref{iformula} it follows that
\begin{equation}
[2]_q(\id \ot \ii)\circ(\id \ot \pp)\circ (\cupa \ot \id  )(v_{(1, 0)}) -[2]_q(\ii\ot \id) (\pp\ot \id)\circ(\id\ot \cupa)(v_{(1, 0)}) 
\end{equation}
is equal to
\begin{equation}\label{assoc1}
\begin{split}
-q^{-3}v_{(1, 0)}\ot&\ii(v_{(0,0)}) +q^{-2}[2]_qv_{(-1, 1)}\ot \ii(v_{(2, -1)}) - [2]_qv_{(1, -1)}\ot \ii(v_{(0,1)}) \\
&+ q^{-3}[2]_q\ii(v_{(0,1)})\ot v_{(1, -1)} -q^{-1}[2]_q\ii(v_{(2, -1)})\ot v_{(-1, 1)} +q \ii(v_{(0,0)})\ot v_{(1, 0)}
\end{split}
\end{equation}
Using \eqref{cup1formula}, we also find that
\begin{equation}
\id \ot \cupa  (v_{(1, 0)}) - \cupa  \ot \id(v_{(1, 0)})
\end{equation}
is equal to
\begin{equation}\label{assoccap1}
\begin{split}
&v_{(1, 0)}\ot \left(-q^{-4}v_{(1,0)}\ot v_{(-1,0)} + q^{-3} v_{(-1,1)} \ot v_{(1, -1)} - q^{-1} v_{(1, -1)} \ot v_{(-1, 1)} + v_{(-1, 0)}\ot v_{(1, 0)}\right) \\
&-\left(q^{-4}v_{(1,0)}\ot v_{(-1,0)} - q^{-3} v_{(-1,1)} \ot v_{(1, -1)} + q^{-1} v_{(1, -1)} \ot v_{(-1, 1)} - v_{(-1, 0)}\ot v_{(1, 0)}\right)\ot v_{(1, 0)}
\end{split}
\end{equation}
Using \eqref{iformula} to show that \eqref{assoc1} = \eqref{assoccap1} is left as an exercise.
\end{proof}

%===========
\subsection{Background on Tilting Modules. }
\label{subsec-objects}
%===========

Let $\ak$ be a field and let $q\in \ak^{\times}$ be such that $q+ q^{-1} \ne 0$. We will write $U_q^{\ak}(\mathfrak{sp}_4) = \ak \ot U_q^{\mathcal{A}}(\mathfrak{sp}_4)$, and $U_q^{\ak}(\mathfrak{sp}_4)-\text{mod}$ for the category of finite dimensional $U_q^{\ak}(\mathfrak{sp}_4)$ modules which are direct sums of their weight spaces and so that $K_{\alpha}$ acts on the $\mu$ weight space as $q^{(\mu, \alpha^{\vee})}$. 

Everything we say in this section is well-known to experts, but the results are essential for our arguments so we include some discussion for completeness. Two excellent references are Jantzen's book \cite{JantzenAgps} (only the second edition contains the appendix on representations of quantum groups and the appendix on tilting modules) and the eprint \cite{tiltnotes}. To deal with specializations when $q$ is an even root of unity we will also need some results from \cite{HansenqKempf} and \cite{Kanedabased}.

For each $\lambda\in X_+$ there is a \textbf{dual Weyl module} of highest weight $\lambda$, denoted $\nabla^{\ak}(\lambda)$, which is defined as an induced module \cite[H.11]{JantzenAgps}. The dual Weyl modules are a direct sum of their weight spaces and therefore have formal characters. Recall that we wrote $\VV(\lambda)$ for the irreducible module $\mathfrak{sp}_4(\mathbb{C})$ module of highest weight $\lambda$. We will write $[\VV(\lambda)]$ for the formal character of $\VV(\lambda)$ in $\mathbb{Z}[X]$, the group algebra of the weight lattice. It is known that a $q$-analogue of Kempf's vanishing holds for any $\ak$ \cite{HansenqKempf}. This implies that dual Weyl modules have formal character $[\VV(\lambda)]$ \cite[Theorem 5.12]{Andersen:1991wl}.

The dual Weyl module always has a unique simple submodule with highest weight $\lambda$. We will denote this module by $L^{\ak}(\lambda)$. The module $L^{\ak}(\lambda)$ should not be thought of as a base change of $\VV(\lambda)$. In fact quite often the two modules will have distinct formal characters. 

Since $U_q^{\ak}(\mathfrak{sp}_4)$ is a Hopf-algebra, it acts on the dual vector space of any finite dimensional representation. Then we define the \text{Weyl module} of highest weight $\lambda$ by $V^{\ak}(\lambda) = \nabla^{\ak}(-w_0\lambda)^*$ \cite[H.15]{JantzenAgps}. The dual Weyl module $V^{\ak}(\lambda)$ has the same formal character as $\nabla^{\ak}(\lambda)$, i.e. $[\VV(\lambda)]$, and $V^{\ak}(\lambda)$ has a unique simple quotient isomorphic to $L^{\ak}(\lambda)$. 

\begin{remark}\label{weylgroupduality}
In type $C_2$ the longest element $w_0$ acts on the weight lattice as $-1$. Therefore $V^{\ak}(\lambda) = \nabla^{\ak}(\lambda)^*$. 
\end{remark}

\begin{defn}
A \textbf{tilting module} is a module which has a (finite) filtration by Weyl modules, and a (finite) filtration by dual Weyl modules. The category of tilting modules, denoted $\textbf{Tilt}(U_q^{\ak}(\mathfrak{sp}_4))$, is the full subcategory of $U_q^{\ak}(\mathfrak{sp}_4)-\text{mod}$ where the objects tilting modules. 
\end{defn}

\begin{prop}\label{weylfilt}
The tensor product of two Weyl modules
\[
V^{\ak}(\lambda_1)\ot V^{\ak}(\lambda_2)
\]
has a filtration by Weyl modules.
\end{prop}
\begin{proof}
That this holds over $\ak$ follows from \cite{Kanedabased} where the result is shown to hold integrally using the theory of crystal bases.
\end{proof}

\begin{cor}
The tensor product of two tilting modules is a tilting module. 
\end{cor}
\begin{proof}
Since $(-)^*$ is exact, it follows from proposition \eqref{weylfilt} that the tensor product of dual Weyl modules 
\[
V^{\ak}(\lambda_1)^*\ot V^{\ak}(\lambda_2)^*
\]
has a filtration by dual Weyl modules. Thus the tensor product of two tilting modules will have a Weyl filtration and a dual Weyl filtration and is therefore a tilting module. 
\end{proof}

\begin{prop}\label{extvanishing}
Let $\lambda, \mu\in X_+$. Then $\dim_{\ak}\Ext^i(V^{\ak}(\lambda), \nabla^{\ak}(\mu))= \delta_{i, 0}\delta_{\lambda, \mu}$ for all $i\ge 0$.
\end{prop}
\begin{proof}
A standard argument \cite[Proof of Claim 3.1]{tiltnotes} shows that the vanishing of higher extension groups follows from Kempf's vanishing \cite{HansenqKempf}.
\end{proof}

\begin{prop}\label{omnibustilt}
The category $\textbf{Tilt}(U_q^{\ak}(\mathfrak{sp}_4))$ is closed under direct sums, direct summands, and tensor products. The isomorphism classes of indecomposable objects in the category are in bijection with $X_+$. We will write $T^{\ak}(\lambda)$ for the indecomposable tilting module corresponding to the dominant integral weight $\lambda$. The module $T^{\ak}(\lambda)$ is characterized as the unique indecomposable tilting module with a one dimensional $\lambda$ highest weight space.
\end{prop}
\begin{proof}
\cite[E.3-E.6]{JantzenAgps}.
\end{proof}

\begin{lemma}\label{grothendieckbasis}
Weyl modules or dual Weyl modules give a basis for the Grothendieck group of $U_q^{\ak}(\mathfrak{sp}_4)-\text{mod}$.
\end{lemma}
\begin{proof}
Both $V^{\ak}(\lambda)$ and $\nabla^{\ak}(\lambda)$ have the same formal character: $[\VV(\lambda)]$. In particular $V^{\ak}(\lambda)$ and $\nabla^{\ak}(\lambda)$ both have one dimensional $\lambda$ weight spaces. 
\end{proof}

For a tilting module $T$, we will write $(T: V^{\ak}(\lambda))$ to denote the filtration multiplicity. Formal character considerations also imply that $(T: V^{\ak}(\lambda))= (T: V^{\ak}(\lambda)^*)$ \cite[E.10]{JantzenAgps}.

\begin{lemma}
The following are equivalent.
\begin{enumerate}\label{weylistilting}
\item The Weyl module $V^{\ak}(\lambda)$ is simple.
\item $V^{\ak}(\lambda)\cong \nabla^{\ak}(\lambda)$
\item The Weyl module $V^{\ak}(\lambda)$ is a tilting module.
\end{enumerate}
\end{lemma}
\begin{proof}
It is not hard to see (1) implies (2) implies (3) \cite[E.1]{JantzenAgps}. That (3) implies (2) follows from Lemma \eqref{grothendieckbasis}, along with the equality of formal characters $[V^{\ak}(\lambda)]= [\nabla^{\ak}(\lambda)]$. To see that (2) implies (1), observe that the composition 
\[
L^{\ak}(\lambda) \rightarrow \nabla^{\ak}(\lambda) \xrightarrow{\sim} V^{\ak}(\lambda) \rightarrow L^{\ak}(\lambda)
\]
is non-zero on the $\lambda$ weight space. So the composition is a non-zero endomorphism of a simple module and therefore is an isomorphism. Thus, $L^{\ak}(\lambda)$ is a direct summand of $\nabla^{\ak}(\lambda)$. Since $\nabla^{\ak}(\lambda)$ has a simple socle, we may conclude that $\nabla^{\ak}(\lambda) \cong L^{\ak}(\lambda)$. 
\end{proof}

\begin{lemma}\label{extvanlemma}
\begin{enumerate} 
\item If $X$ has a filtration by Weyl modules, then for all $\lambda\in X_+$
\[
\text{dim}\Hom_{U_q^{\ak}(\mathfrak{sp}_4)}(X, \nabla^{\ak}(\lambda)) = (X:V^{\ak}(\lambda)).
\] 
\item If $Y$ has a filtration by dual Weyl modules, then for all $\lambda \in X_+$
\[
\text{dim}\Hom_{U_q^{\ak}(\mathfrak{sp}_4)}(V^{\ak}(\lambda) ,Y) = (Y:\nabla^{\ak}(\lambda)).
\] 
\end{enumerate}
\end{lemma}

\begin{proof}
Both claims follow from \eqref{extvanishing} and a long exact sequence argument.
\end{proof}

\begin{prop}\label{dimtilt}
If $T$ and $T'$ are tilting modules, then 
\begin{equation}
\dim \Hom_{U_q^{\ak}(\mathfrak{sp}_4)}(T, T') = \sum_{\lambda\in X_+} (T: V^{\ak}(\lambda))(T': V^{\ak}(\lambda)).
\end{equation}
\end{prop}

\begin{proof}
Since $T$ has both Weyl and dual Weyl filtrations, this follows from \ref{extvanlemma} and the fact that $(T': \nabla^{\ak}(\lambda))= (T': V^{\ak}(\lambda))$. 
\end{proof}

%===========
\subsection{The Image of the Evaluation Functor and Tilting Modules. }
\label{subsec-objects}
%===========

We continue with our assumption that $\ak$ is a field and $q\in \ak^{\times}$ so that $q+q^{-1} \ne 0$.

\begin{defn}
The category $\Fund(U_q^{\ak}(\mathfrak{sp}_4))$, is defined to be the full subcategory of $\text{Rep}(U_q^{\ak}(\mathfrak{sp}_4))$ with objects $V^{\ak}(\underline{w})= V^{\ak}(w_1)\ot V^{\ak}(w_2)\ot \ldots\ot V^{\ak}(w_n)$, where $\underline{w} = w_1w_2\ldots w_n$ and $w_i \in \lbrace \blues, \greent\rbrace$. 
\end{defn}

After changing coefficients to $\ak$, the functor from Theorem \eqref{evalthm} becomes
\begin{equation}
\ak \ot \eval: \DDk \longrightarrow \Fund(U_q^{\ak}(\mathfrak{sp}_4))
\end{equation}
We will abuse notation and write $\eval$ for $\ak \ot \eval$. 

\begin{lemma}\label{fundaretilt}
The modules $V^{\ak}(\underline{w})$ are tilting modules.
\end{lemma}
\begin{proof}
From the description of the integral forms of the modules in \eqref{fund1} and \eqref{fund2}, it is easy to see that $V^{\ak}(\varpi_1)$ and $V^{\ak}(\varpi_2)$ are irreducible with highest weight $\varpi_1$ and $\varpi_2$. They also have the same formal character as $[\VV(\varpi_1)]$ and $[\VV(\varpi_2)]$ respectively. So \eqref{weylistilting} implies that $V^{\ak}(\underline{w})$ is a tensor products of tilting modules and therefore is a tilting module. 
\end{proof}

\begin{remark}
If $q+ q^{-1} = 0$, then the Weyl module $V^{\ak}(\varpi_1)$ is still simple and therefore tilting but the Weyl module $V^{\ak}(\varpi_2)$ is not. In particular, $V^{\ak}(\varpi_2)$ has two Jordan--H\"{o}lder factors, a simple socle isomorphic to $L^{\ak}(0)$ and the simple quotient $L^{\ak}(\varpi_2)$.
\end{remark}

\begin{lemma}
For all $\underline{w}$ and $\underline{u}$ 
\begin{equation}
\dim_{\ak} \Hom_{U_q^{\ak}(\mathfrak{sp}_4)}(V^{\ak}(\underline{w}), V^{\ak}(\underline{u}))= \dim_{\mathbb{C}}\Hom_{\mathfrak{sp}_4(\mathbb{C})}(\VV(\underline{w}), \VV(\underline{u})). 
\end{equation}
\end{lemma}
\begin{proof}
Suppose that
\begin{equation}
\VV(\underline{w}) \cong \bigoplus_{\lambda}\VV(\lambda)^{m_{\lambda}},
\end{equation}
so we have an equality of formal characters $[\VV(\underline{w})]= \sum m_{\lambda}[\VV(\lambda)]$. Since $[V^{\ak}(\underline{w})]= [\VV(\underline{w})]$ and $[V^{\ak}(\lambda)] = [\VV(\lambda)]$ it follows that $(V^{\ak}(\underline{w}): V^{\ak}(\lambda)) = m_{\lambda}$. The claim then follows from proposition \eqref{dimtilt} and \eqref{LLsdimhom}
\end{proof}

\begin{thm}\label{mainthm}
The functor 
\[
\eval: \DDk \longrightarrow \Fund(U_q^{\ak}(\mathfrak{sp}_4)).
\]
is a monoidal equivalence.
\end{thm}

\begin{proof}
The functor $\eval$ is monoidal and essentially surjective, so it suffices to prove $\eval$ is full and faithful. 

Let $\underline{w}$ and $\underline{u}$ be objects in $\DD$. In the next section we will prove that $\eval(\mathbb{LL}_{\underline{w}}^{\underline{u}})$ is a linearly independent set of homomorphisms in $\Fund(U_q^{\ak}(\mathfrak{sp}_4))$.

Since
\begin{equation}
\# \mathbb{LL}_{\underline{w}}^{\underline{u}} = \dim_{\mathbb{C}}\Hom_{\mathfrak{sp}_4(\mathbb{C})}(\VV(\underline{w}), \VV(\underline{u})) = \dim_{\ak} \Hom_{U_q^{\ak}(\mathfrak{sp}_4)}(V^{\ak}(\underline{w}), V^{\ak}(\underline{u})),
\end{equation}
the linear independence of $\eval(\mathbb{LL}_{\underline{w}}^{\underline{u}})$ implies that $\eval$ maps $\mathbb{LL}_{\underline{w}}^{\underline{u}}$ to a basis in $\Fund(U_q^{\ak}(\mathfrak{sp}_4))$. These observations imply that $\mathbb{LL}_{\underline{w}}^{\underline{u}}$ is a linearly independent set of homomorphisms in $\DDk$. From the inequality in Lemma \eqref{hominequality} we deduce that $\mathbb{LL}_{\underline{w}}^{\underline{u}}$ is a basis. So $\eval$ maps a basis to a basis and $\Hom_{\DD}(\underline{w}, \underline{u})\xrightarrow{\eval} \Hom_{U_q^{\ak}(\mathfrak{sp}_4)}(V^{\ak}(\underline{w}), V^{\ak}(\underline{u}))$ is an isomorphism. 
\end{proof}

\begin{cor}
The functor $\eval$ induces a monoidal equivalence between the Karoubi envelope of $\DDk$ and the category $\textbf{Tilt}(U_q^{\ak}(\mathfrak{sp}_4))$. 
\end{cor}
\begin{proof}
Tensor products and direct summands of tilting modules are tilting modules. Therefore, Lemma \eqref{fundaretilt} implies that every direct summand of $V^{\ak}(\underline{w})$ is a tilting module. 

Let $\lambda\in X_+$, so $\lambda = a\varpi_1 + b\varpi_2$ for $a, b\in \mathbb{Z}_{\ge 0}$. The module $V^{\ak}(\blues^{\ot a} \ot \greent^{\ot b})$ has a one dimensional $\lambda$ highest weight space and all other non-zero weight spaces in $X_+$ are less than $\lambda$. From \eqref{omnibustilt} we deduce that $V^{\ak}(\blues^{\ot a} \ot \greent^{\ot b})$ must contain $T^{\ak}(\lambda)$ as a direct summand. Therefore every indecomposable tilting module is a direct summand of some $V^{\ak}(\underline{w})$. 
\end{proof}

\begin{remark}
If we take $\ak$ to be an algebraically closed field of characteristic $p$ and let $q= 1$, then $\textbf{Tilt}(U_q^{\ak}(\mathfrak{sp}_4))$ is equivalent to the category of tilting modules for the reductive algebraic group $\text{Sp}_4(\ak)$ \cite[H.6]{JantzenAgps}. Very little is known about tilting modules for reductive groups in characteristic $p>0$, and our results apply in this setting as well for all $p> 2$. 
\end{remark}

%%%%%%%%%%%%%%%%%%%%%%%%%
\section{Double Ladders are Linearly Independent}
\label{sec-doubleladders}
%%%%%%%%%%%%%%%%%%%%%%%%%

%===========
\subsection{Outline of the Argument}
\label{subsec-outline}
%===========

In this section we will finish the proof of Theorem \eqref{mainthm} by arguing that the set $\eval(\mathbb{LL}_{\underline{w}}^{\underline{u}})$ is linearly independent for all words $\underline{w}$ and $\underline{u}$. 

The idea of the proof is best illustrated as follows. Suppose we just wanted to prove that the image of light ladder diagrams from $\underline{w}$ to $\emptyset$ are linearly independent. Recall that $E(\underline{w}, 0)$ is the set of dominant weight subsequences $\vec{\mu} = (\mu_1, \mu_2, \ldots, \mu_n)$, so that $\sum \mu_i = 0$. Assume that for each dominant weight subsequence in $E(\underline{w}, 0)$, we have fixed a choice of light ladder $LL_{\vec{\mu}}$ and a vector $v_{\vec{\mu}}\in V^{\ak}(\underline{w})$. Consider the following matrix of elements in $\ak$. 
\begin{equation}\label{uppertriangularmatrix}
\left( \eval (LL_{\vec{\mu}})(v_{\vec{\nu}})\right)_{\vec{\mu}, \vec{\nu} \in E(\underline{w}, 0)}
\end{equation}
If \eqref{uppertriangularmatrix} is upper triangular with invertible elements of $\ak$ on the diagonal, then a non-trivial linear dependence among the maps $\eval(LL_{\vec{\mu}})$ will give rise to a non-zero vector in the kernel of the matrix \eqref{uppertriangularmatrix}. 

In the following subsections we will fix a choice of vectors associated to dominant weight subsequences. Then, since we want to argue double ladder diagrams are linearly independent, we must consider the image of the dominant weight subsequence vectors under both light ladders and upside down light ladders. The inductive construction of light ladders allows us to reduce these calculations to elementary light ladders, neutral ladders, and upside down elementary light ladders. In the end we still deduce linear independence of double ladder diagrams from an upper triangularity argument. 

%===========
\subsection{Subsequence Basis}
\label{subsec-subseq}
%===========

Recall that the modules $V^{\ak}(\blues)$ \eqref{fund1basis} and $V^{\ak}(\greent)$ \eqref{fund2basis} both have a fixed basis of weight vectors $v_{\nu}$ for $\nu\in \wt V^{\ak}(\blues)\cup \wt V^{\ak}(\greent)$. 

\begin{defn} 
Fix $\underline{w}= (w_1, \ldots, w_n)$, a word in the alphabet $\lbrace \blues, \greent\rbrace$, and let
\begin{equation}
S(\underline{w}) := \lbrace (\nu_1, \ldots \nu_n) \ : \  \nu_i\in \wt V^{\ak}(w_i)\rbrace.
\end{equation}
We set
\begin{equation}
v_{\underline{w}, +} := v_{w_1}\ot v_{w_2}\ot \ldots \ot v_{w_n}
\end{equation}
where $v_{\blues}= v_{(1, 0)}$ and $v_{\greent}= v_{(0,1)}$. Also, for any sequence of weights  $\vec{\nu} = (\nu_1, ..., \nu_n) \in S(\underline{w})$, we define
\begin{equation}
v_{\underline{w}, \vec{\nu}} := v_{\nu_1}\ot \ldots \ot v_{\nu_n}\in V^{\ak}(\underline{w}). 
\end{equation}
The \textbf{subsequence basis} of $V^{\ak}(\underline{w})$ is the set 
\begin{equation}
\lbrace v_{\vec{\nu}} \ : \ \vec{\nu} \in S(\underline{w})\rbrace.
\end{equation}
\end{defn}

\begin{lemma}\label{subbasis}
The subsequence basis of $V^{\ak}(\underline{w})$ is a basis of $V^{\ak}(\underline{w})$. 
\end{lemma}
\begin{proof}
This is clear.
\end{proof}

\begin{defn}
Let $\chi\in X_+$. The $\chi$ \textbf{weight space} of $V^{\ak}(\underline{w})$, denoted $V^{\ak}(\underline{w})[\chi]$, is the $\ak$-span of the subsequence basis vectors $v_{\vec{\nu}}$ such that $\sum \nu_i = \chi$. 
\end{defn}

Note that $E(\underline{w})\subset S(\underline{w})$. In particular, for each $\vec{\nu}\in E(\underline{w})$ we get a subsequence basis vector $v_{\underline{w}, \vec{\nu}}$. In the special case that the dominant weight subsequence is such that $\nu_i = \wt w_i$ for all $i$, then $v_{\underline{w}, \vec{\nu}}=v_{\underline{w}, +}$. Also, there is a partition of the set of dominant weight subsequences of $\underline{w}$:
\begin{equation}
E(\underline{w}) = \bigcup_{\lambda\in X_+} E(\underline{w}, \lambda),
\end{equation}
where $\vec{\nu}\in E(\underline{w})$ is in $E(\underline{w}, \lambda)$ whenever $\sum \nu_i = \lambda$ or equivalently $v_{\underline{w}, \vec{\nu}}\in V^{\ak}(\underline{w})[\lambda]$.

\begin{defn}\label{orderdefn}
Recall that our choice of simple roots was $\Delta = \lbrace \alpha_s, \alpha_t\rbrace$. There is a partial order on the set of weights defined by $\mu \le \nu$ if $\nu- \mu \in \mathbb{Z}_{\ge 0} \cdot  \Delta$. If we restrict this partial order to the set $\wt V^{\mathcal{A}}(\blues)\cup \wt V^{\mathcal{A}}(\greent)$, the resulting order is:
\begin{equation}
(-1, 0) < (1, -1) < (-1, 1) < (1, 0) 
\end{equation}
\begin{equation}
(0,-1)< (-2, 1)< (0,0) < (2, -1)< (0,1). 
\end{equation}

The lexicographic order gives a total order on the set $S(\underline{w})$. We will transport this total order to give a total order on the subsequence basis. 
\end{defn}

\begin{example}
In the image of $E(\greent\blues\greent\blues, (2, 0))\longrightarrow V^{\ak}(\greent\blues\greent\blues)[(2, 0)]$ we have, 
\[
v_{((0, 1), (1, 0), (2, -1), (-1, 0))} > v_{((0, 1), (1, 0), (0, -1), (1, 0))} > v_{((0, 1), (1, -1), (0, 0), (1, 0))}.
\]
\end{example}

\begin{lemma}\label{wtiszerolemma}
If $\wt \underline{w}\ngeq \chi$, then $V^{\ak}(\underline{w})[\chi] =0$. 
\end{lemma}
\begin{proof}
If $\vec{\nu}\in S(\underline{w})$ is such that $\nu_i\in \wt V^{\ak}(w_i)$, then $\sum \nu_i \le \wt \underline{w}$. The subsequence basis spans $V^{\ak}(\underline{w})$, so whenever $V^{\ak}(\underline{w})[\chi]\ne 0$, we must have $\chi \le \wt\underline{w}$. 
\end{proof}

%===========
\subsection{The Evaluation Functor and Elementary Diagrams}
\label{subsec-lightladdercalc}
%===========

\begin{notation}
In the remainder of the section, we will use the same notation for diagrammatic morphisms and their image under the functor $\eval$. But instead of saying diagram we will say map, for example the image of a light ladder diagram under $\eval$ will be referred to as a light ladder map.

To further simplify some of the statements below, our convention is that $\underline{w}$ and $\underline{u}$ are words in the alphabet $\lbrace \blues, \greent \rbrace$ and $\xi$ represents an invertible element of $\ak$.
\end{notation}

Recall that for each weight $\mu\in \wt V^{\ak}(\blues)\cup \wt V^{\ak}(\greent)$ there is an elementary light ladder diagram. The images of the elementary light ladder diagrams under the evaluation functor are the following elementary light ladder maps:
\begin{align*}
L_{(1, 0)} = \id &: V^{\ak}(\blues)\rightarrow V^{\ak}(\blues) \\
L_{(-1, 1)} = \pp&: V^{\ak}(\blues\blues)\rightarrow V^{\ak}(\greent)\\
L_{(1, -1)} =(\id \ot \capa   )\circ (\ii\ot \id )&: V^{\ak}(\greent\blues)\rightarrow V^{\ak}(\blues) \\
L_{(-1, 0)} = \capa   &: V^{\ak}(\blues\blues)\rightarrow \ak \\
L_{(0,1)} = \id &: V^{\ak}(\greent)\rightarrow V^{\ak}(\greent) \\
L_{(2, -1)} = (\id \ot \capb  \ot \id)\circ (\ii\ot \ii)&: V^{\ak}(\greent\greent)\rightarrow V^{\ak}(\blues\blues) \\
L_{(0,0)} = (\capa   \ot \id )\circ (\id \ot \ii)&: V^{\ak}(\blues\greent)\rightarrow V^{\ak}(\blues)\\
L_{(-2, 1)} = \pp\circ (\id \ot \capa   \ot\id )\circ (\id  \ot \id \ot \ii) &: V^{\ak}(\blues\blues\greent) \rightarrow V^{\ak}(\greent) \\
L_{(0,-1)} = \capb  &: V^{\ak}(\greent\greent)\rightarrow\ak. 
\end{align*}
\noindent There are two simple neutral diagrams, and their images under the evaluation functor are the simple neutral maps:
\begin{align*}
N_{\blues\greent}^{\greent\blues}= (\pp\ot \id )\circ (\id \ot \ii)&:V^{\ak}(\blues\greent)\rightarrow V^{\ak}(\greent \blues) \\
N_{\greent\blues}^{\blues\greent}= (\id \ot \pp)\circ (\ii\ot \id )&:V^{\ak}(\greent\blues)\rightarrow V^{\ak}(\blues\greent).
\end{align*}

\begin{lemma}\label{preserveweight}
If $f: V^{\ak}(\underline{w})\longrightarrow V^{\ak}(\underline{u})$ is a morphism which is in the image of the functor $\eval$, then $f: V^{\ak}(\underline{w})[\chi]\longrightarrow V^{\ak}(\underline{u})[\chi]$, for all $\chi\in X$. 
\end{lemma}
\begin{proof}
It is well known that every $U_q^{\ak}(\mathfrak{sp}_4)$ module homomorphism between finite dimensional modules will preserve weight spaces. But we could also deduce this from observing that the maps $\id$, $\ii$, and the cap and cup maps all preserve weight spaces and that any map in the image of $\eval$ is a linear combination of vertical and horizontal compositions of these basic maps.
\end{proof}

Recall that to construct light ladder diagrams and double ladder diagrams we need to fix a word $\underline{x}_{\lambda}$ in $\blues$ and $\greent$ for all $\lambda\in X_+$, and we need to make choices of neutral diagrams in the algorithmic construction. We now fix an $\underline{x}_{\lambda}$ for all $\lambda\in X_+$ and fix a light ladder diagram $LL_{\underline{w}, (\mu_1, \ldots, \mu_m)}$ for all $\underline{w}$ and all $(\mu_1, \ldots, \mu_m)\in E(\underline{w})$. This allows us to construct double ladder diagrams. The double ladder maps are the image of these double ladder diagrams under the evaluation functor. 

\begin{rmk}
The form of the arguments below do not depend on our choice of light ladder maps. 
\end{rmk}

%===========
\subsection{Pairing Vectors and Neutral Maps}
\label{subsec-neutralmappair}
%===========

\begin{lemma}\label{neutralladders}
If $N: V^{\ak}(\underline{w})\rightarrow V^{\ak}(\underline{u})$ is a neutral map, then $N(v_{\underline{w}, +}) = \xi \cdot v_{\underline{u}, +}$. Furthermore, if $(\mu_1, \ldots, \mu_n)$ is a sequence of weights so that $\mu_i\in \wt V^{\ak}(w_i)$, and $N(v_{\underline{w}, (\mu_1, \ldots, \mu_n)})$ has a non-zero coefficient for $v_{\underline{u}, +}$ after being written in the subsequence basis, then $v_{\underline{w}, (\mu_1, \ldots, \mu_n)}= v_{\underline{w}, +}$. 
\end{lemma}

\begin{proof} Neutral maps are vertical and horizontal compositions of identity maps, and the basic neutral maps $N_{\blues\greent}^{\greent\blues}$ and $N_{\blues\greent}^{\greent\blues}$. The lemma will follow from verifying its validity for the two basic neutral maps. 

The following maps factor through $V^{\ak}(\blues)$:
\begin{equation}
I_{\blues \greent}^{\greent \blues}:= \mathbb{D}(L_{(1, -1)})\circ L_{(0,0)} \ \ \ \ \ \text{and} \ \ \ \ \ I_{\greent \blues}^{\blues \greent}:= \mathbb{D}(L_{(0, 0)})\circ L_{(1,-1)}.
\end{equation}
Since $V^{\ak}(\blues)$ contains no vectors of weight $\varpi_1 + \varpi_2$, it follows that 
\begin{equation}
I_{\blues\greent}^{\greent\blues}(v_{(1, 0)}\ot v_{(0,1)})= 0 \ \ \ \ \ \text{and} \ \ \ \ \  I_{\greent\blues}^{\blues\greent}(v_{(0,1)}\ot v_{(1, 0)})= 0.
\end{equation}

It is easy to use the diagrammatic relations to compute that the maps
\begin{equation}
b_{\blues\greent}^{\greent\blues} = qN_{\blues\greent}^{\greent\blues} + q^{-1}I_{\blues\greent}^{\greent\blues}
\end{equation}
and
\begin{equation}
b_{\greent\blues}^{\blues\greent} = q^{-1}N_{\greent\blues}^{\blues\greent} + qI_{\greent\blues}^{\blues\greent}
\end{equation}
are mutual inverses. 

Both $b_{\blues\greent}^{\greent\blues}$ and $b_{\greent\blues}^{\blues\greent}$ are isomorphisms so they restrict to isomorphisms of weight spaces. Since the $\varpi_1 + \varpi_2$ weight spaces of $V^{\ak}(\blues\greent)$ and $V^{\ak}(\greent \blues)$ are one dimensional, it follows that $N_{\blues\greent}^{\greent\blues}$ sends the vector $v_{(1, 0)}\ot v_{(0,1)}$ to a non-zero scalar multiple of $v_{(0,1)}\ot v_{(1,0)}$ and $N_{\blues\greent}^{\greent\blues}$ sends $v_{(0, 1)}\ot v_{(1, 0)}$ to a non-zero multiple of $v_{(1,0)}\ot v_{(0,1)}$. Furthermore, the only subsequence basis vector which $N_{\blues\greent}^{\greent\blues}$ sends to a non-zero multiple of $v_{(0,1)}\ot v_{(1, 0)}$ is $v_{(1, 0)}\ot v_{(0,1)}$, and the only subsequence basis vector which $N_{\greent\blues}^{\blues\greent}$ sends to a non-zero multiple of $v_{(1, 0)}\ot v_{(0,1)}$ is $v_{(0, 1)}\ot v_{(1, 0)}$. 
\end{proof}

%===========
\subsection{Pairing Vectors and Light Ladders}
\label{subsec-lightladderpair}
%===========

\begin{lemma}\label{elemladders}
Let $\ast \in \lbrace \blues, \greent\rbrace$ and $\mu\in \wt V^{\ak}(\ast)$. Then the map $\id\ot L_{\mu}: V^{\ak}(\underline{w})\ot V^{\ak}(\ast)\rightarrow V^{\ak}(\underline{u})$, is such that for all $\nu \in \wt(V^{\ak}(\ast))$, 
\begin{equation}
\id\ot L_{\mu}(v_{\underline{w}, +}\ot v_{\nu}) =  \begin{cases} 
      0  \ \ \ \ \ \text{if} \  \nu> \mu \\
      \xi \cdot v_{\underline{u}, +}  \ \ \ \ \ \text{if} \ \nu = \mu.
       \end{cases}
\end{equation}
\end{lemma}

\begin{proof} It suffices to check the claim for $L_{\mu}$ and not all $\id\ot L_{\mu}$. The claim is obvious for $L_{(1, 0)}$ and $L_{(0,1)}$. For the rest of the cases, the claim follows from the calculation in section \eqref{subsec-lightladdercalc}. Note that in the $L_{\mu}$ step of the calculation, the first non-zero entry is $v_{\mu}\mapsto \xi \cdot v_{\underline{u}, +}$.
\end{proof}

Let $\vec{\mu} = (\mu_1, \ldots, \mu_n)\in E(\underline{w}, \lambda)$. The light ladder map $LL_{\underline{w}, \vec{\mu}}: V^{\ak}(\underline{w})\rightarrow V^{\ak}(\underline{x}_{\lambda})$ restricts to a map
\begin{equation}
LL_{\underline{w}, \vec{\mu}}: V^{\ak}(\underline{w})[\lambda]\longrightarrow V^{\ak}(\underline{x}_{\lambda})[\lambda].
\end{equation}
Moreover, $V^{\ak}(\underline{x}_{\lambda})[\lambda] = \ak \cdot v_{\underline{x}_{\lambda}, +}$. There is also a totally ordered set of linearly independent vectors in $V^{\ak}(\underline{w})[\lambda]$, namely $v_{\underline{w}, \vec{\nu}}$ for all $\vec{\nu}=  (\nu_1, \ldots, \nu_n)\in E(\underline{w}, \lambda)$. 

\begin{prop}\label{lightleavesunitri} 
\begin{equation}
LL_{\underline{w}, \vec{\mu}}(v_{\underline{w}, \vec{\nu}}) =  \begin{cases} 
      0  \ \ \ \ \ \text{if} \  \vec{\nu}> \vec{\mu} \\
      \xi \cdot v_{\underline{x}_{\lambda}, +}  \ \ \ \ \ \text{if} \ \vec{\nu} = \vec{\mu}.
       \end{cases}
\end{equation}
\end{prop}

\begin{proof} 
By the inductive definition of the light ladder map $LL_{\underline{w}, \vec{\mu}}$ and of the vector $v_{\underline{w}, \vec{\nu}}$, this proposition follows from repeated use of Lemmas \eqref{elemladders} and \eqref{neutralladders}.
\end{proof}

%===========
\subsection{Pairing Vectors and Upside Down Light Ladders}
\label{subsec-upsidelightladderpair}
%===========

In the results of the previous subsection we found the lexicographic order on sequences of weights was adapted to light ladders. There is another order on weights which is convenient for upside down light ladders. 

\begin{defn}
Fix $\underline{w}$ and let $\vec{\mu}= (\mu_1, \ldots , \mu_n)$ and $\vec{\nu} = (\nu_1, \ldots , \nu_n)$ be sequences of weights such that $\mu_i, \nu_i\in \wt V^{\ak}(w_i)$. Define a total order $< ^{\mathbb{D}}$ on weight sequences by setting $\vec{\nu}< ^{\mathbb{D}} \vec{\mu}$ if $(\nu_n, \ldots , \nu_1)< (\mu_n, \ldots , \mu_1)$ in the lexicographic order. We may also transport this order to give a total order on the subsequence basis. 
\end{defn}

\begin{lemma}\label{upsideelemladders}
Let $\ast \in \lbrace \blues, \greent\rbrace$ and $\mu\in \wt V^{\ak}(\ast)$. Then the map $\id\ot \mathbb{D}(L_{\mu}): V^{\ak}(\underline{w})\rightarrow V^{\ak}(\underline{u})\ot V^{\ak}(\ast)$ is such that

\begin{equation}
\id \ot\mathbb{D}(L_{\mu})(v_{\underline{w}, +}) = \xi\cdot v_{\underline{u}, +} \ot v_{\mu}+\sum c_{\vec{\tau}} \cdot v_{\underline{u}, \vec{\tau}}\ot v_{\nu}, \ \ \ \ \ c_{\vec{\tau}} \in \ak,
\end{equation}
where $v_{\underline{u}, \vec{\tau}}\ot v_{\nu}$ is a subsequence basis vector, $v_{\nu} > v_{\mu}$, and $v_{\underline{u}, \vec{\tau}}< v_{\underline{u}, +}$. 

%Moreover, if $v_{\vec{\sigma}}$ is a subsequence basis vector less than $v_{\underline{w}, +}$, then
%\begin{equation}
%\mathbb{D}(L_{\mu})(v_{\vec{\sigma}}) = \text{"lower terms"},
%\end{equation}
%where "lower terms" means a linear combination of subsequence basis vectors so that if $v_{\vec{\tau}}\ot v_{\nu}$ appears with non-zero coefficient, then $v_{\vec{\tau}}\ot v_{\nu} < v_{\underline{u}, +} \ot v_{\mu}$ and $v_{\vec{\tau}}\ot v_{\nu}= v_{\underline{u}, +}\ot v_{\nu}$ only when $\nu< \mu$. 
\end{lemma}

\begin{proof} 
It suffices to check the claim for $\mathbb{D}(L_{\mu})$ and not all $\id\ot \mathbb{D}(L_{\mu})$. The claim is obvious for $\mathbb{D}(L_{(1, 0)})$ and $\mathbb{D}(L_{(0,1)})$. The rest of the cases follow from the calculation in section \eqref{subsec-upsidelightladdercalc}. Note that the first line in the $\mathbb{D}(L_{\mu})$ calculation is $v_{\underline{w}, +}\mapsto \xi \cdot v_{\underline{u}, +}\ot v_{\mu}$, while the remaining terms are of the form $v_{\underline{u}, \vec{\tau}} \ot v_{\nu}$ where $\nu> \mu$. 
\end{proof} 

Let $\vec{\mu} = (\mu_1, \ldots, \mu_n) \in E(\underline{w}, \lambda)$. The associated upside down light ladder map $\mathbb{D}(LL_{\underline{w}, \vec{\mu}}): V^{\ak}(\underline{x}_{\lambda})\longrightarrow V^{\ak}(\underline{w})$ restricts to a map
\begin{equation}
\mathbb{D}(LL_{\underline{w}, \vec{\mu}}): V^{\ak}(\underline{x}_{\lambda})[\lambda]\longrightarrow V^{\ak}(\underline{w})[\lambda].
\end{equation}

\begin{prop}\label{upsidedownlowerterms} 
\begin{equation}
\mathbb{D}(LL_{\underline{w}, \vec{\mu}})(v_{\underline{x}_{\lambda}, +}) = \xi \cdot v_{\underline{w}, \vec{\mu}}  + \sum c_{\vec{\tau}} \cdot v_{\underline{w}, \vec{\tau}}, \ \ \ \ \ c_{\vec{\tau}} \in \ak,
\end{equation}
where $v_{\underline{w}, \vec{\mu}} <^{\mathbb{D}}v_{\underline{w}, \vec{\tau}}$.
\end{prop}

\begin{proof} 
By the inductive definition of the light ladder map $LL_{\underline{w}, (\mu_1, \ldots, \mu_n)}$, this proposition follows from repeated use of Lemmas \eqref{upsideelemladders} and \eqref{neutralladders}.
\end{proof}

%===========
\subsection{Proof of Linear Independence}
\label{subsec-doubelladderproof}

\begin{thm} 
The set 
\begin{equation}
\mathbb{LL}_{\underline{w}}^{\underline{u}} = \bigcup_{\lambda\in X_+}\mathbb{LL}_{\underline{w}}^{\underline{u}}(\lambda)
\end{equation}
is a linearly independent subset of $\Hom_{U_q^{\ak}(\mathfrak{sp}_4)}(V^{\ak}(\underline{w}), V^{\ak}(\underline{u}))$. 
\end{thm}

\begin{proof} 
Let 
\begin{equation}
\sum_{\lambda}\sum_{\substack{\vec{\mu}\in E(\underline{w}, \lambda)\\ \vec{\nu}\in E(\underline{u}, \lambda)}} {^{\lambda}c_{\vec{\mu}}^{\vec{\nu}}}\cdot \mathbb{LL}_{\underline{w}, \vec{\mu}}^{\underline{u}, \vec{\nu}} = 0, \ \ \ \ \ ^{\lambda}c_{\vec{\mu}}^{\vec{\nu}} \in \ak
\end{equation}
be a nontrivial linear relation. There is at least one $\lambda_0\in X_+$ with ${^{\lambda_0}c_{\vec{\mu}}^{\vec{\nu}}}\ne 0$ so that if ${^{\lambda}c_{\vec{\mu}}^{\vec{\nu}}}\ne 0$ then $\lambda \ngtr\lambda_0$. Lemma \eqref{wtiszerolemma} implies that for all $\lambda\ne \lambda_0$ with ${^{\lambda}c_{\vec{\mu}}^{\vec{\nu}}}\ne 0$, $V^{\ak}(\underline{x}_{\lambda})[\lambda_0]= 0$. If $v_0\in V^{\ak}(\underline{w})[\lambda_0]$, then since light ladder maps preserve the weight of a vector \eqref{preserveweight} 
\begin{equation}\label{reducetolambdazero}
0 = \sum_{\lambda}\sum_{\vec{\mu}, \vec{\nu}} {^{\lambda}c_{\vec{\mu}}^{\vec{\nu}}}\cdot \mathbb{LL}_{\underline{w}, \vec{\mu}}^{\underline{u}, \vec{\nu}}(v_0) = \sum_{\vec{\mu}, \vec{\nu}} {^{\lambda_0}c_{\vec{\mu}}^{\vec{\nu}}}\cdot \mathbb{LL}_{\underline{w}, \vec{\mu}}^{\underline{u}, \vec{\nu}}(v_0).
\end{equation}
Note that for $\vec{\mu}\in E(\underline{w}, \lambda_0)$, $v_{\underline{w}, \vec{\mu}}\in V^{\ak}(\underline{w})[\lambda_0]$. 

Let $\vec{\mu_0}$ be the largest $\vec{\mu}$, in the lexicographic order, so that ${^{\lambda_0}c_{\vec{\mu}}^{\vec{\nu}}}\ne 0$. Taking $v_0 = v_{\underline{w}, \vec{\mu_0}}$ in \eqref{reducetolambdazero} results in
\begin{equation}
0 = \sum_{\vec{\mu}, \vec{\nu}} {^{\lambda_0}c_{\vec{\mu}}^{\vec{\nu}}}\cdot \mathbb{LL}_{\underline{w}, \vec{\mu}}^{\underline{u}, \vec{\nu}}(v_{\underline{w}, \vec{\mu_0}})= \sum_{\vec{\mu}, \vec{\nu}} {^{\lambda_0}c_{\vec{\mu}}^{\vec{\nu}}}\cdot \mathbb{D}(LL_{\underline{u}, \vec{\nu}})\circ LL_{\underline{w}, \vec{\mu}}(v_{\underline{w}, \vec{\mu_0}}).
\end{equation}
Proposition \eqref{lightleavesunitri} implies
\begin{equation}
0 = \sum_{\vec{\nu}} {^{\lambda_0}c_{\vec{\mu_0}}^{\vec{\nu}}}\cdot \mathbb{D}(LL_{\underline{u}, \vec{\nu}})\circ LL_{\underline{w}, \vec{\mu_0}}(v_{\underline{w}, \vec{\mu_0}}) \\
= \sum_{\vec{\nu}} {^{\lambda_0}c_{\vec{\mu_0}}^{\vec{\nu}}}\xi\cdot \mathbb{D}(LL_{\underline{u}, \vec{\nu}})(v_{\underline{x}_{\lambda}, +}). 
\end{equation}

Let $\vec{\nu_0}$ be the smallest $\vec{\nu}$, in the $<^{\mathbb{D}}$ order, so that ${^{\lambda_0}c_{\vec{\mu_0}}^{\vec{\nu}}}\ne 0$. Proposition \eqref{upsidedownlowerterms} implies
\begin{equation}
\begin{split}
0 &= {^{\lambda}c_{\vec{\mu_0}}^{\vec{\nu_0}}}\xi \cdot \mathbb{D}(LL_{\underline{u}, \vec{\nu_0}})(v_{\underline{x}_{\lambda}, +}) + \sum_{\vec{\nu_0}<^{\mathbb{D}}\vec{\nu}} {^{\lambda_0}c_{\vec{\mu_0}}^{\vec{\nu}}}\xi\cdot \mathbb{D}(LL_{\underline{u}, \vec{\nu}})(v_{\underline{x}_{\lambda}, +}) \\
&= {^{\lambda_0}c_{\vec{\mu_0}}^{\vec{\nu_0}}}\xi\cdot v_{\underline{u}, \vec{\nu_0}} + \text{``higher terms"},
\end{split}
\end{equation}
where ``higher terms" is a linear combination of subsequence basis vectors all of which are greater than $v_{\underline{u}, \vec{\nu_0}}$ in the $<^{\mathbb{D}}$ order. Since the subsequence basis vectors are linearly independent, we must have ${^{\lambda_0}c_{\vec{\mu_0}}^{\vec{\nu_0}}}\xi = 0$, which is a contradiction. 
\end{proof}

%===========
\subsection{Elementary Light Ladder Calculations}
\label{subsec-lightladdercalc}
%===========

\begin{equation}
L_{(-1,1)}(v_{(1,0)}\ot (-)):  \begin{cases} 
      v_{(1, 0)} &\mapsto 0 \\
v_{(-1, 1)}&\mapsto -v_{(0,1)}\\
v_{(1, -1)}&\mapsto -v_{(2, -1)}\\
v_{(-1,0)}&\mapsto \dfrac{-q}{[2]_q}v_{(0,0)},
       \end{cases}
\end{equation}

\begin{equation}
L_{(1, -1)}(v_{(0,1)}\ot (-)): \begin{cases} 
      v_{(1, 0)} &\mapsto 0 \\
v_{(-1, 1)}&\mapsto 0\\
v_{(1, -1)}&\mapsto -v_{(1, 0)}\\
v_{(-1,0)}&\mapsto -v_{(-1, 1)},
       \end{cases}
\end{equation}

\begin{equation}
L_{(-1, 0)}(v_{(1,0)}\ot (-)):  \begin{cases} 
      v_{(1, 0)} &\mapsto 0 \\
v_{(-1, 1)}&\mapsto 0\\
v_{(1, -1)}&\mapsto 0\\
v_{(-1,0)}&\mapsto 1,
       \end{cases}
\end{equation}

\begin{equation}
L_{(2, -1)}(v_{(0,1)}\ot (-)):  \begin{cases} 
      v_{(0, 1)} &\mapsto 0 \\
v_{(2, -1)}&\mapsto v_{(1, 0)}\ot v_{(1, 0)}\\
v_{(0,0)}&\mapsto v_{(1, 0)}\ot v_{(-1,1)} + q^{-1}v_{(-1, 1)}\ot v_{(1, 0)}\\
v_{(-2, 1)}&\mapsto v_{(-1, 1)}\ot v_{(-1, 1)}\\
v_{(0,-1)}&\mapsto -v_{(1, 0)}\ot v_{(-1, 0)} + v_{(-1, 1)}\ot v_{(1, -1)},
       \end{cases}
\end{equation}

\begin{equation}
L_{(0, 0)}(v_{(1, 0)}\ot (-)): \begin{cases} 
      v_{(0, 1)} &\mapsto 0 \\
v_{(2, -1)}&\mapsto 0 \\
v_{(0,0)}&\mapsto -q^{-1}v_{(1, 0)}\\
v_{(-2, 1)}&\mapsto -v_{(-1, 1)}\\
v_{(0,-1)}&\mapsto -v_{(1, -1)},
       \end{cases}
\end{equation}

\begin{equation}
L_{(-2, 1)}(v_{(1, 0)}\ot v_{(1, 0)}\ot (-)): \begin{cases} 
      v_{(0, 1)} &\mapsto 0\\
v_{(2, -1)}&\mapsto  0\\
v_{(0,0)}&\mapsto 0\\
v_{(-2, 1)}&\mapsto v_{(0,1)}\\
v_{(0,-1)}&\mapsto v_{(2, -1)},
       \end{cases}
\end{equation}

\begin{equation}
L_{(0, -1)}(v_{(0, 1)}\ot (-)):  \begin{cases} 
      v_{(0, 1)} &\mapsto 0\\
v_{(2, -1)}&\mapsto  0\\
v_{(0,0)}&\mapsto 0\\
v_{(-2, 1)}&\mapsto 0\\
v_{(0,-1)}&\mapsto 1.
       \end{cases}
\end{equation}

%===========
\subsection{Upside Down Elementary Light Ladder Calculations}
\label{subsec-upsidelightladdercalc}
%===========

\begin{equation}\begin{split}
\mathbb{D}(L_{(-1, 1)}): v_{(0,1)} \mapsto q^{-1}v_{(1,0)} \ot &v_{(-1,1)} \\ - v_{(-1,1)}\ot &v_{(1,0)} 
%v_{(2,-1)} &\mapsto q^{-1}v_{(1, 0)}\ot v_{(1, -1)} - v_{(1, -1)} \ot v_{(1,0)} \\
%v_{(0,0)} &\mapsto q^{-1} v_{(1,0)} \ot v_{(-1,0)} + q^{-2}v_{(-1,1)} \ot v_{(1,-1)} - v_{(1,-1)}\ot v_{(-1,1)}  -q^{-1} v_{(-1,0)} \ot v_{(1,0)} \\
%v_{(-2, 1)} &\mapsto q^{-1} v_{(-1,1)} \ot v_{(-1,0)} - v_{(-1,0)} \ot v_{(-1,1)} \\
%v_{(0,-1)} &\mapsto q^{-1} v_{(1, -1)} \ot v_{(-1,0)} - v_{(-1,0)} \ot v_{(1, -1)},
\end{split}
\end{equation}

\begin{equation}
\begin{split}
\mathbb{D}(L_{(1, -1)}):  v_{(1,0)} \mapsto -q^{-3}v_{(0,1)}\ot &v_{(1, -1)} \\+ q^{-1}v_{(2, -1)}\ot &v_{(-1, 1)} \\- \dfrac{q}{[2]_q}v_{(0,0)}\ot &v_{(1, 0)}
%v_{(-1, 1)} &\mapsto  -q^{-3}v_{(0,1)}\ot v_{(-1, 0)} - \dfrac{q^{-1}}{[2]_q}v_{(0,0)}\ot v_{(-1,1)} -v_{(-2, 1)}\ot v_{(1,0)}\\
%v_{(1,-1)} &\mapsto  -q^{-3}v_{(2, -1)} \ot v_{(-1, 0)} - \dfrac{q^{-1}}{[2]_q}v_{(0,0)} \ot v_{(1, -1)} - v_{(0,-1)}\ot v_{(1, 0)}\\
%v_{(-1,0)} &\mapsto \dfrac{-q^{-3}}{[2]_q}v_{(0,0)} \ot v_{(-1, 0)} +q^{-2}v_{(-2, 1)}\ot v_{(1, -1)} -v_{(0,-1)}\ot v_{(-1, 1)},
\end{split}
\end{equation}

\begin{equation}
\begin{split}
\mathbb{D}(L_{(-1, 0)}): 1\mapsto -q^{-4}v_{(1,0)}\ot &v_{(-1,0)} \\+ q^{-3} v_{(-1,1)} \ot &v_{(1, -1)} \\- q^{-1} v_{(1, -1)} \ot &v_{(-1, 1)} \\+ v_{(-1, 0)}\ot &v_{(1, 0)},
\end{split}
\end{equation}

\begin{equation}
\begin{split}
\mathbb{D}(L_{(2, -1)}):v_{(1, 0)}\ot v_{(1, 0)} \mapsto -q^{-2}v_{(0,1)}\ot &v_{(2, -1)} \\+ v_{(2, -1)}\ot &v_{(0,1)}
\end{split}
\end{equation}

\begin{equation}
\begin{split}
\mathbb{D}(L_{(0,0)}): v_{(1,0)} \mapsto \dfrac{-q^{-3}}{[2]_q}v_{(1, 0)}\ot &v_{(0,0)} \\+ q^{-2}v_{(-1, 1)}\ot &v_{(2, -1)} \\-v_{(1, -1)}\ot &v_{(0,1)}\\
%v_{(-1, 1)} &\mapsto -q^{-3}v_{(1, 0)}\ot v_{(-2, 1)} + \dfrac{q^{-1}}{[2]_q} v_{(-1,1)}\ot v_{(0,0)} - v_{(-1,0)}v_{(0,1)} \\
%v_{(1,-1)} &\mapsto -q^{-3} v_{(1, 0)}\ot v_{(0,-1)} + \dfrac{q^{-1}}{[2]_q}v_{(1, -1)} \ot v_{(0,0)} -v_{(-1,0)}\ot v_{(2, -1)}\\
%v_{(-1,0)} &\mapsto -q^{-3}v_{(-1, 1)}\ot v_{(0, -1)} +q^{-1}v_{(1, -1)}\ot v_{(-2, 1)} - \dfrac{q}{[2]_q}v_{(-1,0)}\ot v_{(0,0)},
\end{split}
\end{equation}

\begin{equation}
\begin{split}
\mathbb{D}(L_{(-2, 1)}): v_{(0, 1)}\mapsto -q^{-4}v_{(1, 0)}\ot v_{(1, 0)}\ot &v_{(-2, 1)} \\+ \dfrac{q^{-2}}{[2]_q}v_{(1, 0)}\ot v_{(-1, 1)}\ot &v_{(0, 0)} + \dfrac{q^{-3}}{[2]_q} v_{(-1, 1)}\ot v_{(1, 0)}\ot v_{(0, 0)} \\
 - q^{-2} v_{(-1, 1)}\ot v_{(-1, 1)} \ot &v_{(2, -1)}\\- q^{-1} v_{(1, 0)}\ot v_{(-1, 0)}\ot &v_{(0, 1)} + v_{(-1, 1)}\ot v_{(1, -1)}\ot v_{(0, 1)}
\end{split}
\end{equation}

and 
\begin{equation}
\begin{split}
\mathbb{D}(L_{(0, -1)}): 1\mapsto q^{-6} v_{(0,1)}\ot &v_{(0,-1)} \\- q^{-4}v_{(2, -1)}\ot &v_{(-2, 1)} \\+\dfrac{q^{-2}}{[2]_q}v_{(0,0)}\ot &v_{(0,0)} \\- q^{-2} v_{(-2,1)}\ot &v_{(2, -1)} \\+ v_{(0,-1)}\ot &v_{(0,1)}.
\end{split}
\end{equation}

%===========
\subsection{Object Adapted Cellular Category Structure}
\label{subsec-objectadapted}
%===========

We refer to \cite[Definition 2.4]{ELauda} for the definition of a strictly object adapted cellular category or SOACC. 

Let $\ak$ be a field and let $q\in \ak^{\times}$ such that $q+ q^{-1} \ne 0$. In this section we will show that $\DDk$ is an SOACC. It follows that the endomorphism algebras in $\DDk$ are cellular algebras. Since we proved that $\DDk$ is equivalent to $\Fund(U_q^{\ak}(\mathfrak{sp}_4))$, the result about cellular algebras also follows from \cite{tiltcellular}. For more discussion about the relation between our work and \cite{tiltcellular} we recommend \cite[p. 6]{elias2015light} (but replace $\mathfrak{sl}_n$ webs with $\DD$). 

For each $\lambda \in X_+$, choose an object $\underline{x}_{\lambda}$ in $\DDk$ so that $\wt \underline{x}_{\lambda} = \lambda$. The set $\Lambda = \lbrace \underline{x}_{\lambda}\rbrace_{\lambda \in X_+}$ is in bijection with $X_+$, and we define a partial order on $\Lambda$ by setting $\underline{x}_{\lambda} \le \underline{x}_{\mu}$ whenever $\lambda \le \mu$ i.e. $\mu - \lambda\in \mathbb{Z}_{\ge 0} \Phi_+$. 

For any object $\underline{w}$ in $\DDk$ and for all $\vec{\mu} \in E(\underline{w}, \lambda)$ we fix a light ladder diagram $LL_{\vec{\mu}}:= LL_{\underline{w}, \vec{\nu}}\in \Hom_{\DDk}(\underline{w}, \underline{x}_{\lambda})$ and an upside down light ladder diagram $\mathbb{D}(LL_{\vec{\nu}}):=\mathbb{D}(LL_{\underline{w}, \vec{\nu}})\in \Hom_{\DDk}(\underline{x}_{\lambda}, \underline{w})$.

If $\underline{x}_{\lambda} = x_1x_2\ldots x_n$ where $x_i\in \lbrace \blues, \greent\rbrace$, then the set $E(\underline{x}_{\lambda}, \lambda)$ contains a single element, $\vec{\lambda} = (\wt x_1, \wt x_2 \ldots \wt x_n)$. Recall that in our definition of double ladder diagrams we choose $LL_{\vec{\lambda}} = \id_{\underline{x}_{\lambda}} = \mathbb{D}(LL_{\vec{\lambda}})$. 

For $\vec{\mu} \in E(\underline{w}, \lambda)$ and $\vec{\nu}\in E(\underline{u}, \lambda)$ we set
\begin{equation}
\mathbb{LL}_{\vec{\mu}, \vec{\nu}}^{\lambda} := \mathbb{D}(LL_{\vec{\nu}})\circ LL_{\vec{\mu}} \in \Hom_{\DDk}(\underline{w}, \underline{u}).
\end{equation}
It follows from our main theorem that $\lbrace \mathbb{LL}_{\vec{\mu}, \vec{\nu}}^{\lambda} \rbrace_{\lambda \in X_+}$ forms a basis for $\Hom_{\DDk}(\underline{w}, \underline{u})$. 

\begin{remark}
In the definition of an SOACC, one fixes the data of two sets, $E(\underline{w}, \lambda)$ and $M(\underline{w}, \lambda)$, which are in a fixed bijection. We are choosing to ignore the set $M(\underline{w}, \lambda)$. 
\end{remark}

\begin{defn}
Fix $\lambda\in X_+$. Let $(\DDk)_{< {\lambda}}$ be the $\ak$-linear subcategory whose morphisms are spanned by $\mathbb{LL}_{\vec{\mu}, \vec{\nu}}^{\chi}$ with $\chi< \lambda$.
\end{defn}

\begin{lemma}
Let $f\in \Hom_{\DDk}(\underline{w}, \underline{u})$ and let $\vec{\mu} \in E(\underline{u}, \lambda)$. Then 
\begin{equation}
LL_{\vec{\mu}}\circ f \equiv \sum_{\vec{\nu}\in E(\underline{w}, \lambda)} \ast \cdot LL_{\vec{\nu}} \ \ \ \ \ \text{modulo} \  (\DDk)_{< {\lambda}},
\end{equation}
where $\ast$ represents an element of $\ak$.
\end{lemma}
\begin{proof}
Writing $LL_{\vec{\mu}} \circ f$ in the double ladder basis, we find that 
\begin{equation}
\begin{split}
LL_{\vec{\mu}}\circ f &=\sum_{\substack {\chi\in X_+ \\ \vec{\nu} \in E(\underline{w}, \chi) \\ \vec{\tau} \in E(\underline{x}_{\lambda}, \chi)}}\ast \cdot \mathbb{LL}_{\vec{\nu}, \vec{\tau}}^{\chi} \\
&\equiv \sum_{\substack {\vec{\mu} \in E(\underline{w}, \lambda) \\ \vec{\tau} \in E(\underline{x}_{\lambda}, \lambda)}}\ast \cdot \mathbb{LL}_{\vec{\nu}, \vec{\tau}}^{\lambda} \ \ \ \ \ \ \text{modulo} \  (\DDk)_{< {\lambda}} \\
&\equiv \sum_{\vec{\nu}\in E(\underline{w}, \lambda)} \ast \cdot LL_{\vec{\nu}}  \ \ \ \ \ \ \text{modulo} \  (\DDk)_{< {\lambda}}
\end{split}
\end{equation}
The second equality follows from the observation that if $\chi \in X_+$ and $E(\underline{x}_{\lambda}, \chi) \ne \emptyset$, then $\chi \le \lambda$. The third equality follows from recalling that $E(\underline{x}_{\lambda}, \lambda)= \lbrace \vec{\lambda}\rbrace$ and $LL_{\vec{\lambda}} = \id_{\underline{x}_{\lambda}}$. 
\end{proof}

\begin{cor}\label{Dissoacc}
The category $\DDk$ with fixed choices of $\underline{x}_{\lambda}$ and light ladder diagrams is an SOACC. 
\end{cor}

\bibliographystyle{plain}

\bibliography{mastercopy}

%\begin{thebibliography}{99}

%\end{thebibliography}

\end{document}